\documentclass{amsart}
\usepackage{amssymb,epsfig,amscd, verbatim}

\newtheorem{theorem}{Theorem}
\newtheorem{lemma}[theorem]{Lemma}
\newtheorem{proposition}[theorem]{Proposition}
\newtheorem{corollary}[theorem]{Corollary}

\theoremstyle{definition}
\newtheorem{definition}[theorem]{Definition}
\newtheorem{example}[theorem]{Example}
\newtheorem{examples}[theorem]{Examples}

\theoremstyle{remark}
\newtheorem{remark}[theorem]{Remark}

\newtheorem{question}[theorem]{Question}

\numberwithin{equation}{section}
\numberwithin{theorem}{section}
\newcommand{\thref}[1]{Theorem~{\rm\ref{#1}}}
\newcommand{\prref}[1]{Proposition~{\rm\ref{#1}}}
\newcommand{\leref}[1]{Lemma~{\rm\ref{#1}}}
\newcommand{\coref}[1]{Corollary~{\rm\ref{#1}}}

\newcommand{\deref}[1]{Definition~{\rm\ref{#1}}}
\newcommand{\exref}[1]{Example~{\rm\ref{#1}}}

\newcommand{\reref}[1]{Remark~{\rm\ref{#1}}}

\newcommand{\firef}[1]{Figure~{\rm\ref{#1}}}
\newcommand{\seref}[1]{Section~{\rm\ref{#1}}}
\newcommand\lbb[1]{\label{#1}}

\newcommand{\fig}[1]
        {\raisebox{-0.5\height}%
                  {\epsfbox{fig/#1}}%
        }
\def\as{associative}        
\def\conalg{conformal algebra}
\def\conalgs{conformal algebras}
\def\difalg{differential algebra}
\def\difalgs{differential algebras}
\def\psalg{pseudo\-algebra}
\def\psalgs{pseudo\-algebras}
\def\psten{pseudo\-tensor}
\def\tt{\otimes}                               
\def\<{\langle}
\def\>{\rangle}
\def\ti{\tilde}
\def\wti{\widetilde}
\def\what{\widehat}
\def\ov{\overline}

\def\d{\partial}
\def\smash{\,\sharp\,}

\def\tsum{{\textstyle\sum}}

\def\surjto{\twoheadrightarrow}                
\def\injto{\hookrightarrow}                    
\def\isoto{\xrightarrow{\sim}}                 
\def\st{\; | \;}                               
\def\symm{S}                                   

\def\cont{{\mathrm{cont}}}

\newcommand\xb[3]{[{#1}_{#2} {#3}]}
\newcommand\xp[3]{{#1}_{#2} {#3}}
    
\def\Cset{\mathbb{C}}       
\def\Zset{\mathbb{Z}}       
       
\def\Kset{{\mathbf{k}}}                
\def\di{{\mathrm{d}}}
\newcommand{\kk}{\mathbf{k}}
\newcommand{\CC}{\mathbb{C}}

\newcommand{\ZZ}{\mathbb{Z}}
\newcommand{\calA}{\mathcal{A}}
\newcommand{\fd}{{\mathfrak d}}
\newcommand{\fa}{\mathfrak{a}}
\newcommand{\fg}{\mathfrak{g}}

\newcommand{\fn}{\mathfrak{n}}
\newcommand{\fp}{\mathfrak{p}}
\def\al{\alpha}                         
\def\be{\beta}
\def\ga{\gamma}
\def\Ga{\Gamma}
\def\de{\delta}
\def\De{\Delta}
\def\ep{\varepsilon}
\def\la{\lambda}

\def\om{\omega}
\def\Om{\Omega}
\def\ph{\varphi}
\def\si{\sigma}
\def\Si{\Sigma}
\def\th{\theta}

\def\g{{\mathfrak{g}}}      
\def\h{{\mathfrak{h}}}      
\def\dd{{\mathfrak{d}}}
\def\ss{{\mathfrak{s}}}
\def\gl{{\mathfrak{gl}}}
\def\sl{{\mathfrak{sl}}}

\def\fo{{\mathfrak{o}}}
\def\sp{{\mathfrak{sp}}}
\def\csp{{\mathfrak{csp}}}

\def\Wd{W(\dd)}    
\def\Sd{S(\dd,\chi)}
\def\Hd{H(\dd,\chi,\om)}
\def\Kd{K(\dd,\th)}
\def\Wdp{W(\dd')}    
\def\Sdp{S(\dd',\chi')}
\def\Hdp{H(\dd',\chi',\om')}
\def\Kdp{K(\dd',\th')}
\def\A{{\mathcal{A}}}
\def\L{{\mathcal{L}}}
\def\V{{\mathcal{V}}}
\def\W{{\mathcal{W}}}
\def\M{\mathcal{M}}          
\def\N{\mathcal{N}}

\def\O{\mathcal{O}}          
\def\C{\mathcal{C}}           
\def\F{\mathcal{F}}           
           
\def\Vec{\mathcal{V}ec}       

\def\ue{U}                 
\DeclareMathOperator{\rank}{rank}
\DeclareMathOperator{\gw}{gw}
\DeclareMathOperator{\Gr}{Gr}
\DeclareMathOperator{\codim}{codim}
\DeclareMathOperator{\Span}{span}

\DeclareMathOperator{\Ind}{Ind}
\DeclareMathOperator{\ad}{ad}
\DeclareMathOperator{\tr}{tr}

\DeclareMathOperator{\hgt}{ht}
\DeclareMathOperator{\sd}{\ltimes}
\DeclareMathOperator{\id}{id}
\DeclareMathOperator{\Id}{Id}
\DeclareMathOperator{\Div}{div}
\DeclareMathOperator{\gplk}{G}     
\DeclareMathOperator{\primt}{P}    
\DeclareMathOperator{\symp}{S}     
\DeclareMathOperator{\fil}{F}      
\DeclareMathOperator{\coh}{H}      

\DeclareMathOperator{\Der}{Der}

\DeclareMathOperator{\Aut}{Aut}
\DeclareMathOperator{\Hom}{Hom}
\DeclareMathOperator{\Lin}{Lin}     
\DeclareMathOperator{\End}{End}
\DeclareMathOperator{\Chom}{Chom}
\DeclareMathOperator{\Cend}{Cend}
\DeclareMathOperator{\Cder}{Der}
\DeclareMathOperator{\gc}{gc}
\DeclareMathOperator{\oc}{oc}
\DeclareMathOperator{\spc}{spc}

\DeclareMathOperator{\Tor}{Tor}

\DeclareMathOperator{\Cur}{Cur}

\DeclareMathOperator{\coker}{coker}

\DeclareMathOperator{\Rad}{Rad}
\DeclareMathOperator{\Vir}{Vir}

\begin{document}

\title{Theory of finite pseudoalgebras}

\author[B.~Bakalov]{Bojko Bakalov}
\address{Department of Mathematics, MIT, Cambridge MA 02139, USA}
\curraddr{Department of Mathematics, University of California,
Berkeley CA 94720}
\email{bakalov@math.berkeley.edu}
\thanks{The first author was supported in part by NSF grant
    DMS-9622870 and A.~P.~Sloan Dissertation Fellowship.
This research was partially conducted by the first author for the 
Clay Mathematics Institute.}
\author[A.~D'Andrea]{Alessandro D'Andrea}
\address{Dipartimento di Matematica,
Istituto ``Guido Castelnuovo'',
Universit\`a di Roma ``La Sapienza'',
00185 Roma, Italia}
\email{dandrea@mat.uniroma1.it}
\thanks{The second author was supported in part by 
European Union TMR grant
ERB FMRX-CT97-0100 and CNR grant 203.01.71}
\author[V.~G.~Kac]{Victor G.~Kac}
\address{Department of Mathematics, MIT, Cambridge MA 02139, USA}
\email{kac@math.mit.edu}
\thanks{The third author was supported in part by NSF grants
DMS-9622870 and DMS-9970007}
\date{July 19, 2000}


\maketitle
\tableofcontents
\pagebreak

\section{Introduction}\lbb{sintro}
Since the seminal papers of Belavin, Polyakov and Zamolodchikov
\cite{BPZ} and of Borcherds \cite{Bo1} there has been a great deal 
of work towards understanding the algebraic structures
underlying the notion of operator product expansion (OPE) of
chiral fields in conformal field theory.

In physics literature the OPE of local chiral fields $\varphi$
and $\psi$ is written in the form \cite{BPZ}:
\begin{equation}
  \lbb{eq:0.1}
  \varphi (z) \psi (w)= \sum_{j \ll \infty}
      \frac{\varphi (w)_{(j)}\psi (w)}{(z-w)^{j+1}} ,
\end{equation}
where $\varphi (w)_{(j)}\psi (w)$ are some new fields,
which may be viewed as bilinear products of fields
$\varphi$ and $\psi$ for all $j \in \ZZ$ (see
e.g.~\cite{K2} for a rigorous interpretation of (\ref{eq:0.1})).  
If now $V$ is a space
of pairwise local chiral fields which contains $1$, is invariant
with respect to the derivative $\partial = \partial_w$, and is
closed under all $j${th} products, $j \in \ZZ$, we obtain an
algebraic structure which physicists (respectively mathematicians) call
a chiral (respectively vertex) algebra.  In more abstract terms, $V$ is a 
module over $\CC [\partial ]$ with a marked element~$1$ and
infinitely many bilinear over $\CC$ products $\varphi_{(j)}
\psi$, $j \in \ZZ$, satisfying a certain system of identities,
first written down by Borcherds \cite{Bo1}.  (An equivalent system
of axioms, which is much easier to verify, may be found in
\cite{K2}.)

One of the important features of the OPE (\ref{eq:0.1}) is that
its singular part encodes the commutation relations of fields,
namely one has (see e.g.~\cite{K2}):
\begin{equation}
  \lbb{eq:0.2}
  [\varphi (z) , \psi (w)]=\sum_{j \geq 0} \,
  \bigl(\varphi (w)_{(j)} \psi(w)\bigr) \,
  \partial^j_w \delta (z-w) /j!, 
\end{equation}
where $\delta (z-w)= \sum_{n\in\ZZ} \, z^n w^{-n-1}$ is the
delta-function.  This leads to the notion of a Lie conformal
algebra, which is a $\CC [\partial]$-module with $\CC$-bilinear
products $\varphi_{(j)}\psi$ for all non-negative integers $j$, subject to
certain identities \cite{K2}.  In order to write down these
identities in a compact form, it is convenient to consider the
formal Fourier transform of (\ref{eq:0.2}), called the
$\lambda$-bracket (where $\lambda$ is an indeterminate):
\begin{displaymath}
  [\varphi_{\lambda}\psi] = \sum_{j \geq 0}
  \frac{\lambda^j}{j!} (\varphi_{(j)}\psi) .
\end{displaymath}
Then a \emph{Lie conformal algebra} $L$ is defined as a $\CC
[\partial]$-module endowed with a $\CC$-linear map
\begin{displaymath}
  L \otimes L \to \CC [\lambda]\otimes L ,\qquad
  a \otimes b \mapsto [a_{\lambda}b]
\end{displaymath}
satisfying the following axioms \cite{DK} $(a,b,c \in L)$:
\begin{align*}
  \begin{array}{ll}
    \text{(sesquilinearity)}\quad &
    [\partial a_{\lambda}b]=-\lambda [a_{\lambda}b], \;\;
    [a_{\lambda}\partial b]=(\partial +\lambda)[a_{\lambda}b],
    \\
    \text{(skew-commutativity)}\quad &
    [b_{\lambda}a]=-[a_{-\lambda -\partial}b] ,\\
    \text{(Jacobi identity)}\quad &
    [a_{\lambda} [b_{\mu}c]] =
    [[a_{\lambda}b]_{\lambda +\mu}c]+
    [b_{\mu}[a_{\lambda}c]] .
  \end{array}
\end{align*}

In the past few years a structure theory \cite{DK},
representation theory \cite{CK, CKW} and cohomology theory 
\cite{BKV} of finite (i.e.,~finitely generated as $\CC
[\partial]$-modules) Lie conformal algebras have been worked out.  For
example, one of the main results of \cite{DK} states that any
finite simple Lie conformal algebra is isomorphic either to the \emph{Virasoro}
 conformal algebra:
 \begin{equation}\lbb{vir1}
   \Vir = \CC [\partial ] \ell,\quad [\ell_{\lambda}\ell]
   =(\partial +2 \lambda)\ell
 \end{equation}
or to the \emph{current} conformal algebra associated to a
simple finite-dimensional Lie algebra $\fg$:
\begin{equation}
  \Cur \fg = \CC [\partial]\otimes \fg , \quad
  [a_{\lambda}b]=[a,b], \qquad a,b \in \fg .
\end{equation}

The objective of the present paper is to develop a theory of
``multi-dimensional'' Lie conformal algebras, i.e.~a theory where
the algebra of polynomials $\CC [\partial]$ is replaced by a
``multi-dimensional'' associative algebra $H$.  In order to
explain the definition, let us return to the singular part
(\ref{eq:0.2}) of the OPE.  Choosing a set of generators $a^i$ of 
the $\CC [\partial]$-module $L$, we can write:
\begin{displaymath}
  [a^i_{\lambda}a^j]=\tsum_k \, Q^{ij}_k 
  (\lambda ,\partial) a^k ,
\end{displaymath}
where $Q^{ij}_k$ are some polynomials in $\lambda$ and
$\partial$.  The corresponding singular part of the OPE is:
\begin{displaymath}
  [a^i(z),a^j(w)] =\tsum_k \, Q^{ij}_k
  (\partial_w,\partial_t) (a^k(t) \delta (z-w))|_{t=w} .
\end{displaymath}
Letting $P^{ij}_k (x,y) = Q^{ij}_k (-x,x+y)$, we can rewrite this 
in a more symmetric form:
\begin{equation}
  \lbb{eq:0.3}
  [a^i(z),a^j(w)] =\tsum_k \, P^{ij}_k
  (\partial_z,\partial_w) (a^k(w) \delta (z-w)) .
\end{equation}
We thus obtain an $H=\CC [\partial]$-bilinear map (i.e.,~a map of 
$H \otimes H$-modules):
\begin{displaymath}
  L \otimes L \to (H \otimes H) \otimes_H L,\quad
  a \otimes b \mapsto [a*b]
\end{displaymath}
(where $H$ acts on $H \otimes H$ via the comultiplication
map $\Delta (\d)= \d \otimes 1 + 1\otimes \d$), defined by
\begin{displaymath}
  [a^i * a^j]=\tsum_k \, P^{ij}_k
  (\partial \otimes 1, 1 \otimes \partial) \otimes_H a^k .
\end{displaymath}
Hence the notion of $\lambda$-bracket $[a_{\lambda}b]$ is
equivalent to the notion of the $*$-bracket $[a*b]$ introduced by 
Beilinson and Drinfeld \cite{BD}, the relation between the two
brackets being given by letting $\lambda =-\partial \otimes 1$.
For example, the Virasoro conformal algebra \eqref{vir1} corresponds
to the {\em Virasoro pseudoalgebra}
\begin{equation}\lbb{vir2}
   \Vir = \CC [\partial ] \ell,\quad 
   [\ell * \ell] = (1\tt\d - \d\tt1) \tt_{\CC[\d]} \ell.
 \end{equation}

It is natural to introduce the general notion of a conformal 
algebra as a $\CC [\partial]$-module $L$ endowed with a
$\CC$-linear map $L \otimes L \to \CC [\lambda] \otimes L$, $a
\otimes b \mapsto a_{\lambda}b$ satisfying the sesquilinearity
property:
\begin{displaymath}
  (\partial a)_{\lambda}b=-\lambda (a_{\lambda}b), \quad 
  a_{\lambda}(\partial b) =(\partial +\lambda)(a_{\lambda}b) .
\end{displaymath}
Such a conformal algebra is called associative
(respectively commutative) if
\begin{displaymath}
  a_{\lambda}(b_{\mu}c)=(a_{\lambda}b)_{\lambda + \mu}
  c \quad (\text{respectively}\;\; b_{\lambda} a=a_{-\lambda-\partial}b) ,
\end{displaymath}
and the $\lambda$-product of an associative conformal algebra
defines a $\lambda$-bracket
\begin{displaymath}
  [a_{\lambda}b]=a_{\lambda}b-b_{-\lambda -\partial}a ,
\end{displaymath}
making it a Lie conformal algebra \cite{K4, DK}.

As above, we have the equivalent notion of a $*$-product on an
$H=\CC [\partial]$-module $L$, which is an $H$-bilinear map
\begin{equation}
  \lbb{eq:0.4}
  L \otimes L \to (H \otimes H) \otimes_H L , \quad
  a \otimes b \mapsto a * b .
\end{equation}
Now it is clear that the notion of a $*$-product can be defined
by (\ref{eq:0.4}) for any Hopf algebra $H$ by making use of the
comultiplication $\Delta\colon H \to H \otimes H$ to define 
$(H\otimes H) \otimes_H L$.  
A \emph{pseudoalgebra} is a (left) $H$-module $L$ endowed
with an $H$-bilinear map (\ref{eq:0.4}).  The name is motivated by 
the fact that this is an algebra in a pseudotensor category
(introduced in \cite{L}, \cite{BD}).  Accordingly, the $*$-product will be 
called a {\em pseudoproduct}.

One is able to define a pseudoproduct as soon as a 
structure of a bialgebra is given on $H$.
However, in order to generalize the equivalence
of a pseudoalgebra and an $H$-conformal algebra structure on an
$H$-module $L$, we need $H$ to be a Hopf algebra.  In this case
any element of $H \otimes H$ can be uniquely written as a
(finite) sum:
\begin{displaymath}
  \tsum_i \, (h_i \otimes 1) \Delta (f_i), \quad\text{where $h_i$ are
    linearly independent.}
\end{displaymath}
Hence the pseudoproduct on $L$ can be written in the form:
\begin{equation}
  \lbb{eq:0.5}
  a*b=\tsum_i \, (h_i \otimes 1) \otimes_H c_i .
\end{equation}
The corresponding {\em $H$-conformal algebra\/} structure is then a
$\CC$-linear map $L \otimes L \to H \otimes L$ given by
\begin{equation}
  \lbb{eq:0.6}
  ab=\tsum_i \, h_i \otimes c_i .
\end{equation}
Every element $x$ of $H^*$ then defines an {\em $x$-product\/} $L \otimes 
L \to L$:
\begin{equation}
  \lbb{eq:0.7}
  a_{x}b=\tsum_i \, \langle x , S(h_i) \rangle\, c_i ,
\end{equation}
where $S$ is the antipode of $H$.

The $H$-bilinearity property of the pseudoproduct (\ref{eq:0.5})
is, of course, easily translated to certain sesquilinearity
properties of the products~(\ref{eq:0.6}) and (\ref{eq:0.7}).  In 
particular, in the case $H=\CC [\partial]$, the product
(\ref{eq:0.6}) is the $\lambda$-product if we let $\lambda
=-\partial$, and the product (\ref{eq:0.7}) for $x=t^j$ is 
the $j${th} product described above, where $H^*\simeq\CC[[t]]$,
$\langle t,\d \rangle = 1$.
The equivalence of these three structures (discussed in
\seref{shconf}) is very useful in the
study of pseudoalgebras.

In order to define associativity of a pseudoproduct, we extend
it from $ L \otimes L \to H^{\otimes 2} \otimes_H L$ to 
$  (H^{\otimes 2} \otimes_H L) \otimes L \to H^{\otimes 3}
\otimes_H L$ and to $L \otimes (H^{\otimes 2}
\otimes_H L) \to H^{\otimes 3}\otimes_H L$ by letting:
\begin{align*}
(f\otimes_{H}a)*b &= \tsum_i \, (f \otimes 1)\, (\Delta\otimes\id)(g_i)
 \otimes_H c_i ,
\\
a*(f\otimes_{H}b) &= \tsum_i \, (1 \otimes f)\, (\id\otimes\Delta)(g_i)
 \otimes_H c_i ,
\quad\text{where}\;\; a*b=\tsum_i \, g_i \otimes_H c_i .  
\end{align*}
Then the associativity
property is given by the usual equality (in $H^{\otimes 3}
\otimes_H L$):
\begin{displaymath}
  (a*b)*c=a*(b*c) .
\end{displaymath}

The easiest example of a pseudoalgebra is a current
pseudoalgebra, defined as follows.  Let $H'$ be a Hopf subalgebra 
of $H$ and let $A$ be an $H'$-pseudoalgebra (for example, if
$H'=\CC$, then $A$ is an ordinary algebra over $\CC$).  Then the
associated \emph{current} $H$-pseudoalgebra is $\Cur A =H
\otimes_{H'}A$ with the pseudoproduct
\begin{displaymath}
  (f \otimes_{H'}a)*(g \otimes_{H'}b)=
  ((f \otimes g) \otimes_H 1)\, (a*b) .
\end{displaymath}
The $H$-pseudoalgebra $\Cur A$ is associative iff the
$H'$-pseudoalgebra $A$ is.

The most important example of an associative $H$-pseudoalgebra is 
the pseudoalgebra of all pseudolinear endomorphisms of a finitely 
generated $H$-module $V$, which is denoted by $\Cend V$ (see \seref{shcla}).  
A {\em pseudolinear endomorphism\/} of $V$ is a $\CC$-linear map $\phi \colon 
V \to (H \otimes H) \otimes_H V$ such that
\begin{align*}
  \phi (hv) &= ((1 \otimes h)\otimes_H 1) \,
  \phi (v),\qquad h \in H, v \in V .
\intertext{The space $\Cend V$ of all such $\phi$ becomes a (left) $H$-module
if we define}
  (h\phi)(v) &= ((h\otimes 1) \otimes_H1)\, \phi (v) .
\end{align*}
The definition of a pseudoproduct on $\Cend V$ is especially
simple when $V$ is a free $H$-module, $V=H \otimes V_0$, where 
$V_0$ is a finite-dimensional vector space over $\CC$
with a trivial action of $H$. Then
$\Cend V$ is isomorphic to $H \otimes H \otimes \End V_0$, with 
$H$ acting by left multiplication on the first factor, with
the following pseudoproduct:
\begin{displaymath}
  (f \otimes a \otimes A)* (g \otimes b \otimes B)
  =\tsum_i \, (f \otimes ga'_i) \otimes_H
  (1 \otimes b a''_i \otimes AB) ,
\end{displaymath}
where $\Delta (a) =\tsum_i \, a'_i \otimes a''_i $.

The main objects of our study are Lie pseudoalgebras.  The
corresponding pseudoproduct is conventionally called
{\em pseudobracket\/} and denoted by $[a*b]$.  Given an associative
pseudoalgebra with pseudoproduct $a*b$ we may give it a
structure of a Lie pseudoalgebra by defining the pseudobracket
\begin{displaymath}
  [a*b]= a*b - (\si\tt_H\id)\, b*a,
\end{displaymath}
where $\sigma\colon H \otimes H \to H \otimes H$ is the permutation
of factors.  
It is immediate to see that this pseudobracket satisfies the
following skew-commutativity and Jacobi identity axioms:
\begin{align}
  \lbb{eq:0.8}
  [b*a] &= -(\sigma \otimes_H \id)\, [a*b] , \\
  \lbb{eq:0.9}
  [a*[b*c]] &= [[a*b]*c]+((\sigma \otimes \id)\otimes_H\id)\, 
  [b*[a*c]] .
\end{align}

It is important to point out here that the above
pseudobracket and both identities are well defined, provided that 
the Hopf algebra $H$ is cocommutative.  A pseudoalgebra with
pseudoproduct $[a*b]$ satisfying identities (\ref{eq:0.8}) and
(\ref{eq:0.9}) is called a \emph{Lie pseudoalgebra}.  We will always
assume that $H$ is cocommutative when talking about Lie
pseudoalgebras.  Of course, the simplest examples of Lie
pseudoalgebras are $\Cur A$, where $A$ is a $H' (\subset H)$ Lie
pseudoalgebra ($=$ Lie algebra if $H'=\CC$).  It is needless to say that
in the case $H=\CC [\partial]$, $\Delta (\partial)=\partial
\otimes 1+1\otimes \partial$, the $H$-conformal algebras
associated to Lie pseudoalgebras are nothing else but the Lie
conformal algebras discussed above.

We will explain now the connection of the notion of a Lie
pseudoalgebra to the more classical notion of a differential Lie
algebra studied in \cite{R1}--\cite{R4}, \cite{C}, \cite{NW} and 
many other papers (see \seref{slalg}).  Let $Y$ 
be a commutative associative algebra over $\CC$ with compatible
left and right actions of the Hopf algebra $H$.  Then, given a
Lie pseudoalgebra $L$, we let $\calA_Y L =Y \otimes_H L$ with the 
obvious left $H$-module structure and the following Lie algebra
(over $\CC$) structure:
\begin{displaymath}
  [(x\otimes_H a), (y \otimes_H b)]=\tsum_i \, (xf_i)(yg_i)
  \otimes_H c_i \quad\text{if}\;\; [a*b] =\tsum_i \, (f_i\otimes g_i)
  \otimes_H c_i .
\end{displaymath}

Provided that $L$ is a free $H$-module, the Lie algebra $\calA_Y
L$ is a free $Y$-module, hence $\calA_Y L$ is a differential Lie
algebra in the sense of \cite{NW}.  The most classical case is
again $H=\CC [\partial]$, when $Y$ is simply a commutative
associative algebra with a (left and right) derivation
$\partial$, and we get the differential Lie algebras of Ritt
\cite{R1}--\cite{R4}.  Thus, the notion of a Lie pseudoalgebra is 
reminiscent of the notion of a group scheme:  each Lie
pseudoalgebra $L$, which is free as an $H$-module, gives rise to a
functor $\calA$ from the category of commutative associative
algebras with compatible left and right actions of $H$ to the
category of differential Lie algebras ($=$ category of formal
differential groups).

For example, the functor $\calA$ corresponding to the Virasoro
pseudoalgebra \eqref{vir2} associates to any commutative associative algebra 
$Y$ with a derivation~$'$ the differential Lie algebra $Y$ with
bracket $[u,v]=uv'-u'v$, called the substitutional Lie algebra by 
Ritt.  The current pseudoalgebra $\Cur \fg$, where $\fg$ is a
Lie algebra over $\CC$, associates to $Y$ the obvious
differential Lie algebra $Y \otimes \fg$.  Thus, a result of
\cite{DK} asserts that any simple finite differential Lie algebra 
with ``constant coefficients'' is isomorphic either to the
substitutional Lie algebra or to $Y \otimes \fg$ where $\fg$ is a 
simple finite-dimensional Lie algebra.  In the rank~$1$ case, but 
without the constant coefficients assumption, this is the main
result of \cite{R1}.

The main tool in the study of pseudoalgebras is the
\emph{annihilation algebra} $\calA_X L$, where $X=H^*$ is the
associative algebra dual to the coalgebra $H$.  We find it
remarkable that the annihilation algebra of the associative
pseudoalgebra $\Cend H=H \otimes H$ is nothing else but the
Drinfeld double (with the obvious comultiplication) of the Hopf
algebra $H$.  Note that in the associative case $Y$ need not be
commutative in order to define the functor $\calA_Y$, but in the
Lie algebra case it must be.  So, in order to construct the
annihilation Lie algebra we again use cocommutativity of $H$.

Recall that, by Kostant's \thref{tkostant}, any cocommutative Hopf
algebra $H$ is a smash product of a group algebra $\CC [\Gamma]$
and the universal enveloping algebra $\ue(\fd)$ of a Lie algebra
$\fd$. In Sections \ref{sgaalg} and \ref{shkga} we show that the theory of
pseudoalgebras over a smash product of $\CC [\Gamma]$ and any
Hopf algebra $H$ reduces to that over $H$.  This allows us in many cases
to assume, without loss of generality, that $H$ is the universal enveloping 
algebra of a Lie algebra $\fd$.

However, for most of our results we have to assume that $\fd$ is
finite dimensional.  In this case the algebra $H=\ue(\fd)$ is
Noetherian, and the annihilation algebra $\calA_X L$ is
linearly compact, provided that $L$ is \emph{finite}  
(i.e.,~finitely generated as an $H$-module).  Recall that a topological Lie
algebra is called {\em linearly compact\/} if its underlying topological
vector space is a topological product of finite-dimensional vector spaces 
with the discrete topology
(see \seref{scartan}).

In \seref{sreconst} we prove ``reconstruction'' theorems, which
claim that, under some mild assumptions, a Lie pseudoalgebra is
completely determined by its annihilation Lie algebra along with
the action of $\fd$.  This reduces the classification of finite
simple Lie pseudoalgebras to the well developed structure theory of
linearly compact Lie algebras, which goes back to E.~Cartan (see
\cite{G1, G2} and \seref{scartan}).

We turn now to examples of finite Lie pseudoalgebras beyond the rather
obvious examples of current Lie pseudoalgebras.  The first example is
the generalization of the Virasoro pseudoalgebra \eqref{vir2} defined for
$H=\CC [\partial]$ (which is the universal enveloping algebra of
the $1$-dimensional Lie algebra) to the case $H=\ue(\fd)$, where $
\fd$ is any finite-dimensional
Lie algebra.  This is the Lie pseudoalgebra $\Wd=
H \otimes \fd$ with pseudobracket
\begin{displaymath}
  [(1 \otimes a) * (1 \otimes b)] 
  = (1 \otimes 1) \otimes_H (1 \otimes [a,b])
  + (b \otimes 1) \otimes_H (1 \otimes a)
  -(1 \otimes a) \otimes_H (1 \otimes b) .
\end{displaymath}
Since the associated annihilation algebra $\A_X \Wd \simeq X\otimes\fd$ is
isomorphic to the Lie algebra of formal vector fields on the Lie
group $D$ with Lie algebra $\fd$, it is natural to call $\Wd$
the pseudoalgebra of all vector fields.  In fact we develop (in
\seref{svect}) a formalism of pseudoforms similar to the usual
formalism of differential forms,
which may be viewed as the beginning of a ``pseudo differential geometry''.

This allows us to define the
remaining three series of finite simple Lie pseudoalgebras:
$\Sd$, $\Hd$ and $\Kd$.  The
annihilation algebras of the simple Lie pseudoalgebras $\Wd$, 
$\Sd$, $\Hd$ and $\Kd$ are
isomorphic to the four series of Lie--Cartan linearly compact Lie
algebras $W_N$, $S_N$, $P_N$ (which is an extension of
$H_N$ by a $1$-dimensional center) and $K_N$, where $N=\dim\fd$.
However the Lie pseudoalgebras $\Sd$, $\Hd$ and 
$\Kd$ depend on certain parameters $\chi,\om$ and
$\th$, due to inequivalent actions of $\fd$ on the annihilation 
algebra.  The parameter $\chi$ is a $1$-dimensional
representation of $\fd$, i.e., $\chi\in\fd^*$ such that $\chi([\fd,\fd])=0$.
The parameter $\om$ is an element of $\fd^* \wedge \fd^*$ such that
$\om^{N/2} \ne 0$ and $\di\om + \chi\wedge\om = 0$ 
in the case $\Hd$, when $N$ is even.
The parameter $\th\in\fd^*$ is such that
$\th\wedge(\di\th)^{(N-1)/2} \ne 0$
in the case $\Kd$, when $N$ is odd.
In the cases $\Hd$, $\Kd$, these parameters are in one-to-one 
correspondence with ``nondegenerate'' skew-symmetric solutions 
$\al=r+s\tt1-1\tt s$ ($r\in\dd\wedge\dd$, $s\in\dd$)
of a modification of the classical Yang--Baxter equation,
which is a special case of the dynamical classical Yang--Baxter equation
(see \cite{Fe, ES}).
%

The central result of the paper is the classification of finite
simple Lie pseudoalgebras over the Hopf algebra $H=\ue(\fd)$.  As
usual, a Lie pseudoalgebra $L$ is called \emph{simple} if it is nonabelian
(i.e.,~$[L*L]\neq 0$) and its only ideals are $0$ and $L$.  Our
Theorem~\ref{classify} states that any such Lie pseudoalgebra is
isomorphic either to a current pseudoalgebra $\Cur \fg
=\Cur^H_{\CC}\fg$ over a simple finite-dimensional Lie algebra
$\fg$, or to a current pseudoalgebra $\Cur^H_{H'}L'$
over one of the
Lie pseudoalgebras $L' = \Wdp$, $\Sdp$, $\Hdp$ or $\Kdp$,
where $H'=\ue(\fd')$ and $\fd'$ is a subalgebra of $\fd$.

A Lie pseudoalgebra $L$ is called \emph{semisimple} if it
contains no nonzero abelian ideals.  One also defines in the
usual way the derived pseudoalgebra, solvable and nilpotent
pseudoalgebras, and for a finite Lie pseudoalgebra $L$ one has the 
solvable radical $\Rad L$ (so that $L/\Rad L$ is semisimple).

Our Theorem~\ref{tsemisim} states that any finite semisimple Lie 
$\ue(\fd)$-pseudoalgebra is a direct sum of finite simple Lie
pseudoalgebras and of Lie pseudoalgebras of the form $A \ltimes \Cur
\fg$, where $A$ is a subalgebra of $\Wd$ and $\fg$ is a
simple finite-dimensional Lie algebra.
In addition, in \thref{tsubwd} we show that any subalgebra of $\Wd$
is simple, and in \coref{csubwd} we give a complete list of all these
subalgebras.
(A more concise formulation of Theorem~\ref{classify} 
is that any finite simple Lie pseudoalgebra over $U(\dd)$ is either a 
current pseudoalgebra $\Cur\g$
over a simple finite-dimensional Lie algebra $\g$, or a nonzero subalgebra
of $\Wd$.)

Note, however, that Levi's theorem on $L$ being a semidirect sum
of $L/\Rad L$ and $\Rad L$ is not true even in the case $\dim \fd 
=1$.  This stems from the fact that the cohomology of simple Lie
pseudoalgebras with nontrivial coefficients is (highly)
nontrivial (see \seref{cohom} and \cite{BKV}), 
in a sharp contrast with the Lie algebra
case.  For example, it follows from \cite{BKV} that there are
precisely five cases (three isolated examples and two
families) of non-split extensions of $\Vir$ by $\Cur \CC$.
Translated into the language of differential Lie algebras, this
result goes back to Ritt \cite{R3}.

Closely related to the present paper are the papers \cite{Ki} and \cite{NW},
where (in our terminology) the annihilation algebras of rank $1$ over $H$
Lie \psalgs, and of simple Lie \psalgs\ of arbitrary finite rank, respectively,
are studied. In fact, our Theorems~\ref{classify} and \ref{tsemisim}
provide a completed form of the classification results of \cite{NW}
(in the ``constant coefficients'' case).

The structural results of the present paper in the simplest case $\dim \fd
=1$ reproduce the results of \cite{DK}.  However, this case is
much easier than the case $\dim \fd >1$, mainly due to the fact
that only in this case is any finite torsionless $H$-module free.

Note also the close connection of our work to Hamiltonian formalism in the 
theory of nonlinear evolution equations (see the review \cite{DN2}, the 
book \cite{Do} and references there,
and also \cite{GD}, \cite{DN1}, \cite{Z}, \cite{M}, \cite{X}, 
and many other papers). In \seref{poisson} we derive, as a corollary of 
Theorems \ref{classify} and \ref{tsemisim}, a classification of simple and 
semisimple linear Poisson brackets in any finite number of 
indeterminates.

In \seref{srepth} we develop a
representation theory of finite Lie pseudoalgebras.  First, we
prove an analogue of Lie's Lemma that any weight space for an ideal
of a Lie pseudoalgebra $L$ acting on a finite module is an
$L$-submodule (Proposition~\ref{lliel}).  This implies an
analogue of Lie's Theorem that a solvable Lie pseudoalgebra has an
eigenvector in any finite module (Theorem~\ref{tliet}), and an
analogue of Cartan--Jacobson Theorem that describes all finite Lie
pseudoalgebras which have a finite faithful irreducible module
(Theorem~\ref{tcj}).  Finally, we reduce the classification
and construction of finite irreducible modules over semisimple Lie
pseudoalgebras to that of irreducible modules over linearly
compact Lie algebras of the type studied by Rudakov \cite{Ru1, Ru2}
(the complete classification will appear in a future publication).
Note that complete reducibility fails
already in the simplest case of Lie pseudoalgebras with $\dim \fd =1$ 
\cite{CKW}.

In \seref{cohom} we define cohomology of Lie \psalgs\ and show that it
describes module extensions, abelian \psalg\ extensions, and \psalg\
deformations. We also relate this cohomology to the Gelfand--Fuchs cohomology
\cite{F}. These results generalize those of \cite{BKV} in the 
$\dim\dd=1$ case.

Note that in the case $\dim \fd =1$ Lie pseudoalgebras are
closely related to vertex algebras in a way similar to the
relation of Lie algebras to universal enveloping algebras \cite{K2}.
We expect that, under certain conditions, there is a similar
relation of ``multi-dimensional''
Lie  pseudoalgebras to ``multi-dimensional'' vertex algebras defined
in \cite{Bo2}.
In the case of a commutative Lie algebra $\fd$ the Lie pseudoalgebras 
encode the OPE between ultralocal fields (as well as the linear Poisson
brackets).
However, it is not 
clear how Lie pseudoalgebras are related to the OPE of realistic
quantum field theories.

In order to end the introduction on a more optimistic note, we
would like to point out that in the definition of a Lie
pseudoalgebra one may replace the permutation $\sigma$ by the map
$f\tt g \mapsto (g\tt f)R$ where $R$ is an R-matrix for $H$, 
hence one can take for $H$ any quasi-triangular Hopf
algebra (defined in \cite{D}).  This observation, the appearance of the
classical Yang--Baxter equation, and the fact that the
annihilation algebra of the 
associative pseudoalgebra $\Cend H$ is the Drinfeld double of $H$, lead
us to believe that there should be a deep connection between the
theories of pseudoalgebras and quantum groups.

Unless otherwise specified, all vector spaces, linear maps and
tensor products are considered over an algebraically closed field 
$\Kset$ of characteristic $0$.

\subsection*{Acknowledgements}
A major part of the present work was done in the fall of 1998 while two
of the authors were visitors at ENS and in the spring of 1999 while they 
were visitors at IHES. In the spring of 2000 one of the authors
was visiting MIT.
The work has been completed in July 2000 
while two of the authors were visiting ESI. We are grateful to these 
institutions for their hospitality.
One of the authors wishes to thank S.~P.~Novikov for valuable discussions
and references, and E.~B.~Vinberg for valuable correspondence.
Finally, we thank the referees for many remarks which led to improvement
of the exposition.

\section{Preliminaries on Hopf Algebras}\lbb{sprelh}
The goal of this section is to gather some facts and notation which will 
be used throughout the paper. 
The material in Sections \ref{subnotid} and \ref{subfiltop}
is standard and can be found, for example, in Sweedler's book \cite{Sw}.
The material in Section~\ref{subfour} seems new.

\subsection{Notation and basic identities}\lbb{subnotid}
Let $H$ be a 
Hopf algebra with a coproduct $\De$, a counit $\ep$, and 
an antipode $S$. 
%
%
%
We will use the following notation (cf.\ \cite{Sw}):
\begin{align}
\lbb{de1}
\De(h) &= h_{(1)} \tt h_{(2)},
\\
\lbb{de2}
(\De\tt\id)\De(h) &= (\id\tt\De)\De(h) = h_{(1)} \tt h_{(2)} \tt h_{(3)},
\\
\lbb{de3}
(S\tt\id)\De(h) &= h_{(-1)} \tt h_{(2)},
\quad\text{etc.}
\end{align}
Note that notation \eqref{de2} uses the coassociativity of $\De$.
The axioms of the antipode and the counit can be written
as follows:
\begin{align}
\lbb{antip}
h_{(-1)} h_{(2)} &= h_{(1)} h_{(-2)} = \ep(h),
\\
\lbb{cou}
\ep(h_{(1)}) h_{(2)} &= h_{(1)} \ep(h_{(2)}) = h,
\end{align}
while the fact that $\De$ is a homomorphism of algebras
translates as:
\begin{equation}
\lbb{deprod}
(fg)_{(1)} \tt (fg)_{(2)} = f_{(1)} g_{(1)} \tt f_{(2)} g_{(2)}.
\end{equation}
Equations \eqref{antip} and \eqref{cou} imply the following
useful relations:
\begin{equation}
\lbb{cou2}
h_{(-1)} h_{(2)} \tt h_{(3)} = 1\tt h
= h_{(1)} h_{(-2)} \tt h_{(3)}.
\end{equation}

Let $\gplk(H)$ be the subset of group-like elements of $H$, i.e.,
$g\in H$ such that $\De(g)=g\tt g$. Then $\gplk(H)$ is a group,
because $S(g)g=gS(g)=1$ for $g\in \gplk(H)$. 
Let $\primt(H)$ be the subspace of primitive elements of $H$, i.e.,
$p\in H$ such that $\De(p)=p\tt 1 + 1\tt p$. This is a Lie algebra
with respect to the commutator $[p,q]=pq-qp$.
Note that $\gplk(H)$ acts on $\primt(H)$ by inner automorphisms:
$gpg^{-1}\in\primt(H)$ for $p\in\primt(H)$, $g\in\gplk(H)$.

The proof of the following theorem may be found in \cite{Sw}.

\begin{theorem}[Kostant]\lbb{tkostant}
Let $H$ be a cocommutative Hopf algebra over $\Kset$ $($an algebraically
closed field of characteristic $0)$. Then $H$ is isomorphic $($as a 
Hopf algebra$)$ to the smash product of the universal enveloping algebra 
$\ue(\primt(H))$ and the group algebra $\Kset[\gplk(H)]$.
\end{theorem}

An \as\ algebra $A$ is called an {\em $H$-\difalg\/}
if it is also a left $H$-module such that the multiplication
$A\tt A\to A$ is a homomorphism of $H$-modules. In other words,
\begin{equation}\lbb{hxy}
h(xy) = (h_{(1)}x) (h_{(2)}y)
\end{equation}
for $h\in H$, $x,y\in A$. 
The {\em smash product\/} $A\smash H$ of an 
$H$-\difalg\
$A$ with $H$ is the tensor product $A\tt H$ of vector spaces
but with a new multiplication:
\begin{equation}\lbb{smash}
(a\smash g)(b\smash h) = a(g_{(1)}b) \smash g_{(2)}h.
\end{equation}
If both $A$ and $H$ are Hopf algebras, then $A\smash H$ is a Hopf algebra
if we consider it as a tensor product of coalgebras.
In the theorem above, $\ue(\primt(H))$ is a $\Kset[\gplk(H)]$-\difalg\
with respect to the adjoint action of $\gplk(H)$ on $\primt(H)$.

It is worth mentioning that as a byproduct of Kostant's Theorem, 
we obtain that the antipode of a cocommutative Hopf algebra
is an involution, i.e., $S^2 = \id$.

We will often be working with the Hopf algebra $H=\ue(\dd)$, where $\dd$
is a finite-dimensional Lie algebra. 
It is well known that this is a Noetherian domain, and 
any two nonzero elements $f,g \in H$ have a nonzero left 
(respectively right) common multiple. In particular, $H=\ue(\dd)$ has a 
skew-field of fractions $K$.

\begin{lemma}\lbb{rmapstofree}
Let $H$ be a Noetherian domain which has a skew-field of fractions $K$,
and let $L$ be a finite $H$-module.
Then there is a homomorphism $i\colon L\to F$ from
$L$ to a free $H$-module $F$, whose kernel is
the torsion submodule of $L$.  
If $L$ is torsion-free, then the module $F$ can be chosen in such a way that
$hF \subset i(L)$ for some nonzero $h \in H$ and $i(L)/hF$ is torsion.
\end{lemma}
\begin{proof}
The kernel of the
natural map $\iota\colon L \to L_K := K \otimes_H L$ 
is the torsion of $L$. 
The image of $L$ under this map is contained inside a free
$H$-submodule of $L_K$. In order to see this, let us
consider a set of $H$-generators $\{l_1,\dots,l_n\}$ of $L$, and a $K$-basis
$\{v_1,\dots,v_k\}$ of $L_K$. We can express the elements
$\iota(l_j)$ as $K$-linear combinations of the $v_i$'s, and by
rescaling elements of this basis by a common multiple of the denominators, 
we can assume the
$\iota(l_j)$'s to be $H$-linear combinations of the $v_i$'s. Hence the image
$\iota(L)$ is contained in the $H$-module $F$ spanned by the $v_i$'s,
which is free by construction.

The fact that $F/L$ is torsion is clear because there exist nonzero elements 
$h_i \in H$ such that $h_i v_i\in L$. If $h$ is a common multiple of the 
$h_i$'s, then $hF$ is contained in $L$. On the other hand,
the inclusion $L \subset F$ implies $hL \subset hF$, 
hence $h(L/hF)=0$ and $L/hF$ is torsion.
\end{proof}

\subsection{Filtration and topology}\lbb{subfiltop}
We define an increasing 
sequence of 
subspaces of a Hopf algebra $H$ inductively by:
\begin{align}
\lbb{fil0h}
&\fil^n H =0 \;\;\text{for $n<0$,} \quad
\fil^0 H = \Kset[\gplk(H)],
\\
\lbb{filh}
&\begin{split}
\fil^n H = \Span_\kk \big\{ h\in H \;\big|\; 
\De(h)\in \fil^0 H\tt h &+ h\tt\fil^0 H
\\
&+ \tsum_{i=1}^{n-1} \fil^i H\tt\fil^{n-i} H \big\}.
\end{split}
\end{align}
It has the following properties (which are immediate from definitions):
\begin{align}
\lbb{filh1}
(\fil^m H) (\fil^n H) &\subset \fil^{m+n} H,
\\
\lbb{filh2}
\De(\fil^n H) &\subset \tsum_{i=0}^n \fil^i H\tt\fil^{n-i} H,
\\
\lbb{filh3}
S(\fil^n H) &\subset \fil^n H.
\end{align}
When $H$ is cocommutative, using \thref{tkostant}, one can show that:
\begin{equation}\lbb{filh0}
\bigcup_n \fil^n H = H.
\end{equation}
(This condition is also satisfied when $H$ is a quantum universal enveloping 
algebra.)
Provided that \eqref{filh0} holds, we say that a nonzero element $a\in H$
has degree $n$ if $a \in \fil^n H \setminus \fil^{n-1} H$.

When $H$ is a 
universal enveloping algebra, we get its canonical filtration.
Later in some instances we will
also impose the following finiteness condition on $H$:
\begin{equation}\lbb{filh4}
\dim \fil^n H < \infty \qquad \forall n.
\end{equation}
It is satisfied when $H$ is a universal enveloping algebra
of a finite-dimensional Lie algebra, or its smash product
with the group algebra of a finite group.

Now let $X=H^* := \Hom_\Kset(H,\Kset)$ be the dual of $H$. 
Recall that $H$ acts on $X$
by the formula ($h,f\in H$, $x\in X$):
\begin{equation}\lbb{hx}
\langle hx, f\rangle = \langle x, S(h)f\rangle,
\end{equation}
so that $X$ is an \as\
$H$-\difalg\ (see \eqref{hxy}).
Moreover, $X$ is commutative when $H$ is cocommutative.
Similarly, one can define a right action of $H$ on $X$ by
\begin{equation}\lbb{xh}
\langle xh, f\rangle = \langle x, f S(h)\rangle,
\end{equation}
and then we have
\begin{equation}\lbb{xyh}
(xy)h = (x h_{(1)}) (y h_{(2)}).
\end{equation}
Associativity of $H$ implies that $X$ is an $H$-bimodule,
i.e.
\begin{equation}\lbb{fxg}
f(xg) = (fx)g, \qquad f,g\in H, \; x\in X.
\end{equation}

Let $X=\fil_{-1}X\supset\fil_0X\supset\dotsm$
be the decreasing sequence of subspaces
of $X$ dual to $\fil^n H$:
$\fil_n X = (\fil^n H)^\perp$.
It has the following properties:
\begin{align}
\lbb{fils1}
(\fil_m X) (\fil_n X) &\subset \fil_{m+n} X,
\\
\lbb{fils5}
(\fil^m H) (\fil_n X) &\subset \fil_{n-m} X,
\end{align}
and
\begin{equation}\lbb{fils0}
\bigcap_n \fil_n X = 0, 
\quad\text{provided that \eqref{filh0} holds}.
\end{equation}
We define a topology of $X$ by considering $\{\fil_n X\}$
as a fundamental system of neighborhoods of $0$.
We will always consider $X$ with this topology, while $H$ with
the discrete topology. 
It follows from \eqref{fils0} that $X$ is Hausdorff,
provided that \eqref{filh0} holds.
By \eqref{fils1} and \eqref{fils5}, 
the multiplication of $X$ and the action of $H$ on
it are continuous; in other words, $X$ is a topological $H$-\difalg.

%

We define an antipode $S\colon X\to X$ as the dual of
that of $H$:
\begin{equation}\lbb{sxhxsh}
\langle S(x), h \rangle = \langle x, S(h) \rangle .
\end{equation}
Then we have:
\begin{equation}\lbb{sabsbsa}
S(ab) = S(b) S(a) \qquad\text{for}\quad a,b \in X \;\text{or}\; H.
\end{equation}

We will also define a
comultiplication $\De\colon X\to X\widehat{\tt}X$
as the dual of the multiplication $H\tt H\to H$,
where $X\widehat{\tt}X := (H\tt H)^*$ is the completed
tensor product. Formally, we will use the same notation for $X$ as for $H$ 
(see \eqref{de1}--\eqref{de3}), writing for example 
$\De(x) = x_{(1)} \tt x_{(2)}$ for $x\in X$.
By definition, for $x,y \in X$, $f,g\in H$, we have:
\begin{align}
\lbb{xy,f}
\langle xy,f \rangle &= \langle x\tt y, \De(f) \rangle
= \langle x,f_{(1)} \rangle \langle y,f_{(2)} \rangle ,
\\
\lbb{x,fg}
\langle x,fg \rangle &= \langle \De(x), f\tt g \rangle
= \langle x_{(1)}, f \rangle \langle x_{(2)}, g \rangle.
\end{align}
We have:
\begin{align}
\lbb{fils3}
S(\fil_n X) &\subset \fil_n X,
\\
\lbb{fils2}
\De(\fil_n X) &\subset \tsum_{i=0}^n \fil_i X\what\tt\fil_{n-i} X.
\end{align}
If $H$ satisfies the finiteness condition \eqref{filh4}, then
the filtration of $X$ satisfies
\begin{equation}\lbb{fils4}
\dim X / \fil_n X < \infty \qquad \forall n,
\end{equation}
which implies that $X$ is linearly compact (see \seref{scartan} below).


By a basis of $X$ we will always mean a topological basis $\{x_i\}$
which tends to $0$, i.e., such that
for any $n$ all but a finite number of $x_i$ belong to $\fil_n X$.
Let $\{h_i\}$ be a basis of $H$ (as a vector space) compatible with the
filtration. Then the set of elements $\{x_i\}$ of $X$ defined by
$\langle x_i, h_j \rangle = \de_{ij}$ is called the dual basis of $X$.
If $H$ satisfies \eqref{filh4},
then $\{x_i\}$ is a basis of $X$ in the above sense, i.e.,
it tends to $0$. We have for $g\in H$, $y\in X$:
\begin{displaymath}
g = \tsum_i \, \langle g, x_i \rangle h_i,
\quad
y = \tsum_i \, \langle y, h_i \rangle x_i,
\end{displaymath}
where the first sum is finite, and the second one is convergent in $X$.

\begin{example}\lbb{eued}
Let $H=\ue(\dd)$ be the universal enveloping algebra of 
an $N$-dimension\-al Lie algebra $\dd$. Fix a basis $\{\d_i\}$
of $\dd$, and for $I = (i_1,\dots,i_N) \in\Zset_+^N$ let
$\d^{(I)} = \d_1^{i_1} \dotsm \d_N^{i_N} / i_1! \dotsm i_N!$.
Then $\{\d^{(I)}\}$ is a basis of $H$ (the Poincar\'e--Birkhoff--Witt basis).
Moreover, it is easy to see that
\begin{equation}\lbb{dedn}
\De(\d^{(I)}) = \sum_{ J+K=I } 
\d^{(J)} \tt \d^{(K)}.
\end{equation}
If $\{t_I\}$ is the dual basis of $X$, defined by
$\langle t_I, \d^{(J)} \rangle = \de_{I,J}$,
then \eqref{dedn} implies $t_J t_K = t_{ J+K }$.
Therefore, $X$ can be identified with the ring $\O_N = \Kset[[t_1,\dots,t_N]]$
of formal power series in $N$ indeterminates. Then the action of $H$ on
$\O_N$ is given by differential operators.
\end{example}
\begin{lemma}\lbb{lxshxy}
If  $\{h_i\}$, $\{x_i\}$ are dual bases in $H$ and $X$, then
\begin{equation}\lbb{xshxy}
\De(x) = \tsum_i \, xS(h_i) \tt x_i  
       = \tsum_i \, x_i \tt S(h_i) x 
\end{equation}
for any $x\in X$.
\end{lemma}
\begin{proof}
 For $f,g\in H$, we have:
\begin{align*}
\langle \tsum_i\, xS(h_i) \tt x_i, f\tt g \rangle
&= \tsum_i\, \langle xS(h_i), f \rangle\langle x_i,g \rangle
\\
= \langle xS(g), f \rangle
&= \langle x, fg \rangle
= \langle \De(x), f\tt g \rangle,
\end{align*}
which proves the first identity. The second one is proved in the same way.
\end{proof}

\subsection{Fourier transform}\lbb{subfour}
For an arbitrary Hopf algebra $H$,
we introduce a map
$\F\colon H\tt H \to H\tt H$, called the {\em Fourier transform},
by the formula
\begin{equation}
\lbb{ftrans}
\F(f\tt g) = (f\tt 1) \, (S\tt\id)\De(g) = f g_{(-1)} \tt g_{(2)}.
\end{equation}
It follows from \eqref{cou2} that $\F$ is a vector space isomorphism with
an inverse given by
\begin{equation}
\lbb{finv}
\F^{-1}(f\tt g) = (f\tt 1) \,\De(g) = f g_{(1)} \tt g_{(2)}.
\end{equation}
Indeed, using the coassociativity of $\De$ and \eqref{cou2}, we compute
\begin{displaymath}
\F^{-1}(f g_{(-1)} \tt g_{(2)})
= f g_{(-1)} (g_{(2)})_{(1)} \tt (g_{(2)})_{(2)}
= f g_{(-1)} g_{(2)} \tt g_{(3)}
= f\tt g.
\end{displaymath}
The significance of $\F$ is in the identity
\begin{equation}\lbb{sigfour}
f\tt g = \F^{-1}\F(f\tt g) = (f g_{(-1)} \tt 1) \, \De(g_{(2)}),
\end{equation}
which, together with properties \eqref{filh1}--\eqref{filh3}
of the filtration of $H$, implies the next result.

\begin{lemma}\lbb{lhhh}
{\rm(i)}
Every element of $H\tt H$ can be uniquely represented in the form
$\sum_i \, (h_i \tt 1)\De(l_i)$,
where $\{h_i\}$ is a fixed $\Kset$-basis of $H$ and $l_i\in H$.
In other words, $H\tt H = (H\tt\Kset)\De(H)$. 

{\rm(ii)}
We have{\rm:}
\begin{equation}\lbb{lhhh1}
(\fil^n H\tt\Kset) \De(H) = \fil^n (H\tt H) \De(H) = 
(\Kset\tt\fil^n H) \De(H),
\end{equation}
where $\fil^n (H\tt H) = \tsum_{i+j=n} \fil^i H\tt\fil^j H$.

In particular, for any $H$-module $W$, we have{\rm:}
\begin{equation}\lbb{hnkw}
(\fil^n H\tt\Kset) \tt_H W = \fil^n (H\tt H) \tt_H W = 
(\Kset\tt\fil^n H) \tt_H W.
\end{equation}
\end{lemma}
\begin{proof}
For $h\in H\tt H$ we have:
\begin{displaymath}
h = \tsum_i \, (h_i \tt 1)\De(l_i) = \F^{-1}(\tsum_i \, h_i \tt l_i)
\quad\;\;\text{iff}\quad\;\;
\tsum_i \, h_i \tt l_i = \F(h). 
\end{displaymath}
This proves (i).

To prove \eqref{lhhh1}, it is enough to show that
$\fil^n (H\tt H) \subset (\fil^n H\tt\Kset) \De(H)$.
This follows from the above equation and the fact that
$\F( \fil^n (H\tt H) ) \subset \fil^n (H\tt H) \subset \fil^n H\tt H$.
\end{proof}

The Fourier transform $\F$ has the following properties
(which are easy to check using \eqref{antip}--\eqref{deprod}):
\begin{align}
\lbb{propf0}
\F\bigl((f\tt g)\De(h)\bigr) &= \F(f\tt g) \, (1\tt h),
\\
\lbb{propf1}
\F(hf\tt g) &= (h\tt1) \, \F(f\tt g),
\\
\lbb{propf2}
\F(f\tt hg) &= (1\tt h_{(2)}) \, \F(f\tt g) \, (h_{(-1)}\tt1),
\\
\lbb{propf3}
\F_{12}\F_{13}\F_{23} &= \F_{23}\F_{12}.
\end{align}
Here in \eqref{propf3}, we use the standard notation
$\F_{12}=\F\tt\id$ acting on $H\tt H\tt H$.

\section{Pseudotensor Categories and Pseudoalgebras}\lbb{slie*}

In this section, we review some definitions 
of Beilinson and Drinfeld \cite{BD};
we also use the exposition in \cite[Section~12]{BKV}.

The theory of conformal algebras \cite{K2}
is in many ways analogous to
the theory of Lie algebras. The reason is that in fact conformal algebras
can be considered as Lie algebras in a certain ``pseudotensor'' category,
instead of the category of vector spaces. 
A pseudotensor category \cite{BD} is a category equipped with 
``polylinear maps'' and a way to compose them (such categories were
first introduced by Lambek \cite{L} under the name multicategories).
This is enough
to define the notions of Lie algebra, representations, cohomology,
etc.

As an example, consider first the category $\Vec$
of vector spaces (over $\Kset$).
For a finite nonempty set $I$ and
a collection of vector spaces $\{L_i\}_{i\in I}$, $M$,
we can define the space of {\em polylinear maps\/} 
{}from $\{L_i\}_{i\in I}$ to $M$ as
\begin{displaymath}\lbb{plin}
\Lin(\{L_i\}_{i\in I}, M) 
= \Hom(\otimes_{i\in I} L_i, M).
\end{displaymath}
The symmetric group $\symm_I$ acts among
these spaces by permuting the factors in
$\otimes_{i\in I} L_i$.

For any surjection of finite sets $\pi\colon J\surjto I$
and a collection $\{N_j\}_{j\in J}$,
we have the obvious compositions of polylinear maps
\begin{align}\lbb{com1}
&\Lin(\{L_i\}_{i\in I}, M) 
\otimes \bigotimes_{i\in I}
\Lin(\{N_j\}_{j\in J_i}, L_i) 
\to \Lin(\{N_j\}_{j\in J},M),    \\
\lbb{com2}
& \phi \times \{\psi_i\}_{i\in I} \mapsto 
\phi\circ(\otimes_{i\in I} \psi_i) \equiv \phi(\{\psi_i\}_{i\in I}),
\end{align}
where $J_i = \pi^{-1}(i)$ for $i\in I$.

The compositions have the following properties:
\begin{description}
\item[Associativity] 
\begin{sloppypar}
If $K\surjto J$, $\{P_k\}_{k\in K}$ is a family of objects and
$\chi_j\in \Lin(\{P_k\}_{k\in K_j}, N_j)$, then
$\phi\bigl(\bigl\{\psi_i(\{\chi_j\}_{j\in J_i})\bigr\}_{i\in I}\bigr)
\linebreak[1] = \linebreak[0]
\bigl(\phi(\{\psi_i\}_{i\in I})\bigr)(\{\chi_j\}_{j\in J}) 
\linebreak[0] \in \Lin(\{P_k\}_{k\in K}, M)$.
\end{sloppypar}

\item[Unit] 
For any object $M$ there is an element $\id_M\in \Lin(\{M\},M)$ such
that for any $\phi\in \Lin(\{L_i\}_{i\in I}, M)$ one has
$\id_M(\phi)=\phi(\{\id_{L_i}\}_{i\in I})=\phi$.

\item[Equivariance] 
The compositions \eqref{com1}
are equivariant with respect to the natural action of the symmetric group.
\end{description}
\begin{definition}[\cite{BD}]\lbb{dptc} 
A {\em pseudotensor category\/} is a class of objects $\M$ 
together with
vector spaces $\Lin(\{L_i\}_{i\in I}, M)$,
equipped with actions of the symmetric groups $\symm_I$ among them
and composition maps \eqref{com1},
satisfying the above three properties.
\end{definition}
\begin{remark}\lbb{rpt}
For a pseudotensor category $\M$ and objects $L,M\in\M$, let
$\Hom(L,M)=\Lin(\{L\},M)$. This gives 
a structure of an ordinary (additive) category on $\M$ and all
$\Lin$ are functors $(\M^\circ)^I\times\M\to\Vec$,
where $\M^\circ$ is the dual category of $\M$.
\end{remark}
\begin{remark}\lbb{roper}
The notion of pseudotensor category is a straightforward generalization
of the notion of operad. By definition, an {\em operad\/} is a
pseudotensor category with only one object.
\end{remark}

\begin{definition}\lbb{dlie*} 
A {\em Lie algebra in a pseudotensor category $\M$\/} 
is an object $L$ equipped with
$\be\in \Lin(\{L,L\},L)$ satisfying the following properties.

\begin{description}
\item[Skew-commutativity]
$\be = -\sigma_{12} \, \be,$
where $\sigma_{12} = (12) \in\symm_2$.

\item[Jacobi identity]
$\be(\be(\cdot, \cdot), \cdot)=\be(\cdot, \be(\cdot, \cdot))-
\sigma_{12} \, \be(\cdot, \be(\cdot, \cdot)),$
where now $\sigma_{12} = (12)$ is viewed as an element of $\symm_3$.
\end{description}
\end{definition}

It is instructive to think of 
a polylinear map $\phi\in \Lin(\{L_i\}_{i=1}^n, M)$
as an operation with $n$ inputs and $1$ output,
as depicted in \firef{Fpolyn}.
The skew-commutativity and Jacobi identity for a Lie algebra
$(L,\be)$ are represented pictorially in Figures \ref{Fsksymm} and
\ref{Fjacid}.

\begin{figure}[h]
\begin{displaymath}
\fig{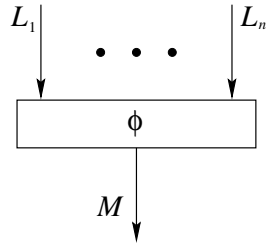}
\end{displaymath}
\caption{A polylinear map from $\{L_i\}_{i=1}^n$ to $M$.}\lbb{Fpolyn}
\end{figure}

\begin{figure}[h]
\begin{displaymath}
\fig{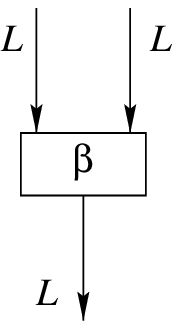}
\;\; = \; - \;
\fig{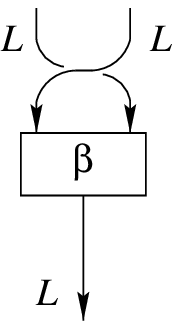}
\end{displaymath}
\caption{Skew-commutativity.}\lbb{Fsksymm}
\end{figure}

\begin{figure}[h]
\begin{displaymath}
\fig{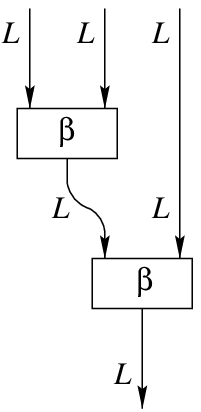}
\quad = \quad
\fig{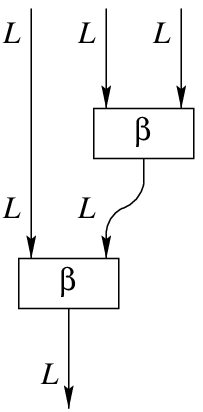}
\quad - \quad
\fig{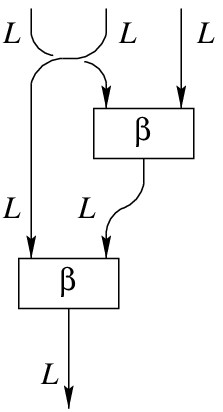}
\end{displaymath}
\caption{Jacobi identity.}\lbb{Fjacid}
\end{figure}

\begin{definition}\lbb{drep}
A {\em representation} of a  Lie algebra $(L,\be)$ is an object $M$
together with $\rho\in \Lin(\{L,M\},M)$ satisfying
\begin{displaymath}\lbb{rep*}
\rho(\be(\cdot, \cdot), \cdot)=\rho(\cdot, \rho(\cdot, \cdot))-
\sigma_{12} \, \rho(\cdot, \rho(\cdot, \cdot)).
\end{displaymath}
\end{definition}
Similarly, one can define cohomology of a  Lie algebra $(L,\be)$
with coefficients in a module $(M,\rho)$ (cf.\ \cite{BKV}).

\begin{definition}\lbb{dcoh}
An {\em $n$-cochain\/} of a Lie algebra $(L,\be)$, with coefficients in
a module $(M,\rho)$ over it, is a polylinear operation 
$\ga\in\Lin(\{\underbrace{L,\dots, L}_n\},M)$
which is skew-symmetric, i.e., satisfying for all $i=1,\dots,n-1$
the identity shown in \firef{Fskscoch}.
The differential $d\ga$ of a cochain $\ga$
is defined by \firef{Fdifcoch}.
The same computation as in the ordinary Lie algebra case shows that
$d^2=0$.  The cohomology of the resulting complex is called the 
{\em cohomology of $L$ with coefficients in $M$\/} and is denoted
by $\coh^\bullet(L,M)$.
\end{definition}

\begin{figure}[h] 
\begin{displaymath}
\fig{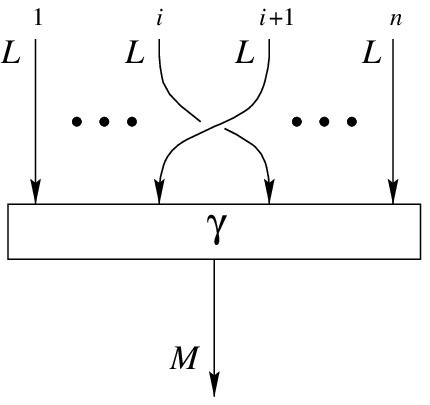}
\quad = - \quad
\fig{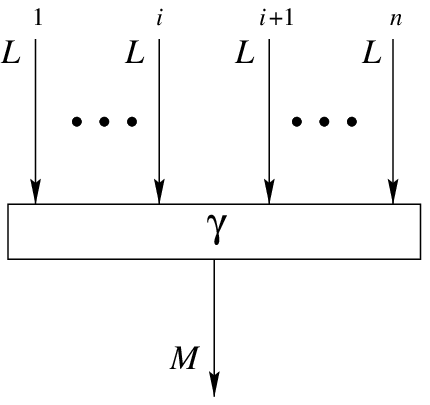}
\end{displaymath}
\caption{Skew-symmetry of a cochain.}\lbb{Fskscoch}
\end{figure}


\begin{figure}[h] 
\begin{align*}
&\fig{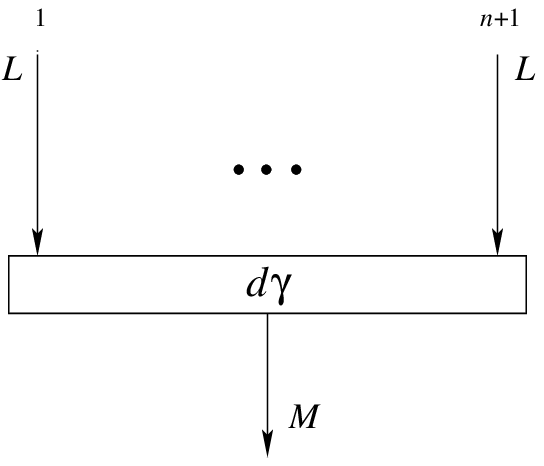}
\\
&\quad
= \sum_{1\le i\le n+1} (-1)^{i+1}
\quad
\fig{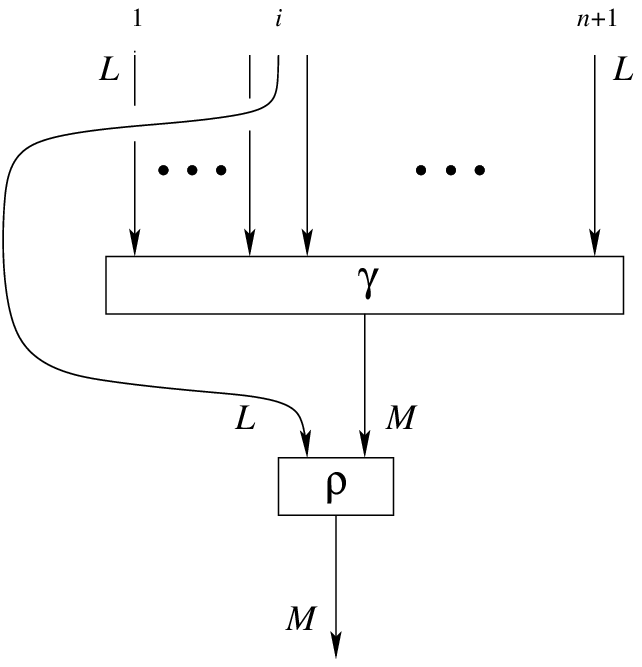}
\\
&\quad
+ \sum_{1\le i<j \le n+1} (-1)^{i+j}
\quad
\fig{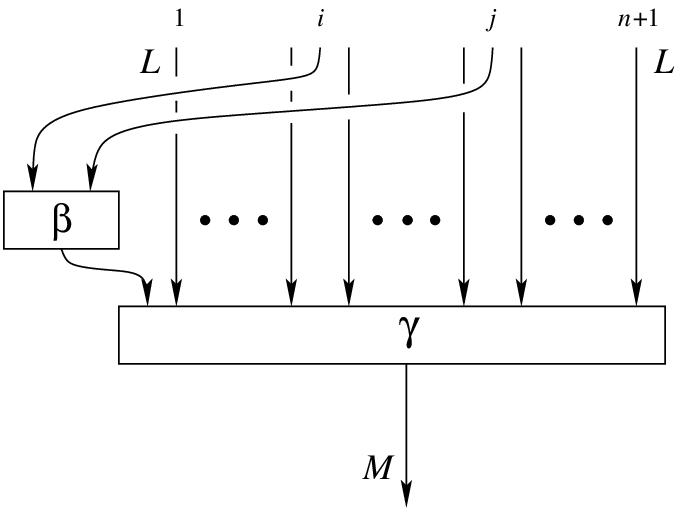}
\end{align*}
\caption{Differential of a cochain.}\lbb{Fdifcoch}
\end{figure}


\begin{example}\lbb{exuslie}
A Lie algebra in the category of vector spaces $\Vec$ is just
an ordinary Lie algebra. The same is true for representations
and cohomology.
\end{example}

\begin{example}\lbb{exmld}
Let $H$ be a cocommutative bialgebra with a comultiplication $\De$
and a counit $\ep$. Then the category $\M^l(H)$ of left $H$-modules
is a symmetric tensor category. Hence, $\M^l(H)$ is a 
pseudotensor category with polylinear maps
\begin{equation}\lbb{mld}
\Lin(\{L_i\}_{i\in I}, M) 
= \Hom_H(\otimes_{i\in I} L_i, M).
\end{equation}
The composition of polylinear maps is given by \eqref{com2}.
An algebra (e.g., Lie or associative)
in the category $\M^l(H)$ will be called an {\em $H$-\difalg}:
this is an ordinary algebra which is also a left $H$-module and such that
the product (or the bracket) is a homomorphism of $H$-modules,
see \eqref{hxy}.
\end{example}

One can also define the notions of associative algebra or 
commutative algebra in a pseudotensor category, 
their representations and 
analogues of the Hochschild, cyclic,  or Harrison cohomology.

\begin{definition}\lbb{dass*} 
An {\em associative algebra in a pseudotensor category $\M$\/} 
is an object $A$ and a product $\mu\in \Lin(\{A,A\},A)$ satisfying
\begin{description}
\item[Associativity]
$\mu(\mu(\cdot, \cdot), \cdot)=\mu(\cdot, \mu(\cdot, \cdot))$,
\end{description}
see \firef{Fass}.
The algebra $(A,\mu)$ is called {\em commutative\/} if, 
in addition, $\mu$ satisfies
\begin{description}
\item[Commutativity]
$\mu = \sigma_{12} \, \mu,$
where $\sigma_{12} = (12) \in\symm_2$.
\end{description}
\end{definition}

\begin{figure}[h]
\begin{displaymath}
\fig{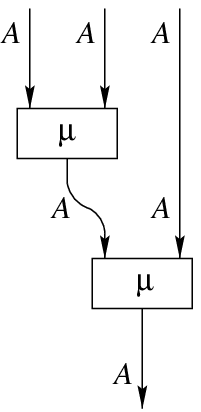}
\quad = \quad
\fig{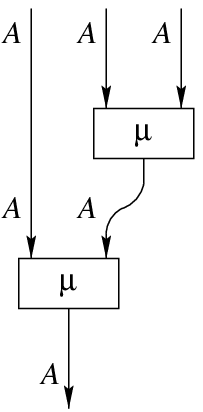}
\end{displaymath}
\caption{Associativity.}\lbb{Fass}
\end{figure}

\begin{remark}\lbb{rsymac}
In order to define the notion of an \as\ algebra in
a pseudotensor category, one does not use the actions
of the symmetric groups among the spaces of polylinear maps.
One can relax the definition of a pseudotensor category
by forgetting these actions. Then what we call a ``pseudotensor category''
should be termed a ``symmetric pseudotensor category'', while there
is a more general notion of a ``braided'' one (cf.\ \cite{So}).
\end{remark}

\begin{proposition}\lbb{assolie}
Let $(A,\mu)$ be an associative algebra in a pseudotensor category $\M$.
Define $\be\in \Lin(\{A,A\},A)$ as the commutator
$\be := \mu - \sigma_{12} \, \mu$, see \firef{Fcomm}.
Then $(A,\be)$ is a Lie algebra in $\M$.
\end{proposition}

\begin{figure}[h]
\begin{displaymath}
\fig{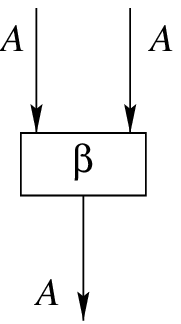}
\;\; = \;\; 
\fig{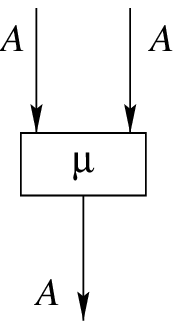}
\;\; - \;\; 
\fig{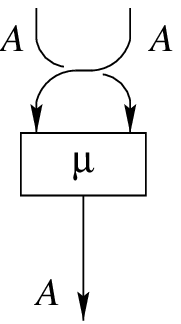}
\end{displaymath}
\caption{Commutator.}\lbb{Fcomm}
\end{figure}


Let $H$ be a cocommutative bialgebra with a comultiplication $\De$.
We introduce a pseudotensor category $\M^*(H)$ with the same objects as 
$\M^l(H)$ (i.e., left $H$-modules)
but with another pseudotensor structure \cite{BD}:
\begin{equation}\lbb{m*d}
\Lin(\{L_i\}_{i\in I}, M) 
= \Hom_{H^{\tt I}} (\boxtimes_{i\in I} L_i, H^{\tt I}\tt_H M).
\end{equation}
Here $\boxtimes_{i\in I}$ is the tensor product functor
$\M^l(H)^I \to \M^l(H^{\tt I})$.
For a surjection $\pi\colon J\surjto I$, the composition of polylinear maps
is defined as follows:
\begin{equation}\lbb{*com}
\phi\bigl(\{\psi_i\}_{i\in I}\bigr) = \De^{(\pi)}\bigl(\phi\bigr) \,
\circ\bigl(\boxtimes_{i\in I}\psi_i\bigr).
\end{equation}
Here $\De^{(\pi)}$ is the functor 
$\M^l(H^{\tt I})\to\M^l(H^{\tt J})$,
$M\mapsto H^{\tt J} \tt_{H^{\tt I}} M$,
where $H^{\tt I}$ acts on $H^{\tt J}$ via the iterated
comultiplication determined by $\pi$.

Explicitly, let $n_j \in N_j$ $(j\in J)$, and write
\begin{equation}\lbb{compoly1}
\psi_i\bigl( \tt_{j\in J_i} n_j \bigr) = \sum_r g_i^r \tt_H l_i^r,
\qquad g_i^r \in H^{\tt J_i}, \; l_i^r \in L_i,
\end{equation}
where, as before, $J_i = \pi^{-1}(i)$ for $i\in I$. Let 
\begin{equation}\lbb{compoly2}
\phi\bigl( \tt_{i\in I} l_i^r \bigr) = \sum_s f^{rs} \tt_H m^{rs},
\qquad f^{rs} \in H^{\tt I}, \; m^{rs} \in M.
\end{equation}
Then, by definition,
\begin{equation}\lbb{compoly3}
\bigl(\phi\bigl(\{\psi_i\}_{i\in I}\bigr)\bigr) \bigl( \tt_{j\in J} n_j \bigr) 
= \sum_{r,s} ( \tt_{i\in I} g_i^r ) \De^{(\pi)}(f^{rs}) \tt_H m^{rs},
\end{equation}
where $\De^{(\pi)} \colon H^{\tt I} \to H^{\tt J}$ is the iterated
comultiplication determined by $\pi$. For example, if 
$\pi\colon\{1,2,3\}\to\{1,2\}$ is given by
$\pi(1)=\pi(2)=1$, $\pi(3)=2$, then $\De^{(\pi)} = \De\tt\id$; if
$\pi(1)=1$, $\pi(2)=\pi(3)=2$, then $\De^{(\pi)} = \id\tt\De$.

The symmetric group $\symm_I$ acts among
the spaces $\Lin(\{L_i\}_{i\in I}, M)$
by simultaneously permuting the factors in $\boxtimes_{i\in I} L_i$
and $H^{\tt I}$. This is the only place where we need the cocommutativity of
$H$; for example, the permutation $\sigma_{12} = (12) \in\symm_2$
acts on $(H\tt H)\tt_H M$ by 
\begin{displaymath}
\si_{12}\bigl( (f\tt g)\tt_H m \bigr) = (g\tt f)\tt_H m,
\end{displaymath}
and this is well defined only when $H$ is cocommutative.

One can generalize the above construction for (quasi)triangular
bialgebras as follows.

\begin{remark}\lbb{rquastr}
Let $H$ be a triangular bialgebra with a universal R-matrix $R$.
Recall that $R$ is an invertible element of $H\tt H$ 
satisfying the following equations:
\begin{align}
\lbb{rma0}
\si(R) &= R^{-1},
\\
\lbb{rma1}
\si(\De(h)) R &= R \De(h) \qquad \forall h\in H,
\\
\lbb{rma2}
(\id\tt\De)R &= R_{13} R_{12},
\\
\lbb{rma3}
(\De\tt\id)R &= R_{13} R_{23},
\end{align}
where $\si$ is the permutation $\si(f\tt g)=g\tt f$, and
we use the standard notation $R_{12}=R\tt\id \in H\tt H\tt H$, etc.
Then we define a pseudotensor category $\M^*(H)$ as above but with 
a modified action of the symmetric groups. It is enough to describe
the action of the transposition $\sigma_{12} = (12) \in\symm_2$
on $(H\tt H)\tt_H M$; it is given by
\begin{displaymath}
\si_{12}\bigl( (f\tt g)\tt_H m \bigr) = (g\tt f)R\tt_H m.
\end{displaymath}
This is well defined because of \eqref{rma1}, and
$\si_{12}^2=\id$ because of \eqref{rma0}.
Since any permutation is a product of transpositions, this can be
extended to an action of the symmetric group among the spaces of
polylinear maps; due to \eqref{rma2}, \eqref{rma3},
this action is compatible with compositions.

If $H$ is quasitriangular, i.e., if we drop relation \eqref{rma0},
we will get an action of the braid group instead of the symmetric one
and a ``braided'' pseudotensor category (cf.\ \reref{rsymac}).
\end{remark}

The following notion will be the main object of our study.

\begin{definition}\lbb{dflie*}
A {\em{Lie $H$-\psalg\/}} (or just a Lie \psalg)
is a Lie algebra $(L,\be)$
in the pseudotensor category $\M^*(H)$ as defined above.
\end{definition}

Examples of Lie \psalgs\ will be given in 
Sections \ref{smainex} and \ref{svect} below.
One can also define {\em{\as\ $H$-\psalgs\/}} as
associative algebras $(A,\mu)$ in the pseudotensor category $\M^*(H)$.
It is convenient to define the general notion of an algebra in $\M^*(H)$
as follows.

\begin{definition}\lbb{dpseu*}
An {\em{$H$-\psalg\/}} (or just a \psalg)
is a left $H$-module $A$ together with an operation
$\mu\in \Hom_{H\tt H}(A\tt A, (H\tt H)\tt_H A)$,
called the {\em{pseudoproduct}}.
\end{definition}

We will denote the pseudoproduct $\mu(a\tt b) \in (H\tt H)\tt_H A$
of two elements $a,b\in A$ by $a*b$. 
It has the following defining property:

\begin{description}
\item[$H$-bilinearity] For $a,b\in A$, $f,g\in H$, one has
\begin{equation}\lbb{bil*2}
fa*gb = ((f\tt g)\tt_H 1) \, (a*b).
\end{equation}
Explicitly, if
\begin{equation}\lbb{ab*'}
a*b = \tsum_i\, (f_i\tt g_i)\tt_H e_i,
\end{equation}
then $fa*gb = \tsum_i\, (f f_i\tt g g_i)\tt_H e_i$.
\end{description}

To describe explicitly the associativity condition 
for a pseudoproduct $\mu$,
we need to compute the compositions $\mu(\mu(\cdot,\cdot),\cdot)$
and $\mu(\cdot,\mu(\cdot,\cdot))$ in $\M^*(H)$.
Let $a*b$ be given by \eqref{ab*'}, and let
\begin{equation}\lbb{abc*1'}
e_i*c = \tsum_{i,j}\, (f_{ij}\tt g_{ij})\tt_H e_{ij}.
\end{equation}
Then $(a*b)*c \equiv \mu(\mu(a\tt b) \tt c)$
is the following element of $H^{\tt3}\tt_H A$ (cf.\ \eqref{compoly3}):
\begin{equation}\lbb{abc*3'}
(a*b)*c = \tsum_{i,j}\, (f_i {f_{ij}}_{(1)}\tt g_i {f_{ij}}_{(2)}\tt g_{ij}) 
\tt_H e_{ij}.
\end{equation}
Similarly, if we write
\begin{align}
\lbb{abc*4'}
b*c &= \tsum_i\, (h_i\tt l_i)\tt_H d_i,
\\
\lbb{abc*5'}
a*d_i &= \tsum_{i,j}\, (h_{ij}\tt l_{ij})\tt_H d_{ij},
\intertext{then}
\lbb{abc*6'}
a*(b*c) &= \tsum_{i,j}\, (h_{ij}\tt h_i {l_{ij}}_{(1)}\tt l_i {l_{ij}}_{(2)})
\tt_H d_{ij}.
\end{align}

Now a pseudoproduct $a*b$ is associative iff it satisfies
\begin{description}
\item[Associativity]
\begin{equation}\lbb{assoc*}
a*(b*c) = (a*b)*c
\end{equation}
in $H^{\tt3}\tt_H A$, where the compositions $(a*b)*c$ and $a*(b*c)$
are given by the above formulas.
\end{description}
The pseudoproduct $a*b$ is commutative iff it satisfies

\begin{description}
\item[Commutativity]
\begin{equation}\lbb{comm*}
b*a = (\si\tt_H\id) \, (a*b),
\end{equation}
where $\si\colon H\tt H\to H\tt H$ is the permutation
$\si(f\tt g)=g\tt f$. Explicitly,
\begin{equation}\lbb{ba*'}
b*a = \tsum_i\, (g_i\tt f_i)\tt_H e_i, 
\end{equation}
if $a*b$ is given by \eqref{ab*'}. Note that the right-hand side
of \eqref{comm*} is well defined due to the cocommutativity of $H$.
\end{description}

In the case of a Lie \psalg\ $(L,\be)$, we will call the pseudoproduct
$\be$ a {\em pseudobracket\/}, and we will denote it by $[a*b]$.
Let us spell out its properties
($a,b,c\in L$, $f,g\in H$):

\begin{description}

\item[$H$-bilinearity]
\begin{equation}\lbb{bil*}
[fa*gb] = ((f\tt g)\tt_H 1) \, [a*b].
\end{equation}

\item[Skew-commutativity]
\begin{equation}\lbb{ssym*}
[b*a] = - (\si\tt_H\id) \, [a*b].
\end{equation}

\item[Jacobi identity]
\begin{equation}\lbb{jac*}
[a*[b*c]] - ((\si\tt\id)\tt_H\id) \, [b*[a*c]]
= [[a*b]*c]
\end{equation}
in $H^{\tt3}\tt_H L$, where the compositions $[[a*b]*c]$ and $[a*[b*c]]$
are defined as above.

\end{description}

\begin{proposition}\lbb{assolie*}
Let $(A,\mu)$ be an \as\ $H$-\psalg.
Define a pseudobracket $\be$ as the commutator
$[a*b] = a*b - (\si\tt_H\id) \, (b*a)$.
Then $(A,\be)$ is a Lie $H$-\psalg\
{\rm(}cf.\ \prref{assolie}{\rm)}.
\end{proposition}

The definitions of representations of Lie \psalgs\ or \as\ \psalgs\
are obvious modifications of the above.

\begin{definition}\lbb{drepas*}
A {\em representation\/} of an  \as\ $H$-\psalg\ $A$
is a left $H$-module $M$ together with an operation
$\rho\in \Lin(\{A,M\},M)$, written as 
$a*c \equiv \rho(a\tt c) \in (H\tt H)\tt_H M$,
which satisfies
\eqref{assoc*} for $a,b\in A$, $c\in M$.
\end{definition}

\begin{definition}\lbb{dreplie*}
A {\em representation\/} of a Lie $H$-\psalg\ $L$
is a left $H$-module $M$ together with an operation
$\rho\in \Lin(\{L,M\},M)$, written as $a*c \equiv \rho(a\tt c)$, 
which satisfies
\begin{equation}\lbb{replie*}
a*(b*c) - ((\si\tt\id)\tt_H\id) \, (b*(a*c))
= [a*b]*c
\end{equation}
for $a,b\in L$, $c\in M$.
\end{definition}

\section{Some Examples of Lie Pseudoalgebras}\lbb{smainex}
In this section we give some examples of Lie \psalgs, 
and discuss their relationship with previously known objects.
Other important examples --- the \psalgs\ of vector fields --- are 
treated in detail in \seref{svect}.

\subsection{Conformal algebras}\lbb{subconfalg}
The (Lie) conformal algebras introduced by Kac \cite{K2} are exactly the 
(Lie) $\Kset[\d]$-\psalgs, where $\Kset[\d]$ is the Hopf
algebra of polynomials in one variable $\d$. The explicit relation
between the $\la$-bracket of \cite{DK} 
and the pseudobracket of \seref{slie*} is:
\begin{displaymath}\lbb{alab}
[a_\la b] = \tsum_i\, p_i(\la)c_i
\;\;\iff\;\;
[a*b] = \tsum_i\, (p_i(-\d)\tt1)\tt_{\Kset[\d]} c_i.
\end{displaymath}
This correspondence has been explained in detail in the introduction.

Similarly, for $H=\Kset[\d_1,\dots,\d_N]$ we get conformal algebras
in $N$ indeterminates, see \cite[Section~10]{BKV}.
We may say that for $N=0$, $H$ is $\Kset$; 
then a $\Kset$-conformal algebra is the 
same as a Lie algebra, cf.\ \exref{exuslie}.

On the other hand, when $H=\Kset[\Ga]$ is the group algebra of a 
group $\Ga$, one obtains the $\Ga$-conformal algebras
studied in \cite{GK}. This is a special case of a more general
construction described in \seref{sgaalg} below.

\subsection{Current \psalgs}\lbb{subcuralg}
Let $H'$ be a Hopf subalgebra of $H$, 
and let $A$ be an $H'$-\psalg.
Then we define the {\em current\/} $H$-\psalg\
$\Cur_{H'}^{H} A \equiv\Cur A$ as ${H}\tt_{H'} A$ by extending the 
pseudoproduct $a*b$ of $A$ using the $H$-bilinearity. 
Explicitly, for $a,b\in A$, we define
\begin{displaymath}\lbb{curalg*}
\begin{split}
(f\tt_{H'}a) * (g\tt_{H'}b)
&= ((f\tt g)\tt_{H}1) \, (a*b)
\\
&= \tsum_i\, (f f_i \tt g g_i)\tt_{H} (1 \tt_{H'} e_i),
\end{split}
\end{displaymath}
if
$a*b = \tsum_i\, (f_i \tt g_i)\tt_{H'} e_i$.
%
%
Then $\Cur_{H'}^{H} A$ is an $H$-\psalg\ which is Lie or \as\
when $A$ is so.

An important special case is when $H'=\Kset$: given a Lie algebra $\g$,
let $\Cur\g=H\tt\g$ with the following pseudobracket
\begin{displaymath}\lbb{curalg2}
[(f \tt a) * (g \tt b)] = (f\tt g)\tt_{H} (1\tt [a,b]).
\end{displaymath}
Then $\Cur\g$ is a Lie $H$-\psalg.

\subsection{$H$-\psalgs\ of rank $1$}\lbb{subrank1}
Let $L=He$ be a Lie \psalg\ which is a free $H$-module
of rank $1$. Then, by $H$-bilinearity, the pseudobracket on $L$
is determined by $[e*e]$, or equivalently, by
an $\al\in H\tt H$ such that  
$[e*e] = \al\tt_H e$.
%

\begin{proposition}\lbb{prank1}
$L=He$ with the pseudobracket $[e*e] = \al\tt_H e$
is a Lie $H$-\psalg\ iff $\al\in H\tt H$ 
satisfies the following equations{\rm:}
\begin{align}
\lbb{al1}
& \al = -\si(\al),
\\
\lbb{al2}
& (\al\tt1)\,(\De\tt\id)(\al) 
= (1\tt\al)\,(\id\tt\De)(\al) 
- (\si\tt\id)\,\bigl((1\tt\al)\,(\id\tt\De)(\al)\bigr).
\end{align}
Similarly, $A=Ha$ with a pseudoproduct $a*a = \al\tt_H a$ 
is an \as\ $H$-\psalg\ iff $\al$ satisfies
\begin{displaymath}\lbb{al3}
(\al\tt1)\,(\De\tt\id)(\al) = (1\tt\al)\,(\id\tt\De)(\al).
\end{displaymath}
\end{proposition}
\begin{proof}
Follows immediately from definitions.
Indeed, if $[e*e] = \al\tt_H e$, then:
\begin{align*}
[[e*e]*e] &= (\al\tt1)\,(\De\tt\id)(\al) \tt_H e,
\\
[e*[e*e]] &= (1\tt\al)\,(\id\tt\De)(\al) \tt_H e.
\end{align*}
\end{proof}

\begin{lemma}\lbb{xcybe1}
Let $H=\ue(\dd)$ be the universal enveloping algebra of a Lie algebra $\dd$.
Then any solution
$\al\in H\tt H$ of equations \eqref{al1}, \eqref{al2} is of the form
$\al = r + s\tt1 - 1\tt s$, where $r\in\dd\wedge\dd$, $s\in\dd$.

In this case \eqref{al2} 
is equivalent to the following system of equations{\rm:}
\begin{align}
\lbb{cybe3}
&[r,\De(s)] = 0,
\\
\lbb{cybe4}
&([r_{12}, r_{13}] + r_{12} s_3) + \text{{\rm{cyclic}}} = 0.
\end{align}
\end{lemma}
$($As usual, $r_{12}=r\tt1$, $s_3=1\tt1\tt s$, etc., and\/ ``{\rm{cyclic}}''
here and further
means applying the two nontrivial cyclic permutations on $H\tt H\tt H.)$
\begin{proof}
Using an argument similar to that of \cite{Ki},
we will show that 
if $\al$ satisfies \eqref{al2} then $\al \in H \tt (\dd + \Kset)$. Then
\eqref{al1} will imply the first claim, that $\al\in(\dd+\Kset)\tt(\dd+\Kset)$.

Let $\{\d_1,\dots, \d_N\}$ be a basis of $\dd$, and let us consider the 
corresponding Poincar\'e--Birkhoff--Witt 
basis of $H = \ue(\dd)$ given by elements
$\d^{(I)} := \d_1^{i_1} \dotsm \d_N^{i_N} / i_1! \dotsm i_N!$, where 
$I = (i_1, \dots, i_N) \in\Zset_+^N$. 
In this basis the comultiplication 
takes the simple form \eqref{dedn}. We can write $\al = \sum_I \al_I \tt
\d^{(I)}$,
$\al_I \in H$.
Equation \eqref{al2} then becomes:
\begin{equation}\lbb{albis}
\sum_I \alpha \Delta(\al_I) \tt \d^{(I)} =
\sum_{I,J,K} (\al_{J+K} \tt \al_I \d^{(J)} - \al_I \d^{(J)} \tt \al_{J+K})
\tt \d^{(I)} \d^{(K)} .
\end{equation}

Let $p$ be the maximal value of $|I| = i_1 + \dots + i_N$ for $I$ such that
$\al_I \neq 0$. We want to show that $p\leq 1$. Among all $I$ such that
$|I| = p$ there will be some (nonzero) $\al_I$ of maximal degree $d$.
Then without
loss of generality we can change the basis $\d_1, \dots, \d_N$ and assume that
the coefficient $\al_{(p, 0, \dots, 0)}$ is nonzero and
of degree $d$.
If $p>1$, then no nonzero term in the left-hand side of \eqref{albis}
has a third tensor factor 
of degree $2p$ or $2p-1$ since $2p-1 > p$.
Hence, terms from the right-hand
side of degree $2p$ (respectively $2p-1$) in the third tensor factor
must cancel against each other.

Terms having degree $2p$ in the third tensor factor cancel, since they
give the following sum:
\begin{equation}\lbb{comm}
\sum_{|I|=|K|=p} \al_K \tt \al_I \tt [\d^{(I)}, \d^{(K)}],
\end{equation}
which in the third tensor factor has degree $2p-1$ and lower.
Note also that their coefficients have total degree at most $2d$.

Terms having third tensor factors of
degree $2p-1$, besides \eqref{comm}, arise when we
choose $|I+K| = 2p-1$. Those with $|I| = p-1, |K| = p$ can be expressed in
terms of commutators as above, and hence only contribute to lower degree. So,
we only need to account for terms with $|I|=p, |K|=p-1$.

Let us focus on such terms having a third tensor factor proportional to
$\d_1^{2p-1}$, whose coefficient must be zero.
They occur in \eqref{albis} only
when $I = (p, 0, \dots, 0)$, 
$K=(p-1, 0, \dots, 0)$. When $J=0$, things cancel as above. The only other
nonzero terms are the following:
\begin{displaymath}
\sum_j (\al_{K +\ep_j} \tt \al_I \d_j -
\al_I \d_j \tt \al_{K + \ep_j}) \tt \d^{(I)}\d^{(K)},
\end{displaymath}
where $\{\ep_j\}$ is the standard basis of $\Zset^N$.

We have seen that all other contributions have coefficients of degree at most
$2d$, so the sum $\sum_j (\al_{K +\ep_j} \tt \al_I \d_j -
\al_I \d_j \tt \al_{K + \ep_j})$ must lie
inside $\fil^{2d}(H\tt H)$.
All $\al_{K +\ep_j}$ are of degree at most $d$ and $\al_I \d_j$
are of degree exactly $d+1$, hence
$\sum_j \al_{K +\ep_j} \tt \al_I \d_j$
must lie in $\fil^{2d}(H\tt H)$ too. 
But this implies that
$\al_{K + \ep_j}\in \fil^{d-1}H$ for all $j$,
so in particular $\al_I\in \fil^{d-1}H$, which is a contradiction.

This proves that $\al\in(\dd+\kk)\tt(\dd+\kk)$.
Now if $\al=r+s_1-s_2$, where $r\in\dd\wedge\dd$, $s\in\dd$, then we have:
\begin{displaymath}
(\De\tt\id)(\al) = r_{13} + r_{23} + s_1 + s_2 - s_3,
\end{displaymath}
and \eqref{al2} becomes
\begin{equation}\lbb{cybe2}
\bigl([r_{12}, r_{13} + s_1 + s_2] + r_{12} s_3\bigr) + \text{cyclic} = 0.
\end{equation}
Comparing the terms in $\dd\tt\dd\tt\Kset$, we see that \eqref{cybe2}
is equivalent to the system (\ref{cybe3}, \ref{cybe4}).
\end{proof}

Note that when $\al=r\in\dd\wedge\dd$, $s=0$,  \eqref{cybe4} is exactly 
the\/ {\em classical Yang--Baxter equation\/}
\begin{equation}\lbb{cybe}
[r_{12},r_{13}] + [r_{12},r_{23}] + [r_{13},r_{23}] =0.
\end{equation}
Eq.~\eqref{cybe4} is a special case of the 
{\em dynamical\/} classical Yang--Baxter equation
(see \cite{Fe, ES}).

\section{$(H\smash\Kset[\Ga])$-Pseudoalgebras}\lbb{sgaalg}

Let again $H$ be a cocommutative Hopf algebra. 
Let $\Ga$ be a group acting on $H$ by automorphisms,
and let $\ti H=H\smash\Kset[\Ga]$
be the smash product of $H$ with the group algebra of $\Ga$. 
As an \as\ algebra this is the semidirect product of $H$ with $\Kset[\Ga]$,
while as a coalgebra it is the tensor product of coalgebras.

We will denote the action of $\Ga$ on $H$ by $g \cdot f$ for $g\in\Ga$,
$f\in H$; then $g \cdot f = gfg^{-1}$. 
Then a left $\ti H$-module $L$ is the same as an $H$-module together
with an action of $\Ga$ on it which is compatible with that of $H$, i.e.,
such that $(g \cdot f)l = g(f(g^{-1}l))$ 
for $g\in\Ga$, $f\in H$, $l\in L$.

In this section we will study the relationship between the \psten\
categories $\M^*(\ti H)$ and $\M^*(H)$. In particular, we will
show that an $\ti H$-\psalg\ is the same as an $H$-\psalg\
on which the group $\Ga$ acts by preserving the pseudoproduct.

Let us start by defining maps 
$\de_I \colon \ti H^{\tt I} \to H^{\tt I} \tt_H \ti H$
for each finite nonempty set $I$. It is enough to define $\de_I$
on elements of the form $\tt_{i\in I} f_i g_i$
where $f_i \in H$, $g_i \in \Ga$, in which case we let
\begin{displaymath}\lbb{dei}
\de_I(\tt_{i\in I} f_i g_i) = \begin{cases}
(\tt_{i\in I} f_i) \tt_H g, \quad &\text{if all $g_i$ are equal to some $g$},
\\
0, &\text{if some of $g_i$ are different}.
\end{cases}
\end{displaymath}
It is easy to see that $\de_I$ is a homomorphism of both
left $H^{\tt I}$-modules and of right $\ti H$-modules.

This allows us to define a \psten\ functor 
$\de \colon \M^*(\ti H) \to \M^*(H)$
as follows. For an object $L$ (a left $\ti H$-module), we let
$\de(L) \equiv L$ be the left $H$-module obtained by restricting
the action of $\ti H$ to $H \subset \ti H$.
For a polylinear map $\phi\in \Lin( \{L_i\}_{i\in I}, M )$ in $\M^*(\ti H)$,
i.e., for a homomorphism of left $\ti H^{\tt I}$-modules
\begin{displaymath}
\phi \colon \boxtimes_{i\in I} L_i \to \ti H^{\tt I}\tt_{\ti H} M,
\end{displaymath}
we let $\de(\phi)$ be the composition
\begin{displaymath}
\de(\phi) \colon \boxtimes_{i\in I} L_i 
\xrightarrow{\phi} \ti H^{\tt I} \tt_{\ti H} M 
\xrightarrow{\de_I \tt_{\ti H} \id} (H^{\tt I} \tt_H \ti H) \tt_{\ti H} M 
\isoto H^{\tt I} \tt_H M.
\end{displaymath}
This is a homomorphism of left $H^{\tt I}$-modules, i.e.,
a polylinear map in $\M^*(H)$.
Moreover, since the maps $\de_I$ are compatible with the actions
of the symmetric groups and with the comultiplication of $\ti H$,
it follows that $\de$ is compatible with the actions of the symmetric groups
and with compositions of polylinear maps, i.e., it is a \psten\ functor.

As usual, the action of $\Ga$ on $H$ can be extended to an action of 
$\Ga$ on $H^{\tt I}$ by using the comultiplication
$\De^{(I)}(g) = \tt_{i\in I} \, g$.
Hence, $\Ga$ also acts on $H^{\tt I} \tt_H M$ by the formula
\begin{displaymath}\lbb{actgahm}
g\cdot\bigl( (\tt_{i\in I} \, f_i) \tt_H m \bigr)
= (\tt_{i\in I} \, g\cdot f_i) \tt_H gm,
\qquad g\in\Ga, \; f_i \in H, \; m\in M.
\end{displaymath}
Then it is easy to see that $\psi=\de(\phi)$ has the following property:
\begin{equation}\lbb{actgaph}
\psi(\tt_{i\in I} \, g l_i) 
= g \cdot \psi(\tt_{i\in I} \, l_i),
\qquad g\in\Ga, \; l_i \in L_i,
\end{equation}
in other words, it commutes with the action of $\Ga$.

We let $\M^*_\Ga(H)$ be the subcategory of $\M^*(H)$ with
objects left $\ti H$-modules, and with polylinear maps
those polylinear maps $\psi$ of $\M^*(H)$ that commute with the action of 
$\Ga$ (see \eqref{actgaph}). This is a \psten\ category, and $\de$
is a \psten\ functor from $\M^*(\ti H)$ to $\M^*_\Ga(H)$.

\begin{theorem}\lbb{tm*gah}
If\/ $\Ga$ is a finite group,
the functor $\de\colon \M^*(\ti H) \to \M^*_\Ga(H)$
constructed above is an equivalence of \psten\ categories.
\end{theorem}
\begin{proof}
We will construct a \psten\ functor $\Si$ from $\M^*_\Ga(H)$
to $\M^*(\ti H)$. On objects $L$ we let $\Si(L)=L$.
In order to define it on polylinear maps, we need to find out how
$\phi$ can be recovered from $\de(\phi)$ and the action of $\Ga$.

Denote by $\iota$ the embedding $H\injto\ti H$, and let $\pi_I$
be the composition
\begin{displaymath}
\pi_I \colon H^{\tt I} \tt_H \ti H
\xrightarrow{\iota^{\tt I} \tt_H \id} \ti H^{\tt I} \tt_H \ti H
\surjto \ti H^{\tt I} \tt_{\ti H} \ti H
\isoto \ti H^{\tt I}.
\end{displaymath}
Explicitly, $\pi_I$ is given by the formula
\begin{displaymath}
\pi_I ( (\tt_{i\in I} f_i) \tt_H g) = \tt_{i\in I} f_i g,
\qquad f_i\in H, \; g\in \Ga.
\end{displaymath}
This is a homomorphism of both
left $H^{\tt I}$-modules and of right $\ti H$-modules.
Moreover, for $f_i \in H$, $g_i \in \Ga$, we have:
\begin{displaymath}\lbb{piidei}
\pi_I\de_I(\tt_{i\in I} f_i g_i) = \begin{cases}
\tt_{i\in I} f_i g_i, \quad &\text{if all $g_i$ are equal},
\\
0, &\text{otherwise}.
\end{cases}
\end{displaymath}
The crucial observation, which will allow us to invert $\de_I$, 
is that for any $h_i \in \ti H$, $g_i \in \Ga$, we have:
\begin{equation}\lbb{sidei}
\sum_{ (g_i)\in \Ga^I/\Ga } 
(\tt_{i\in I} \, g_i) \, (\pi_I\de_I)(\tt_{i\in I} \, g_i^{-1} h_i) =
\tt_{i\in I} h_i.
\end{equation}
Here $\Ga$ acts diagonally on $\Ga^I$ from the right;
the left-hand side of \eqref{sidei} is invariant under
$(g_i) \mapsto (g_i g)$.

Given a polylinear map $\psi\in \Lin( \{L_i\}_{i\in I}, M )$ 
in $\M^*_\Ga(H)$, we can extend it to a map
\begin{displaymath}
\ti\psi \colon \boxtimes_{i\in I} L_i 
\xrightarrow{\psi} H^{\tt I} \tt_H M 
\isoto (H^{\tt I} \tt_H \ti H) \tt_{\ti H} M 
\xrightarrow{\pi_I \tt_{\ti H} \id} \ti H^{\tt I} \tt_{\ti H} M. 
\end{displaymath}
(Note, however, that $\ti\psi$ is not $\ti H^{\tt I}$-linear.)
Now we define 
$\Si\psi\colon \boxtimes_{i\in I} L_i \to \ti H^{\tt I} \tt_{\ti H} M$
by the formula:
\begin{equation}\lbb{sipsi}
(\Si\psi)(\tt_{i\in I} l_i)
= \sum_{ (g_i)\in \Ga^I/\Ga } ((\tt_{i\in I} \, g_i)\tt_{\ti H}1) \, 
\ti\psi(\tt_{i\in I} \, g_i^{-1} l_i).
\end{equation}
It is easy to check that $\Si\psi$ is $\ti H^{\tt I}$-linear, 
so it is a polylinear map in $\M^*(\ti H)$. Moreover, $\de\Si\psi = \psi$.
For a polylinear map $\phi\in \Lin( \{L_i\}_{i\in I}, M )$ 
in $\M^*(\ti H)$, it is immediate from \eqref{sidei}
and the $\ti H^{\tt I}$-linearity of $\phi$ that
$\Si\de\phi = \phi$. 
Therefore, $\Si\colon \M^*_\Ga(H) \to \M^*(\ti H)$
is a \psten\ functor inverse to $\de$.
\end{proof}
\begin{remark}\lbb{rgalfin}
The above theorem holds also for infinite groups $\Ga$ if we
restrict ourselves to polylinear maps $\psi$ of $\M^*_\Ga(H)$ satisfying
the following finiteness condition:
\begin{equation}\lbb{galfin}
\psi(\tt_{i\in I} \, g_i l_i) \ne 0
\quad\text{for only a finite number of $(g_i)\in \Ga \setminus \Ga^I$}
\end{equation}
for any fixed $l_i \in L$. 
(Note that, by \eqref{actgaph}, this condition does not depend on the 
choice of representatives $(g_i)$.)
Indeed, the only place in the proof where we used
the finiteness of $\Ga$ was to insure that the right-hand side of
\eqref{sipsi} is a finite sum.

If $\psi=\de(\phi)$ comes from a polylinear map $\phi$ of $\M^*(\ti H)$, 
then it satisfies \eqref{galfin}, because $\phi$ is $\ti H^{\tt I}$-linear
and for any element $h\in\ti H^{\tt I}$ one has
$\de_I ((\tt_{i\in I} \, g_i)h) \ne 0$ 
for only a finite number of $(g_i)\in \Ga \setminus \Ga^I$.

Therefore, $\de \colon \M^*(\ti H) \to \M^*_{\Ga,\text{fin}}(H)$
is an equivalence of \psten\ categories, where $\M^*_{\Ga,\text{fin}}(H)$
is the subcategory of $\M^*_{\Ga}(H)$ consisting of polylinear maps $\psi$
satisfying \eqref{galfin}.
\end{remark}
\begin{corollary}\lbb{cgapsal}
A Lie $\ti H=(H\smash\Kset[\Ga])$-\psalg\ 
$L$ is the same as a Lie $H$-\psalg\ $L$
on which the group $\Ga$ acts in a way compatible with
the action of $H$, by preserving the $H$-pseudobracket{\rm:}
\begin{equation}\lbb{gapsal1}
[ga*gb] = g \cdot [a*b] 
\qquad\text{for $g\in\Ga$, $a,b\in L$},
\end{equation}
and satisfying the following finiteness condition{\rm:}
\begin{equation}
\text{given $a,b\in L$, $[ga*b]$ is nonzero
for only a finite number of $g\in \Ga$.}
\end{equation}
The $\ti H$-pseudobracket of $L$ is given by the formula{\rm:}
\begin{equation}\lbb{gapsal2}
[a{\,\ti*\,}b] = \sum_{g\in\Ga} \bigl((g^{-1}\tt1)\tt_{\ti H}1\bigr) \, [ga*b],
\qquad a,b\in L.
\end{equation}

A similar statement holds for representations, as well as for \as\ \psalgs.
\end{corollary}

This result, combined with Kostant's \thref{tkostant},
will allow us in many cases to reduce the study of
$H$-\psalgs\ to the case when $H$ is a universal enveloping algebra
(see \seref{shkga}).

\begin{example}\lbb{exgcf}
Let $\Ga$ be a subgroup of $\Kset^*$ and let 
\begin{displaymath}
H = \Kset[\d]\smash\Kset[\Ga] \linebreak[1] = \linebreak[0]
\bigoplus_{m\in\Zset_+, \al\in\Ga}\,\Kset\, \d^m T_\al
\end{displaymath}
with multiplication $T_\al T_\be = T_{\al\be}$, $T_1=1$, $T_\al \d
T_\al^{-1} = \al\d$ and comultiplication $\De(\d)=\d\tt 1 + 1\tt\d$,
$\De(T_\al)=T_\al \tt T_\al$.  
Then the notion of a Lie $H$-\psalg\ 
is equivalent to that of a $\Ga$-conformal algebra
(cf.\ \cite{K4}).
\end{example}
\begin{example}\lbb{exgtwcf}
Let now $H = \Kset[\d]\times F(\Ga)$, where $F(\Ga)$ is the
function algebra of a finite abelian group $\Ga$.  In other words, $H=
\bigoplus_{m\in\Zset_+, \al\in\Ga}\,\Kset\, \d^m \pi_\al$ with
multiplication $\pi_\al \pi_\be = \de_{\al,\be}\pi_\al$,
$\d\pi_\al=\pi_\al\d$ and comultiplication $\De(\d)=\d\tt 1 + 1\tt\d$,
$\De(\pi_\al) = \sum_{\ga\in\Ga} \pi_{\al\ga^{-1}} \tt \pi_\ga$.  Then
one gets the notion of a $\Ga$-twisted conformal algebra 
(cf.\ \cite{K4}).
\end{example}

\section{A Digression to Linearly Compact Lie Algebras}\lbb{scartan}
We will view the base field $\kk$ as a topological field with discrete
topology.
A topological vector space $\L$ over $\kk$
is called \emph{linearly compact} if it is the space of
all linear functionals on a vector space $\V$ with discrete
topology, with the topology on $\L$ defined by taking all subspaces
$\{ U^\perp \subset\L \st U\subset\V, \; \dim U < \infty \}$
as a fundamental system of neighborhoods of $0$ in $\L$.
Here, as usual, $U^\perp$ denotes the subspace of $\L$ consisting of all linear
functionals vanishing on $U$.

In general, given a topological vector space $\W$, we define a topology
on $\W^*$ by taking for the fundamental system of neighborhoods of $0$
the subspaces $U^\perp$ where $U$ is a linearly compact subspace of $\W$.

Several equivalent definitions of linear-compactness
are provided by the next proposition.


\begin{proposition}\lbb{eqdefs}
For a topological vector space $\L$ over the topological field $\kk$ 
the following statements are equivalent{\rm:}
\begin{enumerate}
\item $\L$ is the dual of a discrete vector space.

\item The topological dual $\L^*$ of $\L$ is a discrete topological space.


\item $\L$ is the topological product of finite-dimensional discrete vector 
spaces.

\item $\L$ is the projective limit of finite-dimensional discrete vector 
spaces.

\item $\L$ has a collection of finite-codimensional open subspaces whose 
intersection is $\{0\}$, with respect to which it is complete.
\end{enumerate}
\end{proposition}
\begin{proof}
Can be found in \cite{G1}.
\end{proof}

\begin{remark}\lbb{bidual}
For both discrete and linearly compact vector spaces, the canonical 
map from $\L$ to $\L^{**}$ is an isomorphism.
\end{remark}


A {\em linearly compact\/} (associative or Lie) 
{\em algebra\/} is a topological
(associative or Lie) algebra for which the underlying topological
space is linearly compact. 

The basic example of a linearly
compact associative algebra is the algebra $\O_N = \kk[[t_1,\dots
,t_N]]$ of formal power series over $\kk$ in $N \geq 1$
indeterminates $t_1, \dots ,t_N$, with the usual formal topology
for which $ (t_1 , \dots ,t_N)^j$, the powers 
of the ideal $(t_1,\dots ,t_N)$, form a fundamental system of
neighborhoods of $\O_N$. 

\begin{remark}
The topological vector spaces $\O_N$ 
are isomorphic and characterized among linearly compact
vector spaces by each of the following properties:
\begin{enumerate}
\item 
  $\O^*_N$ is countable-dimensional.
\item 
  $\O_N$ has a filtration by open subspaces.
\end{enumerate}
\end{remark}

\begin{remark}
{\rm(i)}
One defines a {\em completed tensor product\/} of two linearly compact vector 
spaces $\V, \W$ by $\V{\what \tt}\W = (\V^* \tt \W^*)^*$ where we put the 
discrete topology on $\V^* \tt \W^*$. Then $\V{\what \tt}\W$ is
linearly compact.

{\rm(ii)}
With this definition, $\O_{M+N} \simeq \O_M {\what \tt} \O_N$ as topological
algebras.

{\rm(iii)}
Given a commutative associative linearly compact algebra $\O$ and a 
linearly compact Lie algebra $\L$, their 
completed tensor product $\O \what{\otimes} \L$ is again a
linearly compact Lie algebra.
\end{remark}

The basic example of a linearly compact Lie algebra is the Lie
algebra $W_N$ of continuous derivations of the topological
algebra $\O_N$. The filtration 
\begin{displaymath}
\fil_j \O_N = (t_1,\dots,t_N)^{j+1} , \qquad j=-1,0,1,\dots
\end{displaymath}
of $\O_N$
induces the \emph{canonical filtration} $\fil_j W_N$ of $W_N$, where 
\begin{displaymath}
\fil_j W_N = \{ D\in W_N \st D (\fil_i \O_N)
\subset \fil_{i+j} \O_N
\;\;\forall i \}  , \qquad j=-1,0,1,\dots \; .
\end{displaymath}
It is clear that $W_N$ consists of all 
linear differential operators of the form
\begin{displaymath}
  D=\sum^N_{i=1} P_i (t) \frac{\partial}{\partial t_i},
     \quad\text{ where }\; P_i(t) \in \O_N ,
\end{displaymath}
and that $\fil_j W_N$ $(j\ge -1)$ consists of those $D$ for which all $P_i(t)$
lie in $\fil_j \O_N$.

Let $E=\sum^N_{i=1} t_i \frac{\partial}{\partial
    t_i}$ be the Euler operator.  The spectrum of $\ad E$
consists of all integers $j \geq -1$, and, denoting by
$W_{N;j}$ the $j$-th eigenspace of $\ad E$ we obtain the
\emph{canonical} $\ZZ$\emph{-gradation}:
\begin{displaymath}
W_N = \prod_{j \geq -1} W_{N;j} , \quad 
  [W_{N;i}, W_{N;j}] \subset W_{N;i+j} .
\end{displaymath}
The following fact is well known.

\begin{lemma}
  \lbb{lem:1}
$W_{N;0} \simeq \gl_N (\kk)$ and one has the following
isomorphism of $\gl_N(\kk)$-modules{\rm:}
\begin{displaymath}
  W_{N;j} \simeq \kk^N \otimes (\symp^{j+1} \kk^N)^* .
\end{displaymath}
Furthermore, one has a decomposition into a direct sum of
irreducible submodules{\rm:} $W_{N;j}=W'_{N;j} + W''_{N;j}$, where
$W'_{N;j} \simeq (\symp^j\kk^N)^*$ $(=0$ if $j=-1)$ and $W''_{N;j} \simeq$ 
the highest component of
$\kk^N \otimes (\symp^{j+1}\kk^N)^*$.  The subspace $\fp= W_{N;-1} +
W_{N;0} + W'_{N;1}$ is a subalgebra of $W_N$ isomorphic to 
$\sl_{N+1} (\kk)$.
\end{lemma}

Let $\Omega_N = \bigoplus^N_{j=0} \Omega_{N;j}$
denote the algebra of differential forms over $\O_N$.  The
defining representation of $W_N$ on $\O_N$ extends uniquely to a
representation on $\Omega_N$ commuting with the differential
$\di$.

Recall that a \emph{volume form} is a differential $N$-form
$v=f(t_1,\dots ,t_N) \, \di t_1 \wedge \dots \wedge \, \di t_N$ such
that $f(0) \neq 0$, a \emph{symplectic form} is a closed $2$-form
$s= \sum^N_{i,j=1} s_{ij} (t_1,\dots,t_N) \, \di t_i
  \wedge \, \di t_j$ such that $\det (s_{ij} (0))\neq 0$, and a
\emph{contact form} is a $1$-form $c$ such that $c \wedge
(\di c)^{(N-1)/2}$ is a volume form.  The following facts are well
known.

\begin{lemma}\lbb{lem:2}
{\rm(i)}
  Any volume form can be transformed by an automorphism of $\O_N$ 
  to the standard volume form $v_0=\di t_1 \wedge \dots \wedge \,
  \di t_N$.

{\rm(ii)}
  A symplectic form exists iff $N$ is even, $N=2n$, and by an
  automorphism of $\O_N$ it can be transformed to the standard
  symplectic form $s_0 = \sum^N_{i=1} \, \di t_i
    \wedge \, \di t_{n+i}$.

{\rm(iii)}
  A contact form exists iff $N$ is odd, $N=2n+1$, and by an
  automorphism of $\O_N$ it can be brought to the standard
  contact form $c_0=\di t_N+\sum^n_{i=1} t_i \, \di t_{n+i}$.
\end{lemma}

Consider the following (closed) subalgebras of the Lie algebra
$W_N$:
\begin{align*}
  S_N (v) &= \{ D \in W_N \st Dv =0 \} \quad (N \geq 2), \\
  H_N (s) &= \{ D \in W_N \st Ds=0\} \quad 
               (N \text{ even }\geq 2), \\
  K_N (c) &= \{ D \in W_N \st Dc =fc \quad 
               \text{ for some } f \in \O_N \}
               \quad (N \text{ odd } \geq 3).
\end{align*}
Let also $S_N =S_N(v_0)$, $H_N=H_N(s_0)$, $K_N=K_N(c_0)$.
Lemma~\ref{lem:2} implies isomorphisms:  $S_N(v) \simeq S_N$,
$H_N(s) \simeq H_N$, $K_N(c) \simeq K_N$, $S_2 \simeq H_2$.

The  canonical filtration of $W_N$ induces \emph{canonical
  filtrations} $\fil_j S_N (v) := \fil_j W_N \cap S_N (v)$, etc.  
Note that $\dim W_N / \fil_{-1} W_N = N$. A Lie subalgebra $\L$
of $W_N$ is called {\em transitive\/} if $\dim \L / (\L\cap\fil_{-1} W_N) = N$.
It is known that the Lie algebras $W_N$, $S_N$, $H_N$ and $K_N$
are transitive.
In addition, 
the canonical filtrations $\fil_j\L$ of these Lie algebras have the
following transitivity property:
\begin{equation}\lbb{transit}
\fil_{j+1}\L = \{ a \in\fil_j\L \st [a,\L] \subset \fil_j\L \}.
\end{equation}

Noting that $Ev_0=Nv_0$ and
$Es_0=2s_0$, we conclude that $\ad E $ is an
(outer) derivation of $S_N$ and $H_N$, hence the canonical
gradation of $W_N$ induces \emph{canonical}
$\ZZ$\emph{-gradations} ${S_N = \prod_{j \geq -1}
  S_{N;j}}$ and ${H_N=\prod_{j\geq -1} H_{N;j}}$.

Let $E'=2t_N \tfrac{\partial}{\partial t_N} + \sum^{N-1}_{i=1}
t_i \tfrac{\partial}{\partial t_i}$.  Then $E'c_0=2c_0$, hence $
E' \in K_N$ and the eigenspace decomposition of $\ad E'$ defines
the \emph{canonical} $\ZZ$\emph{-gradation}
${K_N=\prod_{j \geq -2} K_{N;j}}$.  The following
facts are well known.

\begin{lemma}\label{lem:3}
{\rm(i)}
  $S_{N;0} \simeq \sl_N (\kk)$, \; $H_{N;0} \simeq \sp_N (\kk)$, \;
  $K_{N;0}\simeq \csp_{N-1} (\kk)$ \break $(\simeq \sp_{N-1}(\kk)\oplus \kk)$.

{\rm(ii)}
The $S_{N;0}$-module $S_{N;j}$ is isomorphic to the highest
component of the $\sl_N(\kk)$-module $\kk^N \otimes (\symp^{j+1}\kk^N)^*$.

{\rm(iii)}
  The $H_{N;0}$-module $H_{N;j}$ is isomorphic to the
  $($irreducible$)$ $\sp_N(\kk)$-module $\symp^{j+2}\kk^N$.

{\rm(iv)}
  $K_{N;0} = \sp_{N-1}(\kk) \oplus \kk E'$ and the $\sp_{N-1}(\kk)$-module 
  $K_{N;j}$ decomposes into the following direct sum of
  irreducible modules{\rm:}
  \begin{displaymath}
    K_{N;j}= \bigoplus^{\left[ \frac{j}{2} \right]+1}_{i=0}
      K^{(i)}_{N;j}, \text{ where }
      K^{(i)}_{N;j} \simeq \symp^{j+2-2i} \kk^{N-1} .
  \end{displaymath}
The subspace $\fp=K_{N;-2} + K_{N;-1} + K_{N;0}+ K^{(1)}_{N;1} +
K^{(2)}_{N;2}$  is a subalgebra of $K_N$ isomorphic to $\sp_{N+1} (\kk)$.
\end{lemma}

The following celebrated theorem goes back to E.~Cartan (see
\cite{G2} for a relatively simple proof).

\begin{theorem}\lbb{tcarg2}
  Any infinite-dimensional simple linearly compact Lie algebra is 
  isomorphic to one of the topological Lie algebras $W_N$, $S_N$, 
  $H_N$ or $K_N$.
\end{theorem}

Let $\g$ be a Lie algebra, and let $\h$ be its subalgebra of codimension $N$.
Then $F = \Hom_{\ue(\h)}(\ue(\g), \Kset)$, with the product
$(f_1 f_2)(u) = f_1(u_{(1)}) f_2(u_{(2)})$, 
is (canonically) isomorphic to the algebra 
of formal power series on $(\g/\h)^*$ \cite{BB},
which is (non-canonically) isomorphic to the linearly compact algebra $\O_N$.
The Lie algebra $D$ of continuous 
derivations of $F$ is then isomorphic to $W_N$. $F$ has
a canonical $\g$-action induced by the left-multiplication
$\g$-action on $\ue(\g)$, which
gives us a homomorphism $\gamma$ of $\g$ to $W_N$. 
(This is non-canonical since the identification of 
$F$ with $\O_N$ is not canonical.)

We will use in the sequel the following theorem of Guillemin and Sternberg 
\cite{GS} (see \cite{BB} for a simple proof).

\begin{proposition}\lbb{existuniq}
Let $\g$ be a Lie algebra, and let $\h$ be its subalgebra of codimension $N$.
Provided that $\h$ contains no nonzero ideals of $\g$,
the above-defined
$\gamma$ is a Lie algebra isomorphism of $\g$ with a subalgebra
of $W_N$ which maps $\h$ into $\fil_0 W_N$.

$($Conversely, if $\g\injto W_N$ maps $\h$ into $\fil_0 W_N$, then
$\h$ doesn't contain nonzero ideals of $\g.)$
Every Lie algebra homomorphism of $\g$ to $W_N$,
which coincides with $\gamma$ modulo
$\fil_0 W_N$, is conjugated to $\gamma$ via a unique automorphism of $\O_N$.
\end{proposition}

We have the following 
important property of the filtrations on $H$ and $X=H^*$,
defined in \seref{subfiltop}.

\begin{lemma}\lbb{minusi}
Let $H=\ue(\g)\smash\kk[\Ga]$ be a cocommutative Hopf algebra, and $X=H^*$.
If $h \in \fil^i \ue(\g) \subset H$ but $h \notin\fil^{i-1}\ue(\g)$,
then $h \fil_n X = \fil_{n-i} X$.
In particular, for any 
$h\in\g\setminus\{0\}$
and for every open
subspace $U\subset X$, there is some $n$ such that $h^n U=X$.
Similar statements hold for the right action of $h$ as well.
\end{lemma}
\begin{proof}
By the construction of the filtrations it is evident
that we can assume $H=\ue(\g)$. 
Then $X\simeq\O_N$ $(N=\dim\g)$, and $\g\injto W_N$ acts on it by linear 
differential operators.
The rest of the proof is clear.
\end{proof}

The following result from \cite{G1, G2} will be essential for our 
purposes.

\begin{proposition}\lbb{guill1}
{\rm(i)}
A linearly compact Lie algebra $\L$ satisfies the 
descending chain condition on 
closed ideals if and only if it has a fundamental subalgebra, i.e., an open 
subalgebra containing no ideals of $\L$.

{\rm(ii)}
When either of the assumptions of {\rm(i)} holds, the 
noncommutative minimal closed ideals of $\L$ are of the form 
$\O_r\what\otimes\ss$ 
where $\ss$ is a simple linearly compact Lie algebra and $r\in\Zset_+$.
\end{proposition}

We will also need the following examples of non-simple linearly
compact Lie algebras:
\begin{align*}
  CS_N (v) &= \{ D \in W_N \st Dv= av, \;\; a \in \kk \} , \\
  CH_N (s) &= \{ D \in W_N \st Ds=as, \;\; a \in \kk \} .
\end{align*}
As before, we have isomorphisms 
$CS_N(v) \simeq CS_N \equiv CS_N (v_0)$
and 
$CH_N (s) \simeq CH_N \equiv CH_N (s_0)$.  Note also
that $CS_N =\kk E \ltimes S_N$ and $CH_N=\kk E \ltimes H_N$.  Another
important example of a non-simple linearly compact Lie algebra is 
the Poisson algebra $P_N$, which is $\O_N$ $(N=2n)$ endowed with the
Poisson bracket:
\begin{displaymath}
  \{ f,g \} = \sum^n_{i=1} \frac{\partial f}{\partial t_i}
  \frac{\partial g}{\partial t_{n+i}} 
  -\frac{\partial f}{\partial t_{n+i}}
  \frac{\partial g}{\partial t_i} .
\end{displaymath}
It is a nontrivial central extension of $H_N$:
\begin{displaymath}
  0 \to \Kset \to P_N \overset{\varphi}{\to} H_N \to 0 ,
\end{displaymath}
where ${\varphi(f)=\sum^n_{i=1} \frac{\partial
    f}{\partial t_i} \frac{\partial}{\partial t_{n+i}} -
  \frac{\partial f}{\partial t_{n+i}} \frac{\partial}{\partial
    t_i}}$.

We can describe also $K_N$ in a more explicit way, similar to the above
description of $P_N$.
For $f,g\in \O_N$, define
\begin{displaymath}
\{f,g\}' = \{f,g\}_{2n}
+ \frac{\d f}{\d t_{2n+1}} (E_{2n}g - 2g)
- (E_{2n}f - 2f) \frac{\d g}{\d t_{2n+1}},
\end{displaymath}
where $\{f,g\}_{2n}$ is the Poisson bracket taken with respect to the
variables $t_1,\dots,t_{2n}$ and $E_{2n}$ is the Euler operator
$\sum_{i=1}^{2n}t_i \frac\d{\d t_i}$.
If we define
\begin{displaymath}
\psi(f) =  \sum_{i=1}^n \Bigl( \frac{\d f}{\d t_i} \frac{\d }{\d t_{n+i}}
- \frac{\d f}{\d t_{n+i}} \frac{\d }{\d t_i} \Bigr)
+ \frac{\d f}{\d t_{2n+1}} E_{2n}
- (E_{2n}f - 2f) \frac{\d}{\d t_{2n+1}},
\end{displaymath}
then we have $\psi(f) g = \{f,g\}' + 2 \frac{\d f}{\d t_{2n+1}} g$.
It is easy to see that $[\psi(f),\psi(g)] = \psi(\{f,g\}')$
and $\psi(f) c_0 = 2\frac{\d f}{\d t_{2n+1}} c_0$.
Thus $K_N$ is isomorphic to $\O_N$ with the bracket $\{,\}'$.

For a linearly compact Lie algebra $\L$ denote by $\Der \L$ the
Lie algebra of its continuous derivations and by $\what{\L}$ the
universal central extension of $\L$. Then we have:

\begin{proposition}\lbb{pder}
{\rm(i)}
      $\Der W_N = W_N$,
      $\Der S_N = CS_N$, $\Der H_N=CH_N$,
      $\Der K_N = K_N$.

{\rm(ii)}
  $\Der (\O_{r} \what{\otimes} \L) =W_{r} \otimes 1 + \O_{r}
  \what{\otimes} \Der \L$ for any simple linearly compact Lie algebra~$\L$.

{\rm(iii)}
  The Lie algebras $\O_{r} \what{\otimes} W_N$, $\O_{r}
  \what{\otimes}S_N$ $($for $N>2)$ and $\O_{r} \what{\otimes} K_N$ have no
  nontrivial central extensions.  The universal central
  extension of $\O_{r} \what{\otimes} H_N$ is $\O_{r}
  \what{\otimes} P_N$.

{\rm(iv)}
  If $\fg$ is a simple finite-dimensional Lie algebra, then
  $(\O_{r} \otimes \fg)^{\what{}} = (\O_{r} \otimes \fg) 
  + (\Omega_{r ;1}/\di \O_{r})$ with the bracket
  \begin{displaymath}
    [f \otimes a , g \otimes b]^{\what{}}
    = fg \otimes [a,b] + (a|b) f\di g \mod \, \di\O_{r} ,
  \end{displaymath}
where $(a|b)$ is the Killing form on $\g$.
\end{proposition}
\begin{proof}
For a proof of~(iv) see \cite{Ka}.

In order to prove (ii), notice that if $d$ is a derivation of 
$\O_r\what\tt \L$,
then its action on $1 \tt \L$ is
given by $d(1\tt x) = \sum_i a_i \tt d_i(x)$ for all $x\in \L$, 
where the $a_i$ form a topological
basis of $\O_r$ and the $d_i$ are continuous derivations of $\L$.
Subtracting $\sum_i a_i\tt d_i$ from $d$, 
we get a derivation $\tilde d$ acting trivially
on $1 \tt \L$. We are going to show that if $\L$ is simple then
$\tilde d$ is of the form $D \tt 1$ where $D \in \Der\O_r = W_r$.

Let us fix $P \in \O_r$. Then $\tilde d (P \tt x)$ can be written as $\sum_i
a_i \tt f_i(x)$, where $f_i$ are continuous $\kk$-endomorphisms of $\L$.
{}From $\tilde d([P\tt x, 1\tt y]) = [\tilde d(P\tt x), 1\tt y]$
we see that $[f_i(x),y] = f_i([x,y])$ for every
$x, y \in \L$. This means that $f_i$ commutes with $\ad y$ for all $y\in\L$.
By a Schur's lemma argument \cite[Proposition~4.4]{G1}
and simplicity of $\L$, we conclude that the $f_i$ are multiples of
the identity map, hence $\tilde d (P \tt x) = a_P \tt x$ 
for some $a_P\in\O_r$ and all $x\in \L$.
It is now immediate to check that the mapping $D\colon P \mapsto a_P$ is indeed
a derivation of $\O_r$, proving (ii).

In order to prove the rest of the statements, denote by $\fa$ the
$0${th} component of the canonical $\ZZ$-gradation of $\L =
W_N$, $S_N$, $H_N$ or $K_N$.  This is a reductive subalgebra of $
\L$, hence $\Der \L =V \oplus \L$, where $[\fa ,V] \subset V$.
But $[V,\L] \subset \L$, hence $[\fa , V]=0$, i.e.,~any element
$D \in V$ defines an endomorphism of $\L$ viewed as an
$\fa$-module.  Since $E \in W_N$ and $E' \in K_N$, we conclude
that $D$ also preserves the canonical gradation of these Lie
algebras and we may assume that $D$ acts trivially on the $(-1)${st} 
component.  Using the
transitivity of $W_N$ and $K_N$, we conclude that $D=0$.  By
Lemma~\ref{lem:3}, all components of the canonical
$\ZZ$-gradation of $S_N$ and $H_N$ are inequivalent
$\fa$-modules, hence $D$ preserves this gradation in this case as 
well.  Subtracting from $D$ a multiple of $E$, we may assume that 
$D$ acts trivially on the $(-1)${st} component and, using
transitivity, we conclude that $D$ is a multiple of $E$.  Thus (i)~is
proved.

Since $\fa$ acts completely reducibly on the space $Z^2$ of
$2$-cocycles on $\O_{r} \what{\otimes} \L$ with values in
$\kk$, and since $\fa$ acts trivially on cohomology, we may choose
a subspace $U$ of $Z^2$, complementary to the space of trivial
$2$-cocycles, on which $\fa$ acts trivially.  Hence for any
$2$-cocycle $\alpha \in U$ we have:  $\alpha (a,b)=0$ if $a \in
M_1$, $b \in M_2$ and $M_i$ are irreducible non-contragredient
$\fa$-submodules of $\overline{\L}:=\O_{r} \what{\otimes} \L$.
Let $\overline{\L}_j=\O_r \otimes \L_j$ for short, where $\L_j$ is the $
j${th} component of the canonical gradation.

It follows from Lemma~\ref{lem:3}(ii) that all pairs of
$\fa$-submodules in $\overline{\L}=\O_{r } \what{\otimes} S_N$ are 
non-contragredient, except for the adjoint $\fa$-submodules in
$\overline{\L}_0 = \O_{r} \otimes S_{N;0}$.  Thus, we have:
\begin{displaymath}
  \alpha (a,b)=0 \;\text{ if } a \in \overline{S}_{N;i}, \,\,
  b \in \overline{S}_{N;j}, \quad i \neq 0 \;\text{ or }\; j \neq 0 .
\end{displaymath}
Taking now $a \in \overline{S}_{N;-1}$, $b \in  \overline{S}_{n;1}$ and $c 
\in \overline{S}_{N;0}$, the cocycle condition
\begin{displaymath}
  \alpha ([a,b],c) + \alpha ([b,c],a) + \alpha ([c,a],b) =0
\end{displaymath}
gives $\alpha ([a,b],c)=0$.  Since
$\overline{S}_{N;0}=[\overline{S}_{N;-1}, \overline{S}_{N;1}]$, we conclude 
that 
$\alpha =0$.  Hence all central extensions of $\overline{S}_N$ are
trivial.

Likewise, $\alpha$ is zero on any pair of subspaces
$\overline{W}_{N;i}, \overline{W}_{N;j}$, unless $i+j=0$, and on the pair
$\overline{W}_{N;-1}, \overline{W}''_{N;1}$ (see Lemma~\ref{lem:1}).
Choosing $a \in \overline{W}_{N;-1}$, $b \in \overline{W}''_{N;1}$, $c \in
\overline{W}_{N;0}$, we obtain, as above, from the cocycle condition, 
that $\alpha$ is zero on the pair $\overline{W}_{N;0},
[\overline{W}_{N;0}, \overline{W}_{N;0}]$.  It follows from (iv) applied to
the subalgebra $\O_{r} \otimes \sl_{N+1} (\kk)$ of
$\overline{W}_N$ (see Lemma~\ref{lem:1}) that $\alpha$ is zero on this 
subalgebra if $N>1$.  Thus any cocycle on $\overline{W}_N$ $(N>1)$
is trivial.  In the case of $\overline{W}_1$ 
the cocycle $\alpha $ is trivial. The case of $\overline{K}_N$ is similar.

In the remaining case of $\overline{H}_N$ we show, as above, that the
cocycle $\alpha$ is trivial on any pair $\overline{H}_{N;i},
\overline{H}_{N;j}$ if $i \neq j$.  Using the cocycle condition for $a 
\in \overline{H}_{N;k}$, $b \in \overline{H}_{N;k+1}$ and $c \in
\overline{H}_{N;-1}$, and the fact that $\overline{H}_{N;k}=
[\overline{H}_{N;k+1}, \overline{H}_{N;-1}]$, we conclude that $\alpha$ is
trivial on any pair $\overline{H}_{N;i}, \overline{H}_{N;i}$ as well,
unless $i =-1$.  It is easy to see that this implies that
$\what{\overline{H}}_N = \overline{P}_N$.
\end{proof}

\section{$H$-Pseudoalgebras and $H$-Differential Algebras}\lbb{slalg}
In this section, $H$ will be a cocommutative Hopf algebra,
and as before, $X=H^*$.

\subsection{The annihilation algebra}\lbb{subanih}
Let $Y$ be an $H$-bimodule which is a
commutative associative $H$-\difalg\ 
both for the left and for the right action of $H$ 
(see \eqref{hxy}, \eqref{xyh}); for example, $Y=X := H^*$.

For a left $H$-module $L$, let
$\A_Y L = Y\tt_H L$. 
We define a left action of $H$ on $\A_Y L$ in the obvious way:
\begin{equation}\lbb{hahx}
h(x\tt_H a) = hx\tt_H a, \qquad h\in H, \; x\in Y, \; a\in L.
\end{equation}
If, in addition, $L$ is an $H$-\psalg\ with a pseudoproduct $a*b$,
we can define a product on $\A_Y L$ by the formula:
\begin{equation}\lbb{ahxbhy}
\begin{split}
(x\tt_H a)(y\tt_H b) &= \tsum_i\, (x f_i)(y g_i) \tt_H e_i,
\\
\text{if}\quad a*b &= \tsum_i\, (f_i\tt g_i) \tt_H e_i.
\end{split}
\end{equation}
By \eqref{xyh} and the $H$-bilinearity \eqref{bil*} of the pseudoproduct,
it is clear that \eqref{ahxbhy} is well defined.

\begin{proposition}\lbb{plieyr}
If $L$ is a Lie $H$-\psalg, then
$\A_Y L$ is a Lie $H$-\difalg, i.e., a Lie algebra which is also
a left $H$-module so that
\begin{equation}\lbb{halbe}
h[\al,\be] = [h_{(1)}\al, h_{(2)}\be], 
\qquad\text{for }\; h\in H, \; \al,\be\in \A_Y L.
\end{equation}
Similarly, if $L$ is an \as\  $H$-\psalg, 
then $\A_Y L$ is an \as\ $H$-\difalg.
A similar statement holds for modules as well\/{\rm:}
if $M$ is an $L$-module, then $\A_Y M$ is an $\A_Y L$-module
with a compatible $H$-action so that
\begin{equation}\lbb{ham12}
h(am)=(h_{(1)}a)(h_{(2)}m)
\quad\text{for}\;\; h\in H, \; a\in \A_Y L, \; m\in \A_Y M.
\end{equation}
\end{proposition}
\begin{proof}
Equation \eqref{halbe} follows from \eqref{hxy}.
The skew-commutativity of the bracket \eqref{ahxbhy} follows immediately
from that of $[a*b]$. The proof of the Jacobi identity is 
straightforward by using \eqref{jac*}. 
Let us check for example that the associativity of $L$ is
equivalent to that of $\A_Y L$; 
the case of the Jacobi identity is similar.

We will use the notation from \eqref{ab*'}--\eqref{abc*6'},
and we will write $a_x \equiv x\tt_H a$ for $a\in L$, $x\in Y$. 
Then we want to compute the products 
$a_x(b_y c_z)$ and $(a_x b_y)c_z$.
By definition, if we have \eqref{abc*4'} and \eqref{abc*5'}, then
\begin{align*}
b_y c_z &= \tsum_i\, (y h_i)(z l_i) \tt_H d_i
\intertext{and}
a_x(b_y c_z) &= \tsum_{i,j}\, 
\bigl(x h_{ij}\bigr) \bigl(\bigl((y h_i)(z l_i)\bigr)l_{ij}\bigr) \tt_H d_i
\\
&= \tsum_{i,j}\, 
(x h_{ij})(y h_i{l_{ij}}_{(1)})(z l_i{l_{ij}}_{(2)}) \tt_H d_i.
\end{align*}
Similarly, if we have \eqref{ab*'} and \eqref{abc*1'}, then
\begin{displaymath}
(a_x b_y)c_z = \tsum_{i,j}\, 
(x f_i{f_{ij}}_{(1)})(y g_i{f_{ij}}_{(2)})(z g_{ij})\tt_H e_i.
\end{displaymath}
Now recalling \eqref{abc*3'} and \eqref{abc*6'}, we see that
the associativity
of $L$ is equivalent to that of $\A_Y L$.
\end{proof}
\begin{definition}\lbb{dannih}
The $H$-\difalg\ $\A(L) \equiv \A_X L := X \tt_H L$ 
is called the {\em annihilation algebra\/} of the \psalg\ $L$. 
We will write $a_x \equiv x\tt_H a$ for $a\in L$, $x\in X$.
\end{definition}
\begin{remark}
When $L$ is an \as\ $H$-\psalg, one does not need the cocommutativity of $H$
or the commutativity of $Y$ in order to define $\A_Y L$ (cf.\ \reref{rsymac}).
\end{remark}
\begin{lemma}\lbb{torzerocoef}
Let $H=\ue(\dd)\smash\kk[\Ga]$, and let $M$ be a left $H$-module.
If an element $a\in M$ is $\ue(\dd)$-torsion, i.e.,
if $ha = 0$ for some $h \in \ue(\dd)\setminus\{0\}$, 
then $X\otimes_H a = 0$.
In particular, for $H=\ue(\dd)$, we have{\rm:}
$\A(M) \simeq \A(M / \Tor M)$, where $\Tor M$ is the torsion
submodule of $M$.
\end{lemma}
\begin{proof}
We have $0 = x \otimes_H ha = xh \otimes_H a$
for every $x \in X$. 
Since the right action of $h$ on $X$ is surjective (see \leref{minusi}),
it follows that $x \otimes_H a = 0$ for any $x\in X$.
\end{proof}

\subsection{The functor $\A_Y$}\lbb{subay}
Analyzing the proof of \prref{plieyr}, one can notice that the definition
of $\A_Y$ is a special case of a more general construction which
we describe below.

First, recall that a commutative associative $H$-\difalg\ $Y$
is the same as a commutative associative algebra in the
\psten\ category $\M^l(H)$ from \exref{exmld}.
We denote by $\M^b(H)$ the category of $H$-bimodules, provided
with a \psten\ structure given by \eqref{mld}, but with
$\Hom_H$ there replaced by $\Hom_{H-H}$ which means homomorphisms
of $H$-bimodules. The composition of polylinear maps in $\M^b(H)$ 
is given again by \eqref{com2}. Then $H$-\difalgs\ $Y$
considered above are exactly the commutative associative algebras
in $\M^b(H)$. 

Instead of one $H$-bimodule $Y$ one can use several: 
for any collections of objects
$Y_i \in \M^b(H)$ and $L_i \in \M^*(H)$ $(i\in I)$
we can consider the left $H$-modules $\A_{Y_i}L_i = Y_i \tt_H L_i$
as objects of $\M^l(H)$. Assume we are given polylinear maps
$f\in \Lin( \{Y_i\}_{i\in I}, Z )$ in $\M^b(H)$ and
$\phi\in \Lin( \{L_i\}_{i\in I}, M )$ in $\M^*(H)$.
Then we define a polylinear map 
$f\tt_H\phi\in 
\Lin( \{Y_i \tt_H L_i\}_{i\in I}, Z \tt_H M )$
in $\M^l(H)$ as the following composition:
\begin{multline*}
\tt_{i\in I} (Y_i \tt_H L_i)
\isoto
(\boxtimes_{i\in I}Y_i) \tt_{H^{\tt I}} (\boxtimes_{i\in I}L_i)
\\
\xrightarrow{\id\tt\phi} (\boxtimes_{i\in I}Y_i) 
\tt_{H^{\tt I}} (H^{\tt I}\tt_H M)
\isoto (\tt_{i\in I}Y_i) \tt_H M
\xrightarrow{f\tt\id} Z \tt_H M.
\end{multline*}

\begin{proposition}\lbb{pfgpp}
The above definition is compatible with compositions of polylinear maps
in $\M^b(H)$, $\M^*(H)$, and $\M^l(H)${\rm:}
\begin{displaymath}\lbb{fgpp}
f(\{g_i\}_{i\in I}) \tt_H \phi(\{\psi_i\}_{i\in I})
= \bigl( f \tt_H \phi \bigr) \bigl( \{g_i \tt_H \psi_i\}_{i\in I} \bigr).
\end{displaymath}
\end{proposition}
The proof of this proposition is straightforward and is left to the reader.

\begin{corollary}\lbb{cfgpp}
Let $(Y,\nu)$ be a commutative associative  algebra in $\M^b(H)$.
For a finite nonempty set $I$, let $\nu^{(I)}\colon Y^{\tt I} \to Y$
be the iterated multiplication
$\nu(\nu\tt\id)\dotsm(\nu\tt\id\tt\dotsm\tt\id)$.
Recall that for an object $L$ in $\M^*(H)$, we define $\A_Y(L) := Y \tt_H L$.
For a polylinear map $\phi\in \Lin( \{L_i\}_{i\in I}, M )$ in $\M^*(H)$,
let $\A_Y(\phi) := \nu^{(I)} \tt_H \phi$.
Then $\A_Y$ is a \psten\ functor from $\M^*(H)$ to $\M^l(H)$.
\end{corollary}

As a special case of this corollary, we obtain \prref{plieyr}.
%

Let us give another application of \prref{pfgpp}.
An instance of an $H$-bimodule is $H$ itself (however, $H$ is
not an $H$-\difalg{}!). The coproduct $\De\colon H\to H\tt H$,
the evaluation map $ev\colon X\tt H\to\Kset$, and the isomorphism
$\Kset\tt H \isoto H$
are all homomorphisms of
$H$-bimodules, so the composition 
\begin{equation}\lbb{etaxh}
\eta\colon X\tt H \xrightarrow{\id\tt\De} X\tt H\tt H 
\xrightarrow{ev\tt\id} \Kset\tt H \isoto H
\end{equation}
is a polylinear map in $\M^b(H)$.
Let again $L$ be a (Lie) \psalg\ and $(M,\rho)$ be an $L$-module, where
$\rho\in\Lin( \{L,M\}, M)$ in $\M^*(H)$.
Then $\eta\tt_H\rho \in \Lin( \{X\tt_H L, H\tt_H M\}, H\tt_H M )$ is a 
polylinear map in $\M^l(H)$. In other words, we get a homomorphism of
$H$-modules $\eta\tt_H\rho\colon \A(L) \tt M \to M$.
\prref{pfgpp} now implies:
\begin{corollary}\lbb{calm}
The above map $\eta\tt_H\rho$ provides $M$
with the structure of an $\A(L)$-module, and this structure
is compatible with that of an $H$-module {\rm(}cf.\ \eqref{ham12}{\rm)}.
\end{corollary}

For $a\in L$, $x\in X$, the action of $a_x \equiv x\tt_H a$
on an element $m\in M$ will be denoted by $a_x \cdot m$.
This defines {\em $x$-products\/} 
$\xp axm := a_x \cdot m \in M$.
When $M=L$ is the Lie \psalg\ with the adjoint action, these
will be called {\em $x$-brackets\/} and denoted as $\xb axb$.
Then all the axioms of (Lie or \as) \psalgs, representations, etc.,
can be reformulated in terms of the properties of the $x$-brackets
or products --- this will be done in \seref{shconf}. 
Although this may seem a mere
tautology, it is more explicit and convenient in some cases.

Finally, let us give two more constructions.

\begin{example}\lbb{eakl}
The base field $\kk$, with the action $h \cdot 1=\ep(h)$ ($h\in H$),
is a commutative associative $H$-\difalg. Then for any
Lie $H$-\psalg\ $L$, $\A_\kk L = \kk \tt_H L$ is a Lie $H$-\difalg.
Explicitly, $\A_\kk L \simeq L / H_+ L$, 
where $H_+ = \{ h\in H \st \ep(h)=0 \}$
is the augmentation ideal. The Lie bracket in $L / H_+ L$ is given by
(cf.\ \eqref{ahxbhy}):
\begin{equation}\lbb{akl1}
[a \;\mathrm{mod}\; H_+ L, b  \;\mathrm{mod}\; H_+ L] 
= \tsum_i\, \ep(f_i) \ep(g_i) e_i  \;\mathrm{mod}\; H_+ L,
\end{equation}
if
\begin{equation}\lbb{psbrab}
[a*b] = \tsum_i\, (f_i\tt g_i) \tt_H e_i.
\end{equation}

In the case when $\dd=\kk\d$ is $1$-dimensional, we recover the usual 
construction $L\mapsto L/\d L$ that assigns a Lie algebra to any
Lie conformal algebra \cite{K2}.
\end{example}
\begin{remark}\lbb{affpsalg}
Let $Y$ be a commutative associative $H$-\difalg\
with a right action of $H$,
and let $L$ be a Lie $H$-\psalg. We provide $Y \tt L$ with the following
structure of a left $H$-module:
\begin{equation}
h(x \tt a) = x h_{(-1)} \tt h_{(2)} a, \qquad 
h\in H, \; x\in Y, \; a\in L.
\end{equation}
Then define a Lie pseudobracket on $Y \tt L$ by the formula:
\begin{equation}\lbb{affps}
[(x \tt a)*(y \tt b)] = \tsum_i\, 
\bigl( {f_i}_{(1)} \tt {g_i}_{(1)} \bigr) \tt_H
\bigl( (x {f_i}_{(2)}) (y {g_i}_{(2)}) \tt e_i \bigr),
\end{equation}
if $[a*b]$ is given by \eqref{psbrab}.
It is easy to check that \eqref{affps} is well defined and provides
$Y \tt L$ with the structure of a Lie $H$-\psalg. Moreover,
$\A_Y L \simeq (Y \tt L) / H_+(Y \tt L)$ as a Lie algebra
(cf.\ \exref{eakl}).

In the case $\dd=\kk\d$, $Y=\kk[t,t^{-1}]$, $\d=-\d_t$,
the Lie $\kk[\d]$-\psalg\ (= conformal algebra) 
$Y \tt L$ is known as an {\em affinization\/}
of the conformal algebra $L$ \cite{K2}.
\end{remark}

\subsection{Relation to differential Lie algebras}\lbb{sublocla}
Fix two positive integers $N,r$ and let $\O_N=\Kset[[t_1,\dots,t_N]]$,
$\L=\O_N\tt\Kset^r$. A structure of a Lie algebra on $\L$ is called
{\em local\/} (and $\L$ is called a {\em local Lie algebra} \cite{Ki})
if the Lie bracket is given by matrix bi-differential operators.
More explicitly, let $\{e^i\}$ be a basis of $\Kset^r$. Then for any
$x,y\in\O_N$, the bracket in $\L$ is given by:
\begin{equation}\lbb{xeiyej}
[x\tt e^i, y\tt e^j] 
= \tsum_{k,l} \, (P^{ij}_{kl} \cdot x) (Q^{ij}_{kl} \cdot y) \tt e^k,
\end{equation}
where $P^{ij}_{kl}, Q^{ij}_{kl}$ are differential operators with coefficients
in $\O_N$. The number $r$ is called the {\em rank\/} of $\L$.

A related notion is that of a {\em differential Lie algebra\/}
\cite{R1}--\cite{R4} (see also \cite{C}). This is a Lie algebra
structure on $\L=Y\tt\Kset^r$, where $Y$ is any commutative associative 
$H=\Kset[\d_1,\dots,\d_N]$-\difalg, given by \eqref{xeiyej} for
$x,y\in Y$, $P^{ij}_{kl}, Q^{ij}_{kl} \in Y\tt H$.
One can allow a universal enveloping algebra $H=\ue(\dd)$ 
($\dim\dd=N$) in place of 
$\Kset[\d_1,\dots,\d_N]$, cf.\ \cite{NW}.

Recall that for $H=\ue(\dd)$, $X=H^*$ is a commutative associative 
$H$-\difalg\ that can be identified with $\O_N$ for $N=\dim\dd$.
Moreover, the action of $H$ (and of $X\tt H$)
on $X$ is given by differential operators
in this identification.
Therefore a differential Lie algebra for $Y=X$ is the same as
a local Lie algebra.

Then the results of \seref{subanih} immediately imply:
\begin{proposition}\lbb{plocla}
Let $L = H\tt\Kset^r$ be a 
Lie \psalg\ which is a free $H$-module of rank $r$.
Let $Y$ be an $H$-bimodule which is a
commutative associative $H$-\difalg\ 
both for the left and for the right action of $H$ 
{\rm(}see \eqref{hxy}, \eqref{xyh}{\rm)}. 
Then $\A_Y L \simeq Y\tt\Kset^r$
is a differential Lie algebra.
In particular, $\A(L) = \A_X L$ is a local Lie algebra.
\end{proposition}

Note that the differential Lie algebras $\A_Y L$ that we get are with
``constant coefficients'': in \eqref{xeiyej} all 
$P^{ij}_{kl}, Q^{ij}_{kl} \in H$.

\subsection{Topology on the annihilation algebra}\lbb{ssubtop}
Now let us discuss the problem
of defining a topology on $\A(M)= X \tt_H M$ where
$M$ is any finite $H$-module.
Recall that $X$ has a decreasing filtration
$X=\fil_{-1}X\supset\fil_0X\supset\dotsm$
defined in \seref{sprelh}. We can use this filtration 
to construct an induced filtration on $\A(M)$
as follows. Choose a finite-dimensional (over $\Kset$) subspace 
$M_0$ of $M$ which generates $M$ over $H$, and set:
\begin{equation}\lbb{filam}
\fil_i \A(M) = \{ x \tt_H m \st x\in \fil_i X, \, m\in M_0\}.
\end{equation}
Note that since $H$ is cocommutative, its filtration satisfies \eqref{filh0},
hence $\bigcap \fil_i X = 0$. This implies:
\begin{equation}\lbb{filam2}
\bigcap_i \fil_i \A(M) = 0.
\end{equation}
The filtration \eqref{filam} will in general depend on the choice of $M_0$, 
but the topology induced by it will not, 
as any two such filtrations are equivalent by the next lemma.

\begin{lemma}\lbb{equivfiltr}
Let $M_0$ and $M'_0$ be two finite-dimensional subspaces of 
$M$ generating it over $H$, 
and let $\{ \fil_i \A(M) \}$, $\{ \fil'_i \A(M) \}$ 
be the corresponding filtrations on $\A(M)$. 
Then there exist integers $a, b$ such that 
$\fil_{i+a} \A(M) \subset \fil'_i \A(M) \subset \fil_{i+b} \A(M)$ 
for all values of $i$.
\end{lemma}
\begin{proof}
Let us choose bases of $M_0$ and $M'_0$, and let us fix expressions of 
elements from the first basis as $H$-linear combinations of elements from the 
second basis. Denote by $a$ the highest degree of the coefficients
of all these expressions.
Using \eqref{fils5}, we see that 
$\fil_i \A(M) \subset \fil'_{i-a} \A(M)$ 
for all $i$.
Repeating the same reasoning after switching the 
roles of $M_0$ and $M'_0$, we get 
$\fil'_i \A(M) \subset \fil_{i-b} \A(M)$ for some $b$ and all $i$.
\end{proof}


\begin{proposition}\lbb{plielc}
Let $H$ be a cocommutative Hopf algebra which satisfies \eqref{filh4}.

{\rm(i)}
If $M$ is a finite $H$-module, then $\A(M)$ is a linearly compact topological
vector space when provided with the filtration \eqref{filam}. 
The action of $H$ on $\A(M)$ is continuous if we endow $H$ with the discrete
topology. 

{\rm(ii)}
If $L$ is a finite Lie $H$-\psalg, then its annihilation algebra
$\A(L)$ is a linearly compact Lie $H$-\difalg, i.e., 
it is a linearly compact topological vector space
and both the Lie bracket and the action of $H$ are continuous.

A similar statement holds for representations and for
\as\ \psalgs\ as well.
\end{proposition}
\begin{proof}
(i) The linear-compactness follows from \prref{eqdefs}, since \eqref{filam}
a filtration by finite-codimensional
subspaces with trivial intersection and $\A(M)$ is complete with 
respect to this filtration. The continuity of the $H$-action follows
from \eqref{fils5}:
\begin{equation}\lbb{hamact}
\fil^i H \cdot \fil_j \A(M) \subset \fil_{j-i} \A(M)
\qquad\text{for all}\;\; i,j.
\end{equation}

(ii) It only remains to check that the Lie bracket of $\A(L)$ is continuous.
Let $L_0$ be a finite-dimensional (over $\Kset$) subspace 
of $L$ which generates it over $H$. For $a,b\in L_0$, we can write
\begin{displaymath}\lbb{psbra1}
[a*b] = \tsum_i\, (f_i\tt g_i) \tt_H e_i
\end{displaymath}
for some $f_i,g_i\in H$ and $e_i\in L_0$. Then the Lie bracket in $\A(L)$,
for $x,y\in X$, is given by:
\begin{displaymath}\lbb{psbra2}
[x\tt_H a, y\tt_H b] = \tsum_i\, (x f_i)(y g_i) \tt_H e_i.
\end{displaymath}

We can find a number $p$ such that
all coefficients $f_i,g_i\in H$ occurring in pseudobrackets of any elements
$a,b\in L_0$ belong to $\fil^p H$. Then equations \eqref{fils1},
\eqref{fils5} imply:
\begin{equation}\lbb{brafil}
[\fil_i \A(L), \fil_j \A(L)] \subset \fil_{i+j-s} \A(L)
\qquad\text{for all}\;\; i,j,
\end{equation}
where $s=2p$. This shows that the Lie bracket is continuous.
%
\end{proof}
\begin{lemma}\lbb{lssurj}
Let $H=\ue(\dd)\smash\kk[\Ga]$.
Then for any nonzero $h\in\dd$ and for every open
subspace $U$ of $\A(M)$ there is some $n$ such that $h^n U=\A(M)$. In
particular, each such $h$ acts surjectively on $\A(M)$.
\end{lemma}
\begin{proof}
Follows immediately from \leref{minusi}.
\end{proof}

\subsection{Growth of the annihilation algebra}\lbb{ssfiltr}
%
Let $M$ be a finite $H$-module. Then any choice of
a finite-dimensional subspace $M_0$ generating $M$ over $H$
provides $\M=\A(M)$ with a filtration
$\M_n := \fil_n X\tt_H M_0$. 

\begin{definition}\lbb{dgrowth}
For a filtered vector space 
$\M = \M_{-1} \supset \M_0\supset\dotsm$
we define its {\em growth\/}
$\gw\M$ to be $d$ if the function $n\mapsto \dim \M/\M_n$ 
can be bounded from above and below by polynomials of degree $d$.
\end{definition}

By \leref{equivfiltr}, a different choice
of $M_0$ would give a uniformly equivalent filtration of the same growth
as $\{\M_n\}$. Hence, we can speak of the growth of $\A(M)$ independently
of the choice of $M_0$. 


\begin{proposition}\lbb{growthisn}
Let $H = \ue(\dd)$ be the universal enveloping algebra of a finite-dimensional 
Lie algebra $\dd$, and $M$ be a finitely generated 
$H$-module. Then the growth of $\A(M)$ is equal to the dimension of $\dd$.
\end{proposition}
\begin{proof}
First of all, notice that we can assume $M$ is torsion-free,
since by \leref{torzerocoef}, $\A(M) \simeq \A(M / \Tor M)$ 
where $\Tor M$ is the torsion submodule of $M$.
The proof of the proposition is then based on \leref{rmapstofree} and
the following two lemmas.

\begin{lemma}\lbb{lemmad}
The map $\A(f)\colon  \A(M) \to \A(F)$ induced by the inclusion 
$f\colon M \injto F$ constructed in
\leref{rmapstofree} is uniformly continuous, i.e., 
for every $i$ we have{\rm:}
\begin{displaymath}
\fil_{i-a} \A(F) \subset \A(f)(\fil_i \A(M)) \subset \fil_{i+b} \A(F),
 \end{displaymath}
where $a$ and $b$ are independent of $i$.

The same is true for $\A(g)\colon  \A(F) \to \A(M)$ 
where $g$ is the embedding $g\colon hF \injto M$ from \leref{rmapstofree}.
\end{lemma}
\begin{proof}
Let us choose finite-dimensional vector subspaces
$F_0$ of $F$ generating $F$ over $H$, and $M_0$ of
$M$ generating $M$ over $H$ and containing $h F_0$. Let us also
choose a second finite-dimensional vector subspace $F'_0$ of $F$
containing $M_0$ and generating $F$ over $H$. We will denote
the filtrations induced by these subspaces by $\{\F_i\}$, $\{\M_i\}$,
and $\{\F'_i\}$, respectively.

Up to identifying $hF$ with $F$, we have constructed injective maps 
$F\xrightarrow{g} M \xrightarrow{f} F$ 
such that the composition $fg$ is a multiplication by $h$.
These maps induce maps 
$\A(F) \xrightarrow{\A(g)} \A(M) \xrightarrow{\A(f)} \A(F)$ which are 
surjective, as one can see by tensoring by $X$
and using that $\A(T)=0$ if $T$ is a torsion $H$-module
(see \leref{torzerocoef}).

The above maps are also continuous with respect to the common topology 
defined by any of the 
above constructed filtrations.
In fact, by construction, one has 
\begin{displaymath}
\A(g)(\F_i)\subset \M_i
\quad\text{and}\quad
\A(f)(\M_i) \subset \F'_i. 
\end{displaymath}

The second inclusion proves that $\A(f)(\M_i) \subset \F_{i+b}$
for some $b$ independent of $i$, because 
the filtrations $\{\F_i\}$ and $\{\F'_i\}$ are uniformly equivalent 
by \leref{equivfiltr}.

Applying $\A(f)$ to the first inclusion, we get
$\A(f) \A(g)(\F_i) \subset \A(f)(\M_i)$.
On the other hand, 
$\A(f) \A(g) = \A(fg) = h\tt_H\id_F$,
and \leref{minusi} implies that 
$\A(f) \A(g)(\F_i) = \F_{i-a}$ where $a$ is such that
$h\in \fil^a H$ but $h\not\in \fil^{a-1} H$.
Therefore $\F_{i-a} \subset \A(f)(\M_i)$ for all $i$.

A similar argument works for $\A(g)$.
\end{proof}
\begin{lemma}\lbb{lemmae}
If $\ph\colon\M\to\N$ is a surjective uniformly continuous map of filtered
modules, then $\gw \M \geq \gw \N$.
\end{lemma}
\begin{proof}
By assumption $\ph(\M_i) \subset \N_{i+b}$
for all $i$ and some $b$ independent of $i$.
The induced map 
$\M/\M_i \to \N/\N_{i+b}$ is a surjective map of finite-dimensional vector 
spaces. Hence $\gw \M \ge \gw \N$.
\end{proof}

Using the above lemmas, now we can complete the proof of \prref{growthisn}.
We have constructed an embedding $f\colon M \injto F$
of $M$ into a free $H$-module $F$, and we have shown that
the induced map $\A(f)\colon  \A(M) \to \A(F)$ is 
surjective and uniformly continuous. This implies $\gw\A(M) \ge \gw\A(F)$.
Similarly, the inclusion $hF \injto M$ gives us 
$\gw\A(F) = \gw\A(hF) \ge \gw\A(M)$.
Therefore, $\gw\A(M) = \gw\A(F)$.
It remains to note that, for any (nonzero) free $H$-module $F$ of finite rank,
one has $\gw \A(F) = \dim \dd$.
This follows from the fact that $\gw X = N =\dim \dd$
because $X\simeq\O_N$.
\end{proof}
\begin{remark}
The growth of a linearly compact Lie algebra $\L$
satisfying the descending chain condition can be defined as follows.
Take a fundamental subalgebra $A\subset\L$, and build a filtration of $\L$ by:
\begin{displaymath}
\L^A_0 = A, \quad \L^A_{i+1} = \{x\in \L^A_i \st [x,\L]\in \L^A_i\},
\;\; i\geq 0. 
\end{displaymath}
Taking $A' = \L_k^A$, we have: $\L^{A'}_i = \L^A_{k+i}$, hence replacing
$A$ by $A'$ does not change the growth. Therefore, by the Chevalley principle
\cite{G1}, the growth of this filtration does not depend on the choice of $A$.
We will denote this common growth by $\gw \L$.

Notice that all simple
linearly compact Lie algebras satisfy the descending chain condition, and
therefore have a well defined growth which equals $N$ for $W_N, S_N, H_N$ and
$K_N$, and $0$ for finite-dimensional Lie algebras.
\end{remark}

\section{Primitive Pseudoalgebras of Vector Fields}\lbb{svect}
In this section, $\dd$ will be a (finite-dimensional)
Lie algebra and $H=\ue(\dd)$ will be its 
universal enveloping algebra. 
As usual, we will identify $\dd$ with its image in $H$.
Then $X:=H^*$ is the 
algebra of formal power series on $\dd^*$, which is isomorphic 
as a topological algebra to $\O_N$ for $N=\dim\dd$.
In this section
we are going to define $H$-\psalg\ analogues of
the primitive linearly compact Lie algebras 
$W_N$, $S_N$, $H_N$, $K_N$, 
which will be called {\em primitive \psalgs\ of vector fields}.

\subsection{$\Wd$}\lbb{subw}
Let $Y$ be a commutative associative algebra on which $\dd$ acts
by derivations from the right (i.e., $Y$ is an $H$-\difalg).
One can define a left action of $Y\tt\dd$ on $Y$ using the right
action of $\dd$ on $Y$:
\begin{equation}\lbb{lieywd2}
(x\tt a)z = -x(za), \qquad x,z\in Y, \; a\in\dd.
\end{equation}
This will define a representation of $Y\tt\dd$ in $Y$ if
the Lie bracket of $\dd$ is extended to $Y\tt\dd$
by the formula
\begin{equation}\lbb{lieywd}
[x\tt a, y\tt b] = xy\tt [a,b] - x(ya)\tt b + (xb)y\tt a.
\end{equation}
In particular, for $Y=X=H^*$, this gives the Lie algebra
of all vector fields on $X$, which is isomorphic to
$W_N$ for $N=\dim\dd$.

Comparing \eqref{lieywd} with \eqref{ahxbhy}, we are led to define
the pseudoalgebra $\Wd=H\tt\dd$ with pseudobracket
\begin{equation}\lbb{wdbr*}
\begin{split}
[(f\tt a)*(g\tt b)] 
&= (f\tt g)\tt_H(1\tt [a,b]) 
\\
&- (f\tt ga)\tt_H(1\tt b) + (fb\tt g)\tt_H(1\tt a) .
\end{split}
\end{equation}
It is easy to check that $\Wd$ is indeed a Lie \psalg, and that
the Lie algebra $\A_Y \Wd$ defined in \seref{slalg} is isomorphic to
$Y\tt\dd$ with bracket defined by \eqref{lieywd}.
In a similar fashion, the module $Y$ over $Y\tt\dd$, 
defined by \eqref{lieywd2}, leads
to a structure of a $\Wd$-module on $H$:
\begin{equation}\lbb{wdac*}
(f\tt a)*g = -(f\tt ga)\tt_H 1.
\end{equation}


\subsection{Differential forms}\lbb{subforms}
We can think of $X=H^*$ as the space of functions on 
$\dd$, and of the elements of $X\tt\dd$ as vector fields.
Then the space of $n$-forms ($n=0,\dots,\dim\dd$) is 
\begin{displaymath}
\Om^n_X := \Hom_\Kset(\textstyle\bigwedge^n \dd, X) 
         \simeq X \tt \textstyle\bigwedge^n \dd^*.
\end{displaymath}
It is convenient to extend the elements
$\om\in\Om^n_X$ to functions from $\bigwedge^n (X\tt\dd)$ to $X$,
polylinear over $X$:
\begin{displaymath}\lbb{omxa}
\om(x_1\tt a_1 \wedge \dots \wedge x_n\tt a_n)
= x_1 \dotsm x_n \, \om(a_1 \wedge \dots \wedge a_n),
\end{displaymath}
so that 
\begin{displaymath}
\Om^n_X = \Hom_X(\textstyle\bigwedge^n (X\tt\dd), X).
\end{displaymath}

We view $X$ as a left ($X\tt\dd$)-module using the right action of $\dd$,
see \eqref{lieywd2}. 
There is a differential
$\di\colon\Om^n_X\to\Om^{n+1}_X$ satisfying $\di^2=0$; this is just the usual 
differential for the cohomology of $\dd$ with coefficients in $X$
where $X$ is viewed as a right $\dd$-module:
\begin{displaymath}\lbb{doma}
\begin{split}
(\di\om)(a_1 \wedge \dots \wedge a_n)
&= \sum_{i<j} (-1)^{i+j} \om([a_i,a_j] \wedge a_1 \wedge \dots \wedge 
\what a_i \wedge \dots \wedge \what a_j \wedge \dots \wedge a_n)
\\
&+ \sum_i (-1)^i \om(a_1 \wedge \dots \wedge \what a_i \wedge \dots 
\wedge a_n) \cdot a_i.
\end{split}
\end{displaymath}
The following analogue of the Poincar\'e Lemma is very useful.

\begin{lemma}\lbb{lpoin}
The complex $(\Om_X^\bullet, \di)$ is acyclic, i.e.,
its $n$-th cohomology is trivial for $n>0$
and $1$-dimensional for $n=0$.
\end{lemma}
\begin{proof}
It is well known that $\coh^n(\dd,U(\dd)^*) \simeq \coh_n(\dd,U(\dd))^*$
is trivial for $n>0$ and $1$-dimensional for $n=0$; see e.g.\ \cite{F}.
\end{proof}

For a vector field $A\in X\tt\dd$, we have the contraction operator
$\iota_A\colon\Om^n_X\to\Om^{n-1}_X$ given by:
\begin{displaymath}\lbb{ioaom}
(\iota_A\om)(a_1 \wedge \dots \wedge a_{n-1})
= \om(A \wedge a_1 \wedge \dots \wedge a_{n-1}).
\end{displaymath}
We define the Lie derivative $L_A\colon\Om^n_X\to\Om^n_X$
by Cartan's formula $L_A = \di \iota_A + \iota_A \di$.
Explicitly, for $x\tt a \in X\tt\dd$, we have:
\begin{equation}\lbb{lxaom}
\begin{split}
(L_{x\tt a}\om)(&a_1 \wedge \dots \wedge a_n)
= -x(\om(a_1 \wedge \dots \wedge a_n) \cdot a)
\\
&+ \sum_i (-1)^i (x \cdot a_i) \, \om(a \wedge a_1 \wedge \dots \wedge 
\what a_i \wedge \dots \wedge a_n)
\\
&+ \sum_i (-1)^i x \, \om([a,a_i] \wedge a_1 \wedge \dots \wedge \what a_i 
\wedge \dots \wedge a_n).
\end{split}
\end{equation}
The Lie derivative provides each $\Om^n_X$ with the structure of a module
over the Lie algebra of vector fields $X\tt\dd$.

For $n=0$, $\Om^0_X = X$ and this is the usual action \eqref{lieywd2}
of $X\tt\dd$ on $X$.
When $n=N=\dim\dd$, we have $\Om^N_X = X v_0$ where 
$v_0\in\bigwedge^N\dd^*$, $v_0\ne0$ is a volume form.
An easy calculation shows that
\begin{equation}\lbb{xav0}
L_{x\tt a}(y v_0) = -\bigl( (xy)(a + \tr\ad a) \bigr) v_0,
\qquad x,y \in X, \; a\in\dd.
\end{equation}

\subsection{Pseudoforms}\lbb{spsforms}
The module $\Om^n_X$ over the Lie algebra $X\tt\dd$ leads
to a module $\Om^n(\dd)$ over the Lie \psalg\ $\Wd$ which we now define.
We let
\begin{displaymath}
\Om^n(\dd) = H\tt \textstyle\bigwedge^n \dd^*,
\quad \Om(\dd) = \textstyle\bigoplus_{n=0}^N \Om^n(\dd)
\quad (N=\dim\dd).
\end{displaymath}
The elements of $\Om(\dd)$ are called {\em pseudoforms}.

$\Om^n(\dd)$ is a free $H$-module, so that 
$\A(\Om^n(\dd)) = X \tt_H \Om^n(\dd) \simeq X \tt \bigwedge^n \dd^* 
= \Om^n_X$.
The action of $\Wd = H\tt\dd$ on $\Om^n(\dd)$
is obtained by comparing \eqref{ahxbhy} with
\eqref{lxaom}. To write an explicit formula, we identify
$\Om^n(\dd)$ with the space of linear maps from $\bigwedge^n \dd$ to $H$,
and $(H\tt H) \tt_H \Om^n(\dd)$ with the space of linear maps from 
$\bigwedge^n \dd$ to $H\tt H$. Then for $f\tt a\in\Wd$, $w\in\Om^n(\dd)$,
and $a_i\in\dd$, we have:
\begin{equation}\lbb{fa*om}
\begin{split}
((f\tt a)*w)(&a_1 \wedge \dots \wedge a_n)
= -f \tt w(a_1 \wedge \dots \wedge a_n) \, a
\\
&+ \sum_i (-1)^i f a_i \tt w(a \wedge a_1 \wedge \dots \wedge \what a_i 
\wedge \dots \wedge a_n)
\\
&+ \sum_i (-1)^i f \tt w([a,a_i] \wedge a_1 \wedge \dots \wedge \what a_i 
\wedge \dots \wedge a_n).
\end{split}
\end{equation}

When $n=0$, $\Om^0(\dd)=H$, and we recover \eqref{wdac*}.
In the other extreme case, when $n=N:=\dim\dd$, $\Om^N(\dd)=Hv_0$
is again a free $H$-module of rank one, where $v_0\in\bigwedge^N \dd^*$,
$v_0 \ne0$
is any volume form on $\dd$. We have (cf.\ \eqref{xav0}):
\begin{equation}\lbb{fa*v}
(f\tt a)*v_0 = -\bigl( f(a + \tr\ad a)\tt1 + f\tt a \bigr) \tt_H v_0.
\end{equation}

Define polylinear maps 
$*_\iota\in\Lin(\{\Wd,\Om^n(\dd)\}, \Om^{n-1}(\dd))$ by
\begin{equation}\lbb{*iota}
((f\tt a)*_\iota w)(a_1 \wedge \dots \wedge a_{n-1})
= f \tt w(a \wedge a_1 \wedge \dots \wedge a_{n-1}).
\end{equation}
Also define a differential $\di\colon\Om^n(\dd)\to\Om^{n+1}(\dd)$ by
\begin{displaymath}
\begin{split}
&\begin{split}
(\di w)(a_1 \wedge \dots \wedge a_{n+1})
&= \sum_{i<j} (-1)^{i+j} w([a_i,a_j] \wedge a_1 \wedge \dots \wedge 
\what a_i \wedge \dots \wedge \what a_j \wedge \dots \wedge a_{n+1})
\\
\lbb{dw}
&+ \sum_i (-1)^i w(a_1 \wedge \dots \wedge \what a_i \wedge \dots 
\wedge a_{n+1}) \, a_i
\qquad\text{if}\;\; n\ge1,
\end{split}
\\
&(\di w)(a) = -wa \qquad\text{if}\;\; n=0,
\end{split}
\end{displaymath}
so that $\di$ is $H$-linear and $\di^2=0$.
For any pseudoform $w$ and $x\in X$, we have a differential form
$x\tt_H w$ and the relation
$\di(x\tt_H w) = x\tt_H \di w$.

\begin{remark}\lbb{rpoin}
The $n$-th cohomology of the complex $(\Om^\bullet(\dd), \di)$ 
is equal to $\coh^n(\dd,U(\dd))$. This space is trivial for $n\ne N=\dim\dd$
and $1$-dimensional for $n=N$. This follows from the Poincar\'e
duality $\coh^n(\dd,U(\dd)) \simeq \coh^{N-n}(\dd,U(\dd)^*)$;
see e.g.\ \cite{F}.
\end{remark}

We have the following analogue of Cartan's formula for the action
of $\Wd$ on $\Om(\dd)$:\begin{equation}\lbb{*cartan}
\al*w = ((\id\tt\id)\tt_H \di)(\al*_\iota w) + \al*_\iota (\di w)
\in (H\tt H)\tt_H \Om(\dd).
\end{equation}
This implies that the action of $\Wd$ commutes with $\di$:
\begin{equation}\lbb{*di*}
\al*(\di w) = ((\id\tt\id)\tt_H \di)(\al* w).
\end{equation}
We note that the maps $\al*_\iota$ anticommute with each other:
\begin{equation}\lbb{i*i*}
\al*_\iota (\be*_\iota w) + ((\si\tt\id)\tt_H \id) \, \be*_\iota (\al*_\iota w)
=0
\end{equation}
for $\al,\be\in\Wd$, $w\in\Om(\dd)$.

The wedge product on $\bigwedge^\bullet \dd^*$ can be extended
to a pseudoproduct $*$ on $\Om(\dd) = H\tt\bigwedge^\bullet \dd^*$,
so that it becomes a current \psalg.
Then it is easy to check that for $\al\in\Wd$, $v\in\bigwedge^m \dd^*$,
$w\in\bigwedge^n \dd^*$, one has:
\begin{align}
\lbb{a*wedge}
\al*(v*w) &= (\al*v)*w + ((\si\tt\id)\tt_H\id) \, v*(\al*w),
\intertext{and similarly,}
\lbb{a*wedge2}
\al*_\iota (v*w) &= (\al*_\iota v)*w 
+ (-1)^m ((\si\tt\id)\tt_H\id) \, v*(\al*_\iota w).
\end{align}
This can be interpreted as
saying that $\al*$ and $\al*_\iota$ are superderivations of $\Om(\dd)$,
see \exref{eomdd} below.

\subsection{$\Sd$}\lbb{subs}
The divergence of a vector field $\tsum_i x_i\tt a_i \in X\tt\dd$
is defined by
$\Div(\tsum_i x_i\tt a_i) = \tsum_i x_i a_i \in X$.
Then one easily checks
\begin{equation}\lbb{div2}
\Div([A,B]) = A \cdot \Div(B) - B \cdot \Div(A), \qquad A,B\in X\tt\dd,
\end{equation}
so that the divergence zero vector fields form a Lie subalgebra 
$S_N$ of $W_N$.
Let $\chi$ be a trace form on $\dd$, i.e., a linear functional
from $\dd$ to $\Kset$ which vanishes on $[\dd,\dd]$. Then we can define
\begin{displaymath}\lbb{divchi}
\Div^\chi(\tsum_i x_i\tt a_i) := \tsum_i x_i (a_i + \chi(a_i)),
\end{displaymath}
which still satisfies \eqref{div2}.

\begin{remark}\lbb{rdivchi}
Let $\chi$ be as above, and let $\psi=\chi-\tr\ad$, which is again
a trace form on $\dd$. We can consider $\psi$ as an element of 
$\Om^1_X = X\tt\dd^*$;
then $\di\psi=0$ and by \leref{lpoin} we have $\psi=-\di z$ for some
$z\in X$. This means that $\psi(a)=za$ for all $a\in\dd$. Let $y=e^z$;
then $ya=y\psi(a)$ for $a\in\dd$. 
Consider the volume form $v=yv_0$, where
$v_0\in\bigwedge^N\dd^*$, $v_0\ne0$. Equation \eqref{xav0} gives
\begin{equation}\lbb{lavdiva}
L_A v = - \Div^\chi(A) v
\qquad\text{for}\quad A \in X\tt\dd.
\end{equation}
Therefore, the Lie algebra of vector fields $A$ with $\Div^\chi(A) = 0$
coincides with the Lie algebra $S_N(v)$ of vector fields 
annihilating the volume form $v$.
\end{remark}

Using the notation
$\al_x \equiv x\tt_H\al \in\A(\Wd) \simeq X\tt\dd$ for $\al\in\Wd$,
$x\in X$, we find for $\al=h\tt a$:
\begin{displaymath}
\al_x = x \tt_H (h\tt a) = xh \tt_H (1\tt a)
\equiv xh \tt a \in X\tt\dd,
\end{displaymath}
hence, $\Div^\chi(\al_x)=xh(a+\chi(a))$.
Define the divergence operator $\Div^\chi\colon\Wd\to H$ by the formula:
\begin{equation}\lbb{div*}
\Div^\chi(\tsum_i \, h_i\tt a_i) = \tsum_i \, h_i (a_i + \chi(a_i)).
\end{equation}
Then we have:
\begin{equation}\lbb{divchial}
\Div^\chi(\al_x) = x \cdot \Div^\chi\al \qquad\text{for}\;\; \al\in\Wd, x\in X.
\end{equation}
Since $\Div^\chi$ is $H$-linear, we can define 
\begin{displaymath}
\Div^\chi_2\colon (H\tt H)\tt_H\Wd
\xrightarrow{\id\tt_H\Div^\chi} (H\tt H)\tt_H H 
\isoto H\tt H.
\end{displaymath}
Similarly to \eqref{div2}, one has:
\begin{equation}\lbb{div2*}
\Div^\chi_2([\al*\be]) = (\Div^\chi\al\tt1)\si(\be) - (1\tt\Div^\chi\be)\al, 
\qquad \al,\be\in\Wd,
\end{equation}
where $\si\colon H\tt H \to H\tt H$ is the transposition.

Equation \eqref{div2*} implies that 
\begin{equation}\lbb{sd}
\Sd := \{ \al\in\Wd \st \Div^\chi\al = 0 \}
\end{equation}
is a subalgebra of the Lie \psalg\ $\Wd$.
By Eq.~\eqref{divchial} and \reref{rdivchi}, its annihilation algebra 
\begin{equation}
\A(\Sd) = \{ A \in W_N \st \Div^\chi A = 0 \} \simeq S_N.
\end{equation}
The rank of $\Sd$ as an $H$-module is $N-1$; 
however, it is free only for $N=2$.

\begin{proposition}\lbb{psd}
$\Sd$ is generated over $H$ by elements
\begin{equation}\lbb{achib}
e_{ab} := (a+\chi(a))\tt b - (b+\chi(b))\tt a - 1\tt [a,b]
\qquad\text{for}\quad a,b\in\dd.
\end{equation}
These elements satisfy $e_{ab} = -e_{ba}$ and the relations 
{\rm(}for $\chi=0${\rm):}
\begin{equation}\lbb{eabrel}
a e_{bc} + b e_{ca} + c e_{ab} = e_{[a,b],c} + e_{[b,c],a} + e_{[c,a],b}.
\end{equation}
For $\chi=0$, their pseudobrackets are given by{\rm:}
\begin{align}
\notag
[e_{ab} * e_{cd}] = (a \tt d) \tt_H e_{bc} &+ (b \tt c) \tt_H e_{ad}
             - (a \tt c) \tt_H e_{bd} - (b \tt d) \tt_H e_{ac}
\\
\notag
             &+ (a \tt 1) \tt_H e_{b,[c,d]} - (b \tt 1) \tt_H e_{a,[c,d]}
\\
\lbb{eab*ecd}
             &- (1 \tt c) \tt_H e_{d,[a,b]} + (1 \tt d) \tt_H e_{c,[a,b]}
\\
\notag
             &- (1 \tt 1) \tt_H e_{[a,b],[c,d]}.
\end{align}
For arbitrary $\chi$, replace everywhere in \eqref{eabrel}, \eqref{eab*ecd}
all $h\in\dd$ with $h+\chi(h)$.
\end{proposition}

\begin{remark}\lbb{reab}
Equation \eqref{eab*ecd} implies that for $\chi=0${\rm:}
\begin{equation}\lbb{eab*eab}
\begin{split}
[e_{ab} * e_{ab}] 
&= (b \tt a - a \tt b \bigr) \tt_H e_{ab}
\\
&+ \bigl( 1 \tt b - b \tt 1 \bigr) \tt_H e_{a,[a,b]}
+ \bigl( a \tt 1 - 1 \tt a \bigr) \tt_H e_{b,[a,b]}.
\end{split}
\end{equation}
{\rm(}Again, for any $\chi$, replace $a,b$ with $a+\chi(a), b+\chi(b)$.{\rm)}
In particular, when the elements $a,b$ span a Lie algebra,
$H e_{ab}$ is a Lie \psalg.
\end{remark}

In the proof of \prref{psd} we are going to use the following lemma.

\begin{lemma}\lbb{decomp}
Let $H = \ue(\dd)$, and let $\{\d_1, \dots, \d_N\}$ be a basis of $\dd$.
If elements $h_i \in \fil^d H$ are such that $\sum_i h_i \d_i \in \fil^d H$, 
then there exist $f_{ij} \in \fil^{d-1}H$ such that
\begin{displaymath}
\sum_i h_i \tt \d_i = \sum_{i,j} (f_{ij} \tt 1)(\d_i\tt\d_j - \d_j\tt\d_i)
\;\mod \fil^{d-1} H \tt\dd.
\end{displaymath}
\end{lemma}
\begin{proof}
The proof is by induction on the number of $h_i$ not contained in 
$\fil^{d-1}H$, the basis of induction being trivial. Consider
 $\sum_{i=1}^n h_i \d_i \in \fil^d H$,
with all $h_i \notin \fil^{d-1}H$. 
We can write
$h_i = f_i \d_1 + k_i$ so that $k_i \in
\fil^d H$ is a linear combination of Poincar\'e--Birkhoff--Witt 
basis elements of $H$ not containing
$\d_1$ in their expression, and $f_i \in \fil^{d-1} H$.
Then:
\begin{align*}
\sum_{i=1}^n h_i \d_i 
&= h_1 \d_1 + \sum_{i=2}^n \bigl( f_i\d_1\d_i + k_i \d_i \bigr)
\\
&= \Bigl( h_1 + \sum_{i=2}^n f_i \d_i \Bigr) \d_1 + \sum_{i=2}^n k_i\d_i
  + \sum_{i=2}^n f_i [\d_1,\d_i] .
\end{align*}
Since the third summand in the right-hand side belongs to $\fil^d H$,
it follows that the first and second
summands lie in $\fil^d H$ too.
This implies:
$h_1 + \sum_{i=2}^n f_i \d_i
\in \fil^{d-1}H$. Hence 
\begin{displaymath}
\sum_{i=1}^n h_i \tt \d_i = 
\sum_{i=2}^n \bigl( f_i\d_1 \tt \d_i - f_i\d_i \tt \d_1 \bigr)
+ \sum_{i=2}^n k_i\tt \d_i 
\;\mod \fil^{d-1} H \tt\dd,
\end{displaymath}
and we can apply the inductive assumption.
\end{proof}

\begin{proof}[Proof of \prref{psd}]
First of all, it is easy to check that the elements \eqref{achib}
indeed belong to $\Sd$. Equation \eqref{eabrel} is easy, and the computation 
of the pseudobrackets is straightforward
using \eqref{wdbr*}, reformulated as
\begin{equation}
\begin{split}
[(1\tt a) * (1\tt b)]
&= \bigl( (a+ \chi(a))\tt 1 \bigr) \tt_H (1 \tt b) 
\\
&- \bigl( 1 \tt (b+\chi(b)) \bigr) \tt_H (1 \tt a) 
 - (1 \tt 1)\tt_H e_{ab}.
\end{split}
\end{equation}

Now let us consider an element $\al =\sum_i h_i\tt \d_i \in \Sd$, $h_i\in H$.
We will prove that $\al$ can be expressed as $H$-linear combination of the
above elements \eqref{achib}
by induction on the maximal degree $d$ of the $h_i$.
Since $\al\in \Sd$, then $\sum_i h_i(\d_i + \chi(\d_i)) = 0$, hence
$\sum_i h_i\d_i \in \fil^d H$.

By \leref{decomp}, 
we can find elements $f_{ij} \in \fil^{d-1}H$ such that 
\begin{displaymath}
\al = \sum_{i,j} (f_{ij}\tt 1)(\d_i \tt \d_j - \d_j \tt \d_i)
\;\mod \fil^{d-1} H \tt\dd.
\end{displaymath}
Therefore the difference
\begin{displaymath}
\al - \sum_{i,j} (f_{ij}\tt 1) \bigl( (\d_i+\chi(\d_i)) \tt \d_j - 
(\d_j+\chi(\d_j)) \tt \d_i - 1\tt[\d_i,\d_j] \bigr)
\end{displaymath}
still lies in $\Sd$ and
its first tensor factor terms have degree strictly less than $d$.
By inductive assumption, we are done.
\end{proof}
\begin{remark}\lbb{rrepwd}
{\rm(i)}
Let, as before, $\chi\in\dd^*$ be such that $\chi([\dd,\dd])=0$.
For any $\la\in\Kset$, let $V_{\la,\chi}=Hv$ be a free $H$-module of rank $1$
with the following action of $\Wd$ on it:
\begin{equation}\lbb{mlachi}
\al*v = (\la\Div^\chi\al \tt 1 - \al) \tt_H v.
\end{equation}
Using \eqref{div2*}, it is easy to check that this indeed
defines a representation of $\Wd$.

For $\la=0$ we get the action \eqref{wdac*}, while for $\la=-1$, $\chi=\tr\ad$
we get \eqref{fa*v}. One can show that all representation of $\Wd$ on a free
$H$-module of rank $1$ are given by 
\begin{equation}
(1 \tt a)* v = \bigl( (\la a + \chi'(a)) \tt 1 - 1 \tt a \bigr) \tt_H v,
\end{equation}
where $a\in \dd$, $\la \in \kk$ and $\chi'$ is a trace form on $\dd$. 
This can be rewritten as in \eqref{mlachi}, for $\chi = \chi'/\la$ 
whenever $\la \neq 0$.

{\rm(ii)}
More generally, let $M$ be any $\Wd$-module, equipped with a
compatible action of $H=\Cur\Kset$. Here $H=\Cur\Kset$ is the \as\
\psalg\ with a pseudoproduct $f*g = (f\tt g)\tt_H 1$, and compatibility
of the actions of $\Wd$ and $H$ means that
\begin{equation}\lbb{wdhcomp}
\al*(h*m) - ((\si\tt\id)\tt_H\id)\, h*(\al*m) = (\al*h)*m
\end{equation}
for $\al\in\Wd$, $h\in H$, $m\in M$, where $\al*h = -(1\tt h)\al\tt_H 1$
is the action \eqref{wdac*} of $\Wd$ on $H$.

Then, for any $\la,\chi$ as above,
\begin{equation}\lbb{mlachi2}
\al *_{\la,\chi} m = \la(\Div^\chi\al)*m + \al*m
\end{equation}
is an action of $\Wd$ on $M$.
\end{remark}

\subsection{Pseudoalgebras of rank $1$}\lbb{subr1wd}
All Lie \psalgs\ that are free of rank one over $H$
were described by \prref{prank1} and \leref{xcybe1}. 
The next lemma implies that all of them are subalgebras of $\Wd$.

\begin{lemma}\lbb{lr1wd}
Let $\al\in H\tt H$ be a solution of equations \eqref{al1}, \eqref{al2}.
Write $\al=r+s\tt1-1\tt s$ with a skew-symmetric $r\in\dd\tt\dd$ and
$s\in\dd$, as in \leref{xcybe1}.
Consider $e=-r+ 1\tt s\in H\tt\dd$ as an element of $\Wd$. 
Then $[e*e]=\al\tt_H e$ in $\Wd$.
\end{lemma}
\begin{proof}
Straightforward computation, using the definition \eqref{wdbr*}
and equations \eqref{cybe3}, \eqref{cybe4}.
\end{proof}

Let us study equations (\ref{cybe3}, \ref{cybe4}) in more detail.
We can write 
\begin{equation}\lbb{raibi}
r=\tsum_i\, (a_i \tt b_i - b_i \tt a_i)
\end{equation}
for some linearly independent $a_i,b_i\in\dd$. 
Denote by $\dd_1$ their linear span, and let $\dd_0 = \dd_1 + \Kset s$.

\begin{lemma}\lbb{lrs1}
$\dd_0$ is a Lie subalgebra of $\dd$, and $\dd_1$ is $\ad s$-invariant.
Moreover, $[a_i,a_j], [b_i,b_j], [a_i,b_j]$ and $[a_i,b_i]+s$
belong to $\dd_1$ for $i\ne j$.
\end{lemma}
\begin{proof}
Similar to that of Proposition 2.2.6 in \cite{CP}.
If $\dd_1=\dd$, there is nothing to prove.
Let $\{c_j\}$ be elements that complement $\{a_i,b_i\}$
to a basis of $\dd$. If $s$ is not in $\dd_1$ we take it to
be one of the $c_j$'s. 

Write out \eqref{cybe3} as
\begin{displaymath}
\tsum_i\, \bigl( [a_i,s] \tt b_i - [b_i,s] \tt a_i 
  + a_i \tt [b_i,s] - b_i \tt [a_i,s] \bigr) = 0.
\end{displaymath}
Now, if $[a_i,s]$ involves some $c_j$'s, there is no way
to cancel out the terms $c_j \tt b_i$. This proves that
$[s,\dd_1] \subset\dd_1$.

Similarly, \eqref{cybe4} reads
\begin{align*}
\tsum_{i,j}\, &\bigl( [a_i,a_j] \tt b_i \tt b_j 
- [b_i,a_j] \tt a_i \tt b_j 
+ [b_i,b_j] \tt a_i \tt a_j 
- [a_i,b_j] \tt b_i \tt a_j +\text{cyclic} \bigr)
\\
&+ \tsum_i\, \bigl( a_i \tt b_i \tt s - b_i \tt a_i \tt s
+\text{cyclic} \bigr) = 0.
\end{align*}
If, for example, $[a_i,a_j]$ involves $c_k$'s, then the
terms $c_k \tt b_i \tt b_j$ cannot be cancelled. 
Therefore $[a_i,a_j] \in\dd_1$.
If $[a_i,b_j]$ involves $c_k$'s, then the terms $c_k \tt b_i \tt a_j$
can be cancelled only with terms coming from $s\tt r$.
This shows that $[a_i,b_j] + \de_{ij} s \in\dd_1$.
\end{proof}

The universal enveloping algebra $H_0=U(\dd_0)$
is a Hopf subalgebra of $H=U(\dd)$. Since $\al\in H_0 \tt H_0$,
we can consider the Lie \psalg\ $H_0 e$ with pseudobracket
$[e*e]=\al\tt_{H_0}e$. Then our \psalg\ $He$ is a current \psalg\
over $H_0 e$.

Clearly, $\dd_1$ is even dimensional. There are two
cases which are treated in detail in the next two subsections:
when $\dd_0=\dd_1$ and when $\dd_0=\dd_1\oplus\Kset s$.
They give rise to Lie \psalgs\ $\Hd,\Kd$
whose annihilation Lie algebras are of
hamiltonian and contact type, respectively.
The following theorem summarizes some of the 
results of Sections \ref{subrank1} and \ref{subr1wd}--\ref{subk}.

\begin{theorem}\lbb{trankone}
Any Lie \psalg\ which is free of rank one is either abelian or
isomorphic to a current
\psalg\ over one of the Lie \psalgs\ $\Hd$, $\Kd$ defined in Sections
{\rm\ref{subh}}, {\rm\ref{subk}}, respectively.
\end{theorem}

\subsection{$\Hd$}\lbb{subh}
This is defined as a Lie $H$-\psalg\ of rank $1$ (see \seref{subr1wd})
corresponding to a solution $(r,s)$ of equations \eqref{cybe3}, \eqref{cybe4}
with a nondegenerate $r\in\dd\wedge\dd$ (i.e., $\dd_1=\dd$),
in which case $N=\dim\dd$ is even. The parameters $\chi$ and $\omega$ are 
defined as follows.


Since $r$ is nondegenerate, the linear map $\dd^*\to\dd$ induced by it
is invertible; its inverse gives rise to a $2$-form
$\om\in\bigwedge^2\dd^*$. 
Explicitly, if $r=\tsum r^{ij} \d_i\tt\d_j$
where $\{\d_i\}$ is a basis of $\dd$, then
$\om(\d_i\wedge\d_j)=\om_{ij}$
is the matrix inverse to $r^{ij}$.
We also define a $1$-form $\chi := \iota_s \om \in\dd^*$.

Conversely, given a nondegenerate skew-symmetric 
$2$-form $\om$ and a $1$-form $\chi$,
we can define uniquely
$r\in\dd\wedge\dd$ as the dual to $\om$ and $s\in\dd$
so that $\chi = \iota_s \om$.

\begin{lemma}\lbb{lomchi}
When $r\in\dd\wedge\dd$ is nondegenerate,
equations \eqref{cybe3}, \eqref{cybe4} are equivalent to the following
identities for the above-defined $\om,\chi:$
\begin{align}
\lbb{omchi1}
\di\om + \chi\wedge\om &= 0,
\\
\lbb{omchi2}
\di\chi &= 0,
\end{align}
which simply mean that $\om$ is a $2$-cocycle for $\dd$ in the $1$-dimensional
$\dd$-module defined by $\chi$.
This establishes a one-to-one correspondence between solutions $(r,s)$
of \eqref{cybe3}, \eqref{cybe4} with nondegenerate $r$
and solutions $(\om,\chi)$ of 
\eqref{omchi1}, \eqref{omchi2} with nondegenerate $\om$.
\end{lemma}
\begin{proof}
Let us write $[\d_i,\d_j] = \tsum c_{ij}^k \d_k$
and $s = \tsum s^k \d_k$ (summation over repeated indices).
Then \eqref{cybe4} is equivalent to
\begin{equation}\lbb{omchi3}
\bigl(\tsum r^{ij} r^{kl} c_{ik}^m + r^{mj} s^l\bigr) +\text{cyclic} = 0,
\end{equation}
where ``cyclic'' means summing over cyclic permutations of the
indices $m,j,l$. Multiply this equation by $\om_{jn}\om_{lp}\om_{mq}$
and sum over $m,j,l$. Using that $\sum r^{ij}\om_{jn} = \de^i_n$,
we get
\begin{equation}\lbb{omchi4}
\bigl(\tsum c_{np}^m \om_{mq} + 
\tsum s^l \om_{lp} \om_{nq}\bigr) +\text{cyclic} = 0,
\end{equation}
where now the cyclic permutations are over $n,p,q$. This is exactly
Eq.~\eqref{omchi1}. Conversely, multiplying \eqref{omchi4} by
$r^{in}r^{jp}r^{kq}$ and summing over $n,p,q$,
we get \eqref{omchi3}.

Similarly, since $[s,\dd_1]\subset\dd_1$, we can write 
$[s,\d_i] = \tsum_k\, c_i^k \d_k$.
Then \eqref{cybe3} is equivalent to
\begin{equation}\lbb{omchi5}
\tsum r^{ij} c_i^k + \tsum r^{kl} c_l^j = 0,
\end{equation}
which after multiplying by $\om_{jm}\om_{kn}$ and summing over $j,k$
becomes
\begin{equation}\lbb{omchi6}
\tsum c_m^k \om_{kn} + \tsum c_n^j \om_{mj} = 0,
\end{equation}
or $L_s\om=0$. Conversely, \eqref{omchi6} gives \eqref{omchi5}
after multiplying by $r^{pm}r^{qn}$ and summing over $m,n$.

Now start with a solution $(r,s)$ of (\ref{cybe3}, \ref{cybe4}).
Above we have deduced \eqref{omchi1} and $L_s\om=0$.
On the other hand, since $\iota_s\chi=0$, we have
$\iota_s(\chi\wedge\om)=0$, and \eqref{omchi1} implies
$\iota_s\di\om=0$. Together with $L_s\om=0$ this gives
$\di\iota_s\om=0$, which is \eqref{omchi2}.

If we start with a solution $(\om,\chi)$ of 
(\ref{omchi1}, \ref{omchi2}),
the above arguments can be inverted to show that $L_s\om=0$,
and we get (\ref{cybe3}, \ref{cybe4}).
\end{proof}

In the basis $\{a_i,b_i\}$ of $\dd$ we have \eqref{raibi}
and $\om(a_i\wedge b_i) = -\om(b_i\wedge a_i) = -1$, 
all other values of $\om$ are zero.
For $e= -r + 1\tt s$ and any $x\in X$, 
the element $e_x := x\tt_H e$ of the annihilation algebra
$\A(\Wd) \simeq X\tt\dd$ is equal to
$-\sum (x a_i\tt b_i - x b_i\tt a_i) + x\tt s$,
and it is easy to check that 
\begin{equation}\lbb{omexa}
\om(e_x \wedge a) = x (-a+\chi(a)), \qquad a\in\dd.
\end{equation}
Since $\di\chi=0$, \leref{lpoin} implies that $\chi=\di y$
for some $y\in\Om^0_X=X$, i.e., $\chi(a)=-ya$.
Then $\ti\om := e^y\om$ satisfies $\ti\om(e_x \wedge a) = -(x e^y)a$
for any $x\in X$, $a\in \dd$. This is equivalent to
$\iota_{e_x} \ti\om = \di(x e^y)$.
Moreover, \eqref{omchi1} implies $\di\ti\om = 0$.
Therefore, $L_{e_x} \ti\om = 0$, and we have the following
proposition.

\begin{proposition}\lbb{phomchi}
Let $\Hd := He$ be a Lie $H$-\psalg\ of rank $1$ 
corresponding to a solution $(r,s)$ of equations \eqref{cybe3}, \eqref{cybe4}
with a nondegenerate $r\in\dd\tt\dd$. Define the $2$-form
$\ti\om$ as above. Then $\ti\om$ is a symplectic form, and
the subalgebra $X\tt_H \Hd$
of $X\tt_H \Wd \simeq X\tt\dd$
is the Lie algebra $H_N(\ti\om)$ of vector fields annihilating $\ti\om$
$($which is isomorphic to $H_N)$.
\end{proposition}
\begin{proof}
It remains to show that, conversely, any vector field that preserves
the form $\ti\om$ is equal to $e_x$ for some $x\in X$. 
Indeed, let $A\in X\tt\dd$ be such that $L_A\ti\om=0$. 
Since $\di\ti\om=0$ and $\ti\om=e^y\om$,
this is equivalent to $\di(e^y\iota_A\om)=0$
which implies $e^y\iota_A\om = \di z$ for some $z\in X$.
In other words, $e^y\om(A\wedge a) = -za$ for any $a\in\dd$.
Using $\chi(a)=-ya$, we get $\om(A\wedge a) = x (-a+\chi(a))$
for $x=e^{-y}z$. This, together with \eqref{omexa}, implies $A=e_x$
since the $2$-form $\om$ is nondegenerate.
\end{proof}
\begin{remark}\lbb{hdsubsd}
Let $r\in\dd\tt\dd$ be given by \eqref{raibi}, and let $x=\sum_i [a_i, b_i]$,
$\phi = -\chi + \iota_x\om = \iota_{x-s}\om$.
Then it is easy to check that $\Div^\phi(-r + 1\tt s) = 0$, so we have:
$\Hd \subset S(\dd, \phi)$.
\end{remark}
\begin{example}\lbb{eybd2}
Let the Lie algebra $\dd$ be $2$-dimensional with basis $\{a,b\}$
and commutation relations $[a,b] = \lambda b$.
Then up to multiplication by a scalar, all nondegenerate solutions $(r,s)$
of \eqref{cybe3} are given by: $r = a\tt b- b\tt a$, any $s$ in case
$\lambda = 0$, and by the same $r$, and $s \in \kk b$ when $\la\ne0$.
It is immediate to see that in both cases $s$ can be written as 
$-\phi(a) b + \phi(b)a + [a,b]$ for some trace form $\phi \in \dd^*$. Then
$r-1\tt s = e_{ab}$ is a 
free generator of $S(\dd, \phi)$, since $\dim \dd = 2$ (see
\prref{psd}).
This shows that the above pairs $(r,s)$ also satisfy \eqref{cybe4}.
We have: $\Hd = S(\dd, \phi)$, where 
$\chi = \iota_s\om = -\phi + \tr\ad$.
(Note that $\tr\ad=\iota_x \om$ for $x=[a,b]=\la b$.)

\end{example}
\begin{example}
When $\dd$ is abelian of dimension $N =2n>2$, then \eqref{omchi1} and
\eqref{omchi2}
become $\chi \wedge\om = 0$, hence $\chi = 0$ and $\om$ is any nondegenerate 
skew-symmetric $2$-form. In this case all solutions of 
\eqref{cybe3}, \eqref{cybe4} are: $s=0$ and $r$ given by \eqref{raibi}
in some basis $\{a_i,b_i\}$ of $\dd$.
\end{example}
\begin{example}\lbb{eybdsim1}
When the Lie algebra $\dd$ is simple, there are no solutions 
$(\om,\chi)$ of \eqref{omchi1}, \eqref{omchi2} with a nondegenerate $\om$.
Indeed, since $[\dd,\dd]=\dd$, we have $\chi=0$, and $\om$ is a $2$-cocycle:
$\di\om=0$. Any $2$-cocycle $\om\in\bigwedge^2\dd^*$ for a simple 
Lie algebra $\dd$ is degenerate, since $\om=\di\al$ for some $\al\in\dd^*$
and the stabilizer $\dd_\al$ of $\al$ is always non-zero.
\end{example}

\subsection{$\Kd$}\lbb{subk}
This is defined as a Lie $H$-\psalg\ of rank $1$  (see \seref{subr1wd})
corresponding to a solution $(r,s)$ of equations \eqref{cybe3}, \eqref{cybe4}
with $\dd=\dd_1\oplus\Kset s$ and nondegenerate $r\in \dd_1 \wedge \dd_1$;
in this
case $N=\dim\dd$ is odd. The
parameter $\th$ is defined below.

Let $\{\d_i\}$ be a basis of $\dd_1$, and $r=\tsum r^{ij} \d_i\tt\d_j$.
As before, we define a $2$-form $\om$ on $\dd_1$ by 
$\om(\d_i\wedge\d_j)=\om_{ij}$, where $(\om_{ij})$ is the matrix inverse 
to $(r^{ij})$. Let us write 
$[\d_i,\d_j] = \tsum c_{ij}^k \d_k + c_{ij} s$
and
$[s,\d_j] = \tsum c_j^k \d_k$.
Then we have:

\begin{lemma}\lbb{lomk}
With the above notation,
equations \eqref{cybe3}, \eqref{cybe4} are equivalent to the following
identities{\rm:}
\begin{align}
\lbb{omk1}
\di\om &= 0 \quad\text{on \; $\textstyle\bigwedge^3\dd_1$},
\\
\lbb{omk2}
c_{ij} &= \om_{ij},
\\
\lbb{omk3}
L_s \om &= 0.
\end{align}
If we extend $\om$ to a $2$-form on $\dd$ by defining $\iota_s\om=0$,
then $\om$ is closed{\rm:} $\di\om=0$.
\end{lemma}
\begin{proof}
The proof is very similar to that of \leref{lomchi}.
There we showed that $L_s \om = 0$ is equivalent to \eqref{cybe3},
and the same argument applies here.
Similarly, \eqref{cybe4} is equivalent to (\ref{omk1}, \ref{omk2}).
Now if $\iota_s\om=0$, then $L_s \om = 0$ implies $\iota_s\di\om=0$, which
together with \eqref{omk1} leads to $\di\om=0$.
\end{proof}

Let $\om$ be extended to a $2$-form on $\dd$ by defining $\iota_s\om=0$,
so that $\di\om=0$. We define a $1$-form $\th\in\dd^*$ by $\th(s):=-1$,
$\th|_{\dd_1} :=0$. Then we have $\di\th=\om$; indeed:
\begin{align*}
(\di\th)(\d_i\wedge\d_j) &= -\th([\d_i,\d_j]) = c_{ij} 
= \om_{ij} = \om(\d_i\wedge\d_j),
\\
(\di\th)(s\wedge\d_j) &= -\th([s,\d_j]) = 0 = \om(s\wedge\d_j),
\end{align*}
using \eqref{omk2} and the fact that $[s,\dd_1]\subset\dd_1$.

\begin{lemma}\lbb{lomk2}
There is a one-to-one correspondence between contact forms $\th$, i.e.\
$1$-forms $\th\in\dd^*$
such that $\th\wedge(\di\th)^{(N-1)/2} \ne 0$ $(N=\dim\dd)$,
and solutions $(r,s)$
of \eqref{cybe3}, \eqref{cybe4}
with $\dd = \dd_1 \oplus \Kset s$ and nondegenerate $r\in\dd_1\tt\dd_1$.
\end{lemma}
\begin{proof}
Given $(r,s)$, above we have defined the $1$-form $\th$ such that
$\th(s)=-1$, $\th|_{\dd_1}=0$ and $\di\th=\om$. Since 
$\om\in\bigwedge^2 \dd_1^*$ is nondegenerate, we have
$\th\wedge\om^{(N-1)/2} \ne 0$.
Conversely, starting with a contact $1$-form $\th\in\dd^*$, we can define
$s$ and $\om$ satisfying \eqref{omk1}--\eqref{omk3}.
\end{proof}
\begin{example}\lbb{eheisenberg}
When $\dd$ is the Heisenberg Lie algebra with a basis $\{a_i,b_i,c\}$
and the only nonzero commutation relations
$[a_i, b_i] = c$ ($1\le i\le n$, $N=2n+1$),
then
\begin{displaymath}\lbb{r1heis}
r = \tsum_{i=1}^n \, (a_i\tt b_i - b_i\tt a_i),
\quad s=-c
\end{displaymath}
is a solution of \eqref{cybe3}, \eqref{cybe4}.
\end{example}
\begin{example}
When $\dd$ is abelian and $\dim \dd = 2n+1 > 1$, 
then equations \eqref{cybe3}, \eqref{cybe4} have
no solutions $(r,s)$ with $\dd = \dd_1 \oplus \Kset s$ 
and a nondegenerate $r\in\dd_1\wedge\dd_1$, because $\di\th = 0$ 
and therefore there are no contact forms.
\end{example}
\begin{example}\lbb{eybdsim2}
When the Lie algebra $\dd$ is simple, a solution $(r,s)$
of \eqref{cybe3}, \eqref{cybe4} with $\dd = \dd_1 \oplus \Kset s$ 
and a nondegenerate $r\in\dd_1\wedge\dd_1$ exists iff $\dd = \sl_2$,
and it is as follows:
\begin{displaymath}
r = e \wedge f := e \tt f - f \tt e, \;\; s=-h.
\end{displaymath}
Only $\dd = \sl_2$ is possible since the dimension of the stabilizer of
$\th$ equals $1$.
\end{example}

Now let us compute $L_{e_x}\th$. Recall that, as in \seref{subh}, 
for any $x\in X$ we identify $e_x:= x\tt_H e$ with 
$-\sum (x a_i \tt b_i - x b_i \tt a_i) + x\tt s$.
Similarly to \eqref{omexa}, it is easy to see that 
$\om(e_x \wedge a) = -xa$ for $a\in\dd_1$ (in this case $\chi=\iota_s\om:=0$).
On the other hand, $\iota_{e_x}\th = \th(e_x) = -x$, and hence
$(\di\iota_{e_x}\th)(a) = -(\di x)(a) = xa$ for any $a\in\dd$.
Therefore $(L_{e_x}\th)(a) = 0$ for $a\in\dd_1$, and $(L_{e_x}\th)(s) = xs$.
In other words, 
\begin{displaymath}\lbb{lexth}
L_{e_x}\th = -(xs)\th, 
\end{displaymath}
and we have the following proposition.

\begin{proposition}\lbb{pthkn}
Let $\Kd := He$ be a Lie $H$-\psalg\ of rank $1$ 
corresponding to a solution $(r,s)$ of equations \eqref{cybe3}, \eqref{cybe4}
with $\dd = \dd_1 \oplus \Kset s$ and a nondegenerate $r\in\dd_1\tt\dd_1$,
where the $1$-form $\th\in\dd^*$ 
is defined by $\th(s)=-1$, $\th|_{\dd_1} =0$.
Then $\th$ is a contact form, and
the subalgebra $X\tt_H \Kd$
of $X\tt_H \Wd \simeq X\tt\dd$
is the Lie algebra $K_N(\th)$ of vector fields that preserve $\th$
up to a multiplication by a function
$($which is isomorphic to $K_N)$.
\end{proposition}
\begin{proof}
It remains to show that, conversely, any vector field from $K_N(\th)$
is equal to $e_x$ for some $x\in X$. 
Indeed, let $A\in X\tt\dd$ be such that $L_A\th=f\th$ for some $f\in X$. 
Let us write $A=\sum_i (x_i\tt a_i + y_i\tt b_i) + x\tt s$
for some $x_i,y_i,x\in X$. Then $\om(A\wedge a_i) = y_i$
and $\om(A\wedge b_i) = -x_i$, while $\th(A)=-x$.
Therefore $(L_A\th)(a) = \om(A\wedge a) + xa$, which implies
$y_i + x a_i=0$, $-x_i + x b_i =0$, and $xs = -f$.
\end{proof}

\begin{remark}
To any $H$-type Lie \psalg, i.e., to any triple
$(\dd, \om, \chi)$ where $\dd$ is a finite-dimensional Lie algebra,
$\om\in \textstyle\bigwedge^2 \dd^*$ is a non-degenerate $2$-form and 
$\chi\in \dd^*$
satisfying \eqref{omchi1} and \eqref{omchi2}, we can associate a $K$-type
Lie \psalg\ as follows.
Set on the vector space $\dd'=\dd\oplus \kk c$ the Lie bracket $[\,\, ,\,]'$
defined as:
\begin{displaymath}
[g,h]' = [g,h]+ \om(g,h)c, \qquad [g,c]' = \chi(g) c,
\end{displaymath}
for $g,h\in \dd$. Then $c+s\in\dd'$ stabilizes $\dd$,
where $s\in \dd$ is the unique element such that $\chi = \iota_s \om$;
indeed,
\begin{displaymath}
[g,s+c]' = [g,s] + \om(g,s) c + \chi(g) c = [g,s]\in \dd.
\end{displaymath}
Define $\th\in (\dd')^*$ as the unique element restricting
to $0$ on $\dd$ such that $\th(c) = -1$.

Note that not all $K$-type data are obtained in this way, since the Lie algebra
$\dd'$ just constructed always has a one dimensional ideal $\kk c$,
and this fails in \exref{eybdsim2}.
\end{remark}

\subsection{Annihilation algebras of \psalgs\ of vector fields}
\lbb{subannih}
To conclude this section, we determine the annihilation algebras of the 
primitive \psalgs\ of vector fields defined above, 
and of current \psalgs\ over them.

\begin{theorem}\lbb{tannihvec}
{\rm(i)}
If $L$ is one of the Lie $H=\ue(\dd)$-\psalgs\ $\Wd$, $\Sd$, $\Hd$ or $\Kd$,
then its annihilation algebra $\A(L)$ is isomorphic to 
$W_N$, $S_N$, $P_N$ or $K_N$, respectively.

{\rm(ii)}
If $L=\Cur L'$ is a current \psalg\ over the Lie $H'$-\psalg\ $L'$,
then its annihilation algebra $\A(L)$ is isomorphic to 
a current Lie algebra $\O_r\what\tt \A(L')$ over $\A(L')$,
where $H'=\ue(\dd')$ and $\dd'$ is a codimension $r$ subalgebra of $\dd$.
\end{theorem}
\begin{proof}
(i) We have seen in \seref{subw} that $\A(\Wd) \simeq W_N$.
Let $L$ be one of the \psalgs\ $\Sd$, $\Hd$ or $\Kd$ and $i$ be its natural
embedding in $\Wd$. We have shown in 
Sections \ref{subs}, \ref{subh} and \ref{subk}
that in this case the image of $\A(L)$ in $W_N$ under $\A(i)$ is
$S_N$, $H_N$ or $K_N$, respectively. 

\leref{linl1} below implies that $\A(L)$ is a central extension of its image
in $W_N$. Moreover, since $L$ is simple, it is equal to its derived
subalgebra, and therefore $\A(L)$ is equal to its derived subalgebra
(see \seref{srad}). Hence, $\A(L)$ is an irreducible central extension
of its image in $W_N$.

Now \prref{pder}(iii) implies that $\A(L) \simeq S_N$, $K_N$ in the cases
$L=\Sd$, $\Kd$ respectively, and $\A(L)$ is a quotient of $P_N$ in the
case $L=\Hd$. However, since $L=\Hd$ is a free $H$-module of rank one,
$\A(L)$ is isomorphic to $X$ as a topological $H$-module.
Therefore, $\A(L) \simeq P_N$.

(ii) Note that $X=H^*$ maps surjectively to $X' = (H')^*$ with kernel
isomorphic to $\O_r$. Moreover $X\simeq\O_N$, $X'\simeq\O_{N'}$ $(N'=N-r)$,
and hence $X\simeq\O_r \what\tt X'$.
We have: $\A(L') := X' \tt_{H'} L'$ and
$\A(L) := X\tt_H L = X\tt_H (H\tt_{H'} L')
\simeq X \tt_{H'} L' \simeq (\O_r \what\tt X')\tt_{H'} L'
\simeq \O_r \what\tt (X' \tt_{H'} L')$.
\end{proof}
\begin{remark}\lbb{rdactwn}
Let us assume that the base field $\Kset=\Cset$, and let $L$ be as
in \thref{tannihvec}(i).
Then the action of $\dd$ on $\A(L)$ can be constructed via the embedding of 
$\dd$ in $W_N$ as follows. 

(i) Any $N$-dimensional Lie algebra $\dd$ can be embedded in $W_N$:
every $a\in\dd$ defines a left-invariant 
vector field on the connected simply-connected Lie group $D$ with Lie algebra 
$\dd$, and we take the corresponding formal vector field in the formal 
neighborhood of the identity element. (See also \prref{existuniq}.)

(ii) If we have a homomorphism of Lie algebras
$\chi\colon\dd\to\Cset$, it defines a
homomorphism $\ti\chi$ of $D$ to $\Cset^\times$. 
Consider a volume form $v$ on $D$  
defined, up to a constant factor, by the property 
$g\cdot v_0=\ti\chi(g)v_0$, $g\in D$,
where $v_0$ is the value of $v$ at the identity 
element. Then we get an embedding of $\dd$ in 
$CS_N(v)=\Der S_N(v) \simeq \Der S_N$.

(iii) Given $\chi$ and $\om\in\bigwedge^2\dd^*$ such that 
$\di\om+\chi\wedge\om=0$,
consider a $2$-form $s$ on $D$ whose value
at the identity element is $\om$ and such that 
$g\cdot s=\ti\chi(g)s$, $g\in D$.
Then $s$ is a symplectic form on $D$, and
we get an embedding of $\dd$ in $CH_N(s)=\Der H_N(s) \simeq \Der P_N$.

(iv) Given a contact form $\th\in\dd^*$, consider the left-invariant 
$1$-form $c$ on $D$ with the value $\th$ at the identity element. 
Then $c$ is a contact form on $D$, and
we get an embedding of $\dd$ in $K_N(c) \simeq K_N$.
\end{remark}

\section{$H$-Conformal Algebras}\lbb{shconf}

The goal of this section is to reformulate the definition
of a Lie (or \as) $H$-\psalg\
in terms of the properties of the $x$-brackets
(or products) introduced in \seref{subay}.
The resulting notion, equivalent to that of an $H$-\psalg,
will be called an $H$-\conalg.

Let us start by deriving explicit formulas for the $x$-brackets.
We will use the notation of \seref{subay}. 
Let $(L,\be)$ be a Lie $H$-\psalg\
with a pseudobracket
\begin{equation}\lbb{psbr1}
[a*b] \equiv \be(a\tt b) = \tsum_i\, (f_i\tt g_i) \tt_H e_i.
\end{equation}
Then for $x\in X$, $h\in H$ we have 
$\eta(x\tt h) = \< x, h_{(1)} \> h_{(2)}$
(see \eqref{etaxh}), and 
\begin{align*}
(\eta\tt_H\be) \bigl( (x\tt_H a) \tt (h\tt_H b) \bigr)
&= \tsum_i\, \eta(x f_i \tt h g_i) \tt_H e_i
\\
&= \tsum_i\, \< x f_i, (h g_i)_{(1)} \> \, (h g_i)_{(2)} \tt_H e_i .
\end{align*}
Taking $h=1$, we get the following expression for the 
{\em $x$-bracket\/} in $L$:
\begin{equation}\lbb{xbrak}
\xb axb
= \tsum_i\, \< S(x), f_i {g_i}_{(-1)} \> \, {g_i}_{(2)} e_i,
\quad\text{if}\quad 
[a*b] = \tsum_i\, (f_i\tt g_i) \tt_H e_i.
\end{equation}

Here we can recognize the Fourier transform $\F$, defined by \eqref{ftrans}:
\begin{displaymath}
\F(f\tt g) = f g_{(-1)} \tt g_{(2)}.
\end{displaymath}
The identity \eqref{sigfour}:
\begin{displaymath}
f\tt g = (f g_{(-1)} \tt 1) \, \De(g_{(2)}),
\end{displaymath}
implies
\begin{equation}\lbb{psbr2}
[a*b] = \tsum_i\, (f_i {g_i}_{(-1)}  \tt 1) \tt_H {g_i}_{(2)} e_i.
\end{equation}
Hence $[a*b]$ can be written uniquely in the form
$\tsum_i\, (h_i \tt 1) \tt_H c_i$,
where $\{h_i\}$ is a fixed $\Kset$-basis of $H$
(cf.\ \leref{lhhh}).

We introduce another bracket $[a,b]\in H\tt L$ 
as the Fourier transform of $[a*b]$:
\begin{equation}\lbb{bracket}
[a,b] = \tsum_i\, \F(f_i\tt g_i) \, (1\tt e_i)
= \tsum_i\, f_i {g_i}_{(-1)} \tt {g_i}_{(2)} e_i.
\end{equation}
In other words,
\begin{equation}\lbb{ab}
[a,b] = \tsum_i\, h_i\tt c_i \quad\text{if}\quad
[a*b] = \tsum_i\, (h_i\tt1)\tt_H c_i.
\end{equation}
Then we have:
\begin{equation}\lbb{xbracket}
\xb{a}{x}{b} = (\langle S(x),\cdot\rangle\tt\id) \, [a,b]
= \tsum_i\, \langle S(x), h_i \rangle c_i.
\end{equation}

Using properties \eqref{propf0}--\eqref{propf3} of the Fourier transform,
it is straightforward to derive the properties of
the bracket \eqref{ab}. 
Then the definition of a Lie \psalg\ can be equivalently reformulated
as follows.

\begin{definition}\lbb{dhconf1}
A {\em Lie $H$-conformal algebra\/} is a left $H$-module $L$
equipped with a bracket $[\cdot,\cdot]\colon L\tt L\to H\tt L$,
satisfying the following properties ($a,b,c\in L$, $h\in H$):
\begin{description}

\item[$H$-sesqui-linearity]
\begin{align}
\lbb{bil1}
[ha,b] &= (h\tt1) \, [a,b],
\\
\lbb{bil2}
[a,hb] &= (1\tt h_{(2)}) \, [a,b] \, (h_{(-1)}\tt1).
\end{align}

\item[Skew-commutativity]
If $[a,b]$ is given by \eqref{ab}, then
\begin{equation}\lbb{ss}
[b,a] = - \tsum_i\, {h_i}_{(-1)} \tt {h_i}_{(2)} c_i.
\end{equation}

\item[Jacobi identity]
\begin{equation}\lbb{jacid}
[a,[b,c]] - (\si\tt\id) \, [b,[a,c]]
= (\F^{-1}\tt\id) \, [[a,b],c]
\end{equation}
in $H\tt H\tt L$, where $\si\colon H\tt H\to H\tt H$ is the permutation
$\si(f\tt g)=g\tt f$, and
\begin{align}
\lbb{abc1}
[a,[b,c]] &= (\si\tt\id) \, (\id\tt[a,\cdot]) \, [b,c],
\\
\lbb{abc2}
[[a,b],c] &= (\id\tt[\cdot,c]) \, [a,b].
\end{align}

\end{description}

\end{definition}
\begin{examples}\lbb{ecwd}
{\rm(i)}
For the current Lie pseudoalgebra
$\Cur\g=H\tt\g$ with the pseudobracket \eqref{curalg2},
the bracket \eqref{ab} is given by:
\begin{displaymath}\lbb{curalg3}
[f \tt a, g \tt b] = f g_{(-1)} \tt (g_{(2)} \tt [a,b]).
\end{displaymath}

{\rm(ii)}
For the Lie pseudoalgebra $\Wd=H\tt\dd$ with pseudobracket defined by
\eqref{wdbr*}, the bracket \eqref{ab} is given by:
\begin{displaymath}\lbb{wdbr}
[1\tt a, 1\tt b] 
= 1\tt(1\tt [a,b]) + a\tt(1\tt b) + b\tt(1\tt a) - 1\tt(a\tt b).
\end{displaymath}
\end{examples}

One can also reformulate \deref{dhconf1} in terms of the 
$x$-brackets \eqref{xbracket}.

\begin{definition}\lbb{dhconf2}
A {\em Lie $H$-conformal algebra\/} is a left $H$-module $L$
equipped with $x$-brackets $\xb axb \in L$ for $a,b\in L$, $x\in X$,
satisfying the following properties:
\begin{description}

\item[Locality]
\begin{equation}\lbb{codimx}
\codim\{ x\in X \st \xb{a}{x}{b} = 0\} < \infty
\quad\text{for any $a,b\in L$}. 
\end{equation}
Equivalently, for any basis
$\{x_i\}$ of $X$,
\begin{equation}\lbb{loc2}
\xb{a}{x_i}{b} \ne 0 \quad\text{for only a finite number of $i$.}
\end{equation}

\item[$H$-sesqui-linearity]
\begin{align}
\lbb{bl1}
\xb{ha}{x}{b} &= \xb{a}{xh}{b},
\\
\lbb{bl2}
\xb{a}{x}{hb} &= h_{(2)} \, \xb{a}{h_{(-1)} x}{b}.
\end{align}

\item[Skew-commutativity]
Choose dual bases $\{h_i\}$, $\{x_i\}$ in $H$ and $X$. Then:
\begin{equation}\lbb{sksy}
\xb{a}{x}{b} = -\tsum_i \langle x, {h_i}_{(-1)} \rangle \,
{h_i}_{(-2)} \, \xb{b}{x_i}{a}.
\end{equation}

\item[Jacobi identity]
\begin{equation}\lbb{jid}
\xb{a}{x}{\xb{b}{y}{c}} - \xb{b}{y}{\xb{a}{x}{c}} 
= [[a_{x_{(2)}} b]_{y x_{(1)}} c].
\end{equation}

\end{description}

\end{definition}

\leref{lxshxy} implies that \eqref{jid} can be rewritten as follows:
\begin{equation}\lbb{jid2}
\xb{a}{x}{\xb{b}{y}{c}} - \xb{b}{y}{\xb{a}{x}{c}} 
= \tsum_i \, [[a_{x_i} b]_{y (xS(h_i))} c].
\end{equation}
In particular, the right-hand side of \eqref{jid} is well defined:
the sum is finite because of \eqref{loc2}.

The definitions of an {\em \as\ $H$-\conalg\/} or of 
{\em representations\/} of $H$-\conalgs\ 
are obvious modifications of the above.
For example, in terms of $x$-products, the associativity looks as
follows (cf.\ \eqref{jid}):
\begin{equation}\lbb{asxpro}
\xp a x {(\xp byc)} = (a_{x_{(2)}} b)_{y x_{(1)}} c.
\end{equation}
The same argument as the one used for $\F$ shows that the map
$x \tt y \mapsto x_{(2)} \tt y x_{(1)}$ has an inverse given by
$x \tt y \mapsto x_{(2)} \tt y x_{(-1)}$. Therefore, \eqref{asxpro}
is equivalent to the following equation:
\begin{equation}\lbb{asxpro2}
\xp a {x_{(2)}} {(\xp b {y x_{(-1)}} c)} = (a_x b)_y c.
\end{equation}
Note that when considering \as\ $H$-\conalgs, $H$ need not be cocommutative,
so $X$ may be noncommutative.

We also note that there is a simple relationship between the
$x$-bracket \eqref{xbracket} of a Lie $H$-\conalg\ 
(or, equivalently, \psalg) $L$
and the commutator in its annihilation Lie algebra $\A(L)$
defined in \seref{slalg}.
Let $\{h_i\}$, $\{x_i\}$ again be dual bases in $H$ and $X$.
Then in \eqref{ab} one has $c_i=\xb a{S^{-1}(x_i)}b$; therefore
\begin{equation}\lbb{ab2}
[a,b] = \tsum_i\, S(h_i) \tt \xb a{x_i}b
\quad\text{and}\quad
[a*b] = \tsum\, ( S(h_i) \tt1)\tt_H \xb a{x_i}b.
\end{equation}
Recall that we denote the element $x \tt_H a$ of $\A(L) := X \tt_H L$
by $a_x$. Then the definition \eqref{ahxbhy} and \eqref{ab2} imply
\begin{equation}\lbb{axby}
[a_x, b_y] = \tsum_i\, [a_{x_i}b]_{(x S(h_i))y}
= [a_{x_{(2)}} b]_{x_{(1)} y},
\end{equation}
using \eqref{xshxy}. This is also equivalent to:
\begin{equation}\lbb{axby2}
[a_xb]_y = [a_{x_{(2)}}, b_{x_{(-1)} y}]
= \tsum_i\, [a_{x_i}, b_{(h_i S(x)) y}].
\end{equation}

Comparing these formulas with Eq.~\eqref{jid}, we obtain the following
important result.

\begin{proposition}\lbb{preplal}
Any module $M$ over a Lie \psalg\ $L$ has a natural structure of
an $\A(L)$-module, given by 
$(x\tt_H a) \cdot m = a_x m$, where
\begin{equation}\lbb{axm2}
\xp axm
= \tsum_i\, \< S(x), f_i {g_i}_{(-1)} \> \, {g_i}_{(2)} v_i,
\quad\text{if}\quad 
a*m = \tsum_i\, (f_i\tt g_i) \tt_H v_i
\end{equation}
for $a\in L, x\in X, m\in M$.
This action is compatible with the action of $H$ $($see \eqref{ham12}$)$
and satisfies the locality condition{\rm:}
\begin{equation}\lbb{codimx2}
\codim\{ x\in X \st a_x m = 0\} < \infty,
\qquad a\in L, \; m\in M,
\end{equation}
or equivalently, for any basis $\{x_i\}$ of $X$,
\begin{equation}\lbb{locc3}
a_{x_i} m \ne 0 \quad\text{for only a finite number of $i$.}
\end{equation}
{\rm(}The above conditions on $M$ mean that, when endowed with the discrete
topology, $M$ is a topological
$\A(L)$-module in the category $\M^l(H)$.{\rm)}

Conversely, any $\A(L)$-module $M$ satisfying the above conditions
has a natural structure of an $L$-module, given by{\rm:}
\begin{equation}\lbb{prpl2}
a*m = \tsum_i\, ( S(h_i) \tt1)\tt_H a_{x_i} \cdot m,
\end{equation}
where $\{h_i\}$, $\{x_i\}$ are dual bases in $H$ and $X$,
and we use the notation $a_x \equiv x\tt_H a$.
\end{proposition}

This proposition provides the main tool for constructing
modules over Lie \psalgs.
Of course, there is an analogous result in the case of \as\ algebras
as well.

Finally, let us give two more expressions for the bracket in $\A(L)$
which will be useful later. Recall that, by \prref{preplal},
we have an action of $\A(L)$
on $L$ given by $a_x \cdot b = \xb axb$. 
Recall also that the action of $H$ on $\A(L)$ is defined by
$h(a_x) = a_{hx}$.
Let $\alpha\in\A(L)$, $b\in L$, $y\in X$. Then:
\begin{align}
\lbb{alby1}
(\alpha \cdot b)_y &= \tsum_i \, [h_i \alpha , b_{x_i y}],
\\
\lbb{alby2}
[\alpha, b_y] &= \tsum_i \, \bigl( (S(h_i) \alpha) \cdot b \bigr)_{x_i y}.
\end{align}
Note that the infinite sums on the right-hand sides make 
sense since they are convergent in the complete topology of $\A(L)$.
It is enough to prove both statements for $\alpha$ of the form 
$a_x = x\tt_H a$ since such elements span $\A(L)$. 
Equation~\eqref{alby2} then follows from 
\eqref{axby} and \eqref{xshxy}.
Analogously, \eqref{alby1} derives from \eqref{axby2} by 
noticing that $x_{(-1)} \tt x_{(2)} = \tsum_i \, x_i \tt h_i x$.

\section{$H$-Pseudolinear Algebra}\lbb{shcla}
The definition of a module over a pseudoalgebra motivates the following 
definition of a pseudolinear map.
\begin{definition}\lbb{dchom}
Let $V$ and $W$ be two $H$-modules. An {\em $H$-pseudolinear map\/}
from $V$ to $W$ is a $\Kset$-linear map 
$\phi\colon V\to (H\tt H)\tt_H W$
such that 
\begin{align}
\lbb{hclm}
\phi(hv) &= ((1\tt h)\tt_H 1) \, \phi(v),
\qquad h\in H, v\in V. 
\intertext{We denote the space of all such $\phi$ by $\Chom(V,W)$.
We define a left action of $H$ on $\Chom(V,W)$ by:}
\lbb{hactchom}
(h\phi)(v) &= ((h\tt1)\tt_H 1)\, \phi(v).
\end{align}
When $V=W$, we set $\Cend V=\Chom(V,V)$.
\end{definition}
\begin{example}\lbb{ecendm}
Let $A$ be an $H$-\psalg\ and $V$ be an $A$-module. Then for any $a\in A$
the map $m_a\colon V\to (H\tt H)\tt_H V$ defined by $m_a(v) = a*v$
is an $H$-pseudolinear map. Moreover, we have $hm_a = m_{ha}$
for $h\in H$.
\end{example}

Consider the map $\rho\colon \Chom(V,W) \tt V \to (H\tt H)\tt_H W$
given by $\rho(\phi\tt v) = \phi(v)$. By definition it is $H$-bilinear,
so it is a polylinear map 
in $\M^*(H)$.
We will also use the notation $\phi * v := \phi(v)$ and consider this
as a pseudoproduct (or rather action, see \prref{pcend} below).

The corresponding $x$-products are called {\em Fourier coefficients\/}
of $\phi$ and are given by a formula analogous to \eqref{xbrak}:
\begin{equation}\lbb{fourcoe}
\phi_x v
= \tsum_i\, \< S(x), f_i {g_i}_{(-1)} \> \, {g_i}_{(2)} w_i,
\qquad
\text{if}\quad \phi(v) = \tsum_i\, (f_i\tt g_i) \tt_H w_i.
\end{equation}
They satisfy a locality relation and an $H$-sesqui-linearity relation 
similar to \eqref{codimx} and \eqref{bl2}:
\begin{equation}\lbb{codimphix}
\codim\{ x\in X \st \phi_x v = 0\} < \infty \quad\text{for any}\;\; v\in V,
\end{equation}
\begin{equation}\lbb{bl2phix}
\phi_x (hv) = h_{(2)} ( \phi_{h_{(-1)} x} v ).
\end{equation}
Conversely, any collection of maps $\phi_x\in\Hom(V,W)$, $x\in X$,
satisfying relations \eqref{codimphix}, \eqref{bl2phix} 
comes from an $H$-pseudolinear map $\phi\in\Chom(V,W)$.
Explicitly (cf.~\eqref{prpl2}):
\begin{equation}\lbb{fourcoe2}
\phi(v) = \tsum_i\, ( S(h_i) \tt1)\tt_H \phi_{x_i} v,
\end{equation}
where $\{h_i\}$, $\{x_i\}$ are dual bases in $H$ and $X$.

\begin{remark}\lbb{rhomchom}
It follows from \eqref{bl2phix} that for $\phi\in\Chom(V,W)$,
the map $\phi_1\colon V \to W$ is $H$-linear, where $1\in X$ is the unit
element. This establishes an isomorphism 
$\Hom_H(V,W) \simeq \kk\tt_H \Chom(V,W) \simeq \Chom(V,W) / H_+ \Chom(V,W)$,
where $H_+ = \{ h\in H \st \ep(h)=0 \}$ is the augmentation ideal.
\end{remark}
\begin{lemma}\lbb{chom*chom}
Let $U,V,W$ be three $H$-modules, and assume that $U$ is finite.
Then there is a unique polylinear map 
\begin{displaymath}
\mu\in\Lin( \{ \Chom(V,W), \Chom(U,V)\}, \Chom(U,W))
\end{displaymath}
in $\M^*(H)$, denoted as $\mu(\phi\tt\psi) = \phi*\psi$, such that
\begin{equation}\lbb{phipsiv}
(\phi*\psi)*u = \phi*(\psi*u)
\end{equation}
in $H^{\tt3}\tt_H W$
for $\phi\in\Chom(V,W)$, $\psi\in\Chom(U,V)$, $u\in U$.
\end{lemma}
\begin{proof}
We define $\phi*\psi$ in terms of its Fourier coefficients ---
the $x$-products $\xp \phi x \psi$.
We have already seen, when we discussed associativity, that
\eqref{phipsiv} is equivalent to the following equation
(cf.\ \eqref{asxpro}):
\begin{displaymath}
\xp \phi x {(\xp \psi y u)} = (\phi_{x_{(2)}} \psi)_{y x_{(1)}} u.
\end{displaymath}
This can be inverted to find (cf.\ \eqref{asxpro2}):
\begin{displaymath}
(\phi_x \psi)_y u = \xp \phi {x_{(2)}} {(\xp \psi {y x_{(-1)}} u)}
= \tsum_i\, \xp \phi {x_i} {(\xp \psi {y (h_i S(x))} u)}.
\end{displaymath}

The $H$-sesqui-linearity properties of $(\phi_x \psi)_y u$ with respect to 
$x$ and $y$ are easy to check by a direct calculation.
By properties \eqref{fils1}, \eqref{fils3}, \eqref{fils2} of the filtration
$\{\fil_n X\}$, and locality of $\psi$,
it follows that for each fixed $x\in X$, $u\in U$ there is an $n$ such that
$(\phi_x \psi)_y u = 0$ for $y\in \fil_n X$.
Therefore, for each $x\in X$ we have defined $\xp \phi x \psi \in \Chom(U,W)$.

In order that $\phi*\psi$ be well defined, we need
to check that $\xp \phi x \psi$
satisfies locality, i.e., that $\xp \phi x \psi = 0$
when $x\in \fil_n X$ with $n\gg0$. 
By the locality of $\phi$ and $\psi$, for each $u\in U$
there is an $n$ such that
$(\phi_x \psi)_y u = 0$ for $x\in \fil_n X$ and all $y\in X$.
Since $U$ is finite, we can choose an $n$ that works
for all $u$ belonging to a set of generators of $U$ over $H$.
Now the $H$-sesqui-linearity of $(\phi_x \psi)_y u$ with respect
to $y$ (for fixed $x$) implies that 
$(\phi_x \psi)_y u = 0$ for all $y$ and $u$.
Hence $\xp \phi x \psi = 0$ for $x\in \fil_n X$.
\end{proof}

Specifying to the case $U=V=W$, we obtain a pseudoproduct $\mu$
on $\Cend V$, and an action $\rho$ of $\Cend V$ on $V$.

\begin{proposition}\lbb{pcend}
{\rm(i)}
For any finite $H$-module $V$, the above pseudoproduct
provides $\Cend V$ with the structure of an \as\ $H$-\psalg.
$V$ has a natural structure of a $\Cend V$-module given by 
$\phi*v \equiv\phi(v)$. 

{\rm(ii)} For an \as\ $H$-\psalg\ $A$,
giving a structure of an 
$A$-module on $V$ is equivalent to giving a homomorphism of
\as\ $H$-\psalgs\
from $A$ to $\Cend V$.
\end{proposition}
\begin{proof}
Part (i) is an immediate consequence of \leref{chom*chom}. Indeed, the only
thing that remains to be checked is the associativity of $\Cend V$,
and it follows from \eqref{phipsiv}:
\begin{align*}
\bigl((\phi*\psi)*\chi\bigr)*v 
&= (\phi*\psi)*(\chi*v)
= \phi*(\psi*(\chi*v)) 
\\
&= \phi*((\psi*\chi)*v)
= \bigl(\phi*(\psi*\chi)\bigr)*v.
\end{align*}

To prove part (ii), we associate with each $a\in A$ the
$H$-pseudolinear map $m_a\in\Cend V$ given by
$m_a(v) = a*v$. Then 
\begin{displaymath}
(m_a*m_b)*v = m_a*(m_b*v) =
a*(b*v) = (a*b)*v = m_{a*b}*v,
\end{displaymath}
which shows that
$m_a*m_b = m_{a*b}$.
\end{proof}

We denote by $\gc V$ the Lie $H$-\psalg\
obtained from the \as\ one $\Cend V$ by the construction of \prref{assolie*}.
Then $V$ is a $\gc V$-module, and one has a statement
analogous to part (ii) above.

\begin{proposition}\lbb{gcreps}
Let $V$ be a finite $H$-module. Then,
for a Lie $H$-\psalg\ $L$,
giving a structure of an $L$-module on $V$ is equivalent to giving a 
homomorphism of Lie $H$-\psalgs\ from $L$ to $\gc V$.
\end{proposition}
\begin{remark}\lbb{ractchom}
Let $L$ be a Lie $H$-\psalg, and $U,V$ be finite $L$-modules.
Then the formula ($a\in L$, $u\in U$, $\phi\in\Chom(U,V)$)
\begin{equation}\lbb{actchom}
(a*\phi)(u) = a*(\phi*u) - ((\si\tt\id)\tt_H\id) \, \phi*(a*u)
\end{equation}
provides $\Chom(U,V)$ with the structure of an $L$-module.
\end{remark}
\begin{definition}\lbb{dcder}
{\rm(i)} Let $A$ be an \as\ $H$-\psalg.
A {\em derivation\/} of $A$ is an
$H$-pseudolinear map $\phi\in\gc A$ which satisfies
\begin{equation}\lbb{cdera}
\phi*(a*b) = (\phi*a)*b + ((\si\tt\id)\tt_H\id) \, a*(\phi*b),
\qquad a,b\in A. 
\end{equation}
We denote the space of all such $\phi$ by $\Cder A$.

{\rm(ii)} Similarly, for a Lie $H$-\psalg\ $L$, let $\Cder L$
be the space of all $\phi\in\gc L$ satisfying
\begin{equation}\lbb{cderl}
\phi*[a*b] = [(\phi*a)*b] + ((\si\tt\id)\tt_H\id) \, [a*(\phi*b)],
\qquad a,b\in L. 
\end{equation}
\end{definition}

The next result follows easily from definitions.

\begin{lemma}\lbb{lcder}
{\rm(i)}
For any $H$-\psalg\ $A$, $\Cder A$ is a 
subalgebra of $\gc A$. 

{\rm(ii)}
When $A$ is \as\ {\rm(}respectively Lie{\rm)}, 
we have a homomorphism of \psalgs\
$i\colon A\to\Cder A$ given by 
$i(a)(b) = a*b - (\si\tt_H\id)\, b*a$ 
$($respectively $i(a)(b) = [a*b])$,
whose kernel is the center of $A$.

{\rm(iii)} 
For any $x\in X$ and $\phi\in\Cder A$, $\phi_x$
is a derivation of the annihilation algebra of $A$.
In other words, we have{\rm:} $\A(\Der A) \subset \Der \A(A)$.

{\rm(iv)}
Let $A$ be an \as\ $H$-\psalg\ and $L$ be the corresponding 
Lie \psalg\ with pseudobracket given by commutator.
Then $\Cder A \subset \Cder L$.
\end{lemma}
\begin{example}\lbb{eomdd}
Recall that in \seref{spsforms} we defined the $\Wd$-module 
of pseudoforms $\Om(\dd) = H\tt\bigwedge^\bullet \dd^*$.
Since $\bigwedge^\bullet \dd^*$ is an \as\ algebra with respect
to the wedge product, we can consider $\Om(\dd)$ as an \as\
\psalg: the current \psalg\ over $\bigwedge^\bullet \dd^*$.
Then, as in the case of usual differential forms, for any $\al\in\Wd$,
$\al*$ and $\al*_\iota$ are {\em superderivations\/} of $\Om(\dd)$,
see (\ref{a*wedge}, \ref{a*wedge2}).
\end{example}

In the case when $V$ is a free $H$-module of finite rank, one can give an
explicit description of $\Cend V$, and hence of $\gc V$, as follows.

\begin{proposition}\lbb{pcendexpl}
Let $V=H\tt V_0$, where $H$ acts trivially on $V_0$ and $\dim V_0 < \infty$.
Then $\Cend V$
is isomorphic to $H\tt H\tt\End V_0$, with $H$ acting by a left
multiplication on the first factor,
and with the following pseudoproduct{\rm:}
\begin{equation}\lbb{prodcend}
(f\tt a\tt A)*(g\tt b\tt B) 
= (f \tt g a_{(1)}) \tt_H (1 \tt b a_{(2)} \tt AB).
\end{equation}
The action of\/ $\Cend V$ on $V=H\tt V_0$ is given by{\rm:}
\begin{equation}\lbb{actcendv}
(f\tt a\tt A)*(h\tt v) = (f \tt ha) \tt_H (1 \tt Av).
\end{equation}
The pseudobracket in $\gc V$ is given by{\rm:}
\begin{equation}\lbb{brackgc}
\begin{split}
[(f\tt a\tt A)*(g\tt b\tt B)]
&= (f \tt g a_{(1)}) \tt_H (1 \tt b a_{(2)} \tt AB)
\\
&- (f b_{(1)} \tt g) \tt_H (1 \tt a b_{(2)} \tt BA).
\end{split}
\end{equation}
\end{proposition}
\begin{proof}
Since $(H \tt H) \tt_H V \simeq H \tt H \tt V_0$,
we can identify $\Cend V$ with $H\tt H\tt\End V_0$ so that
its action on $V$ is given by \eqref{actcendv}.
To prove \eqref{prodcend}, we use \eqref{phipsiv}
and the explicit definition of associativity from \seref{slie*}.
Due to $H$-bilinearity, we can assume that $f=g=h=1$. Then:
\begin{align*}
(1\tt a\tt A)*\bigl((1\tt b\tt B) *(1\tt v)\bigr) 
&=(1\tt a\tt A)*\bigl((1 \tt b) \tt_H (1 \tt Bv)\bigr)
\\
&= (1\tt a_{(1)} \tt b a_{(2)}) \tt_H (1 \tt ABv).
\end{align*}
On the other hand, we have:
\begin{displaymath}
\bigl( (1 \tt a_{(1)}) \tt_H (1 \tt b a_{(2)} \tt AB) \bigr)*(1\tt v)
= (1\tt a_{(1)} \tt b a_{(2)}) \tt_H (1 \tt ABv).
\end{displaymath}
Now \eqref{prodcend} follows from 
the uniqueness from \leref{chom*chom}.
\end{proof}
\begin{remark}\lbb{curglsubgc}
Let $V=H\tt V_0$, where $H$ acts trivially on $V_0$ and $\dim V_0 < \infty$.
Then $\Cur\End V_0$ can be identified with 
$H\tt 1\tt\End V_0 \subset \Cend V$.
Similarly, $\Cur\gl\, V_0$ is a subalgebra of $\gc V$.
\end{remark}

When $V=H\tt\Kset^n$, we will denote $\Cend V$ by $\Cend_n$,
and $\gc V$ by $\gc_n$. Of course, the essential case is when $V=H$ is of 
rank one. Let us describe the \as\ algebra $\A_Y\Cend_1$, where
$\A_Y$ is as in \seref{slalg}. As an $H$-module it is isomorphic to 
$Y\tt_H\Cend_1 \simeq Y\tt H$ with $H$ acting on the first factor.
We have $a_x = x\tt_H (1\tt a) \equiv x\tt a$ for $x\in Y$, $a\in H$.
Comparing \eqref{ahxbhy} with \eqref{prodcend}, we see that
the product in $Y\tt H$ is given by:
\begin{equation}\lbb{aycend1}
(x\tt a)(y\tt b) 
= x (y a_{(1)}) \tt b a_{(2)}.
\end{equation}
Hence $\A_Y\Cend_1$ is isomorphic to the smash product 
$Y\smash H$ (see \seref{sprelh}).
The annihilation algebra $\A(\Cend_1) \equiv \A_X\Cend_1 \simeq X\smash H$
is isomorphic as an \as\ algebra to the Drinfeld double of $H$
(see \cite{D}).
For $H=\ue(\dd)$, 
$\A(\Cend_1)$ can be
identified with the \as\ algebra of all differential operators on $X$,
while $\A(\gc_1)$ with the corresponding Lie algebra.

\begin{example}\lbb{cendud}
Let $H=\ue(\dd)$ be the universal enveloping algebra of a Lie algebra $\dd$.
We identify $\dd$ with its image in $H$, so that $\gc_1 = H\tt H$ contains
$H\tt\dd$. We claim that $f\tt a \mapsto -f\tt a$ ($f\in H$, $a\in\dd$)
is an embedding of Lie \psalgs\ $\Wd\injto\gc_1$,
compatible with their actions on $H$. This is immediate
by comparing \eqref{brackgc} with \eqref{wdbr*} 
and \eqref{actcendv} with \eqref{wdac*}. 

Consider $H$ as $\Cur\Kset$, i.e., as an \as\ $H$-\psalg\
with a pseudoproduct $f*g = (f\tt g)\tt_H 1$.
Then $\Wd = \Cder H \subset\gc_1$.
\end{example}

\begin{example}\lbb{cendcga}
Let $H=\Kset[\Ga]$ be the group algebra of a group $\Ga$.
Then for $V=H$ and $f,g,a,b\in\Ga$, \eqref{prodcend} takes the form:
\begin{displaymath}\lbb{prodcga}
(f\tt a)*(g\tt b) 
= (f \tt g a) \tt_H (1 \tt b a).
\end{displaymath}
\end{example}

We end this section with two lemmas that 
will be useful in representation theory.

\begin{lemma}\lbb{lkey}
For $\phi\in\Chom(V,W)$, let 
\begin{displaymath}\lbb{kernphi}
\ker_n\phi 
= \{ v \in V \st \phi_x v = 0 \;\;\;\forall x\in\fil_n X \},
\end{displaymath}
so that, for example, $\ker_{-1}\phi=\ker\phi$. If $V$ is a finite
$H$-module and $\fil^n H$ is finite dimensional,
then\/
$\ker_n\phi / \ker\phi$ is finite dimensional.
\end{lemma}
\begin{proof}
Since $\ker\phi$ is an $H$-submodule of $V$,
after replacing $V$ with $V/\ker\phi$, we can assume that $\ker\phi=0$.

By definition, 
$\ker_n\phi = \phi^{-1}\bigl( (\fil^n H\tt \Kset) \tt_H W \bigr)$.
Since, by \leref{hnkw}, 
$(\fil^n H\tt\Kset) \tt_H W = (\Kset\tt\fil^n H) \tt_H W$,
we have 
$\phi(\ker_n\phi) \subset (\Kset\tt\fil^n H) \tt_H W$.

On the other hand, since $V$ is finite over $H$
and $\phi$ satisfies \eqref{hclm}, there exists a finite-dimensional
subspace $W'$ of $W$ such that
$\phi(\ker_n\phi) \subset (\Kset\tt H) \tt_H W'$.
It follows that $\phi(\ker_n\phi) \subset (\Kset\tt\fil^n H) \tt_H W'$,
which is finite dimensional.
Since $\phi$ is injective, $\ker_n\phi$ is
finite dimensional. 
\end{proof}
\begin{lemma}\lbb{ltor0}
Let $\phi\in\Chom(V,W)$ and $h\in H$. If $h$ is not a divisor of zero, 
then{\rm:}

{\rm (i)}
$h\phi=0$ implies $\phi=0;$

{\rm (ii)}
$hv\in\ker\phi$ implies $v\in\ker\phi$.
\end{lemma}
\begin{proof}
Part (i) follows from Eq.~\eqref{hactchom}: if 
$\phi(v)=\tsum_i \, (f_i\tt1)\tt_H w_i$
with linearly independent $w_i$, then
$(h\phi)(v)=\tsum_i \, (h f_i\tt1)\tt_H w_i$ can be zero
only if all $h f_i=0$, which implies $f_i=0$.

Similarly, part (ii) follows from \eqref{hclm}, since we can
write $\phi(v)$ uniquely in the form $\tsum_i \, (1\tt g_i)\tt_H w_i$.
\end{proof}
\begin{corollary}\lbb{ctorcen}
Let $L$ be a 
\psalg, and $M$ be an $L$-module.
Then any torsion element from $L$ acts trivially on $M$, and any
torsion element from $M$ is acted on trivially by $L$.
In particular, the torsion of a Lie \psalg\ is central.
\end{corollary}

\section{Reconstruction of an $H$-Pseudoalgebra
from an $H$-Differential Algebra}\lbb{sreconst}
\subsection{The reconstruction functor $\C$}\lbb{subfc}
Let, as before, $H$ be a cocommutative Hopf algebra and $X=H^*$.
Given a topological left $H$-module $\L$
(where $H$ is endowed with the discrete topology),
let
\begin{equation}\lbb{confl}
\C(\L) = \Hom_H^\cont(X,\L)
\end{equation}
be the space of continuous $H$-homomorphisms.
We define a structure of a left $H$-module on $\C(\L)$ by
\begin{equation}\lbb{hax}
(h\al)(x) = \al(xh).
\end{equation}
Then $\C$ is a covariant functor from the category of topological $H$-modules
to the category of $H$-modules.

\begin{lemma}\lbb{lcl}
{\rm(i)}
The functor $\C$ is left exact{\rm:}
$\C(i)$ is injective if $i$ is injective.

{\rm(ii)}
The functor $\C$ preserves direct sums{\rm:}
$\C(\L_1\oplus\L_2) = \C(\L_1)\oplus\C(\L_2)$.

{\rm(iii)}
Assume that the Hopf algebra $H$ contains nonzero primitive elements.
If $\L$ is finite dimensional over $\Kset$ with discrete topology,
then $\C(\L)=0$.

{\rm(iv)}
If $H=\ue(\dd)$,
then $\C(\L)$ has no torsion as an $H$-module.
\end{lemma}
\begin{proof}
Parts (i) and (ii) are obvious.

(iii) By Kostant's \thref{tkostant}, 
$H=\ue(\dd)\smash\Kset[\Ga]$ with $\dd\ne0$.
If $\L$ is finite dimensional, any continuos
homomorphism $\al\colon X\to\L$ must contain some $\fil_n X$ in its kernel.
Let $h\in\fil^{n-1} \ue(\dd)$ but $h\not\in\fil^{n-2} \ue(\dd)$. 
Then, by \leref{minusi},
$h\fil_n X = X$ so that for each $x\in X$, $x=hy$ for some $y\in\fil_n X$.
This implies $\al(x)=\al(hy)=h(\al(y))=0$, since $\al(y)=0$,
proving part (iii).

Similarly, part (iv) follows from the fact that $Xh=X$ for any
nonzero $h\in \ue(\dd)$.
\end{proof}

If, in addition, $\L$ is a topological Lie $H$-\difalg,
we define $x$-brackets in $\C(\L)$ 
by the formula (cf.\ \eqref{axby2}):
\begin{align}
\lbb{xbconfl}
\xb{\al}{x}{\be}(y) &= [\al(x_{(2)}), \be(y x_{(-1)})]
= \tsum_i\, [\al(x_i), \be(y (h_i S(x)))].
\intertext{This is well defined because the infinite sum
in the right-hand side converges in $\L$. Equation~\eqref{xbconfl}
is also equivalent to (cf.\ \eqref{axby}):}
\lbb{xbconfl2}
[\al(x), \be(y)] &= \xb{\al}{x_{(2)}}{\be}(y x_{(1)})
= \tsum_i\, \xb \al{x_i}\be (y (x S(h_i))).
\end{align}

\begin{proposition}\lbb{pconfl}
For any topological Lie $H$-\difalg\ $\L$,
$\C(\L)$ satisfies properties \eqref{bl1}--\eqref{jid}.
\end{proposition}
\begin{proof}
This can be verified by straightforward but rather tedious computations.
To illustrate them, let us check \eqref{bl1}. By definition, we have:
\begin{displaymath}
\xb{h\al}{x}{\be}(y) = \tsum_i\, [(h\al)(x_i), \be(y (h_i S(x)))]
= \tsum_i\, [\al(x_i h), \be(y (h_i S(x)))],
\end{displaymath}
while
\begin{displaymath}
\xb{\al}{xh}{\be}(y) = \tsum_i\, [\al(x_i), \be(y (h_i S(xh)))].
\end{displaymath}
Hence \eqref{bl1} is a consequence of the following identity:
\begin{equation}\lbb{hhixi}
\tsum_i \, x_i h \tt h_i = \tsum_i \, x_i \tt h_i S(h),
\end{equation}
which can be checked by pairing both sides with $f\tt z\in H\tt X$.
Indeed,
\begin{displaymath}
\tsum_i \, \< x_i h, f \> \< h_i, z \> = \< zh, f \> = \< z, f S(h) \>
= \tsum_i \, \< x_i, f \> \< h_i S(h), z \>.
\end{displaymath}

A more conceptual proof can be given by noticing that
formula \eqref{xbconfl} is the same as the formula for the 
commutator of $H$-pseudolinear maps.
For $\al\in\C(\L)$
consider the family $\ad\al(x) \in \Hom(\L,\L)$ indexed by $x\in X$.
It is easy to see that it satisfies \eqref{bl2phix}.
So, if it also satisfies \eqref{codimphix}, it would give
an  $H$-pseudolinear map from $\L$ to itself.
Although this is not true in general, the argument still
works because all infinite series that appear will be convergent.
(In other words, we embed $\C(\L)$ in a certain completion
of $\gc\L$.)
\end{proof}

In order to have the locality \eqref{codimx}, one has to impose
some restrictions on $\L$. 
In particular, the condition that $\L$ be a linearly 
compact topological Lie algebra will often suffice to guarantee locality of 
$\C(\L)$.

In what follows, we explain how to reconstruct an  $H$-\psalg\ $L$ which is
finitely generated over $H$ from its annihilation Lie algebra $\A(L)$. 
Recall that $\A(L)$ is a linearly compact topological Lie algebra 
(\prref{plielc}).
In many of the proofs we never 
exploit the algebra structure on $\A(L)$, so the corresponding 
statements hold for finite $H$-modules in general.

We let $\what L = \C\A(L) := \Hom_H^\cont(X,X \tt_H L)$.
There is an obvious map
\begin{equation}\lbb{rtorhat}
\Phi\colon L\to\what L, \qquad a\mapsto \al(x) = x \tt_H a.
\end{equation}
It is clear by definitions that $\Phi$ is a homomorphism of $H$-\psalgs\
(or $H$-modules if $L$ is only an $H$-module).

\subsection{The case of free modules}\lbb{subfree}
Let $L$ be a Lie $H$-\psalg\ which is free as an $H$-module:
$L=H\tt L_0$ with the trivial action of $H$ on $L_0$.
Then 
$\L= \A(L) := X \tt_H (H\tt L_0) \simeq X\tt L_0$ as an $H$-module.
%

\begin{proposition}\lbb{prrhat}
When $L$ is a Lie $H$-\psalg\ that is a free $H$-module,
the map $\Phi$ defined by \eqref{rtorhat}
is an isomorphism of Lie $H$-\psalgs.
\end{proposition}
\begin{proof}
To construct the inverse of $\Phi$, identify $\L$ with $X\tt L_0$ and consider 
\begin{displaymath}\lbb{rhattor}
\Psi\colon \what L\to L, \qquad \al\mapsto 
\tsum_i \, S(h_i) \tt (\ep\tt\id)\al(x_i).
\end{displaymath}
Here, as before, $\{h_i\}$, $\{x_i\}$ are dual bases in $H$ and $X$,
and $\ep(x)=\langle 1,x\rangle$ for $x\in X$.
This is well defined, i.e., the sum is finite, because
$\al(x_i) \in \fil_1 X \tt L_0$ for all but a finite number of $x_i$
and $\ep(\fil_1 X) =0$.
Using identity \eqref{hhixi},
it is easy to see that $\Psi$ is $H$-linear.
Next, we have for $a\in L_0$:
\begin{displaymath}
\Psi\Phi(1\tt a) = \tsum_i \, S(h_i) \tt \langle 1,x_i\rangle a 
= S(1)\tt a = 1\tt a,
\end{displaymath}
showing that $\Psi\Phi=\id$. In particular, $\Psi$ is surjective.

Assume that $\Psi(\al)=0$ for some $\al\in\what L$.
This means that $(\langle 1,\cdot\rangle \tt\id)\al(x) = 0$
for any $x\in X$. But then for any $h\in H$, we have:
\begin{align*}
\bigl(\langle S(h),\cdot\rangle \tt\id\bigr) \al(x) 
&= \bigl(\langle 1,\cdot\rangle \tt\id\bigr) \bigl( (h\tt1)\al(x) \bigr)
\\
&= \bigl(\langle 1,\cdot\rangle \tt\id\bigr) \bigl( h(\al(x)) \bigr)
= \bigl(\langle 1,\cdot\rangle \tt\id\bigr)\al(hx)
= 0,
\end{align*}
which implies $\al=0$.
Hence $\Psi$ is injective.
\end{proof}
\begin{remark}\lbb{rrrhat}
If $L$ is only a free $H$-module, then $\Phi$ is an isomorphism of
$H$-modules. Analogous results hold for representations, or
for \as\ $H$-\psalgs.
\end{remark}

\subsection{Reconstruction of a non-free module}\lbb{subnonfree}
Throughout this subsection $L$ will be a (possibly non-free) finitely
generated $H$-module, and $H$ will be the universal enveloping algebra
$\ue(\dd)$ of a finite-dimensional Lie algebra $\dd$.

The natural map $\Phi\colon L \to \what L$ (see \eqref{rtorhat})
is in general neither injective nor surjective. 
However, the induced map $\ph = \A(\Phi)\colon \A(L) \to
\A(\what L)$ has a left inverse 
$\psi\colon x\tt_H \al \mapsto \al(x)$. This
shows that $\A(\Phi)$ is injective, and that $\psi$ is surjective. 

We want to figure out to what extent injectivity and surjectivity of $\Phi$
fail. First of all let us remark that, by \leref{torzerocoef},
every torsion element $a \in L$ has all zero Fourier coefficients, i.e.,
it belongs to the kernel of $\Phi$.
In fact, we have:
%
%
\begin{proposition}\lbb{keristor}
For any finite $H$-module $L$,
the kernel of $\Phi\colon  L \to \what L$ equals the torsion of $L$.
\end{proposition}
\begin{proof}
It remains to show that a non-torsion element $a\in L$
does not lie in the kernel of $\Phi$.
Consider the map $i\colon L \to F$ constructed in
\leref{rmapstofree}. 
Then $i(a)\ne0$.
The map $\A(i)$ induced by $i$ maps the
Fourier coefficient $x \otimes_H a$ of $a$
to the corresponding Fourier coefficient $x \otimes_H i(a)$ of a
nonzero element in the free $H$-module $F$. Now, it is clear from
\prref{prrhat} that $x \otimes_H i(a) \ne 0$ for some
$x \in X$, hence $x\otimes_H a$ must be nonzero too. 
\end{proof}
\begin{corollary}\lbb{tor0coef}
A finite $H$-module $L$ is torsion iff $X \tt_H L = 0$.
\end{corollary}
\begin{corollary}\lbb{a=iff=}
Let $M,N$ be finite $H$-modules, $f\colon M \to N$ be an $H$-linear map, 
and assume $N$ to be torsionless. Then $\A(f) = 0$ if and only if $f=0$.
\end{corollary}
\begin{proof}
$\A(f) = 0$ means that $X\tt_H f(M) = 0$, hence $f(M)\subset N$ is torsion.
\end{proof}
\begin{remark}\lbb{rcenttor}
By \coref{ctorcen}, the torsion of a Lie $H$-\psalg\ $L$ is always central, 
hence the map $\Phi$ is injective if $L$ is centerless.
\end{remark}
\begin{remark}\lbb{phiinverse}
If $\Phi$ is an isomorphism, then $\Phi^{-1}$ induces $\psi$, i.e.,
$\psi = \A(\Phi^{-1})$. \coref{a=iff=} tells us that if $L$ is torsionless and
$\psi$ is induced by some map $\Psi$, then $\Phi$ is an isomorphism and $\Psi =
\Phi^{-1}$.
\end{remark}
\begin{proposition}
For any map of finite $H$-modules $f\colon M \to N$,
the following conditions are equivalent{\rm:}
\begin{enumerate}
\item $N/f(M)$ is torsion.
\item $\A(f)\colon \A(M) \to \A(N)$ is surjective.
\item $\gw \coker \A(f) < \dim \dd$.
\end{enumerate}
\end{proposition}
\begin{proof}
To show the equivalence of (1) and (2), 
it is enough to tensor the exact sequence 
$M \xrightarrow{f} N \to N/f(M) \to 0$ 
with $X$, getting the exact sequence 
$\A(M) \to \A(N) \to X \tt_H N/f(M) \to 0$, and to apply
\coref{tor0coef}.

Assume that (3) holds but $N/f(M)$ is not torsion. Then it contains a nonzero
element $a$ which generates a free $H$-module. Then $X\tt_H a \simeq X$
has growth $\dim \dd$, which is a contradiction.
\end{proof}

\subsection{Reconstruction of a Lie \psalg}\lbb{subrecalg}
Now let $L$ be a Lie $H$-\psalg\ which is finite as an $H$-module.
Again, $H=\ue(\dd)$ will be the universal enveloping algebra
of a finite-dimensional Lie algebra $\dd$.
Let $\L = \A(L)$ be the  annihilation Lie algebra of $L$, and
$\what L=\C(\L)$, as before.

By \prref{pconfl}, $\what L$ satisfies all the properties
of a Lie $H$-\psalg\ except the locality \eqref{codimx}.
An indirect way to establish the locality property for
$\what L$ is by embedding it in the bigger
(local) Lie \psalg\ $\gc L$.
In order to map $\what L$ to $\gc L$, we need to assign to each element of
$\what L$ a pseudolinear map from $L$ to itself.

This can be done as follows.
Recall that $\L$ acts on $L$ by $(x\tt_H a) \cdot b = \xb axb$ for $a,b\in L$,
$x\in X$ (see \prref{preplal}).
Now in terms of $x$-products the action of $\what L$ on $L$
is given by $\xp \al x b = \al(x) \cdot b$.
The locality condition $\xp \al x b = 0$ for $x\in\fil_n X$, $n\gg0$
is satisfied because $\al$ is continuous and $L$ is a discrete topological
$\L$-module (see \prref{preplal}).
All the other axioms of a Lie \psalg\ representation follow easily from
definitions. 

We now need to find conditions for the above-defined
$j\colon\what L \to \gc L$ to be injective.
Then $\what L$ would embed into $\gc L$, which will show locality.

\begin{lemma}\lbb{lkerj}
If $L$ is torsionless,
the kernel of the above-defined $j\colon \what L \to \gc L$
consists of all elements $\al$ such that $\al(X)$ 
is contained in the center of $\L$.
\end{lemma}
\begin{proof}
Since $L$ is torsionless, $\Phi$ is injective by \prref{keristor}.
Hence, for $a,b\in L$, $x\in X$, one has 
$\xb axb = 0$ iff $[a_xb]_y = 0$ for all $y\in X$.
By \eqref{axby}, \eqref{axby2}, this is equivalent to
$[a_x, b_y] = 0$. Hence, for $l\in\L$, $l\cdot b =0$ for all $b$ iff
$l$ lies in the center of $\L$.
Now $\al\in\what L$ is in the kernel of $j$ iff $\al(x)\cdot b =0$
for all $x$ and $b$, which means that $\al(x)$ is central for all $x$.
\end{proof}
\begin{lemma}\lbb{finitectrembeds}
If $L$ is torsionless and $\L=\A(L)$ has a finite-dimensional center, then
$j\colon \what L \to \gc L$ is injective.
\end{lemma}
\begin{proof}
Let $\al\in\what L$ be in the kernel of $j$;
then by the previous lemma $\al(X)$
is contained in the center of $\L$. The latter is finite dimensional
by assumption, so the kernel $N$ of $\al$ is of finite codimension
in $X$. This implies that $N$ is open in $X$, and it contains
$\fil_i X$ for some $i$.
Let $h \in \fil^{i+1} H$ but $h \notin \fil^i H$; then 
by \leref{minusi}, $h \fil_i X = \fil_{-1} X = X$.
Since $\al$ is $H$-linear, $N$ is an $H$-submodule of $X$.
Then $X = h \fil_i X \subset h N \subset N$, therefore $N=X$ and $\al=0$.
\end{proof}
\begin{proposition}\lbb{rislocal}
Let $L$ be a Lie $H$-\psalg\ which is finite 
and torsionless as an $H$-module.
If its annihilation Lie algebra $\L = \A(L)$ has a finite-dimensional center, 
then $\what L=\C(\L)$ is a Lie $H$-\psalg\ containing $L$ as an ideal.
\end{proposition}
\begin{proof}
The only thing that remains to be checked is the locality property for
$\what L$. It follows from that of $\gc L$, since in this case
$j\colon \what L \to \gc L$ is injective.
\end{proof}

In the proof of \leref{lkerj} we have shown that, if $L$ is finite 
and torsionless, the kernel of the action of $\A(L)$ on $L$
is exactly the center of $\A(L)$. This implies the following
result which was used in the proof of \thref{tannihvec}.

\begin{lemma}\lbb{linl1}
Let $i\colon L\injto L_1$ be an injective map of Lie $H$-\psalgs,
and assume that $L$ is finite and torsionless. Then the kernel
of the induced map $\A(i)\colon \A(L)\to \A(L_1)$ is contained
in the center of $\A(L)$.
\end{lemma}
\begin{proof}
The kernel of $\A(i)$ acts trivially on $L_1$ and hence on $L$.
\end{proof}

In \seref{scders} we will need the following lemma.

\begin{lemma}\lbb{lsubgcv}
Let $L$ be a subalgebra of $\gc V$ for some finite $H$-module $V$
$(L$ may be infinite$)$. Then the map 
$\Phi \colon L \to \widehat L$ is injective.
\end{lemma}
\begin{proof}
Assume that $a$ belongs to $\ker\Phi$; then all $x \otimes_H a = 0$, 
$x \in X$.
This implies that all Fourier coefficients $a_x \in \End V$ of $a\in\gc V$
are zero, hence $a=0$.
\end{proof}


For any topological Lie $H$-\difalg\ $\L$, we have a natural homomorphism
$\psi\colon \A\C(\L) \to \L$, given by $x\tt_H a \mapsto a(x)$
for $a\in \C(\L)$, $x\in X$.
The map $\psi$ does not need to be surjective, but we have a good control on 
injectivity, which can sometime prove useful.
\begin{lemma}
The kernel of $\psi\colon \A\C(\L) \to \L$ lies in the center of $\A\C(\L)$.
\end{lemma}
\begin{proof}
Follows easily from \eqref{alby1} and \eqref{alby2}. Say that $\alpha$ lies
in the kernel of $\psi$. Since $\psi$ is a homomorphism of Lie $H$-\difalgs, 
its kernel is an $H$-stable ideal of $\A\C(\L)$.
Then by \eqref{alby1}, $(\alpha\cdot b)_y\in \ker\psi$ for all 
$b\in \C(\L)$, $y\in X$, because
in the right-hand side all elements $h_i\alpha$ lie in $\ker \psi$. This
means that $(\alpha\cdot b)(y) = 0$ for all $y\in X$, hence
$\alpha\cdot b = 0$ for every $b \in \C(\L)$.
Now, use this in \eqref{alby2} to obtain that $\alpha$ is central.
\end{proof}

\section{Reconstruction of Pseudoalgebras of Vector Fields}\lbb{srecvec}
In this section, we show that the reconstruction procedure
of \seref{sreconst}, when applied to the primitive Lie algebras
of vector fields (or current algebras over them), gives 
the primitive pseudoalgebras of vector fields defined in \seref{svect}
(or current \psalgs\ over them).

As before, $\dd$ will be an $N$-dimensional Lie algebra, and
$H=U(\dd)$ its universal enveloping algebra. 
$\L$ will be a Lie algebra provided with an action of $\dd$ and
a filtration by subspaces $\L=\L_{-1}\supset\L_0\supset\dotsm$.
When $\L$ is a subalgebra of $W_N$, it will always be considered
with the filtration induced by the canonical filtration of $W_N$.

The Lie algebra $\Der\L$ of derivations of $\L$ has the
induced filtration: 
\begin{displaymath}
(\Der\L)_i := \{ d\in\Der\L \st d (\L_j) \subset \L_{i+j}
\quad\forall j \}.
\end{displaymath}
The action of $\dd$ is called {\em transitive\/} if the composition of the
homomorphism $\dd\to\Der\L$ and the projection 
$\Der\L \to \Gr_{-1}(\Der\L) := \Der\L / (\Der\L)_0$
is a linear isomorphism. This is equivalent to the following two
conditions: $\dd\injto\Der\L$ intersects $(\Der\L)_0$ trivially and
$\dim \Gr_{-1}(\Der\L) = N$.

\begin{lemma}\lbb{ldatrans}
Let $L$ be a current Lie $H$-\psalg\ over a finite-dimensional simple Lie 
algebra or over one of the primitive \psalgs\ of vector fields. Then the
action of $\dd$ on its annihilation Lie algebra $\L=\A(L)$ is transitive.
\end{lemma}
\begin{proof}
By \thref{tannihvec}, $\L=\O_r\what\tt \L'$ is a current Lie algebra over 
$\L'$, where 
$\L'$ is either a finite-dimensional simple Lie algebra
$\g$ (for $r=N=\dim\dd$),
or one of the Lie algebras of vector fields $W_{N'}$, $S_{N'}$, 
$P_{N'}$ or $K_{N'}$ ($N'=N-r$).
In particular, we know that $\dim \Gr_{-1}(\Der\L) = N$.
By \leref{lssurj}, a sufficiently high power of any nonzero element $a\in\dd$
maps any given open subspace of $\L$ surjectively onto $\L$. This cannot
hold if $a$ belongs to $(\Der\L)_0$, therefore
$\dd\to\Der\L$ is injective and the image of $\dd$
intersects $(\Der\L)_0$ trivially.
Comparing the dimensions,
we get that $\dd\to\Gr_{-1}(\Der\L)$ is an isomorphism.
\end{proof}

The main results of this section can be summarized by \thref{trecvec} below.
Its proof follows from Sections~\ref{swn}--\ref{sknpn}.

\begin{theorem}\lbb{trecvec}
Let $\L=\O_r\what\tt \L'$ be a current Lie algebra over $\L'$,
where $\L'$ is a simple linearly compact Lie algebra of growth
$N'=N-r$. 
Assume that $\dd$ acts transitively on $\L$.
%
%
Then there is a codimension $r$ subalgebra $\dd'$ of $\dd$,
acting transitively on $\L'$,
such that the $H$-\psalg\ $\C(\L)$ is isomorphic to a current
\psalg\ over the $H'$-\psalg\ $\C(\L')$, where $H'=\ue(\dd')$.
Moreover, $\C(\L')$ is either a finite-dimensional simple Lie 
algebra {\rm(}and $\dd'=0${\rm)} or one of the primitive $H'$-\psalgs\
of vector fields $\Wdp$, $\Sdp$, $\Hdp$ or $\Kdp$.
%
\end{theorem}
\subsection{Reconstruction from $W_N$}\lbb{swn}
Recall that $X=H^*$ can  be identified with $\O_N=\Kset[[t_1,\dots,t_N]]$.
The action of $\dd$ on $X$ gives an action on $\O_N$ in terms of linear
differential operators, i.e., an embedding $\dd\injto W_N = \Der\O_N$
which we call the {\em canonical embedding\/} of $\dd$ in $W_N$.
Note that this embedding is transitive, i.e., $\dd\subset W_N$
is complementary to $\fil_0 W_N$.

A structure of an $H$-\difalg\ on $W_N$ is equivalent to a
transitive action of
$\dd$ on $W_N$ by derivations. Since $\Der W_N = W_N$, this is the same
as a transitive embedding of $\dd$ in $W_N$. By \prref{existuniq},
any two such embeddings are
equivalent, i.e., conjugate by an automorphism of $W_N$. 
With the canonical action of $\dd$,
$W_N$ becomes isomorphic to the annihilation
algebra of the Lie \psalg\ $\Wd$ defined in \seref{subw}. Since $\Wd$ is
a free $H$-module, \prref{prrhat} shows that the
reconstruction of $W_N$ is $\Wd$, i.e., $\C(W_N)=\Wd$.

\subsection{Reconstruction from subalgebras of $W_N$}\lbb{ssubwn}
Let $\L$ be a linearly compact Lie subalgebra of $W_N$, with the induced
filtration and with a transitive action of $\dd$ on it.
After an automorphism of $W_N$, we can
assume that the action of $\dd$ is the canonical one. Then $\C(\L)$ is a
subalgebra of $\Wd=\C(W_N)$, because the functor $\C$ is left exact.
Below we will be concerned with the case when $\L$ is the subalgebra
consisting of vector fields annihilating some differential form.

Let $w\in\Om^n(\dd)$ be a pseudoform, and $I\subset H$ be a right
ideal. We denote by $W(\dd,w,I)$ the set of all elements
$\al\in\Wd=H\tt\dd$ such that
\begin{equation}\lbb{wdwi}
\al*w \in (H\tt I)\tt_H \Om^n(\dd).
\end{equation}
It is easy to check that $W(\dd,w,I)$ is a subalgebra of $\Wd$.
\begin{lemma}\lbb{lwdwi}
Let $\om\in\Om^n_X$ be a differential form, 
and $W_N(\om)$ be the Lie subalgebra of $W_N$ consisting of vector fields
annihilating $\om$. 
If $\om=y\tt_H w$ for some $y\in X$, $w\in\Om^n(\dd)$, then
$\C(W_N(\om))$ is isomorphic to the Lie \psalg\ $W(\dd,w,I)$ where
$I=\{h\in H \st y h =0\}$.
\end{lemma}
\begin{proof}
As was already remarked, $\C(W_N(\om))$ is a subalgebra of $\Wd$.
Since $\Om^n(\dd) = H\tt\bigwedge^n \dd^*$ is a free $H$-module,
we have $(H\tt H)\tt_H \Om^n(\dd) \simeq H\tt H\tt\bigwedge^n \dd^*$.
For $\al\in\Wd$, write $\al*w = \sum_i (f_i\tt g_i)\tt_H w_i$
with $f_i,g_i\in H$ and linearly independent $w_i\in\bigwedge^n \dd^*$.
Then for any $x\in X$ we have (cf.\ \eqref{ahxbhy}):
\begin{displaymath}\lbb{al*om}
L_{x\tt_H\al} \om = \tsum_i (x f_i)(y g_i) w_i.
\end{displaymath}
This is zero for any $x$ iff $y g_i = 0$ for all $i$,
which means $g_i \in I$.
\end{proof}
\subsection{Reconstruction from current algebras over $W_{N'}$}\lbb{scurwn}
Let now $\L=\O_r\what\tt W_{N'}$ be a current algebra over $W_{N'}$,
and $\dd$ be an $N=N'+r$ dimensional Lie algebra acting transitively on $\L$.
Then, by \prref{pder}, 
$\dd\injto\Der\L = W_r\tt1 + \O_r\what\tt W_{N'} \subset W_N$.
The Lie algebra $\L$ is described as the subalgebra of $W_N$
consisting of vector fields annihilating the functions
$t_{N'+1},\dots,t_N$, hence it is an intersection of algebras of the
form $W_N(f)$ ($f\in\Om^0_X = X$), see \seref{ssubwn}.

After an automorphism of $W_N$, we can
assume that the action of $\dd$ on it is the canonical one. Then $\L$
becomes the intersection of $W_N(f_i)$ ($i=N'+1,\dots,N$) where $f_i\in X$ is
the image of $t_i$. 
Now \leref{lwdwi} implies that $\C(\L)=W(\dd,1,I)$ 
where $1\in\Om_0(\dd)=H$ and 
$I=\{h\in H \st f_i  h =0 \;\; (i=N'+1,\dots,N) \}$.

Recall that for $\al\in\Wd=H\tt\dd$, its action on $1\in H$ is given by 
$\al*1 = -\al\tt_H 1 \equiv -\al$. Therefore $\al\in W(\dd,1,I)$
iff $\al$ belongs to $(H\tt\dd) \cap (H\tt I) = H\tt\dd'$, where
the intersection $\dd'=\dd\cap I$ is a Lie subalgebra of $\dd$.
Then $H'=U(\dd')$ is a Hopf subalgebra of $H$, and 
$H\tt\dd' \simeq H\tt_{H'} (H'\tt\dd')$ is a current \psalg\ over
$H'\tt\dd' = \Wdp$. We have thus proved the following lemma.

\begin{lemma}\lbb{lreccurw}
The reconstruction of a current Lie algebra over $W_{N'}$, 
provided with a transitive action of a Lie algebra $\dd$, is a current Lie
$H$-\psalg\ over $\Wdp$ where $\dd'$ is an $N'$-dimensional
Lie subalgebra of $\dd$.
\end{lemma}

This result is a special case of \leref{lreccurl} below.

\subsection{Solving compatible systems of linear differential equations}
\lbb{seqs2}
Let $A$ be any associative $\Kset$-algebra, and let
$\O_r=\Kset[[t_1,\dots,t_r]]$, $W_r=\Der\O_r$, as before.
For fixed $n\ge0$,
let $f_i(t) \in A[[t_1,\dots,t_r]]$ $(i=1,\dots,r+n)$
be formal power series with coefficients in $A$,
where $t=(t_1,\dots,t_r)$.
Note that $W_r$ acts on $A[[t_1,\dots,t_r]]$ by derivations.

Given $r+n$ linear differential operators 
$D_1,\dots,D_{r+n}\in W_r$,
consider the following system of differential equations
for an unknown $y(t)\in A[[t_1,\dots,t_r]]$:
\begin{equation}\lbb{diggf1}
D_i(y(t)) = y(t) f_i(t), \qquad i=1,\dots,r+n.
\end{equation}
We assume that the operators $D_i$ satisfy
\begin{equation}\lbb{di2}
[D_i,D_j] = \tsum_k \, c_{ij}^k(t) D_k 
\quad\text{with}\;\; c_{ij}^k(t) \in\O_r;
\end{equation}
in other words, the space of all operators of the form $\tsum_i \, p_i(t)D_i$
with $p_i(t) \in\O_r$ is a Lie algebra.

Suppose we have found a solution to the system \eqref{diggf1}.
Combining equations \eqref{diggf1} and \eqref{di2}, we get:
\begin{align*}
[D_i,D_j] (y) &= D_i D_j (y) - D_j D_i (y)
= D_i (y f_j) - D_j (y f_i)
\\
&= y f_i f_j + y D_i(f_j)
- y f_j f_i - y D_j(f_i),
\intertext{and}
[D_i,D_j] (y) &= \tsum_k \, c_{ij}^k D_k(y)
= \tsum_k \, c_{ij}^k y f_k.
\end{align*}
The system \eqref{diggf1} is called {\em compatible\/} if
\begin{equation}\lbb{di3}
[f_i(t), f_j(t)] + D_i(f_j(t)) - D_j(f_i(t))
= \tsum_k \, c_{ij}^k(t) f_k(t)
\quad\text{for all}\;\; i,j.
\end{equation}
When $y(t)$ is not a divisor of zero in $A[[t_1,\dots,t_r]]$
the compatibility of the system is a necessary condition for having a solution.
The compatibility \eqref{di3} is equivalent to saying
that $\tsum_i \, p_i(t)D_i \mapsto \tsum_i \, p_i(t)(D_i + f_i(t))$
is a homomorphism of Lie algebras.

We will be interested in solving a more general system of
equations than \eqref{diggf1}. Before formulating it, let us note
that the above remarks have obvious analogues for systems of the form
\begin{equation}\lbb{diggf2}
D_i(z(t)) = -h_i(t) z(t), \qquad i=1,\dots,r+n
\end{equation}
with $z(t), h_i(t) \in A[[t_1,\dots,t_r]]$.
The compatibility of \eqref{diggf2} is equivalent to \eqref{di3}
with $f_i$ replaced by $h_i$.

Now consider the system
\begin{equation}\lbb{diggf}
D_i(g(t)) = g(t) f_i(t) - h_i(t) g(t), \qquad i=1,\dots,r+n
\end{equation}
for an unknown $g(t)\in A[[t_1,\dots,t_r]]$.
We will show it has a solution, provided that both
\eqref{diggf1} and \eqref{diggf2} are compatible and
some initial conditions at $t=0$ are satisfied.
(The compatibility of \eqref{diggf1} and \eqref{diggf2}
implies the compatibility of \eqref{diggf}.)

\begin{proposition}\lbb{psolve}
In the above notation, let the operators $D_i \in W_r$ satisfy
\eqref{di2} and
\begin{equation}\lbb{di1}
D_i|_{t=0} = \begin{cases}
\d_{t_i}, \quad &1\le i\le r,
\\
0, \quad &r+1\le i\le r+n.
\end{cases}
\end{equation}
Assume that the systems \eqref{diggf1} and \eqref{diggf2} are compatible
{\rm(}cf.\ \eqref{di3}{\rm)}, and that
\begin{equation}\lbb{di4}
f_i(0) = h_i(0), \qquad r+1\le i\le r+n.
\end{equation}
Then the system \eqref{diggf} has a unique
solution $g(t)\in A[[t_1,\dots,t_r]]$ for any given initial
condition $g(0)\in A$ which commutes with $f_i(0)$ $(r+1\le i\le r+n)$.
\end{proposition}
\begin{proof}
For $r=0$, both sides of \eqref{diggf} are trivial.
For $r\ge1$, we will proceed by induction on $r$.

First of all, note that the compatibility or solvability of the systems
\eqref{diggf1}, \eqref{diggf2} or \eqref{diggf}
does not change when we apply an automorphism of $\O_r$.
The same is true when we make an elementary transformation:
multiply one equation by a function (an element of $\O_r$) 
and add it to another equation. For example, we can replace all
$D_i$ ($i\ne r$) by $D_i-p_i(t)D_r$, 
and correspondingly
$f_i(t)$ by $f_i(t)-p_i(t)f_r(t)$
and $h_i(t)$ by $h_i(t)-p_i(t)h_r(t)$,
as long as we do not violate (\ref{di1}, \ref{di4}).

Any vector field $D_r\in W_r$ satisfying $D_r|_{t=0} = \d_{t_r}$
can be brought to $\d_{t_r}$
after an automorphism of $\O_r$, so we will assume that
$D_r = \d_{t_r}$. Replacing $D_i$ ($i\ne r$) by 
$D_i-D_i(t_r)D_r$, we can assume in addition that $D_i(t_r) = 0$
for $i\ne r$.

Now it makes sense to put $t_r=0$ in the equations with $i\ne r$
in \eqref{diggf}. Let us denote $\bar D_i = D_i|_{t_r=0}$,
$\bar f_i(\bar t) = f_i(t_1,\dots,t_{r-1},0)$,
$\bar h_i(\bar t) = h_i(t_1,\dots,t_{r-1},0)$,
$\bar t = (t_1,\dots,t_{r-1})$.
Consider the reduced system 
\begin{equation}\lbb{digred}
\bar D_i(\bar g(\bar t)) 
= \bar g(\bar t) \bar f_i(\bar t)
- \bar h_i(\bar t) \bar g(\bar t), 
\qquad i=1,\dots,r-1,r+1,\dots,r+n
\end{equation}
for an unknown $\bar g(\bar t) \in A[[t_1,\dots,t_{r-1}]]$.
Note that, since $D_i(t_r) = \de_{ir}$, we have:
$[D_i,D_j](t_r) = 0$ for any $i,j$, hence
$[D_i,D_j]$ does not contain $D_r$. In particular,
putting $t_r=0$ we see that the operators $\bar D_i$
satisfy \eqref{di2}. The other assumptions of the proposition
are also easy to check, so by induction
the system \eqref{digred} has a solution $\bar g(\bar t)$. 

The equation
\begin{equation}\lbb{digr1}
\d_{t_r} g(t) = g(t) f_r(t) - h_r(t) g(t)
\end{equation}
has a unique solution $g(t)$ satisfying the initial condition
\begin{equation}\lbb{digr2}
g(t_1,\dots,t_{r-1},0) = \bar g(\bar t).
\end{equation}
We claim that this $g(t)$ is then a solution of the system \eqref{diggf}.
Indeed, it satisfies \eqref{diggf} for $t_r=0$.
Next, we compute for $i\ne r$
(using \eqref{digred}, \eqref{digr1}, and the compatibility
of \eqref{diggf1}, \eqref{diggf2}):
\begin{align*}
D_r &D_i(g) |_{t_r=0}
= [D_r,D_i](g) |_{t_r=0} + D_i D_r(g) |_{t_r=0}
\\
&= \tsum_j \, \bar c_{ri}^j \bar D_j \bar g
+ \bar D_i (\bar g \bar f_r - \bar h_r \bar g)
\\
&= \tsum_j \, \bar c_{ri}^j (\bar g \bar f_j - \bar h_j \bar g)
+ (\bar g \bar f_i - \bar h_i \bar g) \bar f_r
+ \bar g \bar D_i(\bar f_r) 
- \bar D_i(\bar h_r) \bar g - \bar h_r (\bar g \bar f_i - \bar h_i \bar g)
\\
&= g \bigl( D_r(f_i) + f_r f_i \bigr) |_{t_r=0}
 - \bigl( D_r(h_i) - h_i h_r \bigr) g |_{t_r=0}
 - \bar h_i \bar g \bar f_r - \bar h_r \bar g \bar f_i
\\
&= D_r (g f_i - h_i g) |_{t_r=0}.
\end{align*}
This shows that 
$\d_{t_r} (D_i(g) - g f_i + h_i g) |_{t_r=0} = 0$.
We can apply the same argument with $[D_r,D_i]$ instead of $D_i$,
and so on, to show that all derivatives with respect 
to $t_r$ vanish at $t_r=0$.
\end{proof}
\begin{remark}\lbb{rinv}
Any solution $g(t)$ of the system \eqref{diggf}, such that
$g(0)$ is invertible in $A$, is invertible in $A[[t_1,\dots,t_r]]$.
Its inverse $g(t)^{-1}$ satisfies \eqref{diggf} 
with $f_i \leftrightarrow -h_i$.
\end{remark}

\subsection{Reconstruction from a current Lie algebra}\lbb{sreccur}
Let $\L'$ be a simple linearly compact Lie algebra,
and let $\L=\O_r\what\tt \L'$ be a current algebra over $\L'$.
The filtration by subspaces $\L'=\L'_{-1}\supset\L'_0\supset\dotsm$
and the canonical filtration of $\O_r$ give rise to the product
filtration of $\L$.
Assume that $\dd$ acts on $\L$
transitively by derivations. By \prref{pder}, we have
$\Der\L = W_r\tt1 + \O_r\what\tt \Der\L'$.

Denote by $j$ the embedding $\dd\injto\Der\L$, and by $p$ the projection
$\Der\L\to W_r$. The preimage $\dd' := (pj)^{-1}(\fil_0 W_r)$ 
is a Lie subalgebra
of $\dd$ of codimension $r$. We have 
$\dd' \injto \fil_0 W_r\tt1 + \O_r\what\tt \Der\L'$. The latter contains
$\fil_0 W_r\tt1 + \fil_0\O_r \what\tt \Der\L'$ as an ideal, hence we get
a Lie algebra homomorphism $j'\colon\dd'\to\Der\L'$. 
It leads to a transitive action of $\dd'$ on $\L'$,
because the action of $\dd$ on $\L$ is transitive.

\begin{lemma}\lbb{lunj}
Any two transitive embeddings $j\colon\dd\injto\Der\L$, that induce the same
subalgebra $\dd'$ and the same $j'\colon\dd'\to\Der\L'$, are equivalent
up to an automorphism of\/ $\Der\L$.
\end{lemma}
\begin{proof}
Let us choose a basis $\{\d_i\}$ of $\dd$ and write 
$j(\d_i) = D_i + f_i(t)$ ($i=1,\dots,N=r+n=\dim\dd$)
where $D_i \in W_r$ and $f_i(t) \in \O_r\what\tt \Der\L'$,
$t=(t_1,\dots,t_r)$. Note that $D_i=(pj)(\d_i)$ and $pj\colon\dd\to W_r$
is a Lie algebra homomorphism. We can choose the basis $\{\d_i\}$ in such a way
that $D_i|_{t=0} = \d_{t_i}$ for $1\le i\le r$, and $D_i|_{t=0} = 0$ for
$r+1\le i\le r+n$. Then $\{\d_i\}_{i=r+1,\dots,r+n}$ is a basis of $\dd'$.
Moreover, note that $j'\colon\dd'\to\Der\L'$ is given by 
$j'(\d_i) = f_i(0)$.

Let $\ti j$ be another transitive embedding of $\dd$ into $\Der\L$.
Since, by \prref{existuniq},
the homomorphism $pj$ is uniquely determined by the choice of 
$\dd'$, we can assume that $(p\ti j)(\d_i) = D_i$. Then
$\ti j(\d_i) = D_i + h_i(t)$ for some
$h_i(t) \in \O_r\what\tt \Der\L'$. By assumption, $j' = \ti j'$, hence
$f_i(0) = h_i(0)$ for $r+1\le i\le r+n$.

Now we want to find an automorphism $g(t)\in\O_r\what\tt \Aut\L'$
such that $g(0)=\id$ and
$g(t) \circ (D_i + f_i(t)) = (D_i + h_i(t)) \circ g(t)$.
This equation is equivalent to \eqref{diggf}, and it is easy to see that
all conditions of \prref{psolve} are satisfied: for example, the system
\eqref{diggf1} is compatible because $j$ is a homomorphism.
This completes the proof.
\end{proof}

Now given the embedding $j'\colon\dd'\to\Der\L'$ we can consider
the reconstruction $L':=\C(\L')$ which is a Lie $H'$-\psalg, where
$H'=\ue(\dd')$. Given $L'$ we can take the current $H$-\psalg\
$L := \Cur L' = H\tt_{H'}L'$. Since its annihilation Lie algebra 
$\A(L)$ is isomorphic to $\L$, we get an embedding 
$\ti j\colon\dd\injto\Der\L$. It induces the same embedding $j'$
as our initial $j$,
so by the previous lemma $j$ and $\ti j$ are equivalent. But then the
reconstruction $\C(\L)$ of $\L$ provided with $j$ is isomorphic
to the reconstruction of $\L$ provided with $\ti j$, which is $L$.
This can be summarized as follows.

\begin{lemma}\lbb{lreccurl}
The reconstruction $\C(\L)$ of a current Lie algebra $\L=\O_r\what\tt \L'$
over a simple linearly compact Lie algebra $\L'$, 
provided with a transitive action of a Lie algebra $\dd$, 
is a current Lie $H$-\psalg\ over the $H'$-\psalg\
$\C(\L')$, where $H'=\ue(\dd')$ and $\dd'$ is a
Lie subalgebra of $\dd$ of codimension $r$.
\end{lemma}

\subsection{Reconstruction from $S_N$}\lbb{ssn}
Now consider $S_N$ with a transitive action of $\dd$ on it. Since 
$\Der S_N = CS_N \subset W_N$, we have $\dd\injto W_N$.
After an automorphism of $W_N$, we can assume $\dd\injto W_N$
is the canonical embedding, while $S_N$ becomes 
$W_N(\om) \, (\equiv S_N(\om))$ where 
$\om\in\Om^N_X$ is a volume form. We can write $\om=y\tt_H v$ with 
$y\in X$ and $v\in\bigwedge^N\dd^*$. Then, by \leref{lwdwi},
the reconstruction of $W_N(\om)$ is $W(\dd,v,I)$ where 
$I=\{h\in H \st y h =0\}$ is as before.

The action of $\Wd$ on $v$ is given by \eqref{fa*v}. In the notation of
\seref{subs}, we have for $\al\in\Wd$:
\begin{displaymath}\lbb{al*v1}
\al*v = -(\Div^{\tr\ad}(\al)\tt1 + \al) \tt_H v.
\end{displaymath}
This shows that $\al\in W(\dd,v,I)$ iff 
$\Div^{\tr\ad}(\al)\tt1 + \al \in H\tt I$. 

Note that, since $\om\ne0$, we have $I\cap\Kset=0$.
The intersection $I\cap(\dd+\Kset)$ is a Lie algebra. 
The projection $\pi\colon (\dd+\Kset) \to \dd$ is a Lie algebra
homomorphism, which maps $I\cap(\dd+\Kset)$ isomorphically
onto a subalgebra $\dd'$ of $\dd$. 
The inverse isomorphism $\dd' \to I\cap(\dd+\Kset)$
is given by $a\mapsto a+\chi(a)$ for some linear functional
$\chi\colon\dd'\to\Kset$ which vanishes on $[\dd',\dd']$.
Conversely, any such $\chi$ gives rise to an isomorphism as above.

For $\be\in H\tt (\dd+\Kset)$, the equation
$\be\in H\tt I$
is equivalent to the following two conditions:
$(\id\tt\pi)(\be) \in H\tt\dd'$
and
$(\id\tt\pi + \id\tt\chi\pi)(\be) = \be$.
Applying this for $\be=\Div^{\tr\ad}(\al)\tt1 + \al$,
we get $(\id\tt\pi)(\be) = \al \in H\tt\dd'$
and $(\id\tt\pi + \id\tt\chi\pi)(\be) = \al + (\id\tt\chi)(\al) = \be$.
The latter equation is equivalent to
$(\id\tt\chi)(\al) = \Div^{\tr\ad}(\al)\tt1$,
i.e.\ to $\Div^{\tr\ad - \chi}(\al) = 0$.
We have proved:

\begin{lemma}\lbb{lrecs}
The reconstruction of the Lie algebra  $S_N$, 
provided with a transitive action of a Lie algebra $\dd$, is a current Lie
$H$-\psalg\ over $\Sdp$ where $\dd'$ is a Lie subalgebra of $\dd$
and $\chi'$ is a linear functional $\dd'\to\Kset$ 
which vanishes on $[\dd',\dd']$.
\end{lemma}

In fact, one can show that in this case $\dd'=\dd$, i.e.,
$\dim I\cap(\dd+\Kset) = N$, but the above statement is sufficient
for our purposes.

\subsection{Reconstruction from $K_N$ and $H_N$}\lbb{sknpn}
Now let $\L$ be one of the Lie algebras $K_N$ or $P_N$, together 
with a transitive action of $\dd$ on it. We know from \seref{scartan} that
as a topological vector space $\L$ is homeomorphic to $X$. 
Since, by \prref{existuniq}, all transitive
actions of $\dd$ on $X$ are equivalent, $\L$ is isomorphic to the ``canonical''
$H$-module $X$, i.e., we may assume that
the embedding $\dd\injto W_N = \Der X$ is the canonical one
(\seref{swn}).
Then, by \prref{prrhat}, the reconstruction of $\L$ is
isomorphic to $H$ as an $H$-module. In other words, $\C(\L)$ is a free
$H$-module of rank one.

\begin{lemma}\lbb{lrech}
The reconstruction of the Lie algebras  $K_N$ and $H_N$,
provided with a transitive action of a Lie algebra $\dd$, 
is a free $H$-module of rank one.
\end{lemma}
\begin{proof}
It is enough to show that the reconstruction functor $\C$ gives the
same result on the topological $H$-modules $X = P_N$ and 
$X/\kk = H_N$. In order to do so, we must show that every $H$-linear
continuous homomorphism of $X$ to $X/\kk$ can be obtained from a unique
$H$-linear continuous homomorphism of $X$ to itself by composing with
the canonical projection $X \to X/\kk$.

Since $X/\kk$ is linearly compact, by \reref{bidual}
there is a bijection between
$\Hom_H^\cont(X, X/\kk)$ and $\Hom_H((X/\kk)^*, H)$. The space
$(X/\kk)^*$ is nothing but the augmentation ideal $H_+ = \ker \ep \subset H$.
Therefore we are
reduced to show that every $H$-linear map 
$\phi\colon H_+ \to H$ is
a restriction of a unique $H$-linear map $H \to H$.

An $H$-linear $\phi\colon H_+ \to H$ is determined by its value on
$\dd \subset H_+$. If $a,b \in \dd$, then $ab-ba = [a,b]$, hence
$a\phi(b) - b\phi(a) = \phi([a,b])$.
Let $d$ be the maximal degree of $\phi(a)$ for $a\in\dd$.
Then $a\phi(b) = b\phi(a)$ modulo
$\fil^d H$. This means that there exists some $\alpha \in \fil^{d-1} H$
such that for every $a\in \dd$, $\phi(a) = a \alpha$ modulo $\fil^{d-1} H$. 
Then the difference between $\phi$ and right multiplication by
$\alpha$ is still $H$-linear, and its maximal degree on elements from
$\dd$ is strictly less than $d$. The proof now follows by induction.
\end{proof}


All Lie \psalgs\ that are free $H$-module of rank one
are classified in \thref{trankone}:
they are isomorphic to current \psalgs\ over 
$\Kdp$ 
or
$\Hdp$.
%

\section{Structure Theory of Lie Pseudoalgebras}\lbb{sstrth}

\subsection{Structural correspondence between a Lie pseudo\-algebra
and its annihilation algebra}\lbb{srad}
Recall that a Lie $H$-\psalg\ $L$ is called {\em finite\/} if it is
finitely generated as an $H$-module. If $H$ is Noetherian
(e.g., $H=\ue(\dd)$ for a finite-dimensional Lie algebra $\dd$)
and $L$ is finite, then
$L$ is a Noetherian $H$-module, i.e., every increasing sequence
of $H$-submodules of $L$ stabilizes.

For any two subspaces $A$ and $B$ of $L$, let 
\begin{equation}
[A,B] = \{ \xb{a}{x}{b} \st a\in A, b\in B, x\in X \}.
\end{equation}
Define the {\em derived series\/} of $L$ by
$L^{(0)} = L$, $L^{(1)} =[L,L]$,
$L^{(n+1)} = [L^{(n)}, L^{(n)}]$. 
A Lie \psalg\ $L$ is called {\em solvable\/} if
$L^{(n)}=0$ for some $n$.
Similarly, define the {\em central series\/} of $L$ by
$L^0 = L$, $L^1 =[L,L]$,
$L^{n+1} = [L^n, L]$. 
The Lie \psalg\ $L$ is called {\em nilpotent\/} if
$L^n=0$ for some $n$.
As usual, $L$ is called {\em abelian\/} if $[L,L]=0$, i.e.,
if $[a*b]=0$ for all $a,b \in L$.

A Lie \psalg\ $L$ is called {\em simple\/} if it 
contains no nontrivial ideals and is not abelian.
Note that $[L,L]$ is an ideal of $L$, so in particular, $[L,L]=L$
if $L$ is simple. $L$ is called {\em semisimple\/} if it 
contains no nonzero abelian ideals. 

We will show that, as in the Lie algebra case,
$L$ is semisimple if and only if its radical is zero.
Provided that it exists,
we define the {\em radical\/} of $L$, $\Rad L$, to be its maximal
solvable ideal. When $H$ is Noetherian and $L$ is finite,
$\Rad L$ exists because of the
Noetherianity of $L$ and part (ii) of the next lemma.

\begin{lemma}\lbb{lsolvable}
{\rm(i)}
If $S$ is a solvable ideal in $L$ and $L/S$ is solvable, then $L$ is solvable.

{\rm(ii)}
If $S_1,S_2$ are solvable ideals in $L$, then their sum $S_1+S_2$
is a solvable ideal.

{\rm(iii)}
$L/\Rad L$ is semisimple. $L$ is semisimple iff\/ $\Rad L=0$.
\end{lemma}
\begin{proof}
(i) is standard.

(ii) follows from (i) and the fact that 
$(S_1+S_2)/S_1 \simeq S_2/(S_1 \cap S_2)$.

(iii) If $L/\Rad L$ has an abelian ideal $I$, then the preimage of $I$ 
under the natural projection $L \to L/\Rad L$ must be solvable and strictly 
bigger than $\Rad L$, which is a contradiction.
\end{proof}

It is easy to see, using
(\ref{axby}, \ref{axby2}), that for any two subspaces
$A,B\subset L$, we have: 
\begin{equation}\lbb{xhab}
[X\tt_H A, X\tt_H B] = X\tt_H [A,B]
\end{equation}
as subspaces of $\A(L)= X\tt_H L$. 
In particular, if $I$ is an ideal of $L$, then $X\tt_H I$
is an ideal of $\A(L)$. We will call an ideal of $\A(L)$ {\em regular\/}
if it is of the form $X\tt_H I$ for some ideal $I$ of $L$.

\begin{lemma}\lbb{lreg0}
Let $L$ be a Lie $H$-\psalg\ and $I\subset L$ be an ideal. Then{\rm:}

{\rm(i)} 
$X\tt_H I = 0$ only if $I$ is central.

{\rm(ii)} 
$X\tt_H I = \A(L)$ only if\/ $[L,L] \subset I$.
\end{lemma}
\begin{proof}
(i) has already been proved, when $L$ is finite, in 
\coref{tor0coef} and \reref{rcenttor}. In the general case, 
it can be deduced from \prref{preplal}. Let $a\in I$, then 
$a_x \equiv x\tt_H a = 0$ for any $x\in X$.
Hence the action of $a_x$ on $L$ 
is trivial, and by \eqref{ab2}, $[a*b]=0$ for any $b\in L$.

In order to prove (ii), notice that
$X\tt_H L/I=0$. Then build a Lie $H$-\psalg\ structure on
$\wti L = L \oplus L/I$
by letting $L$ act on the abelian ideal $L/I$ via the adjoint action.
Then by part (i),
$L/I$ is central in $\wti L$, hence $L$ acts trivially on $L/I$. This
means $[L,L] \subset I$.
\end{proof}

Using this lemma and \eqref{xhab},
it is easy to prove the next two results.

\begin{proposition}\lbb{solann}
A Lie \psalg\ $L$ is solvable {\rm(}respectively nilpotent{\rm)}
if and only if its annihilation Lie algebra $\A(L)$ is.
\end{proposition}
\begin{proposition}\lbb{ssimann}
Let $L$ be a centerless Lie $H$-\psalg\ which is equal
to its derived subalgebra $[L,L]$. Then $L$ is simple if
either of the following conditions holds{\rm:}

{\rm(i)} $\A(L)$ has no nontrivial $H$-invariant ideals.

{\rm(ii)} $L$ is finite or free, and $\A(L)$ has no non-central
$H$-invariant ideals.
\end{proposition}
\begin{proof}
(i) is immediate from \leref{lreg0}.

Assume that (ii) holds but $L$ is not simple. Then $\A(L)$ has a
nontrivial central regular ideal.
If $a_x \equiv x\tt_H a$ is central in $\A(L)$ for every $x\in X$,
then by \eqref{axby2} $[a_x b]_y = 0$ for every $b\in L$, $x,y\in X$. 
When $L$ is either finite or free, $l_y = 0$ for all $y\in X$ if and only 
if $l = 0$ (cf.\ \coref{tor0coef}).
Therefore $[a_x b] = 0$ for all $b\in L$, $x\in X$,
and by \eqref{ab2} we get $[a*b]=0$ for any $b\in L$. Hence $a=0$.
\end{proof}

As an immediate consequence we obtain:

\begin{corollary}\lbb{issimple}
Let $L$ be a current Lie $H$-\psalg\ over a finite-dimensional simple Lie 
algebra or over one of the primitive \psalgs\ of vector fields.
Then $L$ is simple.
\end{corollary}
\begin{proof}
It is easy to check that $L$ satisfies the assumptions of \prref{ssimann}
(see \thref{tannihvec}).
\end{proof}

The following proposition will play an important role in the
classification of finite simple Lie \psalgs.

\begin{proposition}\lbb{irregid}
For any Lie $H$-\psalg\ $L$, any non-central $H$-invariant ideal $J$ of $\A(L)$
contains a nonzero regular ideal.
\end{proposition}
\begin{proof}
Let $\alpha\in J$ be non-central. Assume that
$X\tt_H \al\cdot l = 0$ for all $l\in L$.
Note that by \prref{preplal}, we have: 
$h(\al\cdot l) = (h_{(1)} \al) \cdot (h_{(2)} l)$ 
for $h\in H$. This implies:
$(h\al)\cdot l = h_{(1)} \bigl( \al\cdot (h_{(-2)} l) \bigr)$,
which gives
$X\tt_H (h\al)\cdot l = 0$ for any $h\in H$, $l\in L$.
Then we can use \eqref{alby2} to show that $\alpha$ is central in $\A(L)$,
which is a contradiction. 

Therefore, there is 
some $l\in L$ such that $\alpha\cdot l = a$ 
has a nonzero Fourier coefficient, i.e., $X\tt_H a \ne 0$.
Since 
$a_y = (\alpha\cdot l)_y = \sum_i [h_i\alpha, l_{y x_i}]$, 
and $J$ is $H$-stable, we see that all Fourier 
coefficients of $a$ lie in $J$. Then, due to \eqref{axby2}, 
all elements in the ideal $(a)$ of $L$ generated by $a$ have all of 
their Fourier coefficients in $J$, i.e., 
$0\ne X\tt_H (a) \subset J$.
\end{proof}

\subsection{Annihilation algebras of finite simple Lie $\ue(\dd)$-\psalgs}
\lbb{sslclie}
We will now approach the problem of classification of all 
finite simple Lie $H$-pseudo\-algebras. 
In view of Kostant's \thref{tkostant}
and the results of \seref{sgaalg}, we will first restrict ourselves to
the case when $H$ is the universal enveloping algebra of a Lie algebra
$\dd$. Moreover, we will assume that $\dd$ is finite dimensional;
in this case $H=\ue(\dd)$ is filtered by finite-dimensional subspaces.
The classification is done in two steps: the first one (done in this 
subsection) is classifying all Lie algebras
that can arise as $\A(L)$ for some finite simple Lie $H$-\psalg\ $L$, 
the second step (done in the next subsection) involves a reconstruction of $L$ 
from its annihilation Lie algebra $\A(L)$
and the action of $H$ on it.



\begin{theorem}\lbb{tminquotients}
If $L$ is a finite simple Lie $H=\ue(\dd)$-\psalg, then its annihilation Lie
algebra $\A(L)$ is isomorphic {\rm(}as a topological Lie algebra{\rm)}
to an irreducible 
central extension of a current Lie algebra $\O_r\what\otimes\ss$
where $\ss$ is
a simple linearly compact Lie algebra of growth $\gw\ss = \dim\dd - r$.
\end{theorem}
\begin{proof}
First of all, we observe that $\L=\A(L)$ 
is a linearly compact Lie algebra with 
respect to the topology defined in \seref{ssubtop}, see \prref{plielc}(ii).
%
%
Consider the {\em extended annihilation algebra\/} $\L^e = \dd\ltimes\L$,
obtained by letting $\dd$ act on $\L=\A(L)$ according to its 
$H=\ue(\dd)$-module structure.

\begin{lemma}\lbb{lle}
$\L^e = \dd\ltimes\L$ is a linearly compact Lie algebra possessing a
fundamental subalgebra, i.e., 
an open subalgebra containing no ideals of $\L^e$.
\end{lemma}
\begin{proof}
Indeed, if $L_0$ is a finite-dimensional subspace of $L$ 
generating it over $H$, then because of \eqref{brafil},
$\L_i = \fil_i X \tt_H L_0$ is a subalgebra of $\L$ for $i\ge s$.
None of $\L_i$ contains ideals of $\dd\ltimes\L$, 
since every such ideal is stable under the action of $H$
and $H\cdot\fil_i X = X$, which implies $H\cdot\L_i = \L$.
\end{proof}

The center $Z$ of $\L$ is an $H$-stable closed ideal.
The quotient $\L^e/Z = \dd\ltimes(\L/Z)$ 
is a linearly compact Lie algebra possessing a fundamental subalgebra
$\L_s/(Z\cap\L_s)$.
\thref{tminquotients} will be deduced from 
\prref{guill1} applied for $\ov\L^e := \L^e/Z$.



By \prref{irregid}, any nonzero $H$-stable ideal of $\ov\L := \L/Z$ 
contains the image of a nonzero regular ideal of $\L$.
Since $L$ is simple, this means that the only 
nonzero $H$-stable ideal of $\ov\L$ is the whole $\ov\L$.
Then every nonzero ideal of $\ov\L^e$ 
contained in $\ov\L$ must equal $\ov\L$.
Hence $\ov\L$ is a minimal closed ideal of a linearly compact Lie algebra 
satisfying the assumptions of 
\prref{guill1}(i), and is therefore of the form stated in 
part (ii) of this proposition. 

Therefore, $\L$ is a central extension of a current Lie algebra over
a simple linearly compact Lie algebra. Moreover, $\L$ equals its derived 
subalgebra (otherwise we would have a proper nonabelian subideal of $\L$).
Hence it is an irreducible central extension. 

Consider the canonical filtration 
$\fil_n (\O_r\what\otimes\ss) := 
\tsum_i\, \fil_{n-i}\O_r \what\otimes \fil_i\ss$,
where $\fil_i\ss$ is the canonical filtration of $\ss$ defined
in \seref{scartan} (if $\dim\ss < \infty$ we put $\fil_i\ss = 0$ for $i\ge0$).
Then the growth of $\O_r\what\otimes\ss$ (with respect to this filtration)
equals $\gw\O_r + \gw\ss = r+\gw\ss$.
It is clear from \prref{pder} that any irreducible central extension
of $\O_r\what\otimes\ss$ has the same growth.
On the other hand, with respect to the filtration defined by \eqref{filam},
the growth of $\L$ is equal to $N=\dim\dd$ (see \prref{growthisn}).
We have to show that the two different filtrations give the same growth.

Recall that by \leref{lssurj}, a sufficiently high power of any nonzero 
element $a\in\dd$ maps any given open subspace of $\L$ surjectively onto $\L$.
Then the same argument as in the proof of \leref{ldatrans} shows that
$\dd\injto\Der\L$ intersects $\fil_0(\Der\L)$ trivially, where
$\fil_0(\Der\L)$ is induced by the canonical filtration on 
$\O_r\what\otimes\ss$. This implies $N \le r+\gw\ss$.

To show the inverse inequality, note that since 
$\fil_0 (\O_r\what\otimes\ss)$ is open in 
$\ov\L = \L/Z \simeq \O_r\what\otimes\ss$,
it contains some $\ov\L_m := \L_m/(Z\cap\L_m)$.
Now \eqref{brafil} implies $[\ov\L_i,\ov\L] \subset \ov\L_{i-s-1}$,
which together with \eqref{transit} leads to 
$\ov\L_{m+n(s+1)} \subset \fil_n (\O_r\what\otimes\ss)$
for all $n\ge0$. This implies $N \ge r+\gw\ss$.

This completes the proof of \thref{tminquotients}. 
\end{proof}

In fact, the above arguments can be used to
prove a stronger statement than \thref{tminquotients}.

\begin{corollary}\lbb{cminid1}
Let $L$ be a finite Lie $H$-\psalg\ and $M$ be a minimal 
nonabelian ideal of $L$. 
Then the annihilation algebra of $M$ is one of the Lie
algebras described in \thref{tminquotients}. 
\end{corollary}
\begin{proof}
The only place in the proof of \thref{tminquotients} where we used the
simplicity of $L$ was where we deduced that any 
nonzero regular ideal of $\A(L)$ 
must equal the whole $\A(L)$. This argument is modified
as follows. Let $J = X\tt_H I$ be a nonabelian
regular ideal of $\A(L)$ contained in $\A(M)$. 
Then the minimality of $M$ implies that $I=M$ and $J=\A(M)$.
The proof is concluded again by applying \prref{guill1}.
\end{proof}

\subsection{Classification of finite simple Lie $\ue(\dd)$-\psalgs}
\lbb{ssclass}
We will call a 
{\em \psalg\ of vector fields\/} any subalgebra of the Lie \psalg\ $\Wd$.
As in \seref{svect}, a
\psalg\ of vector fields is called {\em primitive\/} if it is one of
the following: $\Wd$, $\Sd$, $\Hd$ or $\Kd$ (then its annihilation algebra 
$\A(L)$ is isomorphic to one of the primitive Lie algebras
$W_N$, $S_N$, $P_N$ or $K_N$).

The following is the main theorem of this section.

\begin{theorem}\lbb{classify}
Let $H = U(\dd)$ be the universal enveloping algebra of a 
finite-dimensional Lie algebra $\dd$.
Then any finite simple Lie $H$-\psalg\ $L$ is isomorphic to 
a current \psalg\ over a finite-dimensional simple Lie algebra
or over one of the primitive \psalgs\ of vector fields.

Explicitly, $L\simeq\Cur_{H'}^H L'$, where 
$H'=U(\dd')$, $\dd'$ is a subalgebra of $\dd$, and
$L'$ is one of the following{\rm:}

{\rm(a)}
$L'$ is a finite-dimensional simple Lie algebra and $\dd'=0${\rm;}

{\rm(b)}
$L'=\Wdp$, $\dd'$ is arbitrary{\rm;}

{\rm(c)}
$L'=\Sdp$, where $\dd'$ is arbitrary and $\chi'\in(\dd')^*$ is such that
$\chi'([\dd',\dd']) =0${\rm;}

{\rm(d)}
$L'=\Hdp$, where $N'=\dim\dd'$ is even, $\chi'$ is as in {\rm(c)},
and $\om'\in\bigwedge^2(\dd')^*$ is such that 
$(\om')^{N'/2} \ne 0$ and $\di\om' + \chi'\wedge\om' = 0${\rm;}

{\rm(e)}
$L'=\Kdp$,  where $N'=\dim\dd'$ is odd and $\th'\in(\dd')^*$ is such that
$\th'\wedge(\di\th')^{(N'-1)/2} \ne 0$.
%
\end{theorem}
\begin{proof}
By \thref{tminquotients},
the annihilation algebra $\L$ of $L$ is an irreducible central extension of 
a current Lie algebra $\ov\L = \O_r\what\otimes\ov\ss$,
where $\ov\ss$ is a simple linearly compact Lie algebra of growth $N'=N-r$.
We have surjective maps
\begin{equation}\lbb{surj}
\O_r\what\otimes \ss \to \L \to \O_r \what\otimes \ov\ss,
\end{equation}
where $\ss$ is the universal central extension of $\ov\ss$. 
By \thref{tcarg2}, $\ov\ss$ is either finite dimensional (when $N'=0$)
or one of the Lie algebras $W_{N'}$, $S_{N'}$, $H_{N'}$ or $K_{N'}$.
By \prref{pder}, we have $\ss=\ov\ss$ in all cases, except $\ov\ss=H_{N'}$
in which case the center of $\ss=P_{N'}$ is 1-dimensional. 

Note that $\Der\ss = \Der\ov\ss$, and therefore, by \prref{pder}, we have
$\Der(\O_r\what\otimes\ss) = \Der(\O_r\what\otimes\ov\ss)$.
This implies 
$\Der(\O_r\what\otimes\ss) = \Der\L = \Der(\O_r\what\otimes\ov\ss)$.
Then the action of $\dd$ on $\L$ induces actions on $\O_r\what\otimes\ss$
and $\O_r\what\otimes\ov\ss$. The argument from the proof of \leref{ldatrans}
shows that these actions are transitive.

Now, let us apply the reconstruction functor $\C$ to 
the maps in \eqref{surj}. By \thref{trecvec}, 
$\C(\O_r\what\otimes \ss) \simeq \Cur_{H'}^H \C(\ss)$,
and $S:=\C(\ss)$ is one of the Lie \psalgs\ described in (a)--(e) above.
Moreover, by \leref{lrech}, we have
$\C(\ov\ss) \simeq \C(\ss) = S$,
and hence $\C(\O_r\what\otimes \ov\ss) \simeq \Cur_{H'}^H S$.
We therefore obtain $H$-linear maps
$\Cur_{H'}^H S \to \what L \to \Cur_{H'}^H S$ 
whose composition is the identity.
Hence $\what L := \C(\L)$ is isomorphic to $\Cur_{H'}^H S$, which is a 
simple Lie \psalg\ (\coref{issimple}).

The homomorphism $\Phi\colon L \to \what L$ given by \eqref{rtorhat}
is injective because $L$ is centerless (\reref{rcenttor}).
The action of $\what L$ on $L$ built in \seref{subrecalg} shows that the
image of $\Phi$ is an ideal of $\what L$.
Since $\what L$ is simple, it follows that $\Phi$ is an isomorphism.
\end{proof}

\coref{cminid1} and the above proof imply the following result.

\begin{corollary}\lbb{cminid2}
Let $L$ be a finite Lie \psalg\ and $M$ be a minimal 
nonabelian ideal of $L$. Then $M$ is a simple Lie \psalg.
\end{corollary}
\begin{lemma}\lbb{lminidex}
If $L$ is a centerless Lie \psalg, then any nonzero finite ideal of $L$
contains a nonzero minimal ideal.
\end{lemma}
\begin{proof}
By Zorn's Lemma, it is enough to show that $\bigcap I_\alpha \ne 0$ 
for any collection of finite ideals $\{I_\alpha\}_{\alpha \in A}$ 
of $L$ such that $I_\alpha \subset I_\beta$ for $\alpha < \beta$, where $A$
is a totally ordered index set. Assume that $\bigcap I_\alpha = 0$.
Then there is a chain of ideals 
$\{I_\alpha\}_{\alpha \in A'}$ $(A'\subset A)$ of the same rank whose
intersection is zero. Fix some $\al_0 \in A'$. Then for any $\be\in A'$,
$\be < \al_0$, the module $I_{\al_0} / I_\be$ is torsion, 
so by \coref{ctorcen}, $L$ acts trivially on it. 
This implies $[L,I_{\al_0}] \subset I_\be$
for each such $\be$, hence $I_{\al_0}$ is central. 
\end{proof}

\subsection{Derivations of finite simple Lie $\ue(\dd)$-\psalgs}\lbb{scders}
We will determine all derivations of a finite
simple Lie $H=\ue(\dd)$-\psalg\ $L$ (see \deref{dcder}).

First let us consider the case when $L=\Cur\g := H\tt\g$ 
is a current \psalg\ over a finite-dimensional Lie algebra $\g$. 
The Lie \psalg\ $\Wd$ acts on $L$
by just acting on the first factor in $H\tt\g$ (cf.\ \eqref{wdac*}):
\begin{equation}\lbb{wdcurg}
(f\tt a)*(g\tt b) = -(f\tt ga)\tt_H (1\tt b),
\qquad f,g \in H, \; a\in\dd, \; b\in\g.
\end{equation}
We also have an embedding $\Cur\Der\g\subset\Cder L$. 
The image of $\Cur\Der\g$ in $\Cder L$ is normalized by that
of $\Wd$, and the two form a semidirect sum
$\Wd\ltimes\Cur\Der\g$ which as an $H$-module is isomorphic to
$H\tt(\dd\oplus\Der\g)$.

\begin{proposition}\lbb{lcdercurg}
For any simple finite-dimensional Lie algebra $\g$, we have
$\Cder\Cur\g = \Wd\ltimes\Cur\g$.
\end{proposition}
\begin{proof}
By \leref{lcder}(iii), the annihilation algebra 
$\A(\Cder\Cur\g) \subset \Der \A(\Cur\g) = \Der (X\tt\g)$.
By \prref{pder}(ii), the latter is isomorphic to
$W_N \tt 1 + \O_N\tt\g$, since $X\simeq\O_N$.
Then $\C\A(\Cder\Cur\g) \subset \C(\Der (X\tt\g)) = \Wd + \Cur\g$
(see \thref{trecvec}).
Now by \leref{lsubgcv}, 
$\Cder\Cur\g \subset \C\A(\Cder\Cur\g) \subset \Wd + \Cur\g$.
\end{proof}

A similar argument as in the proof of the proposition shows that $\Cder L = L$
when $L$ is one of the primitive \psalgs\ of vector fields
$\Wd$, $\Sd$, $\Hd$ or $\Kd$.
In fact, the same holds when $L$ is a current \psalg\ over one of them.

\begin{proposition}\lbb{lcdercurs}
Let $L$ be a simple \psalg\ of vector fields
$($i.e., $L$ is a current \psalg\ over one of the primitive ones$)$. 
Then $\Cder L = L$.
\end{proposition}
\begin{proof}
Let $L$ be a current \psalg\ over $L'$, and $L'\subset \Wdp$
be one of the primitive \psalgs\ of vector fields, where $\dd'$ is a
Lie subalgebra of $\dd$. Then, by \thref{tannihvec},
the annihilation algebra $\L=\A(L)$
is a current Lie algebra over $\L'=\A(L')$: $\L=\O_r\what\tt\L'$,
and $\L'$ is isomorphic to $W_{N'}$, $S_{N'}$, $P_{N'}$ or $K_{N'}$,
where $N'=\dim\dd'=N-r$, $N=\dim\dd$.

As in the proof of \prref{lcdercurg}, we have
$\Cder L \subset \C\A(\Cder L) \subset \C(\Der\L)$.
By \prref{pder}, we have:
$\Der\L = W_r\tt1 + \O_r\what\tt \Der\L'$,
and
$\Der\L' = $ $W_{N'}$, $CS_{N'}$, $CH_{N'}$ or $K_{N'}$
is a Lie subalgebra of $W_{N'}$.
In particular, we see that $\Der\L\subset W_N$, and hence
$\C(\Der\L)$ is a subalgebra of $\Wd$.

So, we have: $\Cder L \subset \Wd$,
$L = \Cur L' \subset \Cur\Wdp \simeq H\tt\dd' \subset\Wd$.
Take any two nonzero elements $a\in\Wd$, $b\in H\tt\dd'$.
Then we claim that $[a*b]\in H\tt\dd'$ implies $a\in H\tt\dd'$.
This follows easily from the definition \eqref{wdbr*},
using that $\dd'$ is a subalgebra of $\dd$
(see the proof of \prref{psubwd} below for a similar argument).

Therefore $\Cder L \subset \Cur\Wdp$.
For $a\in\Cur\Wdp$, we can write $a = \sum h_i \tt_{H'} a_i$
for some $a_i\in\Wdp$ and $h_i\in H$ such that the classes
$h_i H'$ are linearly independent in $H/H'$.
Then if $a\in \Cder L$, it is easy to see that all $a_i$ must belong to
$\Cder L'$. Hence $\Cder L = \Cur\Cder L'$.
But $\Cder L' \subset \C(\Der\L') = L'$,
so $\Cder L = \Cur L' = L$.
\end{proof}

\subsection{Finite semisimple Lie $\ue(\dd)$-\psalgs}\lbb{sclss}
Recall that a Lie $H$-\psalg\ $L$ is called {\em semisimple\/} if it contains
no nonzero abelian ideals. Let $H=\ue(\dd)$, for a finite-dimensional Lie 
algebra $\dd$.

If $\g$ is a simple finite-dimensional Lie algebra, then
by \prref{lcdercurg}, we have
$\Cder\Cur\g = \Wd\ltimes\Cur\g$. 
It is easy to see that for any subalgebra $A$ of the Lie \psalg\ $\Wd$,
the Lie \psalg\ $A\ltimes\Cur\g$ is semisimple. 
Indeed, assume that $I\subset A\ltimes\Cur\g$ is an abelian ideal.
Then $I\cap\Cur\g$ is an abelian ideal in $\Cur\g$,
hence $I\cap\Cur\g = 0$. But this is impossible unless $I=0$
because the pseudobracket of any element from 
$(\Wd+\Cur\g) \setminus \Cur\g$
with elements from $\Cur\g$ gives nonzero elements from $\Cur\g$
(see \eqref{wdcurg}).
Note that this argument implies that any nonzero ideal of $A\ltimes\Cur\g$
contains $\Cur\g$.

Now we can classify all finite semisimple Lie $\ue(\dd)$-\psalgs.

\begin{theorem}\lbb{tsemisim}
Any finite semisimple Lie $\ue(\dd)$-\psalg\ $L$ is a direct sum of finite
simple Lie \psalgs\ $($described by \thref{classify}$)$ and of 
\psalgs\ of the form $A\ltimes\Cur\g$, where $A$ is a
subalgebra of $\Wd$ and $\g$ is a simple finite-dimensional Lie algebra.
\end{theorem}
\begin{proof}
Consider the set $\{M_i\}$ of all minimal nonzero ideals of $L$.
This set is nonempty by \leref{lminidex}, and finite because
$L$ is a Noetherian $H$-module.
The adjoint action of $L$ on $M_i$ gives a homomorphism of Lie \psalgs\ 
$L\to\Cder M_i$, cf.\ \leref{lcder}(ii). 

We claim that the direct sum of
these homomorphisms is an injective map. Indeed, let $N\subset L$ be the
set of all elements that act trivially on all $M_i$. This set is an
ideal of $L$. If it is nonzero it must contain some minimal ideal
$M_i$. But then this $M_i$ is abelian, which contradicts the
semisimplicity of $L$.

Therefore we have embeddings $\bigoplus M_i\subset L\subset\bigoplus\Cder M_i$.
By \coref{cminid2} all $M_i$ are simple Lie \psalgs. 
If $M_i$ is not a current
\psalg\ over a finite-dimensional Lie algebra, then by \prref{lcdercurs},
$\Cder M_i = M_i$. For $M_i = \Cur\g$, 
we have $\Cder\Cur\g = \Wd\ltimes\Cur\Der\g$.
Any subalgebra of $\Wd\ltimes\Cur\g$ containing $\Cur\g$ is
of the form $A\ltimes\Cur\g$, where $A$ is a
subalgebra of $\Wd$.
\end{proof}


%
%

Recall that a \psalg\ of vector fields is 
any subalgebra of the Lie \psalg\ $\Wd$.

\begin{proposition}\lbb{psubwd}
For any two nonzero elements $a,b\in\Wd$, we have $[a*b]\ne0$.
In particular, $\Wd$ does not contain nonzero abelian
elements, i.e.,  elements $a$ such that $[a*a]=0$.
\end{proposition}
\begin{proof}
%
%
Let us write
\begin{displaymath}
a = \tsum_i \, h_i \tt \d_i, 
\quad
b = \tsum_j \, k_j \tt \d_j,
\end{displaymath}
where $h_i,k_j \in H$ and $\{\d_i\}$ is a basis of $\dd$. Denote
by $m$ (respectively $n$) the maximal degree of the $h_i$ 
(respectively $k_j$).
We have (cf.\ \eqref{wdbr*}):
\begin{align*}
[a * b] &= \tsum_{i,j} \, (h_i \tt k_j) \tt_H (1 \tt [\d_i, \d_j])
\\
&- \tsum_{i,j} \, (h_i \tt k_j\d_i) \tt_H (1 \tt \d_j) 
+ \tsum_{i,j} \, (h_i\d_j \tt k_j) \tt_H (1 \tt \d_i).
\end{align*}
Assume that $[a*b] = 0$. Notice that only the third summand contains
coefficients from $H \tt H$ of degree
$(m+1,n)$, hence it must be zero modulo $\fil^m H \tt \fil^n H$.
Since the $\d_i$ are linearly independent, the same is true for each term
$\sum_j h_i\d_j \tt k_j$. If we choose $i$ such that $h_i$
is of degree exactly $m$, we get a contradiction.
\end{proof}

\begin{corollary}\lbb{pabel}
%
A finite Lie $\ue(\dd)$-\psalg\ $L$ contains no nonzero abelian elements
iff it is a direct sum of \psalgs\ of vector fields.
\end{corollary}
\begin{proof}
Assume that $L$ is not a direct sum of \psalgs\ of vector fields.
If $L$ is not semisimple, then
$\Rad L$ contains nonzero abelian elements. If $L$ is semisimple,
\thref{tsemisim} implies that
$L$ contains a subalgebra of the form $A\ltimes\Cur\g$ with $\g \neq 0$,
and therefore contains nonzero abelian elements
(for example $1\tt a$ for any $a\in\g$).

The converse statement follows from \prref{psubwd}.
\end{proof}

\begin{theorem}\lbb{tsubwd}
Any \psalg\ of vector fields is simple.
\end{theorem}
\begin{proof}
By \prref{psubwd}, a \psalg\ $L$ of vector fields does not contain nonzero
abelian elements, and hence is semisimple.
Then, by \thref{tsemisim}, $L$ is a direct sum of finite
simple Lie \psalgs\ and of \psalgs\ of the form $A\ltimes\Cur\g$.
In fact, there is only one summand, as $[a*b]\ne0$ for any 
two nonzero elements $a,b\in\Wd$. Furthermore, $L$ cannot be of the
form $A\ltimes\Cur\g$ with $\g \neq 0$,
 because $\Cur\g$ contains nonzero abelian elements.
\end{proof}

\begin{corollary}\lbb{csemis}
Any finite semisimple Lie $\ue(\dd)$-\psalg\ $L$ is a direct sum of
\psalgs\ of the form $A\ltimes\Cur\g$, where $A$ is either $0$ or
one of the simple \psalgs\ of vector fields
$($described by \thref{classify}$)$, and
$\g$ is either $0$ or a simple finite-dimensional Lie algebra.
\end{corollary}

We can also describe all ideals of a finite semisimple Lie \psalg\ $L$.
By the above corollary, it is enough to consider the case $L=A\ltimes\Cur\g$
with $A\ne0$, $\g\ne0$.

\begin{proposition}\lbb{idsemis}
Let $L=A\ltimes\Cur\g$ where $A$ is a \psalg\ of vector fields
and $\g$ is a simple finite-dimensional Lie algebra.
Then the only nonzero proper ideal of $L$ is $\Cur\g$.
\end{proposition}
\begin{proof}
We have already noticed (see the paragraph before \thref{tsemisim})
that any nonzero ideal $I$ of $L$ contains $\Cur\g$.
Then $I/\Cur\g$ is an ideal of $L/\Cur\g \simeq A$, but $A$ is simple
by \thref{tsubwd}.
\end{proof}

\subsection{Homomorphisms between finite simple Lie $\ue(\dd)$-\psalgs}
\lbb{sembaut}
In this subsection, $H=\ue(\dd)$ is again the universal enveloping algebra of 
a finite-dimensional Lie algebra $\dd$.

\begin{theorem}\lbb{tcurvec0}
For any finite-dimensional Lie algebra $\g$ and any
\psalg\ of vector fields $L$, 
there are no nontrivial homomorphisms between $L$ and $\Cur\g$.
\end{theorem}
\begin{proof}
Any homomorphism $\Cur\g\to L$ leads to abelian elements in $L$,
and therefore is zero (see \prref{psubwd}).

Let $f$ be a homomorphism from $L$ to $\Cur\g$.
Then $f$ induces a homomorphism
of Lie algebras $\A(f)\colon \A(L) \to \A(\Cur\g)$.
By \thref{tsubwd},  $L$ is simple, so
$L=\Cur_{H'}^H L'$ where $L'$ is a primitive $H'$-\psalg\ of vector fields
($H'=\ue(\dd')$ and $\dd'$ is a subalgebra of $\dd$).
By \thref{tannihvec}, the annihilation algebra $\L=\A(L)$ is isomorphic
to a current Lie algebra $\O_r \what\tt\L'$ over $\L'=\A(L')$.
Moreover, the quotient of $\L'$ by its center
is infinite dimensional and simple.
On the other hand, the annihilation algebra $\A(\Cur\g) \simeq X\tt\g$
is a current Lie algebra over $\g$, which is a projective limit
of finite-dimensional Lie algebras $(X / \fil_n X) \tt\g$.
Hence the adjoint action
of $\L' \equiv 1\tt\L'$ on $\L$ maps trivially to 
each of them via $\A(f)$.
But since $[\L',\L] = \L$, this implies that
each $\L \to (X / \fil_n X) \tt\g$ is trivial.
Therefore $\A(f)=0$, and by \coref{a=iff=}, we get $f=0$.
\end{proof}

\begin{theorem}\lbb{homocur}
Let $\g$ and $\h$ be finite-dimensional simple Lie algebras.
Then any isomorphism $f\colon \Cur\g\isoto\Cur\h$ maps 
$1 \tt \g$ onto $1 \tt \h$,
and thus is induced by some isomorphism of Lie algebras $\g\isoto\h$.
In particular, $\Aut\Cur\g \simeq \Aut\g$.
\end{theorem}

Recall that $\A(\Cur\g)\simeq X\tt\g$ is a current Lie algebra.
In the proof of the theorem we are going to use the following lemma.

\begin{lemma}\lbb{curideal}
Let $\g$ be a finite-dimensional simple Lie algebra, and $R$ be a
commutative \as\ algebra. Then all ideals of $R\tt\g$
are of the form $I\tt\g$ where $I$ is an ideal of $R$. 
\end{lemma}
\begin{proof}
As a $\g$-module, $R\tt\g$ is isomorphic to a direct sum of several
copies of $\g$. Any ideal $J$ of $R\tt\g$ is in particular a $\g$-module,
hence it is spanned over $\kk$ by elements of the form $r\tt a \in J$
where $r\in R$ and $a$ is a root vector in $\g$.
If $r\ne0$ is such that $r\tt a\in J$
for some nonzero $a\in\g$, then $r\tt\g\subset J$,
since $\{a \in \g \st r\tt a \in J\}$ is an ideal of $\g$ and $\g$ is simple.
Setting $I = \{ r \in R \st r\tt\g\subset J \}$,
we see that $I$ is an ideal of $R$ and $J = I\tt\g$.
\end{proof}
\begin{proof}[Proof of \thref{homocur}]
Define a map $\rho\colon \A(\Cur\h)\to\h$ by the formula:
\begin{equation}\lbb{xbrcur4}
\rho(x\tt_H (1\tt a)) = \langle x, 1 \rangle a,
\qquad x \in X, \; a \in \h.
\end{equation}
Then for $a = \tsum_i \, h_i\tt a_i \in \Cur\h = H \tt \h$, we have:
\begin{equation}\lbb{xbrcur3}
\rho(x\tt_H a) = \tsum_i\, \langle S(x), h_i \rangle a_i.
\end{equation}
It is easy to see that $\rho$ is a surjective Lie algebra homomorphism.

Any isomorphism $f\colon \Cur\g\isoto\Cur\h$
induces an isomorphism of Lie algebras 
$\ph=\A(f)\colon \A(\Cur\g) \isoto \A(\Cur\h)$.
By \leref{curideal}, $\ker\rho\ph = I\tt\g$
for some proper ideal $I$ of $X$.
Recall that $X$ is isomorphic as a topological algebra to
$\O_N = \kk[[t_1,\dots, t_N]]$ ($N=\dim\dd$), 
and $\O_N$ has a unique
maximal ideal $M_0=(t_1,\dots, t_N)$. Noting that $M_0$ corresponds to 
$\fil_0 X := \{ x\in X \st \langle x, 1\rangle = 0 \}$ 
via the isomorphism $X\simeq\O_N$, 
we deduce that $I\subset\fil_0 X$.
If $I\ne\fil_0 X$, then $(\fil_0 X/I) \tt\g$ is a nontrivial ideal of 
$(X/I) \tt\g \simeq \h$,
which is impossible because
$\h$ is simple. It follows that $\ker\rho\ph =\fil_0 X \tt\g$.

Now fix $a\in\g$ and write $f(a)=\tsum_i \, h_i\tt a_i$
for some $h_i\in H$ and linearly independent $a_i\in\h$.
Assume that, say, $h_1\not\in\kk=\fil^0 H$. Then we can find
$x\in\fil_0 X$ such that $\langle S(x), h_1\rangle \ne 0$. 
Then, by \eqref{xbrcur3}, the element 
$x\tt a\in \fil_0 X \tt\g$ is mapped by $\rho\ph$ to
$\tsum_i \, \langle S(x), h_i\rangle \tt a_i \ne 0$,
which is a contradiction. This shows that $f(a)\in 1\tt\h$, 
completing the proof.
\end{proof}

We turn now to the description of subalgebras of $\Wd$.
Recall that the Lie \psalg\ $\Wd$ acts on $H=\ue(\dd)$ by \eqref{wdac*}.
Hence any homomorphism of Lie \psalgs\ $L\to\Wd$ gives rise to a structure
of an $L$-module on $H$. 

Let us first consider the case when $L$
is a free $H$-module of rank one: $L = He$ with a pseudobracket
$[e * e] = \alpha \tt_H e$, $\alpha\in H\tt H$. Let $M = Hm$
be an $L$-module, with action $e *m = \beta\tt_H m$, $\beta \in H\tt H$.
We already know (\leref{xcybe1})
 that $\alpha$ must be of the form $\alpha = r + s\tt 1
- 1 \tt s$ where $r \in \dd \wedge \dd, s\in \dd$. Moreover $r$ and $s$
satisfy equations
\eqref{cybe3} and \eqref{cybe4}.
Furthermore, $\beta$ defines a representation of $L$ if and only if
it satisfies the 
following equation in $H\tt H\tt H$
(cf.\ \prref{prank1}):
\begin{equation}\lbb{repeq}
(1 \tt \beta)(\id \tt \Delta)(\beta) 
- (\sigma \tt \id)\, \bigl( (1 \tt \beta)(\id \tt \Delta)(\beta) \bigr)
= (\alpha \tt 1)(\Delta \tt \id)(\beta).
\end{equation}

\begin{proposition}\lbb{definingaction}
If $L = He$ with $[e*e] = \alpha\tt_H e$, $\alpha = r+s\tt 1 - 1\tt s$, then 
the only nonzero homomorphism $L \to \Wd$ is given by
$e \mapsto -r + 1\tt s$.
\end{proposition}
\begin{proof}
The statement of the proposition is equivalent to saying that all 
solutions $\be$ of \eqref{repeq} with
$\beta \in H\tt \dd$ are either trivial
or of the form $\beta = r -1\tt s$.
It is easy to check that the latter is indeed a solution (cf.\ \leref{lr1wd}).

Let us choose a basis $\{\d_i\}$ of $\dd$, and 
write $\beta = \sum h^i \tt \d_i$ and $r = \sum_{i,j} r^{ij}
\d_i \tt \d_j$ for some $h^i\in H$,
$r^{ij}\in\Kset$.
We will assume throughout the
proof that $\beta \neq 0$, and denote by $d$ the maximal degree of 
the $h^i$. Substituting the above expressions for
$\al$ and $\beta$ in \eqref{repeq}, we get:
\begin{displaymath}
\begin{split}
\sum_{i,j} h^j \tt &h^i \tt [\d_i, \d_j] 
+ \sum_{i,j}(h^j \tt h^i \d_j - h^i \d_j  \tt h^j)\tt \d_i
\\
&= \sum_{i,j,k} r^{ij} \d_i h^k_{(1)} \tt \d_j h^k_{(2)} \tt \d_k 
+ \sum_k (s h^k_{(1)} \tt h^k_{(2)} - h^k_{(1)} \tt s h^k_{(2)}) \tt \d_k.
\end{split}
\end{displaymath}

If $d>1$, expressing all $h^k$ in the Poincar\'e--Birkhoff--Witt 
basis relative to the basis $\{\d_i\}$, we
see that $H\tt H$-coefficients of degree $2d+1$ in the second summation
in the left-hand side cannot cancel with terms from other summations,
which only contribute lower degree terms. Therefore
\begin{displaymath}
\sum_{i,j} (h^j \tt h^i \d_j - h^i \d_j  \tt h^j)\tt \d_i = 0
\mod \fil^{2d}(H\tt H) \tt\dd.
\end{displaymath}
This implies that 
$\sum_j h^j \tt h^i \d_j = \sum_j h^i \d_j  \tt h^j \mod \fil^{2d}(H\tt H)$
for every $i$, which gives a contradiction, since
we can choose $h^i$ to have degree exactly $d$.
So $d \leq 1$, and we can write $\beta = \sum_{i,j}
\beta^{ij} \d_i \tt \d_j + 1\tt t$ with $\beta^{ij}\in \Kset$, $t\in \dd$.

Substituting this into \eqref{repeq} and comparing degree four terms
we get $\beta^{ij}\beta^{kl} = r^{ij}\beta^{kl}$ for all
$i,j,k,l$. Since $\beta \neq 0$ we conclude that $\beta^{ij} = r^{ij}$
for all $i,j$.
We are only left with showing that $t = -s$. 
Substitute $\al = r + s_1 - s_2$ and $\be = r + t_2$ into \eqref{repeq},
and then use \eqref{cybe4} to obtain:
\begin{equation}\lbb{strange}
\begin{split}
r_{12}(s_3+t_3) + [t_3, r_{13} - r_{23}] + r_{23}t_2 - r_{13}t_1 - &s_1 r_{13}
+ s_2 r_{23} \\
&= t_3(t_1 -t_2 + s_1 - s_2).
\end{split}
\end{equation}
Notice that
$r_{12}(s_3+t_3)$ is the only term lying in $\dd \tt \dd \tt \dd$ and
everything else belongs to 
$H\tt H\tt\Kset + H \tt\Kset\tt H + \Kset\tt H \tt H$. 
Hence
$r_{12}(s_3+t_3) = 0$, which is only possible if $r=0$ or $s+t = 0$. In the
latter case we are done. In the former, the left-hand 
side of \eqref{strange} becomes zero, and $t\neq 0$ since $\beta \neq 0$.
Thus $t_1+s_1 - t_2-s_2=0$ and
$t+s = 0$.
\end{proof}

\prref{definingaction} shows that nonabelian Lie \psalgs\ that are free of
rank one over $H$ embed uniquely in $\Wd$. We will show that
the other simple \psalgs\ of vector fields are spanned as Lie \psalgs\
by subalgebras of rank one, and therefore also embed uniquely in $\Wd$.
Recall that any \psalg\ of vector fields is in fact simple 
(\thref{tsubwd}).

\begin{theorem}\lbb{uniquely}
{\rm(i)}
For any subalgebra $L$ of $\Wd$, there is a unique nonzero
homomorphism $L\to\Wd$.

{\rm(ii)}
There is at most one nonzero homomorphism between any two
\psalgs\ of vector fields.
\end{theorem}
\begin{proof}
Part (ii) is an immediate consequence of (i) and \thref{tsubwd}.

By the above remarks, it remains to prove (i) in the cases when
$L$ is a current \psalg\ over either $\Wdp$ 
or $\Sdp$, where $\dd'$ is a subalgebra of $\dd$.

In the former case ($L=\Cur\Wdp:=H\tt_{H'}(H'\tt\dd')$, $H' = \ue(\dd')$), 
note that $L$ is spanned over 
$H = \ue(\dd)$ by elements
$\tilde a = 1 \tt_{H'} (1 \tt a)$ for $a\in \dd'$. 
Then $[\tilde a * \tilde a] = (a\tt 1 - 1\tt a)\tt_H \tilde a$, and by
\prref{definingaction} we know that the only nonzero homomorphism of
the Lie $H$-\psalg\ $H\tilde a$ to $\Wd$
maps $\tilde a$ to $1 \tt a$. Hence any embedding
of $L$ in $\Wd$ maps each $\tilde a$ to the corresponding element $1 \tt a$ of
$\Wd$.

Now let $L$ be a current \psalg\ over $\Sdp$. We will give the proof in the 
case when $L = S(\dd, \chi)$, the case of
currents being completely analogous. We are going to
make use of the following lemma.

\begin{lemma}\lbb{sub2}
If $\dd$ is a finite-dimensional Lie algebra of dimension $N>1$, then there
exist $2$-dimensional subalgebras $\dd_i$ $(i = 1,\dots,N-1)$ such that
$\dim\sum_{i=1}^r \dd_i = r+1$ for every $r = 1,\dots, N-1$.
\end{lemma}
\begin{proof}
If $\dd$ has a semisimple element $h$, we complement it to a basis
of $\ad h$ eigenvectors $\{h, h_1,\dots, h_{N-1}\}$. The subalgebras
$\dd_i = \kk h + \kk h_i$
then satisfy the statement of the lemma.

If $\dd$ has no semisimple elements, then from Levi's theorem we know that
$\dd$ must be solvable. In this case it has a $1$-dimensional ideal $\kk h$.
Complementing $h$ to a basis $\{h, h_1,\dots, h_{N-1}\}$, 
we conclude as before.
\end{proof}

Now consider a $2$-dimensional subalgebra of $\dd$ with basis $\{a,b\}$. 
Then the element $e_{ab}\in\Sd$ from \prref{psd}
depends on the choice of basis only up to
multiplication by a nonzero
element of $\kk$. Moreover, the $H$-span of this element is a (free) rank one 
subalgebra of $S(\dd, \chi)$, as can be easily checked (cf.\ \reref{reab}).

Let $S_i$ be the rank one subalgebras of $S(\dd, \chi)$ associated 
as above with the 2-dimensional
subalgebras $\dd_i$ of $\dd$ constructed in \leref{sub2}. Then by
comparing second tensor factors, we see that
$S_{r+1} \cap \sum_{i=1}^r S_i = 0$ for each $r = 1,\dots, N-2$.
Therefore
the sum of all $S_i$ is a free $H$-submodule $F$ of $S(\dd,\chi)$ of
rank $N-1$. Since the rank of $S(\dd,\chi)$ is also $N-1$, we see that
$S(\dd, \chi)/F$ is a torsion $H$-module.

Denote by $S$ the subalgebra of $S(\dd, \chi)$ generated by $F$. Since
$S(\dd, \chi)/S$ is a torsion $H$-module, we conclude, by \coref{ctorcen},
that $S$
is an ideal of $S(\dd, \chi)$. Hence $S(\dd, \chi) = S$ by simplicity of
$S(\dd, \chi)$.
Now by \prref{definingaction}, each subalgebra $S_i$ embeds uniquely in
$\Wd$. Hence $S = S(\dd, \chi)$ embeds uniquely in $\Wd$.

This completes the proof of \thref{uniquely}.
\end{proof}

Combining a number of previous results, we get
an explicit description of all subalgebras of $\Wd$, and of all isomorphisms
between the simple Lie pseudoalgebras listed in \thref{classify}.

\begin{corollary}\lbb{csubwd}
A complete list of all subalgebras $L$ of $\Wd=H\tt\dd$ is{\rm:}

{\rm(a)} 
$L=H\tt\dd'\simeq\Cur_{H'}^H \Wdp$, 
where $\dd'$ is any subalgebra of $\dd$ and $H'=\ue(\dd')${\rm;}

{\rm(b)}
$L=\{ \sum_i h_i\tt a_i \in  H\tt\dd' \st \sum_i h_i(a_i+\chi'(a_i)) = 0 \}
\simeq\Cur_{H'}^H \Sdp$, where $\dd'$ is any subalgebra of $\dd$ and
$\chi'\in(\dd')^*$ is such that
$\chi'([\dd',\dd']) =0${\rm;}

{\rm(c)}
$L=\{ (h\tt1)(r-1\tt s) \st h\in H\}$, where
$r \in \dd \wedge \dd$ and $s\in \dd$ satisfy 
\eqref{cybe3}, \eqref{cybe4}.
In this case, $L$ is isomorphic to a current \psalg\ over 
$\Hdp$ or $\Kdp$ {\rm(}see Sections {\rm \ref{subh}} and {\rm \ref{subk})}.
\end{corollary}

\begin{corollary}\lbb{isosimple}
All nontrivial isomorphisms among
the simple Lie $H = \ue(\dd)$-\psalgs\ listed in \thref{classify} are
the following $(H' = \ue(\dd'))${\rm:}

{\rm(i)}
$\Cur \g' \simeq \Cur \g''$ when $\g'$ and $\g''$ are isomorphic
Lie algebras.

{\rm(ii)}
$\Cur_{H'}^H H(\dd', \chi', \om') \simeq \Cur_{H'}^H H(\dd', \chi',\om'')$
when $\om'' = c\om'$ 
for some nonzero $c\in \kk$.

{\rm(iii)}
$\Cur_{H'}^H K(\dd', \th') \simeq \Cur_{H'}^H K(\dd', \th'')$
when $\th'' = c\th'$ for some nonzero $c\in \kk$.

{\rm(iv)}
$\Cur_{H'}^H W(\dd') \simeq \Cur_{H'}^H K(\dd',\th')$ when $\dim \dd' = 1$.

{\rm(v)}
$\Cur_{H'}^H H(\dd', \chi', \om) \simeq \Cur_{H'}^H S(\dd', \chi'')$ 
when $\dim \dd' = 2$ and $\chi'' = -\chi' + \tr\ad$.
\end{corollary}

\subsection{Finite simple and semisimple
Lie $(\ue(\dd)\smash\Kset[\Ga])$-\psalgs}\lbb{shkga}
Let, as before, $H=\ue(\dd)$ be the universal enveloping algebra of a 
finite-dimensional Lie algebra $\dd$. Let $\Ga$ be a (not necessarily
finite) group acting on $\dd$ by automorphisms.
The action of $\Ga$ on $\dd$ can be extended
to an action on $H$ which we denote by $g \cdot f$ for $g\in\Ga$, $f\in H$.
Recall that the smash product $\ti H = H\smash\Kset[\Ga]$
is a Hopf algebra, with the product determined by
$g \cdot f = gfg^{-1}$, and coproduct $\De(fg)=\De(f)\De(g)$
($g\in\Ga$, $f\in H$). 

A left $\ti H$-module $L$ is the same as an $H$-module together
with an action of $\Ga$ on it which is compatible with that of $H$.
An $\ti H$-module $L$ will be called {\em finite\/} if it is finite as an
$H$-module.

Let $\ti L$ be a Lie $\ti H$-\psalg\ with a pseudobracket denoted as
$[a{\,\ti*\,}b]$. By \coref{cgapsal}, $\ti L$ is also a Lie $H$-\psalg,
which we denote as $L$ with a pseudobracket $[a*b]$. $L$ is equipped
with an action of $\Ga$, and $[a*b]$ is $\Ga$-equivariant, 
see \eqref{gapsal1}.
As an $\ti H$-module, $L=\ti L$. 
The relationship between the two pseudobrackets is given by
\eqref{gapsal2}.

Then the following statements are easy to check.

\begin{lemma}\lbb{lhkga1}
{\rm(i)}
$[a{\,\ti*\,}b] = 0$ iff\/ $[ga*b]=0$ for all $g\in\Ga$.

{\rm(ii)}
$I\subset L=\ti L$ is an ideal of the Lie $\ti H$-\psalg\ $\ti L$
iff it is a $\Ga$-invariant ideal of the Lie $H$-\psalg\ $L$.

{\rm(iii)}
If $I$ is as in {\rm(ii)}, then its derived \psalg\ $[I,I]$ is the same
with respect to both pseudobrackets $[a{\,\ti*\,}b]$ and $[a*b]$.

{\rm(iv)}
$\Rad\ti L = \Rad L$.
\end{lemma}
\begin{proof}
(i) If $[a{\,\ti*\,}b] = 0$ then all its 
coefficients in front of $(g H\tt \Kset)\tt_{\ti H} L$
are zero for different $g\in\Ga$. Since 
$[a*b]\in (H\tt H)\tt_H L$, it follows that all $[ga*b]=0$.

(ii) and (iii) are clear by \eqref{gapsal2}.

(iv) follows from (i)--(iii) and the fact that $\Rad L$ is $\Ga$-invariant.
($\Rad L$ is $\Ga$-invariant because $[a*b]$ is $\Ga$-equivariant, 
see \eqref{gapsal1}.)
\end{proof}
\begin{proposition}\lbb{phkga2}
The Lie $\ti H$-\psalg\ $\ti L$ is solvable {\rm(}respectively 
semisimple{\rm)} if and only if the Lie $H$-\psalg\ $L$ is.
\end{proposition}
\begin{proof}
Follows from Lemmas \ref{lhkga1}(iv) and \ref{lsolvable}(iii).
\end{proof}
\begin{proposition}\lbb{phkga3}
The Lie $\ti H$-\psalg\ $\ti L$ is finite and simple 
if and only if the Lie $H$-\psalg\ $L$ is
a finite direct sum of 
isomorphic finite simple Lie $H$-\psalgs\
and $\Ga$ acts on them transitively.
\end{proposition}
\begin{proof}
By \leref{lhkga1}, $\ti L$ is simple iff $L$ is not abelian
and has no nontrivial $\Ga$-invariant ideals. In particular,
$L$ is semisimple. Using \thref{tsemisim} and the fact that
$\Kset[\Ga]I$ is an ideal of $L$ if $I$ is an ideal,
we see that $L$ is a direct sum of isomorphic finite simple Lie $H$-\psalgs.
\end{proof}

\subsection{Examples of infinite simple subalgebras of $\gc_n$}\lbb{sinfsubgc}
In this subsection, $H$ is an arbitrary cocommutative Hopf algebra.
Let us define a map $\om\colon H\tt H \to H\tt H$ by the formula:
\begin{equation}\lbb{omfa}
\om(f\tt a) = f a_{(-1)} \tt a_{(-2)} = (f\tt1)\,\De(S(a)).
\end{equation}
It is easy to check that $\om^2=\id$; this also follows from the identities
$\om=\F^{-1}(\id\tt S)=(\id\tt S)\F$ where $\F$ is the Fourier transform
defined by \eqref{ftrans}.

\begin{lemma}\lbb{lomcend}
The above
$\om$ is an anti-involution of $\Cend_1 = H\tt H$, i.e., it is an $H$-linear
map satisfying $\om^2=\id$ and 
\begin{equation}\lbb{omaomb}
\om(a)*\om(b) = (\si\tt_H\om)\,(b*a),
\qquad a,b\in\Cend_1,
\end{equation}
where, as before, $\si\colon H\tt H \to H\tt H$ is the permutation of 
the factors.
\end{lemma}
\begin{proof}
It only remains to check \eqref{omaomb}, which is straightforward.
\end{proof}

When $H=\ue(\dd)$, the annihilation algebra $\A(\Cend_1)$ is isomorphic
to the \as\ algebra of all differential operators on $X$, and $\om$
induces its standard anti-involution $*$ (formal adjoint).

Let $\ga\colon\End(\Kset^n)\to\End(\Kset^n)$ be an anti-involution,
i.e., $\ga^2=\id$ and $\ga(A)\ga(B) = \ga(BA)$.
Then we can define
an anti-involution $\om$ of $\Cend_n = H\tt H\tt\End(\Kset^n)$ by the formula
(cf.\ \eqref{omfa}):
\begin{equation}\lbb{omfan}
\om(f\tt a\tt A) = f a_{(-1)} \tt a_{(-2)} \tt \ga(A).
\end{equation}
Let $\gc_n(\om)$ be the set of all $a\in\Cend_n$ such that $\om(a)=-a$.
This is a subalgebra of the Lie \psalg\ $\gc_n$. Indeed, it is an
$H$-submodule because $\om$ is $H$-linear. If $\om(a)=-a$, $\om(b)=-b$,
then:
\begin{align*}
(\id\tt_H\om)\,[a*b] 
&= (\id\tt_H\om)\, \bigl( a*b - (\si\tt_H\id)\,(b*a) \bigr)
\\
&= (\si\tt_H\id)\, \om(b)*\om(a) - \om(a)*\om(b)
= - [a*b].
\end{align*}


Two important examples of Lie \psalgs\ $\gc_n(\om)$ are obtained when
$\Kset^n$ is endowed with a symmetric or skew-symmetric nondegenerate
bilinear form, and $\ga(A)$ is the adjoint of $A$ with respect to this form.
In these cases, we denote $\gc_n(\om)$ by $\oc_n$ and $\spc_n$, respectively.

\begin{proposition}\lbb{ocspc}
Let $H=\ue(\dd)$, $\dd\ne0$.
Then $\oc_n$ and $\spc_n$ are infinite subalgebras
of $\gc_n$ that act irreducibly on $H\tt\Kset^n$. We have{\rm:}
$\oc_n \cap \Cur\,\gl_n = \Cur\,\fo_n$ and
$\spc_n \cap \Cur\,\gl_n = \Cur\,\sp_n$.
\end{proposition}
\begin{proof}
The second statement is obvious by the definitions.
Since $\fo_n$ ($n\ge 3$) and $\sp_n$ ($n\ge 2$) act irreducibly on
$\Kset^n$, we only have to check that the action of $\oc_n$
on $H\tt\Kset^n$ is irreducible for $n=1,2$.
Using diagonal matrices, we see that
it suffices to check that $\oc_1$ acts irreducibly on $H$.

Recall that this action is given by (see \eqref{actcendv}):
\begin{displaymath}
\al * h = (1 \tt h)\al \tt_H 1
\qquad\text{for}\quad \al\in\gc_1 = H \tt H, \; h\in H.
\end{displaymath}
For $a\in\dd$, let 
$\al = 1\tt a - \om(1\tt a) = 2 \tt a + a \tt 1 \in \oc_1$.
Then $\al*h = (1 \tt ha) \tt_H 1 + (1 \tt h) \tt_H a$.
If $M\subset H$ is an $\oc_1$-submodule, and $h\in M$, $h\ne0$,
then the previous formula implies $1\in M$. Therefore $M=H$.
\end{proof}
\begin{remark}\lbb{rgcsimple}
It follows from \thref{tcj} below that in the case $H=\ue(\dd)$, $\dd\ne0$,
the Lie \psalgs\ $\gc_n$, $\oc_n$ and $\spc_n$ are semisimple. 
In fact, one can show that in this case they are simple.
In the case $H=\Kset[\Ga]$ with a finite group $\Ga$, 
the Lie \psalg\ $\gc_n$ has a center that
is a free $H$-module of rank $1$, the quotient by which is simple.
\end{remark}

If $I$ is a left or right ideal of the \as\ \psalg\ $\Cend_n$ and $L$ is a
subalgebra of the Lie \psalg\ $\gc_n$, then their intersection
$I\cap L$ is again a subalgebra of $\gc_n$. All ideals of $\Cend_n$
are described in the next proposition.

\begin{proposition}\lbb{plidcend}
{\rm(i)}
Any left ideal of the \as\ \psalg\ $\Cend_n$ is a sum of ideals
of the form $H\tt R\tt E$ where $R\subset H$ is a right ideal
and $E\subset\End(\Kset^n)$ is a left ideal. 

{\rm(ii)}
Any right ideal of\/ $\Cend_n$ is of the form $\om(I)$ for a unique
left ideal $I$.

{\rm(iii)}
$\Cend_n$ has no two-sided ideals, i.e., it is
a simple \as\ \psalg.
\end{proposition}
\begin{proof}
Let $I\subset\Cend_n$ be a left ideal, $\al=1\tt a\tt A \in\Cend_n$,
and $\be = \sum_i \, g_i\tt b_i\tt B_i \in I$ with linearly independent $g_i$.
Then:
\begin{displaymath}
\al*\be = \sum_i \, (1\tt g_i a_{(1)}) \tt_H (1\tt b_i a_{(2)} \tt A B_i)
\in (H\tt H)\tt_H I. 
\end{displaymath}
Taking $a=1$, we see that all 
$1\tt b_i \tt A B_i \in I$. 
In particular, $1\tt b_i \tt B_i \in I$, and hence
each element from $I$ is an $H$-linear
combination of elements of the form $\be=1\tt b\tt B$. 
For such $\be$, we have
$\al*\be = (1\tt a_{(1)}) \tt_H (1\tt b a_{(2)} \tt A B)$.
For $a\in\dd$, $A=\Id$,  we get that $1\tt b a \tt B \in I$. 
This proves the first part of the proposition.

Part (ii) is obvious, and part (iii) follows easily from (i) and (ii).
%
\end{proof}


\section{Representation Theory of Lie Pseudoalgebras}\lbb{srepth}

\subsection{Conformal version of the Lie Lemma}\lbb{subliel}
%
Let $L$ be a Lie $H$-\psalg\ and $V$ be an $L$-module.
In this subsection, $H=\ue(\dd)$ will be the universal enveloping
algebra of a finite-dimensional Lie algebra $\dd$.
In particular, $H$ is a Noetherian ring with no divisors of zero.

Let $I\subset L$ be an ideal and $\ph\in\Hom_H(I,H)$ be such that
\begin{equation}\lbb{vph}
V_\ph := \{ v\in V \st a*v = (\ph(a)\tt1)\tt_H v \;\;\;\forall a\in I \}
\end{equation}
is nonzero.
We will call the elements of $V_\ph$ {\em eigenvectors\/} for $I$
with an {\em eigenvalue\/} $\ph$. 
Note that every $v\in V_\ph$ is an eigenvector for the action of
$X \tt_H I \subset \A(L)$ on $V$.
By abuse of notation, we will also
write $a_x v = \ph(a_x)v$ for $a\in I$, $x\in X$, $v\in V_\ph$,
where $\ph(a_x)=\langle S(x), \ph(a)\rangle$, cf.\ \eqref{xbracket}.

Clearly, if $\ph=0$, then $V_0$ is an $L$-submodule of $V$.

\begin{lemma}\lbb{lvph}
If $\ph\ne0$, then $H V_\ph$ is a free $H$-module,
isomorphic to $H \tt V_\ph$ with $H$ acting on the first tensor factor.
\end{lemma}
\begin{proof}
Assume that
\begin{equation}\lbb{sumhivi}
\tsum_i \, f_i v_i = 0
\end{equation}
for some $f_i \in H$, $v_i \in V_\ph$. Let \eqref{sumhivi} be a relation
of this form with $f_i \in\fil^{n_i} H$ so that $\sum_i n_i$ is minimal.
We call $\sum_i n_i$ the degree of the relation \eqref{sumhivi}.
Assume that $v_i$'s are linearly independent, so that the degree of
\eqref{sumhivi} is positive.

We can find $a\in I$, $x\in X$ such that
$\ph(a_x) \ne0$. Applying $a_x$ to \eqref{sumhivi}
and using \eqref{bl2}, we obtain
\begin{displaymath}
\tsum_i \, {f_i}_{(2)} \, \ph(a_{{f_i}_{(-1)} x}) v_i = 0.
\end{displaymath}
Subtracting this from \eqref{sumhivi}, we get 
a relation of lower degree than \eqref{sumhivi}, because
$\De(f) \in 1\tt f + \tsum_{j=1}^n \fil^j H \tt \fil^{n-j} H$
for $f\in\fil^n H$.
\end{proof}

The following result is an analogue of Lie's Lemma.

\begin{proposition}\lbb{lliel}
If $V$ is finite as an $H$-module, then{\rm:}
\begin{equation}
L*V_\ph \subset (H\tt\Kset) \tt_H (\dd V_\ph + V_\ph).
\end{equation}
In other words, for every $\be\in \A(L)$, there exist
$\d_\be \in\dd$ and $A_\be \in\End V_\ph$ such that
\begin{equation}\lbb{eqliel}
\be v = (\d_\be + A_\be) v \qquad\text{for any}\;\; v \in V_\ph.
\end{equation}
In particular, $HV_\ph$ is an $L$-submodule of $V$.
\end{proposition}
\begin{proof}
Fix nonzero elements $w\in V_\ph$, $\be\in \A(L)$, and let $w_n = \be^n w$.
Let $W_n$ be the linear span of $w_0,\dots,w_n$;
we set $W_n=0$ for $n<0$.
For $a\in I$, $x\in X$, we have: $a_x w = \ph(a_x) w$, and by induction,
\begin{equation}\lbb{axwn}
a_x w_n \in \ph(a_x) w_n + n \ph([a_x,\be]) w_{n-1} + W_{n-2}.
\end{equation}
In particular, all $H W_n$ are $I$-modules.

Since $V$ is a Noetherian $H$-module, there exists $N\ge0$ such that
$H W_{N-1} \ne H W_N = H W_{N+1}$. In particular,
\begin{equation}\lbb{wn+1}
w_{N+1} \in (N+1)h w_N + H W_{N-1}
\end{equation}
for some $h\in H$. 

Writing \eqref{axwn} for $n=N+1$ and using \eqref{wn+1}, we get
\begin{equation}\lbb{axwn1}
a_x w_{N+1} \in \ph(a_x) (N+1)h w_N + (N+1)\ph([a_x,\be]) w_N + H W_{N-1}.
\end{equation}
On the other hand, applying $a_x$ to both sides of \eqref{wn+1} 
and using the $H$-sesqui-linearity gives
\begin{equation}\lbb{axwn2}
a_x w_{N+1} \in \ph(a_{h_{(-1)} x}) (N+1)h_{(2)} w_N  + H W_{N-1}.
\end{equation}
Comparison of the last two equations gives
\begin{equation}\lbb{fwn}
f w_N \in H W_{N-1} \quad\text{for}\quad
f = \ph(a_x) h + \ph([a_x,\be]) - \ph(a_{h_{(-1)} x}) h_{(2)}.
\end{equation}

If $f\ne0$, then the module $H W_N / H W_{N-1}$ is torsion, hence
$I$ acts on it as zero by \coref{ctorcen}.
This gives $a_x w_N \in H W_{N-1}$
for all $a\in I$, $x\in X$. Then \eqref{axwn} implies 
$\ph(a_x) w_N \in H W_{N-1}$. Since $H W_{N-1} \ne H W_N$, it follows
that $\ph=0$, which contradicts to the assumption $f\ne0$.

Therefore $f=0$. This is possible only when $h\in\fil^1 H = \dd + \Kset$.
Then for any $v\in V_\ph$, one has:
\begin{displaymath}
0 = fv = h a_x v + [a_x,\be] v - h_{(2)} a_{h_{(-1)} x} v
= [a_x, \be - h] v.
\end{displaymath}
This implies that $(\be - h)v \in V_\ph$, proving \eqref{eqliel}.
\end{proof}

\subsection{Conformal version of the Lie Theorem}\lbb{ssolv}


\begin{theorem}\lbb{tliet}
Let $H=\ue(\dd)\smash\Kset[\Ga]$ with $\dim\dd < \infty$.
Let $L$ be a solvable Lie $H$-\psalg\ and $V$ be an $L$-module
which is finite over $\ue(\dd)$.
Then there exists an eigenvector for the action of $L$ on $V$,
i.e., $v\in V\setminus\{0\}$ and $\ph\in\Hom_H(L,H)$
such that $a*v=(\ph(a)\tt1) \tt_H v$ for all $a\in L$.
\end{theorem}
\begin{proof}
Using \coref{cgapsal} and \prref{phkga2}, we can assume that $H=\ue(\dd)$.
The proof will be by
induction on the length of the derived series of $L$.

First consider the case when $L$ is abelian. By a Zorn's Lemma
argument, it is enough to find an eigenvector when $L=Ha$
is abelian generated by one element $a$. 
We may assume that $\ker a = 0$; then by \leref{lkey}
all $\ker_n a$ are finite dimensional.
Let $n$ be such that $\ker_n a \ne0$.
Then the statement follows from the usual Lie Theorem
applied to the $\A(L)$-module $\ker_n a$.

Now let $L$ be nonabelian, $I=[L,L]\ne0$. By the inductive
assumption, $I$ has a space of eigenvectors $V_\ph \ne0$.
If $\ph=0$, then $V_0$ is an $L$-submodule of $V$ on which $I$ acts as zero.
The abelian $H$-\psalg\ $L/I$ has an eigenvector in $V_0$,
which is also an eigenvector for $L$. 

Now assume that $\ph\ne0$. By \prref{lliel}, we have
for $\al,\be\in \A(L)$, $v\in V_\ph$:
\begin{align*}
\al v &= (\d_\al + A_\al) v,
\\
\be v &= (\d_\be + A_\be) v,
\\
[\al,\be] v &= \ph([\al,\be]) v.
\end{align*}
On the other hand, we can compute:
\begin{align*}
\al\be v &= \al (\d_\be + A_\be) v 
          = \d_\be(\al v) - (\d_\be \al) v + \al (A_\be v)
\\
&= \d_\be(\d_\al + A_\al) v - (\d_{\d_\be \al} + A_{\d_\be \al}) v
   + (\d_\al + A_\al) A_\be v
\\
&= \d_\be \d_\al v - \d_{\d_\be \al} v + \d_\be A_\al v + \d_\al A_\be v
   - A_{\d_\be \al} v + A_\al A_\be v.
\end{align*}
It follows that 
\begin{displaymath}
[\d_\al, \d_\be] = \d_{\d_\al \be} - \d_{\d_\be \al}.
\end{displaymath}

Assume that $\d_{a_x} \ne 0$ for some $a\in L$, $x\in X$, and write
$\d_x = \d_{a_x}$ for short. For $\al=a_x$, $\be=a_y$,
the above equation becomes:
\begin{displaymath}
[\d_x, \d_y]  = \d_{ \d_x y } - \d_{ \d_y x }
\end{displaymath}
(recall that $h a_x = a_{hx}$ for $h\in H$).
Note that $\d_y=0$ if $y\in\fil_n X$ for sufficiently large $n$.
Take the minimal such $n$, and let $x\in\fil_{n-1} X$
be such that $\d_x \ne 0$. 
By \leref{minusi}, there exists $y\in\fil_n X$ such that
$x=\d_x y$. Then $\d_y=0$ and
$\d_x = \d_{ \d_x y } - \d_{ \d_y x } = [\d_x, \d_y] =0$,
which is a contradiction. 

It follows that all $\d_{a_x} = 0$,
hence $L$ preserves $V_\ph$. By \leref{lvph}, $\dim V_\ph < 0$,
and therefore $L$ has an eigenvector by the usual Lie Theorem for $\A(L)$.
\end{proof}

\begin{corollary}\lbb{cliet}
Let $L$ be a solvable Lie $H$-\psalg\ and $V$ be a finite $L$-module
$($i.e., finite over $\ue(\dd))$.
Then $V$ has a filtration by $L$-submodules 
$0 = V_0 \subset V_1 \subset\dots \subset V_n = V$
such that for any $i$ 
the $L$-module $V_{i+1}/ V_i$ is generated over $H$ by eigenvectors
of some given eigenvalue $\ph_i\in\Hom_H(L,H)$.
\end{corollary}

\subsection{Conformal version of the Cartan--Jacobson Theorem}
\lbb{subcjthm}

\begin{theorem}\lbb{tcj}
Let $H=\ue(\dd)$ be the universal enveloping
algebra of a finite-dimensional Lie algebra $\dd$.
Let $L$ be a Lie $H$-\psalg\ acting faithfully and irreducibly on
the finite $H$-module $V$. Then one of the following two possibilities
holds{\rm:}

{\rm(i)}
$L$ is semisimple, either finite or infinite.

{\rm(ii)}
$L$ is finite, $\Rad L$ is abelian and of rank one as an $H$-module.
In this case, there is a subspace $\ov V$ of $V$ such that 
$V\simeq H\tt \ov V$ is a free $H$-module and $L$ can be identified with
$(A\ltimes\Cur\g)\ltimes (R\tt\id_{\ov V}) \subset \gc V$, where
$A$ is a subalgebra of $\Wd$, $\g$ is zero or a semisimple subalgebra of
$\sl\, \ov V$, and $R$ is a nonzero left ideal of $H$.
\end{theorem}
\begin{proof}
Assume that $L$ is not semisimple, i.e., it has a nonzero abelian ideal $I$.
Then, by \thref{tliet}, $I$ has an eigenvector in $V$.
If $\ov V= V_\ph$ is the corresponding eigenspace in $V$, then,
by \prref{lliel}, $H \ov V$ is an $L$-submodule of $V$.
The irreducibility of $V$ implies that $V=H \ov V$.
Now, by \leref{lvph}, $V\simeq H\tt \ov V$ is a free $H$-module,
since $\ph\ne0$ by the faithfulness of $V$.

\prref{lliel} and the faithfulness of $V$ 
 also show that $L$ embeds into $\Wd\ltimes\Cur\gl\, \ov V \subset \gc V$.
In particular, $L$ is finite. Then $\Rad L$ exists, and we can assume that 
$\ov V$ is an eigenspace for $\Rad L$. 
For each element $a\in\Rad L$ and $v\in \ov V$
we have $a*v = (\ph(a)\tt1)\tt_H v$, which means that
$\Rad L$ is identified with $R\tt\id_{\ov V} \subset\Cur\gl\, \ov V$
for $R=\ph(\Rad L)$.
Note that $R$ is of rank one, because $R\ne0$ and $H$ has no zero divisors.

Let $L_1 = L \cap (\Wd\ltimes\Cur\sl\, \ov V)$.
Notice that $L_1$ is a subalgebra of $L_2 := L/\Rad L$ and
$L_1 + \Rad L$ is a semidirect sum, because
$\Wd\ltimes\Cur\gl\, \ov V
= (\Wd\ltimes\Cur\sl\, \ov V)\ltimes(H\tt\id_{\ov V})$.
Since $\Rad L = R\tt\id_{\ov V}$, this also implies that
$RL$ is contained in $L_1 + \Rad L$.
Then $RL_2$ embeds in $L_1$, i.e., $R(L_2/L_1) = 0$.
Hence $L_2/L_1$ is torsion, and by \coref{ctorcen},
$L_1$ is an ideal of $L_2$. But $L_2$ is semisimple,
and by \prref{idsemis} a finite semisimple Lie \psalg\ does not have 
proper ideals of the same rank. Therefore, $L_1 = L_2$ and 
$L=L_1\ltimes\Rad L$ is a semidirect sum of \psalgs.

Finally, to show that $L_1$ is of the form $A\ltimes\Cur\g$,
notice that $L_1 \cap \Cur\sl\, \ov V$ is an ideal of $L_1$,
since $\Cur\sl\, \ov V$ is an ideal of $\Wd\ltimes\Cur\sl\, \ov V$.
This ideal is generated over $H$ by abelian elements, so by
Propositions \ref{psubwd} and \ref{idsemis} if it is nonzero it is of the form
$\Cur\g$ for some semisimple subalgebra $\g$ of $\sl\, \ov V$.
Hence $\Cur\g \subset L_1 \subset \Wd\ltimes\Cur\g$.
But any subalgebra of $\Wd\ltimes\Cur\g$ containing $\Cur\g$
is equal to $A\ltimes\Cur\g$ for some subalgebra $A$ of $\Wd$.
This completes the proof.
\end{proof}
\subsection{Conformal version of Engel's Theorem}\lbb{subengel}
%
%
As an application of the results of \seref{ssolv},
we can prove a conformal analogue of Engel's Theorem.

\begin{theorem}\lbb{tcengt}
Let $H=\ue(\dd)\smash\Kset[\Ga]$ with $\dim\dd < \infty$, and
let $L$ be a finite Lie $H$-\psalg\ 
$($i.e., finite over $\ue(\dd))$.
Assume that the action of any element $\al\in\A(L)$ on $L$ is nilpotent. 
Then $L$ is a nilpotent Lie \psalg.
\end{theorem}
\begin{proof}
First of all, note that the property that any element $\al$ of
$\A(L)$ acts nilpotently on $L$ remains valid when we replace $L$
by any quotient of $L$ by an ideal. In particular, $L/\Rad L$ will
have that property. However, $L/\Rad L$ is semisimple, and from the 
classification of finite semisimple Lie \psalgs\ we see that this is 
impossible, unless $L/\Rad L = 0$.

Therefore $L$ is solvable. The nilpotence of all $\al\in\A(L)$ imply that
all eigenvalues for $L$ are zero. Now \coref{cliet}
implies that $L$ is a nilpotent Lie \psalg.
\end{proof}

\subsection{Generalized weight
decomposition for nilpotent Lie \psalgs}\lbb{nilpopsalg}
Let $L$ be a (not necessarily finite) Lie $H$-\psalg,
and $V$ be a finite $L$-module,
where $H=\ue(\dd)$ for a finite-dimensional Lie algebra $\dd$.

Recall that for any $\ph\in\Hom_H(L,H)$, the {\em eigenspace\/}
$V_\ph$ of $V$ is defined by:
\begin{equation}\lbb{vph2}
V_\ph = \{ v\in V \st a*v = (\ph(a)\tt1)\tt_H v \;\;\;\forall a\in L \}.
\end{equation}
Let $V^\ph_{-1} = 0$ and set inductively
\begin{equation}\lbb{vph3}
V^\ph_{i+1} = H \{v\in V \st a*v - (\ph(a)\tt1)\tt_H v \in 
(H\tt H)\tt_H V^\ph_i  \;\;\;\forall a\in L \}.
\end{equation}
Then $V^\ph_0 = HV_\ph$ and $V^\ph_{i+1} / V^\ph_i = H(V/V^\ph_i)_\ph$.
The $V^\ph_i$ form an increasing sequence of
$H$-submodules of $V$ which stabilizes (because of noetherianity) to
some $H$-submodule of $V$ denoted $V^\ph$. If $V^\ph_{n-1} \neq
V^\ph_n = V^\ph$, then we set the {\em depth} of $V^\ph$ to be
$n$. We call $V^\ph$ the {\em generalized weight submodule\/} of $V$
relative to the {\em weight}~$\ph$.

When $L$ is nilpotent, it is obviously solvable, and, by \coref{cliet},
any finite $L$-module $V$ has a filtration by $L$-submodules 
so that the successive quotients are generalized weight modules. 


%
The main result of this subsection is the following theorem.

\begin{theorem}\lbb{tnilp}
Let $L$ be a nilpotent Lie $H$-\psalg\ and $V$ be a finite $L$-module.
Then $V$ decomposes as a direct sum of generalized weight modules.
\end{theorem}
\begin{proof}
In order to prove the statement, it is enough to show that all
$L$-module extensions between generalized weight
modules relative to distinct weights are trivial.

The strategy is to consider first the case when
$L = \langle T \rangle$ is the Lie \psalg\ generated by one element 
$T\in \gc V$. Then in the general case, we show that the generalized weight 
spaces $V^\ph$ relative to some element $T\in L$ are $L$-invariant.

\begin{lemma}\lbb{trivialextensionT}
Let $V$ be a finite $H$-module,
$T \in \gc V$, and $L=\langle T \rangle$ be a nilpotent Lie \psalg.
If $V$ contains a $T$-generalized weight module $V^\ph$ and
$V/V^\ph = W = W^\psi$ with $\psi \neq \ph$, then $V \simeq V^\ph \oplus W$
as $L$-modules.
\end{lemma}
\begin{proof}
%
Since $W^\psi_{i+1} / W^\psi_i = H(W / W^\psi_i)_\psi$
for any $i$, it suffices to prove the statement when $W = W^\psi_0 = H W_\psi$.

Let us first consider the case when
$W = H \bar v$ is a cyclic $H$-module. In order to prove that
the extension is trivial, we need to find a lifting $v\in V$ of $\bar v$
such that $T * v = (\psi \tt 1)\tt_H v$ 
and to show that $Hv + V^\ph$ is a direct sum
of $H$-modules (here and below, we write just $\psi$ instead of $\psi(T)$).
We will prove this by induction on the depth of $V^\ph$, the basis
of induction being trivial.

Let thus the statement be true for all $T$-generalized weight modules of 
depth $\leq n$ and consider a module $V^\ph$ of depth $n+1$. 
Fix an arbitrary lifting $v\in V$ of $\bar v$; then:
\begin{equation}\lbb{azioneT}
T * v = (\psi\tt 1)\tt_H v \mod (H\tt H)\tt_H V^\ph.
\end{equation}
Set 
\begin{displaymath}
T_1 = T, \qquad
T_{m+1} = [T_m * T] \in H^{\tt(m+1)} \tt_H L \quad\text{for}\;\; m\ge1.
\end{displaymath}
Then we claim that for $m>1$,
$T_{m+1} *v \in H^{\tt(m+2)} \tt_H V^\ph_n$ 
implies $T_m * v \in H^{\tt(m+1)} \tt_H V^\ph_n$.
We are going to show this first in the case when $\ph\neq 0$, 
the proof for $\ph=0$ only requiring minor changes.

So, let $\ph\neq 0$. Then 
$V^\ph / V^\ph_n = V^\ph_{n+1} / V^\ph_n = H(V^\ph / V^\ph_n)_\ph$ 
is a free $H$-module, because it is
generated by its $\ph$-eigenspace and we can apply \leref{lvph}.
We pick some $H$-basis $\{w^j\}$ for $V^\ph$ modulo $V^\ph_n$. 
If $\{h^i\}$ is some
$\Kset$-basis of $H$ compatible with its filtration, we write
\begin{equation}\lbb{azioneX}
T_m * v = \tsum_{i,j} (\alpha^i_j\tt h^i)\tt_H w^j 
\mod H^{\tt (m+1)}\tt_H V^\ph_n,
\end{equation}
where $\alpha^i_j \in H^{\tt m}$.

Notice that for $m>1$,
$T_m$ belongs to $H^{\tt m}\tt_H [L,L]$ where $[L,L]$ is the 
derived algebra of $L$, hence all weights are zero on it. 
This means that 
\begin{equation}\lbb{azioneX2}
T_m * V^\ph \subset H^{\tt (m+1)}\tt_H V^\ph_n 
\quad\text{for}\;\; m>1.
\end{equation}
We have:
\begin{displaymath}
T_{m+1} * v = [T_m * T] * v =
T_m *(T * v) - ((\sigma \tt \id)\tt_H \id)\, T * (T_m * v).
\end{displaymath}
We compute the right-hand side, using \eqref{abc*3'}, \eqref{abc*6'} and
\eqref{azioneT}--\eqref{azioneX2}, and obtain:
\begin{displaymath}
T_{m+1} * v = \tsum_{i,j} (\alpha^i_j \tt \psi h^i_{(1)} \tt
h^i_{(2)} - \alpha^i_j \tt \ph \tt h^i) \tt_H w^j \mod
H^{\tt (m+2)}\tt_H V^\ph_n.
\end{displaymath}
Now the assumption $T_{m+1} *v \in H^{\tt(m+2)} \tt_H V^\ph_n$ 
implies that coefficients of all $w^j$ must be zero. Let us
choose the highest degree $d$ for which there is some $h^i$ of degree
$d$ such that $\alpha^i_j \neq 0$ for some $j$. Then we get
$\alpha^i_j \tt (\psi-\ph) = 0$ for all $j$ and all $h^i$ of degree
$d$, hence $\alpha^i_j = 0$, giving a contradiction. 
This proves that all $\alpha^i_j = 0$, and therefore
$T_m * v \in H^{\tt(m+1)} \tt_H V^\ph_n$.

Now, because of nilpotence of $L$, $T_N = 0$ for $N\gg0$, and obviously
$0 * v \in H^{\tt (N+1)}\tt_H V^\ph_n$.
Thus we can pull the statement back to $m=2$ to obtain that $[T*T]$ maps any
lifting $v$ of $\bar v$ inside $H^{\tt3}\tt_H V^\ph_n$.

Now we can choose the lifting $v$ of $\bar v$ so that 
$T * v - (\psi\tt 1)\tt_H v \in H^{\tt2} \tt_H V^\ph_n$.
Indeed, performing the same computation as above, using
instead of \eqref{azioneT}
\begin{displaymath}
T * v = (\psi\tt 1)\tt_H v + \sum_{i,j} (\alpha^i_j\tt h^i)\tt_H w^j 
\mod H^{\tt2} \tt_H V^\ph_n
\end{displaymath}
for some $\alpha^i_j \in H$, 
we get $\alpha^i_j \tt(\ph-\psi) - (\ph-\psi) \tt \alpha^i_j = 0$.
This shows that $\alpha^i_j = c^i_j (\ph-\psi)$ for some choice of 
$c^i_j\in\kk$.
Now choose $v$ to be the lifting of $\bar v$ minimizing
the top degree $d$ of $h^i$ such that some $\alpha^i_j$ is nonzero. Then
if we replace $v$ by $v' = v + \sum c^i_j h^i w^j$,
all coefficients $\alpha^i_j$ in degree $d$ vanish, 
against minimality of $v$. This contradiction shows that the lifting $v$ can be chosen in such a way
that $\alpha^i_j=0$
for all $i,j$, and $T*v = (\psi\tt 1)\tt_H v$ modulo 
$H^{\tt2} \tt_H V^\ph_n$.

This shows that $Hv + V^\ph_n$ is indeed a submodule of $V$, and it satisfies
the hypotheses of our claim. Moreover, $V^\ph_n$ is of depth $n$ and we can
apply the inductive assumption to show that $Hv + V^\ph_n$ decomposes as a
direct sum of $L$-submodules. This means that we can find a 
lifting $\tilde v$ of $v+ V^\ph_n$
for which $T*\tilde v =(\psi\tt 1)\tt_H \tilde v$ holds exactly.

We have found a lifting $\tilde v$ of $\bar v$ proving that $V=H\tilde v +
V^\ph$. We are left with showing that this is a direct sum of $H$-modules.
This is clear if $\psi\neq 0$ since in this case $H\tilde v$ is free, hence
projective. If instead $\psi = 0$, assume the sum not to be free. This means
that some multiple $h\tilde v$ of $\tilde v$
lies in $V^\ph$. Since $\tilde v$ is
killed by $T$, so is $h\tilde v$, showing $h\tilde v = 0$ as no other
vector in a generalized weight module of nonzero weight $\ph$ is
killed by $T$.
This concludes the proof in case $\ph\neq 0$.



If $\ph = 0$, then we choose a $\Kset$-basis of $V^\ph$ modulo
$V^\ph_n$, and use in (\ref{azioneX}) coefficients of the form
$\alpha_j \tt 1$. The rest of the proof is the same.

Finally, consider the general case of
a non-cyclic $H$-module $W$.
We distinguish two cases. If $\psi\neq 0$, 
then $W=HW_\psi$ is free by \leref{lvph},
and it decomposes as a direct sum of cyclic modules to
which we can apply the above argument independently.
If $\psi = 0$, then we choose generators $\bar v^i$ of $W$ over $H$,
lift them to elements $v^i$ of $V$ in such a way that each of them is mapped
by $T$ to zero, and then argue that if $\sum h_i \bar v^i = 0$
then $\sum h_i v^i$ is an element of $V^\ph$ killed by $T$, hence is
zero. Therefore the extension of $H$-modules splits, and so does that
of $L$-modules, by the above computation.
\end{proof}

Now let $L$ be any nilpotent Lie $H$-\psalg, $V$ be a finite $L$-module,
and $T\in L$, $T\ne0$.

\begin{lemma}
Every $T$-generalized weight submodule of $V$ is stabilized by the
action of $L$.
\end{lemma}
\begin{proof}
We set 
\begin{equation}\lbb{filtraL}
L_{(-1)} = 0, \;\; L_{(i+1)} = \{a\in L \st [T*a]\in (H\tt H)\tt_H L_{(i)}\}
\end{equation}
and
\begin{equation}\lbb{filtraV}
V_{(-1)} = 0, \;\; V_{(i+1)} = \{v \in V \st T*v - (\ph(T)\tt 1)\tt_H v
\in (H \tt H)\tt_H V_{(i)} \}.
\end{equation}
Then the $L_{(i)}$ are $H$-submodules of $L$ whose union is $L$
(since $L$ is nilpotent),
and the $V_{(i)}$ are vector subspaces of $V$ whose $H$-span is
the $T$-generalized weight space $V^\ph$ 
(because $V_i^\ph = H V_{(i)}$ for all $i$).

It is easy to show by induction on $n = i+j$ that 
\begin{equation}\lbb{filtraLV}
L_{(i)} * V_{(j)} \subset (H\tt H)\tt_H V_{(i+j)}.
\end{equation}
Indeed, the basis of induction (say $n = -1$) is trivial, and
the inductive step follows from \eqref{filtraL}, \eqref{filtraV}
and the identity 
$[T*a]*v = T*(a*v) - \bigl((\sigma\tt\id)\tt_H \id\bigr) a*(T*v)$.
Equation \eqref{filtraLV} implies that $L*V^\ph \subset (H\tt H)\tt_H V^\ph$,
as desired.
\end{proof}

We are now able to complete the proof of \thref{tnilp}.
Let $V = \bigoplus V_i$ be finest among all
decompositions into direct sum of $L$-submodules of $V$ such that all of
the $H$-torsion of $V$
is contained in one of the $V_i$.
Note that such a finest decomposition always exists,
because any decomposition defines a partition of $\rank V$ into
non-negative integers, and finer decompositions define finer partitions.

We claim that each
$V_i$ is a generalized weight module for $L$. Otherwise, there must be
some element $T\in L$ for which some of the $V_i$ is not a
$T$-generalized weight module. But if so, then $V_i$ decomposes into a direct
sum of its $T$-generalized  weight submodules, and all torsion elements lie 
in the $T$-eigenspace of eigenvalue $0$. Since all 
$T$-generalized  weight submodules are $L$-invariant, we obtain 
a contradiction. Therefore $V$ is a direct sum of its 
generalized  weight submodules.
\end{proof}

\subsection{Representations of a Lie \psalg\ and of its annihilation algebra}
\lbb{srepannih}
Let $H=\ue(\dd)$ be the universal enveloping
algebra of a finite-dimensional Lie algebra $\dd$, 
and $L$ be a finite Lie $H$-\psalg.

Recall that the annihilation algebra $\L=\A(L)$ of $L$ possesses a filtration
by subspaces $\L = \L_{-1} \supset \L_0 \supset\dotsm$
satisfying \eqref{brafil}:
\begin{displaymath}\lbb{brafil2}
[\L_i, \L_j] \subset \L_{i+j-s}
\qquad\text{for all $i,j$ and some fixed $s$},
\end{displaymath}
that make $\L$ a linearly compact Lie algebra (\prref{plielc}).
Moreover, $\L$ is an $H$-\difalg, i.e., $\dd$ acts on it by derivations.
The semidirect sum $\L^e := \dd\ltimes\L$ is called the {\em extended
annihilation algebra}. Letting $\L^e_n = \L_n$ for all $n$
makes $\L^e$ a topological Lie algebra as well.

An $\L^e$-module (or $\L$-module) $V$ is called {\em conformal\/}
if any $v\in V$ is killed by some $\L_n$; in other words, if $V$ is 
a topological $\L^e$-module when endowed with the discrete topology.
Now \prref{preplal} can be reformulated as follows.

\begin{proposition}\lbb{preplal2}
Any module $V$ over the Lie \psalg\ $L$ has a natural structure of a
conformal $\L^e$-module, and vice versa. Moreover, $V$ is irreducible
as an $L$-module iff it is irreducible as an $\L^e$-module.
\end{proposition}

Together with the next two lemmas, this proposition is an
important tool in the study of representation theory of Lie \psalgs.

\begin{lemma}\lbb{lkey2}
Let $L$ be a finite Lie \psalg\
and $V$ be a finite $L$-module. For $n\ge -1$, let
\begin{displaymath}\lbb{kernv}
\ker_n V
= \{ v \in V \st \L_n v = 0 \},
\end{displaymath}
so that, for example, $\ker_{-1} V = \ker V$
and $V = \bigcup_n \ker_n V$.
Then all vector spaces
$\ker_n V / \ker V$ are finite dimensional.
\end{lemma}
\begin{proof}
The proof is an application of \leref{lkey}, using the following fact.
Let $A$ be a vector space and $A_i\supset B_i$ ($i=1,\dots,k$)
be subspaces of $A$ such that all $A_i/B_i$ are finite dimensional,
then $\bigcap A_i / \bigcap B_i$ is finite dimensional.
It is enough to show this for $k=2$, in which case it follows from
the isomorphism
\begin{displaymath}
\frac {(A_1 \cap A_2) / (B_1 \cap B_2)} {(A_1 \cap B_2) / (B_1 \cap B_2)}
\simeq \frac{A_1 \cap A_2}{A_1 \cap B_2} .
\end{displaymath}
\end{proof}

Note that $[\L_s, \L_n] \subset \L_n$ for any $n$, and in particular
$\L_s$ is a Lie algebra.

\begin{lemma}\lbb{lkey21}
Let $L$ be a finite Lie \psalg\
and $V$ be a finite $L$-module such that $\ker V=0$.
Then $V$ is locally finite as an $\L_s$-module, i.e., 
any vector $v\in V$ is contained in a finite-dimensional subspace
invariant under $\L_s$.
\end{lemma}
\begin{proof}
Any $v\in V$ is contained in
some $\ker_n V$, which is finite dimensional by \leref{lkey2},
and $\L_s$-invariant because $[\L_s, \L_n] \subset \L_n$.
\end{proof}

Let $V$ be a finite irreducible $L$-module. Then $\ker V=0$.
Take some $n$ such that $\ker_n V \ne 0$. This space is finite
dimensional and
$\L_s$-invariant; let $U$ be an irreducible $\L_s$-submodule of $\ker_n V$.
The $\L^e$-submodule of $V$ generated by $U$ is a factor of the induced
module $\Ind_{\L_s}^{\L^e} U$. 
Thefore, $V$ is a factor module of $\Ind_{\L_s}^{\L^e} U$.

In many cases $\dd$ acts on $\L$ by inner derivations so that we have an
injective homomorphism $\dd\injto\L$. In this case, $\L^e$ is isomorphic
to the direct sum of $\dd$ and $\L$, and we have 
$\Ind_{\L_s}^{\L^e} U \simeq H\tt\Ind_{\L_s}^{\L} U$.

The above results, combined with the results of Rudakov \cite{Ru1, Ru2}
and \cite{Ko}, will allow us to classify all finite
irreducible representations of all
finite semisimple Lie \psalgs\ (work in progress).

\section{Cohomology of Lie Pseudoalgebras}\lbb{cohom}

\subsection{The complexes $C^\bullet(L,M)$ and 
$\wti C^\bullet(L,M)$}
\lbb{sclm}
Recall that in \seref{slie*} we defined cohomology of a Lie
algebra in any pseudotensor category (\deref{dcoh}). 
Now we will spell out this definition for the case of Lie $H$-\psalgs,
i.e., for the pseudotensor category $\M^*(H)$ (see \eqref{m*d}).
As before, $H$ is a cocommutative Hopf algebra.
Let $L$ be a Lie $H$-\psalg\ and $M$ be an $L$-module.

By definition, $C^n(L,M)$, $n\ge1$, consists of all 
\begin{equation}\lbb{nga}
\ga\in\Lin(\{\underbrace{L,\dots, L}_n\},M) 
:= \Hom_{H^{\tt n}} (L^{\tt n}, H^{\tt n}\tt_H M)
\end{equation}
that are skew-symmetric (see \firef{Fskscoch}). Explicitly, $\ga$ has the
following defining properties (cf.\ \eqref{bil*}, \eqref{ssym*}):
\begin{description}
\item[$H$-polylinearity]
\begin{equation}\lbb{gapolin*}
\ga(h_1 a_1 \tt\dotsm\tt h_n a_n) = 
((h_1 \tt\dotsm\tt h_n)\tt_H 1)\,  \ga(a_1 \tt\dotsm\tt a_n)
\end{equation}
for $h_i\in H$, $a_i\in L$.

\item[Skew-symmetry]
\begin{equation}\lbb{ssym*ga}
\begin{split}
\ga(a_1 \tt\dotsm &\tt a_{i+1}\tt a_i \tt\dotsm\tt a_n)
\\
&= -(\si_{i,i+1} \tt_H\id)\, 
\ga(a_1 \tt\dotsm\tt a_i\tt a_{i+1} \tt\dotsm\tt a_n),
\end{split}
\end{equation}
where $\si_{i,i+1}\colon H^{\tt n} \to H^{\tt n}$
is the transposition of the $i$th and $(i+1)$st factors.
\end{description}

For $n=0$, we put $C^0(L,M) = \kk\tt_H M \simeq M/H_+M$,
where $H_+ = \{ h\in H \st \ep(h)=0 \}$ is the augmentation ideal.
The differential 
$d\colon C^0(L,M) = \kk\tt_H M \to C^1(L,M) = \Hom_H(L,M)$ is given by:
\begin{equation}\lbb{diga0}
\begin{split}
\bigl( d(1\tt_H m) \bigr)(a) = &\tsum_i \, (\id\tt\ep)(h_i) \, m_i \in M
\\
&\;\text{if}\quad
a*m = \tsum_i \, h_i \tt_H m_i \in H^{\tt2} \tt_H M
\end{split}
\end{equation}
for $a\in L$, $m\in M$.

For $n\ge1$, the differential $d\colon C^n(L,M) \to C^{n+1}(L,M)$ is given by 
\firef{Fdifcoch}. Explicitly:
\begin{equation}\lbb{diga}
\begin{split}
(d\ga)&(a_1\tt\dotsm\tt a_{n+1})
\\
& = \;\; \sum_{1\le i\le n+1} (-1)^{i+1}
(\si_{1 \to i} \tt_H\id)\, 
a_i * \ga(a_1 \tt\dotsm\tt {\what a}_i \tt\dotsm\tt a_{n+1})
\\
& + \sum_{1\le i<j \le n+1} (-1)^{i+j}
(\si_{1 \to i, \, 2 \to j} \tt_H\id)
\\
&\qquad\qquad\qquad\times
\ga([a_i * a_j] \tt a_1 \tt\dotsm\tt 
        {\what a}_i \tt\dotsm\tt {\what a}_j \tt\dotsm\tt a_{n+1}),
\end{split}
\end{equation}
where $\si_{1 \to i}$ is the permutation
$h_i \tt h_1 \tt\dotsm\tt h_{i-1}\tt h_{i+1} \tt\dotsm\tt h_{n+1}
\mapsto h_1 \tt\dotsm\tt h_{n+1}$,
and $\si_{1 \to i, \, 2\to j}$ is the permutation
$h_i \tt h_j \tt h_1 \tt\dotsm\tt h_{i-1}\tt h_{i+1} \tt\dotsm\tt
h_{j-1}\tt h_{j+1} \tt\dotsm\tt h_{n+1}
\mapsto h_1 \tt\dotsm\tt h_{n+1}$.

In \eqref{diga} we also use the following conventions.
If $a*b = \tsum_i\, f_i \tt_H c_i \in H^{\tt 2}\tt_H M$
for $a\in L$, $b\in M$, then for any $f\in H^{\tt n}$ we set:
\begin{displaymath}
a*(f\tt_H b) = \tsum_i\, (1\tt f) \, (\id\tt\De^{(n-1)})(f_i) \tt_H c_i
\in H^{\tt (n+1)} \tt_H M,
\end{displaymath}
where $\De^{(n-1)} = (\id\tt\dotsm\tt\id\tt\De)\dotsm(\id\tt\De)\De
\colon H \to H^{\tt n}$ is the iterated comultiplication 
($\De^{(0)} := \id$). Similarly, if 
$\ga(a_1 \tt\dotsm\tt a_n) = \tsum_i\, g_i \tt_H v_i \in H^{\tt n} \tt_H M$,
then for $g\in H^{\tt 2}$ we set:
\begin{multline*}
\ga((g\tt_H a_1) \tt a_2 \tt\dotsm\tt a_n)
\\
= \tsum_i\, (g\tt 1^{\tt(n-1)}) \, (\De\tt\id^{\tt(n-1)})(g_i) \tt_H v_i
\in H^{\tt (n+1)} \tt_H M.
\end{multline*}
These conventions reflect the compositions of
polylinear maps in $\M^*(H)$, see \eqref{compoly3}.
Note that \eqref{diga} holds also for $n=0$ if we define 
$\De^{(-1)} := \ep$.

The fact that $d^2=0$ is most easily checked using \firef{Fdifcoch}
and the same argument as in the usual Lie algebra case.
The cohomology of the resulting complex $C^\bullet(L,M)$ is called the 
{\em cohomology of $L$ with coefficients in $M$\/} and is denoted
by $\coh^\bullet(L,M)$.

One can also modify the above definition by replacing everywhere
$\tt_H$ with $\tt$. Let ${\wti{C}}^n(L,M)$ consist of all skew-symmetric
$\ga\in\Hom_{H^{\tt n}} (L^{\tt n}, H^{\tt n}\tt M)$,
cf.\ \eqref{gapolin*}, \eqref{ssym*ga}. Then we can define a differential 
$d\colon {\wti{C}}^n(L,M) \to {\wti{C}}^{n+1}(L,M)$ by \eqref{diga}
with $\tt_H$ replaced everywhere by $\tt$; then again $d^2=0$.
(In fact, one can define a pseudotensor category $\wti\M^*(H)$
by replacing $\tt_H$ with $\tt$ everywhere in the definition of $\M^*(H)$.)
The corresponding cohomology ${\mathrm{\wti H}}^\bullet(L,M)$ will be called
the {\em basic cohomology\/} of $L$ with coefficients in $M$.
In contrast, $\coh^\bullet(L,M)$ is sometimes called the 
{\em reduced\/} cohomology (cf.\ \cite{BKV}).

\subsection{Extensions and deformations}\lbb{sextdefo}
We will show that the cohomology theory of Lie \psalgs\ defined
in \seref{sclm} describes extensions and deformations, 
just as any cohomology theory.
This result is a straightforward generalization of
Theorem~3.1 from \cite{BKV}.

\begin{theorem}\lbb{exts}
{\rm(i)}
The isomorphism classes of $H$-split extensions
\begin{displaymath}
0 \to M \to E \to N \to 0
\end{displaymath}
of finite modules over a Lie $H$-\psalg\ $L$ are in one-to-one 
correspondence with elements of $\coh^1(L, \Chom(N,M))$.

{\rm(ii)}
Let $C$ be an $L$-module, considered as a Lie $H$-\psalg\
with respect to the zero pseudobracket. Then the equivalence
classes of $H$-split ``abelian'' extensions
\begin{displaymath}
0 \to C \to \what L \to L \to 0
\end{displaymath}
of the Lie $H$-\psalg\ $L$ correspond bijectively to
$\coh^2(L,C)$.

{\rm(iii)}
The equivalence classes of first-order deformations of a Lie $H$-\psalg\ $L$
$($leaving the $H$-action intact$)$ correspond bijectively to $\coh^2(L,L)$.
\end{theorem}
\begin{proof}
(i) Let
\begin{displaymath}
0 \to M \stackrel{i}\to E \stackrel{p}\to N \to 0
\end{displaymath}
be an extension of $L$-modules, which is split over $H$.
Choose a splitting $E=M\oplus N = \{ m+n \st m\in M, n\in N\}$
as $H$-modules. The fact that $i$ and $p$ are homomorphisms of
$L$-modules implies ($a\in L$, $m\in M$, $n\in N$):
\begin{align}
\lbb{aemamm}
a*_E m &= a*_M m,
\\
\lbb{aenann}
a*_E n - a*_N n &=: \ga(a)(n) \in H^{\tt 2} \tt_H M.
\end{align}
It is clear that $\ga(a)\in\Chom(N,M)$ and $\ga\colon L\to \Chom(N,M)$
is $H$-linear; in other words, 
$\ga\in C^1(L,\Chom(N,M)) = \Hom_H(L,\Chom(N,M))$.

For $a,b\in L$, $n\in N$, we have (cf.\ \eqref{replie*}):
\begin{align*}
[a*b]*_E n &= a*_E (b*_E n) - ((\si\tt\id)\tt_H\id) \, (b*_E (a*_E n)),
\\
[a*b]*_N n &= a*_N (b*_N n) - ((\si\tt\id)\tt_H\id) \, (b*_N (a*_N n)).
\end{align*}
Subtracting these two equations and using \eqref{aemamm}, \eqref{aenann}, 
we get:
\begin{align*}
\ga([a*b])(n) &= a*_M \ga(b)(n) - ((\si\tt\id)\tt_H\id) \, \ga(b)(a*_N n)
\\
& \qquad\qquad\quad\;\;\,
- ((\si\tt\id)\tt_H\id) \, b*_M \ga(a)(n) + \ga(a)(b*_N n)
\\
&= \bigl( (a*\ga)(b) \bigr) (n) 
  - ((\si\tt\id)\tt_H\id) \, \bigl( (b*\ga)(a) \bigr) (n)
\end{align*}
(recall that the action of $L$ on $\Chom(N,M)$ was defined in 
\reref{ractchom}). The last equation means that $d\ga=0$.

If we choose another splitting of $H$-modules
$E=M\oplus' N = \{ m+'n \st m\in M, n\in N\}$,
then it will differ by an element $\ph$ of $\Hom_H(N,M)$:
$m+n = (m+\ph(n)) +' n$. Then the corresponding 
\begin{displaymath}
\ga(a)(n) = a*_M \ph(n) - (\id_{H\tt H} \tt_H \ph) \, (a*_N n) + \ga'(a)(n). 
\end{displaymath}
Since $\Hom_H(N,M) \simeq \kk\tt_H \Chom(N,M) = C^0(L,\Chom(N,M))$ 
(see \reref{rhomchom}),
we get $\ga(a) = a*\ph + \ga'(a)$, i.e.,
$\ga=d\ph+\ga'$. 

Conversely, given an element of $\coh^1(L, \Chom(N,M))$, we can
choose a representative $\ga\in C^1(L,\Chom(N,M))$ and define
an action $*_E$ of $L$ on $E=M\oplus N$ by \eqref{aemamm}, \eqref{aenann},
which will depend only on the cohomology class of $\ga$.
This proves~(i).

The proof of (ii) is similar. Write 
$\what L = L\oplus C = \{ a+c \st a\in L, c\in C \}$
as $H$-modules. Denoting the pseudobracket of $\what L$ by
$[a{\,\hat*\,}b]$, we have for $a,b\in L$, $c,c_1\in C$:
\begin{align*}
[a{\,\hat*\,}c] &= a*c,
\\
[c{\,\hat*\,}c_1] &= 0,
\\
[a{\,\hat*\,}b] - [a*b] &=: \ga(a\tt b) \in H^{\tt 2} \tt_H C.
\end{align*}
It is clear that $\ga\in C^2(L,C)$, and the Jacobi identity
for $\what L$ implies $d\ga=0$.

(iii) A first-order deformation of $L$ is the structure of a
Lie $H$-\psalg\ on $\what L = L[\epsilon]/(\epsilon^2) = L\oplus L\epsilon$, 
where $H$ acts trivially on $\epsilon$, such that the map $\what L\to L$ given
by putting $\epsilon=0$ is a homomorphism of Lie \psalgs.
This means that
\begin{displaymath}
0 \to L\epsilon \to \what L \to L \to 0
\end{displaymath}
is an abelian extension of Lie \psalgs, so (iii) follows from (ii).
\end{proof}

\subsection{Relation to Gelfand--Fuchs cohomology}\lbb{srelgf}
Let again $L$ be a Lie $H$-\psalg\ and $\L=\A(L):=X\tt_H L$
be its annihilation Lie algebra. Recall that (by \prref{preplal})
any $L$-module $M$ has a natural structure of an $\L$-module, given by 
$(x\tt_H a) \cdot m = a_x m$ $(a\in L, x\in X, m\in M)$,
where $a_x m$ is the $x$-product defined by (cf.\ \eqref{xbracket}):
\begin{displaymath}
a_x m = \tsum_i\, \langle S(x), g_i \rangle\, v_i
\quad\hbox{if}\quad
a*m = \tsum_i\, (g_i\tt1)\tt_H v_i \in H^{\tt 2}\tt_H M.
\end{displaymath}
Similarly, for $\ga\in {\wti{C}}^n(L,M)$ and $x_1,\dots,x_n \in X$, we define
\begin{displaymath}
\ga_{x_1,\dots,x_n}(a_1 \tt\dotsm\tt a_n)
= \tsum_i\, \langle S(x_1), g_{i,1} \rangle \dotsm 
            \langle S(x_n), g_{i,n} \rangle\, v_i
\end{displaymath}
if
\begin{displaymath}
\ga(a_1 \tt\dotsm\tt a_n)
= \tsum_i\, (g_{i,1} \tt\dotsm\tt g_{i,n}) \tt v_i \in H^{\tt n}\tt M.
\end{displaymath}
The $H$-polylinearity \eqref{gapolin*} of $\ga$ implies that
the map $\A\ga \colon \L^{\tt n} \to M$, given by 
\begin{displaymath}
(\A\ga)\bigl( (x_1\tt_H a_1) \tt\dotsm\tt (x_n\tt_H a_n) \bigr)
:= \ga_{x_1,\dots,x_n}(a_1 \tt\dotsm\tt a_n),
\end{displaymath}
is well defined. Moreover, $\A\ga$ is skew-symmetric (i.e., it is 
map from $\bigwedge^n \L$ to $M$) because of skew-symmetry \eqref{ssym*ga}
of $\ga$. 

Therefore, we can consider $\A\ga$ as an $n$-cochain for the
Lie algebra $\L$ with coefficients in $M$.
It is not difficult to check that the map 
$\A\colon {\wti{C}}^n(L,M) \to C^n(\L,M)$
commutes with the differentials (this also follows from the results
of \seref{subay}).
The following result is proved in the same way as \prref{preplal}.

\begin{proposition}\lbb{pcohlal}
The above map $\A\colon {\wti{C}}^\bullet(L,M) \to C^\bullet(\L,M)$ is
an isomorphism from the complex ${\wti{C}}^\bullet(L,M)$ to the subcomplex
$C_{\rm{GF}}^\bullet(\L,M)$ of $C^\bullet(\L,M)$ consisting of local
cochains, i.e., cochains $\A\ga$ satisfying
\begin{equation}\lbb{xprga4}
(\A\ga)\bigl( (x_1\tt_H a_1) \tt\dotsm\tt (x_n\tt_H a_n) \bigr) = 0,
\qquad\hbox{for}\;\; x_1\in\fil_k X, \; k\gg0,
\end{equation}
for any fixed $x_2,\dots,x_n$ and  $a_1,\dots,a_n$.
\end{proposition}

Note that the locality condition \eqref{xprga4} means that
$\A\ga$ is continuous when $M$ is endowed with the discrete
topology and $\L$ with the topology defined in \seref{ssubtop}.
Therefore we have:

\begin{corollary}\lbb{cohgf}
The basic cohomology ${\mathrm{\wti H}}^\bullet(L,M)$ of a Lie \psalg\ $L$
is isomor\-phic to the Gelfand--Fuchs cohomology 
$\coh_{\rm{GF}}^\bullet(\L,M)$ of its annihilation Lie algebra~$\L$.
\end{corollary}

Recall that $H$ acts on $\L = X\tt_H L$ via its left action on $X$:
$h(x\tt_H a) = hx\tt_H a$ ($h\in H$, $x\in X$, $a\in L$).
Using the comultiplication 
$\De^{(n-1)}(h) = \tsum h_{(1)} \tt\dotsm\tt h_{(n)}$,
we also get an action of $H$ on $\L^{\tt n}$. 
It follows from \eqref{xh}, \eqref{sabsbsa}
that for $h\in H$, $\al\in\L^{\tt n}$,
$\ga \in{\wti{C}}^n(L,M)$, one has:
\begin{displaymath}
(\A\ga)(h\al) = (\A(\ga\cdot h))(\al),
\end{displaymath}
where $\ga\cdot h \in{\wti{C}}^n(L,M)$ is defined by:
\begin{align*}
(\ga\cdot h)(a_1 \tt\dotsm\tt a_n) &= \tsum_i\, g_i\,\De^{(n-1)}(h) \tt v_i 
\intertext{if}
\ga(a_1 \tt\dotsm\tt a_n) &= \tsum_i\, g_i \tt v_i \in H^{\tt n}\tt M.
\end{align*}
Considering $C^n(L,M)$ instead of ${\wti{C}}^n(L,M)$ amounts to replacing
$\tt$ with $\tt_H$, i.e., to factoring by the relations
\begin{displaymath}
(\ga\cdot h)(a_1 \tt\dotsm\tt a_n) 
- (1^{\tt n} \tt h) \, \ga(a_1 \tt\dotsm\tt a_n),
\qquad h\in H.
\end{displaymath}
In terms of $\A\ga$, this corresponds to factoring by
\begin{displaymath}
h \bigl( (\A\ga)(\al) \bigr) - (\A\ga)(h\al)
= \bigl( h \circ (\A\ga) - (\A\ga) \circ h \bigr)(\al).
\end{displaymath}
This implies the next result.

\begin{proposition}\lbb{pcohlal2}
The isomorphism 
$\A\colon {\wti{C}}^\bullet(L,M) \isoto C_{\rm{GF}}^\bullet(\L,M)$
induces an isomorphism from $C^\bullet(L,M)$ to the quotient complex
of $C_{\rm{GF}}^\bullet(\L,M)$ with respect to the subcomplex
$\{ h \circ c - c \circ h \st c\in C_{\rm{GF}}^\bullet(\L,M), 
\; h\in H \}$.
\end{proposition}

When $H=\ue(\dd)$, we can define an action of $H$ on 
$C_{\rm{GF}}^\bullet \equiv C_{\rm{GF}}^\bullet(\L,M)$ 
by $h \cdot c := h \circ c + c \circ S(h)$.
This action commutes with the differential $d$, and 
$\A$ induces an isomorphism from $C^\bullet(L,M)$ to the quotient complex
$C_{\rm{GF}}^\bullet / H \cdot C_{\rm{GF}}^\bullet$.
The Lie algebra $\dd$ acts on 
$C_{\rm{GF}}^\bullet$, and clearly
$C_{\rm{GF}}^\bullet / H \cdot C_{\rm{GF}}^\bullet 
= C_{\rm{GF}}^\bullet / \dd \cdot C_{\rm{GF}}^\bullet$.
We have an exact sequence of complexes
\begin{displaymath}
0 \to \dd \cdot C_{\rm{GF}}^\bullet \to C_{\rm{GF}}^\bullet \to 
C_{\rm{GF}}^\bullet / \dd \cdot C_{\rm{GF}}^\bullet \to 0,
\end{displaymath}
which gives a long exact sequence for cohomology
\begin{equation}\lbb{exactseq}
\begin{split}
\dotsm \to \coh^i(\dd \cdot C_{\rm{GF}}^\bullet) 
\to \coh^i(C_{\rm{GF}}^\bullet) 
&\to \coh^i(C_{\rm{GF}}^\bullet / \dd \cdot C_{\rm{GF}}^\bullet) 
\\
\to \coh^{i+1}(\dd \cdot C_{\rm{GF}}^\bullet) 
&\to \coh^{i+1}(C_{\rm{GF}}^\bullet) 
\to\dotsm \, .
\end{split}
\end{equation}
\begin{remark}[\cite{BKV}] \lbb{bkv}
If $\dim\dd = 1$, then $\dd$ acts freely on $C_{\rm{GF}}^i$ for $i>0$,
and we have 
$\coh^i(\dd \cdot C_{\rm{GF}}^\bullet) \simeq \coh^i(C_{\rm{GF}}^\bullet)$
for $i>0$. When $M$ is a free $H$-module, this is also true for $i=0$.
\end{remark}

\begin{proposition}\lbb{pgfcoh}
Assume that $\dd$ acts on $\L$ by inner derivations and that the
action of $\dd$
on $M$ coincides with that of its image in $\L$.
Then for any $i\ge0$, we have isomorphisms
\begin{align}
\lbb{gf1}
\coh^i(L,M) &\simeq  \coh^i_{\rm{GF}}(\L,M) \oplus
\coh^{i+1}(\dd\cdot C_{\rm{GF}}^\bullet).
\intertext{If, in addition, $\dim\dd =1$, then we have for $i>0$}
\lbb{formula}
\coh^i(L,M) &\simeq  \coh^i_{\rm{GF}}(\L,M) \oplus
\coh^{i+1}_{\rm{GF}}(\L,M).
\end{align}
This also holds for $i=0$, provided that $M$ is a free $H$-module.
\end{proposition}
\begin{proof}
Since the adjoint
action of $\L$ on $\coh^\bullet_{\rm{GF}}(\L,M)$ is trivial,
we obtain that \mbox{$\coh^i(\dd \cdot C_{\rm{GF}}^\bullet)$} maps to zero
in the exact sequence \eqref{exactseq}.
Therefore we have exact sequences
\begin{displaymath}
0 \to \coh^i(C_{\rm{GF}}^\bullet) 
\to \coh^i(C_{\rm{GF}}^\bullet / \dd \cdot C_{\rm{GF}}^\bullet) 
\to \coh^{i+1}(\dd \cdot C_{\rm{GF}}^\bullet) \to 0,
\end{displaymath}
which lead to isomorphisms \eqref{gf1}. Formula \eqref{formula} follows from
\reref{bkv}.
\end{proof}

Note that in general 
we have:
\begin{displaymath}
\dim\coh^i(L,M) \le  \dim\coh^i_{\rm{GF}}(\L,M) + 
\dim\coh^{i+1}(\dd \cdot C_{\rm{GF}}^\bullet).
\end{displaymath}

The above results provide a tool for computing the cohomology
of Lie \psalgs, by making use of the known results on
Gelfand--Fuchs cohomology of Lie algebras of vector fields 
\cite{F}.

\subsection{Central extensions of finite simple Lie \psalgs}\lbb{centrexts}
In this section we determine by a direct computation
all nontrivial central extensions of a finite simple Lie
\psalg\ $L$ with trivial coefficients (see \thref{centext} below).

Such a central extension of $L$ is isomorphic as an $H$-module
to $\what L = L \oplus \kk 1$, where the action of $H$ on $1$ is given 
by $h \cdot 1 = \ep(h)1$. The pseudobracket is then
\begin{equation}\lbb{centr1}
[a{\,{\hat *}\,}b] = [a * b] +
\ga(a,b) \tt_H 1, \qquad a,b \in L,
\end{equation}
where $\ga(a,b) \in H \tt H$.

Notice that a tensor product 
$(h^1 \tt \dotsm\tt h^n)\tt_H 1 \in H^{\tt n} \tt_H \kk$ 
can always be re-expressed as 
$(h^1 h^n_{(-1)} \tt \dotsm\tt h^{n-1} h^n_{(-(n-1))} \tt 1) \tt_H 1$, 
and this coefficient is unique in $H^{\tt(n-1)}\tt 1$
(see \leref{lhhh}).

Therefore, the above bracket is uniquely determined by the unique $\beta(a,b)
\in H$ such that 
\begin{equation}\lbb{centr2}
\ga(a,b) \tt_H 1 = (\beta(a,b) \tt 1) \tt_H 1;
\end{equation}
we will call this map $\beta\colon L\tt L \to H$
the {\em cocycle} representing the central extension. Then $H$-bilinearity 
and skew-symmetry of the pseudobracket give the following properties
of this cocycle:
\begin{equation}\lbb{bilinskew}
\beta(ha,b) = h\beta(a,b), \quad \beta(a,hb) = \beta(a,b) S(h),
\quad \beta(a,b) = -S\bigl( \beta(b,a) \bigr),
\end{equation}
for all $a,b \in L$, $h\in H$.

We consider two central extensions equivalent if they are isomorphic 
Lie \psalgs.
An isomorphism of $\what L$ is given by an embedding $L \injto \what L$ 
projecting to the identity on $L$. All such embeddings are uniquely determined
by $H$-linear maps $\phi\colon L \to \kk$. 
Then the cocycles representing the two
equivalent central extensions differ by $\tau_\phi(a,b)$ such that
\begin{equation}\lbb{centr3}
(\tau_\phi(a,b) \tt 1)\tt_H 1 = (\id_{H\tt H} \tt_H \phi) ([a*b]).
\end{equation}
This is called a {\em trivial cocycle}.

If $L=He$ is an $H=\ue(\dd)$-module which is free
on the generator $e$, such that 
\begin{displaymath}
[e*e] = \alpha\tt_H e, \qquad
\alpha = r + s\tt 1 - 1\tt s, \quad
r \in \dd\wedge\dd, \; s\in\dd,
\end{displaymath}
then a cocycle
$\beta(a,b)$ is completely determined by its value 
$\beta = \beta(e, e) \in H$.
Trivial cocycles are of the form
\begin{displaymath}
\tau = \tau(e,e) = \phi(e)(2s - x), \qquad \phi(e) \in \kk,
\end{displaymath}
where
\begin{equation}\lbb{xr}
x = \frac{1}{2} \sum_{i,j} r^{ij} [\d_i, \d_j]
\qquad\text{if}\quad
r = \sum_{i,j} r^{ij} \d_i \tt \d_j.
\end{equation}

\begin{lemma}\lbb{centrextension}
Let $L = He$ be a Lie \psalg\ as above. Then
$\coh^2(L, \kk) \simeq B/\kk(2s-x)$, where $B$ is the space of
elements $\be \in H$ satisfying the following two conditions{\rm:}
\begin{align}
\lbb{cocycle1}
\be &= -S(\be),
\\
\lbb{cocycle}
\al\Delta(\be) &= (\be\tt 1 + 1\tt\be)\al + \be\tt(3s-x) - (3s-x)\tt\be.
\end{align}
Moreover, when $r\neq 0$, then $\be\in\dd$, and \eqref{cocycle} becomes 
equivalent to the following system of equations{\rm:}
\begin{align}
\lbb{cocycle2}
[s,\be] &= 0,
\\
\lbb{cocycle3}
[r,\Delta(\be)] &= \be\tt(3s-x) - (3s-x)\tt\be.
\end{align}
\end{lemma}
\begin{proof}
Let $\what L = He + \kk 1$ be a central extension of $L$ with
a pseudobracket 
\begin{displaymath}
[e{\,\hat*\,}e] = \al\tt_H e + (\be \tt 1)\tt_H 1, 
\end{displaymath}
where $h \cdot 1 = \ep(h) 1$ for $h\in H$.

The skew-symmetry of $[e{\,\hat*\,}e]$ is equivalent to \eqref{cocycle1}.
The Jacobi identity is equivalent to Jacobi identity for $[e*e]$
together with the following cocycle condition for $\ga=\be \tt 1$
(cf.\ \prref{prank1}):
\begin{equation}\lbb{cocycle4}
\begin{split}
(\al\tt1)\,(\De\tt\id)(\ga) \tt_H 1
&= (1\tt\al)\,(\id\tt\De)(\ga) \tt_H 1
\\
&- (\si\tt\id)\,\bigl((1\tt\al)\,(\id\tt\De)(\ga)\bigr) \tt_H 1.
\end{split}
\end{equation}
With the usual notation $r_{12}=r\tt1$, $s_1=s\tt1\tt1$, etc., we have:
\begin{align*}
(\al\tt1)\,(\De\tt\id)(\ga) \tt_H 1 
&= (\al\De(\be) \tt 1) \tt_H 1,
\\
(1\tt\al)\,(\id\tt\De)(\ga) \tt_H 1
&= (r_{23} + s_2 - s_3) \be_1 \tt_H 1
= \be_1 (r_{23} + s_2 - s_3) \tt_H 1
\\
= \be_1 (-& r_{21} - x_2 + s_1 + 2 s_2) \tt_H 1
= \be_1 (\al_{12} + 3 s_2 - x_2) \tt_H 1.
\end{align*}
{}From here it is easy to see that \eqref{cocycle4} is 
equivalent to \eqref{cocycle}.

Let now $r$ be nonzero. Rewrite \eqref{cocycle} in the form
\begin{displaymath}
\al \bigl( \Delta(\be) - \be\tt 1 - 1\tt\be \bigr) 
= [\be\tt 1 + 1\tt\be, \al] + \be\tt(3s-x) - (3s-x)\tt\be.
\end{displaymath}
If $\be \notin \dd + \kk$,
then the degree of the left-hand side equals $\deg \be + 2$ 
while that of the right-hand side is at most $\deg \be+1$,
giving a contradiction. So $\be \in \dd + \kk$, and \eqref{cocycle1}
shows that $\be \in \dd$.
\end{proof}
\begin{proposition}\lbb{virasoro}
Let $\dd' \subset \dd$ be finite-dimensional
Lie algebras, $H = \ue(\dd)$, $H' = \ue(\dd')$, 
and let $L = \Cur_{H'}^H W(\dd')$.

{\rm(i)}
If\/ $\dim \dd' = 1$, then $\coh^2(L, \kk)$ is $1$-dimensional.

{\rm(ii)}
If\/ $\dd$ is abelian and $\dim \dd' > 1$, then $\coh^2(L, \kk) = 0$.
\end{proposition}
\begin{proof}
(i) The Lie \psalg\ $L=He$ is free of rank one, 
with $e=1\tt s$, $s\in \dd' \setminus\{0\}$, hence we can
use \leref{centrextension}. In this case $\al = s\tt 1 - 1\tt s$, 
and equation \eqref{cocycle} becomes 
\begin{displaymath}
(s\tt 1 - 1 \tt s) (\Delta (\be) - \be \tt 1 - 1\tt\be)
= 3(\be \tt s - s\tt\be) 
\end{displaymath}
for $\be \in H$. We choose
a basis $\{ \d_i\}$ of $\dd$
such that $\d_1 = s$, and express $\be$ in a Poincar\'e--Birkhoff--Witt
basis as $\be = \sum_I \be_I \d^{(I)}$, $\be_I \in\kk$ (see \exref{eued}).
Then the above equation becomes:
\begin{displaymath}
\sum_{I,J\neq0}
\be_{I+J} (\d_1\d^{(I)}\tt\d^{(J)} - \d^{(I)}\tt \d_1\d^{(J)}) =
\sum_I 3\be_I ( \d^{(I)} \tt \d_1 - \d_1 \tt\d^{(I)}).
\end{displaymath}
Comparing terms of the form $h\tt\d_j$ ($j\ne1$) we find that $\be_I$ is zero
unless $I = (i, 0, \dots, 0)$ for some $i$. 
Hence $\be = \sum_i \be_i s^i$, $\be_i\in\kk$.
Substituting and comparing coefficients, we obtain that
$\be = \be_1 s + \be_3 s^3$. This obviously satisfies \eqref{cocycle1}.
The trivial cocycles are multiples of
$2s$, hence $s^3$ is the unique central extension up to scalar multiples. This
is the well-known Virasoro central extension.

(ii) Choose a basis of $\dd'$ and
let $\be$ be a cocycle representing a central extension of $L\simeq H\tt\dd'$. 
Then for each basis element $a$,
$\be$ restricts to a cocycle of $H\tt a \subset L$,
which is a current Lie \psalg\ over $W(\kk a)$. 
By part (i) we can then add to $\be$ a trivial cocycle as to make
$\be (1\tt a,1\tt a) = c_a a^3$, $c_a \in \kk$,
for every such basis element $a\in\dd'$. Denoting
$\be = \be(1\tt a, 1\tt b)$,
the Jacobi identity for elements $1\tt a, 1\tt a,1\tt b$ then gives:
\begin{equation}\lbb{caa3}
(a\tt 1 - 1\tt a) \Delta(\be) 
= c_a (a^3 \tt b - b \tt a^3)  + (\be\tt 1- 1\tt \be)\Delta(a).
\end{equation}
Let $a,b$ be distinct elements in the above basis, which we extend
to a basis $\{\d_i\}$ of $\dd$ with $\d_1 = a$, $\d_2 = b$.
We substitute the Poincar\'e--Birkhoff--Witt 
basis expression $\be = \sum_I \be_I \d^{(I)}$ in \eqref{caa3}, to get:
\begin{align*}
c_a( \d_1^3\tt\d_2 - \d_2\tt\d_1^3) &= \sum_{I,J} \be_{I+J}
( \d_1 \d^{(I)} \tt \d^{(J)} - \d^{(I)} \tt \d_1 \d^{(J)}) 
\\
&- \sum_J \be_J (\d_1\d^{(J)}\tt 1 + \d^{(J)}\tt\d_1 - \d_1\tt\d^{(J)} -
1\tt\d_1\d^{(J)}).
\end{align*}
Comparing coefficients of the form $h\tt \d_j$ for $j \neq 1$,
we find that $\be_I$ can be nonzero only
when $I = (2,1,0, \dots, 0)$, in which case $\be_I = 2c_a$, and when $I =
(i,0,\dots, 0)$ for some $i$.
This means that $\be = \be(1\tt a, 1\tt b) = f(a) + c_a a^2 b$
for some polynomial $f$.

We can repeat the same argument after switching the roles of $a$ and $b$, 
to get:
$\be(1\tt b, 1\tt a) = g(b) + c_b b^2 a$. Then the skew-symmetry 
$\be(1\tt a, 1\tt b) = -S \bigl( \be(1\tt b, 1\tt a) \bigr)$ implies:
$f(a) + c_a a^2 b = -g(-b) + c_b a b^2$.
This is possible only when $f=0$, $c_a=0$.
Therefore $\be$ is identically zero.
\end{proof}

\begin{proposition}\lbb{centrH}
Let $\dd' \subset \dd$ be abelian finite-dimensional Lie
algebras, $H = \ue(\dd), H' = \ue(\dd')$, and let
$L = \Cur_{H'}^H H(\dd', \chi,\om)$.
Then $\coh^2(L, \kk)$ is
isomorphic to $\dd$ if $\chi = 0$, and is trivial otherwise.
\end{proposition}
\begin{proof}
$L$ is free of rank one and $r \neq 0$, hence
\eqref{cocycle} becomes
$3(\be \tt s - s\tt\be) = 0$, $\be\in \dd$. This is
satisfied only by multiples of $s$ if $s\neq 0$ and by all elements of $\dd$
otherwise. Since $\dd$ is abelian, then $x=0$ and trivial cocycles are
multiples of $s$.
\end{proof}

\begin{proposition}\lbb{centrK}
Let $\dd'$ be the Heisenberg Lie algebra of dimension $N = 2n+1\geq 3$, 
and $\dd =\dd' \oplus \dd_0$ be the direct sum of $\dd'$ and an abelian
Lie algebra $\dd_0$.
Let $H = \ue(\dd),
H' = \ue(\dd')$, and $L = \Cur_{H'}^H K(\dd', \th)$. Then $\coh^2(L,\kk) = 0$.
\end{proposition}
\begin{proof}
$L$ is free of rank one, and 
\begin{displaymath}
\al = \sum_{i=1}^n (a_i\tt b_i - b_i\tt a_i) - c\tt 1 + 1\tt c, 
\end{displaymath}
where $\{a_i,b_i,c\}$ is a basis of $\dd'$ with 
the only nonzero commutation relations
$[a_i, b_i] = c$, $1\le i\le n$ (see \exref{eheisenberg}).

It is immediate to check that $[r,d\tt 1 + 1\tt d] = 
c\tt d - d\tt c$ for all $d\in \dd'$. Moreover, the element $x$ from
\eqref{xr} equals $nc$. Then, if $\be = \be'+\be_0$ with
$\be' \in \dd'$, $\be_0 \in \dd_0$, equation \eqref{cocycle3} becomes:
\begin{displaymath}
\be' \tt c - c\tt \be' = (n+3) (\be \tt c - c\tt \be).
\end{displaymath}
All solutions $\be$ of this equation are multiples of $c$.
Trivial cocycles are multiples of $2s-x = -(n+2)c$, hence all cocycles
are trivial.
\end{proof}
\begin{proposition}\lbb{centrS}
Let $\dd' \subset \dd$ be abelian finite-dimensional Lie algebras
such that $\dim\dd'>2$, let
$H = \ue(\dd), H' = \ue(\dd')$, and let $L = \Cur_{H'}^H S(\dd', 0)$. Then 
$\coh^2(L, \kk) = 0$.
\end{proposition}
\begin{proof}
By \prref{psd}, $L$ is spanned over $H$ by elements 
\begin{displaymath}
e_{ab} = a\tt b - b\tt a, \qquad a,b\in \dd',
\end{displaymath}
satisfying the relations $e_{ab} = -e_{ba}$ and
\begin{equation}\lbb{eabrel2}
a e_{bc} + b e_{ca} + c e_{ab} = 0.
\end{equation}
The pseudobrackets are (see \eqref{eab*ecd}):
\begin{align*}
[e_{ab} * e_{cd}] &= (a \tt d) \tt_H e_{bc} + (b \tt c) \tt_H e_{ad}
             - (a \tt c) \tt_H e_{bd} - (b \tt d) \tt_H e_{ac},
\intertext{and in particular}
[e_{ab} * e_{ac}] &= - (a \tt c) \tt_H e_{ab}
+ (b \tt a) \tt_H e_{ac} -(a \tt a) \tt_H e_{bc},
\\
[e_{ab} * e_{ab}] &= (b \tt a - a \tt b) \tt_H e_{ab}.
\end{align*}

Trivial cocycles $\tau_\phi$ are determined by the identity
(see \eqref{centr3}):
\begin{multline*}
(\tau_\phi(e_{ab} , e_{cd}) \tt 1) \tt_H 1
\\
             = (a \tt d) \tt_H \phi_{bc} + (b \tt c) \tt_H \phi_{ad}
             - (a \tt c) \tt_H \phi_{bd} - (b \tt d) \tt_H \phi_{ac},
\end{multline*}
where $\phi_{ab} = \phi(e_{ab}) = -\phi_{ba} \in\kk$, which gives:
\begin{equation}\lbb{trivphi}
\tau_\phi(e_{ab},e_{cd}) = -ad \phi_{bc} - bc \phi_{ad}
            + ac \phi_{bd} + bd \phi_{ac}.
\end{equation}

Let $\be$ be a cocycle for $L$ representing a central extension.
Write $\be_{ab,cd} = \be(e_{ab}, e_{cd})$ for short.
Equations \eqref{bilinskew}, \eqref{eabrel2} give the identities:
\begin{gather}
\lbb{berel3}
\be_{ab,cd} = -\be_{ba,cd} = -\be_{ab,dc} = -S(\be_{cd,ab}),
\\
\lbb{berel2}
a \be_{bc,cd} + b \be_{ca,cd} + c \be_{ab,cd} = 0.
\end{gather}
Using this, Jacobi identity for the
elements $e_{ab}, e_{ab}, e_{ac}$ gives the following
equation for $\be$:
\begin{equation}\lbb{Scocycle}
\begin{split}
(b \tt a - a \tt b) \bigl( \Delta(\be_{ab,ac}) - \be_{ab,ac}\tt 1
&- 1 \tt \be_{ab,ac} \bigr)
\\
= ab \tt \be_{ab,ac} - \be_{ab,ac} \tt ab
&+ \be_{ab,bc} \tt a^2 - a^2 \tt \be_{ab,bc} 
\\
&+ \be_{ab,ab} \tt ac - ac \tt \be_{ab,ab} .
\end{split}
\end{equation}

This is a homogeneous equation and can be solved degree by degree. If $\be$
is homogeneous of degree one, then the left-hand side is zero,
and we immediately see that $\be_{ab,ac} = \be_{ab,bc} = \be_{ab,ab} = 0$.
Then, by \eqref{berel2}, $\be_{ab,cd} = 0$.

If $\be$ is homogeneous of degree other than one, 
then $\be_{ab,ab} = 0$,
since $\be$ restricts to a cocycle of the free rank one Lie \psalg\
$H e_{ab}$, which has been shown in \prref{centrH} to take values in $\dd$.
Then equations \eqref{berel3}, \eqref{berel2} give 
$a\be_{ab,bc} = b\be_{ab,ac}$.
Hence if $a$ and $b$ are linearly independent, we can find some 
$p = p^a_{bc} \in H$ such that
\begin{equation}\lbb{compa}
\be_{ab,ac} = ap^a_{bc}, \qquad \be_{ab,bc} = bp^a_{bc}.
\end{equation}
We substitute this into \eqref{Scocycle} and 
after simplification obtain
$\Delta(p) = p\tt 1 + 1 \tt p$. Therefore,
$p\in\dd$, hence the only nonzero solutions to \eqref{Scocycle} occur
in degree two.

Now using \eqref{berel2} and \eqref{compa}, we get:
\begin{equation}\lbb{bsab1}
\be_{ab,cd} = a p^c_{bd} - b p^c_{ad}.
\end{equation}
The skew-symmetry $\be_{ab,cd} = -\be_{cd,ab}$ gives the equations
$p^a_{bc} = -p^a_{cb}$ and
\begin{equation}\lbb{bsab2}
a p^c_{bd} - b p^c_{ad} = c p^a_{bd} - d p^a_{bc}.
\end{equation}
{}From this we obtain that $p^a_{bd}$ lies in the linear span of $a,b,d$.
Comparing the coefficients in front of $ac$ in \eqref{bsab2}, we see
that the coefficient of $a$ in $p^a_{bd}$ is equal to the coefficient
of $c$ in $p^c_{bd}$. Call this coefficient $\phi_{bd}$; then 
$\phi_{bd} = -\phi_{db}$. Then comparison of other coefficients in 
\eqref{bsab2} shows that
\begin{equation}
p^a_{bc} = a \phi_{bc} + b \phi_{ca} + c \phi_{ab}
\end{equation}
for all $a,b,c\in\dd'$. Substitute this in \eqref{bsab1} to obtain
that $\be=\tau_\phi$ is trivial (cf.~\eqref{trivphi}).
\end{proof}
\begin{proposition}\lbb{cursimple}
Let $H = \ue(\dd)$,
and let $\g$ be a simple finite-dimensional Lie
algebra. If $L = \Cur \g$, then $\coh^2(L, \kk)\simeq \dd$.
\end{proposition}
\begin{proof}
Let $\be$ be a cocycle for $L$. We
will write $\be(a,b) = \be(1\tt a, 1\tt b)$ for $a,b\in\g$.
Then Jacobi identity leads to the equation
\begin{equation}\lbb{curcocycle}
\be(a,[b,c]) \tt 1 - 1 \tt
\be(b,[a,c]) = \Delta\bigl(\be ([a,b],c)\bigr).
\end{equation}
This immediately implies: $\be([a,b],c)
\in \dd+ \kk$. Since $[\g,\g] =
\g$ this shows that $\be(a,b) \in \dd+\kk$ for all $a,b\in \g$.

We can now solve the homogeneous equation \eqref{curcocycle} degree by degree.
Solutions of degree zero are cocycles of the Lie algebra $\g$, hence they are
all trivial. Solutions of degree one satisfy $\be(a, [b,c]) = \be([a,b],c)$,
and skew-symmetry implies $\be(a,b) = -S\bigl(\be(b,a)\bigr) = \be(b,a)$.
Therefore every such $\be$ is an invariant symmetric
bilinear form on $\g$ with values
in $\dd$. Any such bilinear form can be written as $\be(a,b) = (a|b)d$ where
$(\cdot|\cdot)$ is the Killing form on $\g$ and $d$ is some
element of $\dd$. Such cocycles $\be$ are nontrivial,  hence inequivalent
central extensions are in one-to-one correspondence
with elements of $\dd$.
\end{proof}

\begin{theorem}\lbb{centext}
Let $H=\ue(\dd)$ and $L$ be a simple Lie $H$-\psalg. Then $L$ may have a
nontrivial central extension 
{\rm(}given by \eqref{centr1}, \eqref{centr2}{\rm)}
only if\/{\rm:}

{\rm(i)} 
$L = \Cur \g$, in which case $\coh^2(L, \kk) \simeq \dd$ and
cocycles $\be$ are given by 
$\be_d(1\tt a, 1\tt b) = (a|b) d$ for $a,b\in\g$,
where $d\in \dd$ and
$(\cdot|\cdot)$ is the Killing form.

{\rm(ii)}
$L = \Cur_{H'}^H W(\dd')$ with $\kk s = \dd' \subset \dd$, $\dim \dd' =1$, 
in which
case $\coh^2(L, \kk)$ is of dimension one, generated by the Virasoro
cocycle $\be(1\tt s,1\tt s) = s^3$.

{\rm(iii)}
$L = \Cur_{H'}^H H(\dd', \chi', \om')$ with $\dd' \subset \dd$, in which
case $\coh^2(L, \kk)$ is isomorphic to the quotient of 
the space of all solutions 
$\be\in\dd$ to equations \eqref{cocycle2}, \eqref{cocycle3}
by the subspace $\kk(2s-x)$,
where $r \in \dd'\wedge\dd'$ is dual to $\om'$, $s\in\dd'$ is such that
$\chi' = \iota_s \om'$, and $x$ is given by \eqref{xr}.
\end{theorem}
\begin{proof}
(i) and (ii) follow from Propositions \ref{cursimple} and \ref{virasoro}(i),
and (iii) from a direct application of \leref{centrextension}.

For any other simple \psalg\ $L$, the strategy is to construct a continuous
family of \psalgs\ $L_t$,
indexed by $t\in \kk$ endowed with the Zariski topology, that are all
isomorphic to $L$ when $t \neq 0$, and whose fiber at $t=0$ is
one of the \psalgs\ already considered in Propositions 
\ref{virasoro}(ii), \ref{centrK} and \ref{centrS}.
Then, since $\coh^2(L_t, \kk)=0$ for $t=0$,
it will follow that $\coh^2(L_t, \kk)=0$
whenever $t$ lies in a neighborhood of $0$, hence for all $t\in\kk$.

In the case of a current pseudoalgebra over a $W$ or $S$ type Lie \psalg,
choose a basis $\{\d_i\}$ of $\dd$ that contains a basis of $\dd'$,
and construct the family $\dd'_t \subset \dd_t$ of Lie algebras indexed by 
$t \in \kk$ generated by elements
$\{\d^t_i\}$ with Lie bracket $[\d^t_i, \d^t_j] = t [\d_i, \d_j]^t$. 
Then for $t \neq 0$ we
have an isomorphism $\dd_t \ni \d^t_i \mapsto t \d_i \in \dd$, 
whereas $\dd_0$ is an abelian Lie algebra.
Then $\{\Cur_{H'_t}^{H_t} W(\dd'_t) \}_{t\in\kk}$,
where $H_t = \ue(\dd_t)$, $H'_t = \ue(\dd'_t)$,
is a family of \psalgs\ all isomorphic to $\Cur_{H'}^H W(\dd')$ for $t\ne0$.
The fiber of this family at $t = 0$ has been
shown in \prref{virasoro}(ii) to have no nontrivial central extensions.
In the same way, if we set $\chi_t(\d^t_i)= t\,\chi(\d_i)$, then
$\{\Cur_{H'_t}^{H_t} S(\dd'_t, \chi_t) \}_{t\in\kk}$ is a family of
\psalgs\ all isomorphic to $\Cur_{H'}^H S(\dd', \chi)$ for $t\ne0$,
and the fiber at $t=0$ is $\Cur_{H'_0}^{H_0} S(\dd'_0, 0)$ 
where $\dd_0' \subset \dd_0$ are
isomorphic to $\dd' \subset \dd$ as vector spaces but have trivial bracket. 

If $L$ is a current \psalg\ over $K(\dd', \th)$, for finite-dimensional
Lie algebras $\dd' \subset \dd$, choose a basis $\{a_i, b_i, s\}$ of $\dd'$
as in \leref{lrs1}, and complete it with $\{d_1, \dots, d_r\}$
to a basis of $\dd$. Then a continuous family $\{\dd_t\}$ 
of Lie algebras can be constructed for $t \neq 0$ 
as $\dd_t \simeq \dd$ spanned by $a^t_i = t a_i$, $b^t_i = t b_i$,
$s^t = t^2 s$, $d^t_i = t^2 d_i$, 
and by setting $a^0_i, b^0_i, c^0=-s^0$ to span a
Heisenberg algebra, and all brackets involving $d^0_i$ to be zero.
Define $\th_t \in (\dd'_t)^*$ by $\th_t(a^t_i) = \th_t(b^t_i) = 0$,
$\th_t(s^t) = -1$.
Then $\Cur_{H'_0}^{H_0} K(\dd'_0, \th_0)$ is the limit of the Lie \psalgs\
$\{\Cur_{H'_t}^{H_t} K(\dd'_t, \th_t) \}_{ t\neq 0}$, 
which are all isomorphic to 
$\Cur_{H'}^H K(\dd', \th)$, and the former Lie
\psalg\ is of the type treated in \prref{centrK}.
\end{proof}

\section{Application to the Classification of Poisson Brackets in Calculus of 
  Variations}\lbb{poisson}						      

In calculus of variations the phase space consists of $C^\infty$ vector 
functions
 $u = (u_1(x), \dots, u_r(x))$ where $u_i(x)$ are, for example, 
functions 
with
 compact support on a closed $N$-dimensional manifold. We will consider
 {\em linear} local Poisson brackets:
\begin{equation}\lbb{16.1}
\{ u_i(x), u_j(y)\} = \sum_\al B_{\al ij}(y) \, \d_y^\al \delta(x-y)
\end{equation}
where the sum runs over a finite set of multi-indices 
$\al = (\al_1, \dots, \al_N) \in\Zset_+^N$, the $B_{\al ij}$ 
are linear combinations  of the
 $u_k$ and of their derivatives 
$u_k^{(\ga)}: = \d_x^{\ga} u_k$, where $\d_x^\ga:= 
\bigl(\frac{\d}{\d x_1}\bigr)^{\ga_1}
\dotsm\bigl(\frac{\d}{\d x_N}\bigr)^{\ga_N}$, 
and $\delta(x-y)$ is the delta-function 
(defined by $\int f(x)\delta(x-y) \di x = f(y)$). By Leibniz rule and 
bilinearity,
the Poisson bracket \eqref{16.1} extends to arbitrary polynomials $P$ 
and $Q$ in 
the $u_i$ and their derivatives. Explicitly:
\begin{equation}\lbb{16.2}
\{P(x), Q(y)\} = \sum_{\al, \be, i,j} \frac{\d P(x)}{\d u_i^{(\al)}}
\frac{\d Q(y)}{\d u_j^{(\be)}} \,\d_x^\al \d_y^\be \{u_i(x), u_j(y)\}.
\end{equation}
This bracket, apart from bilinearity and Leibniz rule, should satisfy 
skew-commu\-tativity and the Jacobi identity. 

The basic quantities in calculus of 
variations are local functionals (Hamiltonians) $I_P = \int P(x)\di x$. Using 
bilinearity and integration by parts 
($ \int \frac{\d P}{\d x_i} Q \di x = -\int P \frac{\d Q}{\d x_i} \di x$), 
we get from
 \eqref{16.2} the following well-known formula:
\begin{equation}\lbb{16.3}
\{I_P, I_Q\} = \sum_{i,j} \iint \frac{\delta P(x)}{\delta u_i} 
\frac{\delta Q(y)}{\delta u_j} \{u_i(x), u_j(y)\} \di x \di y,
\end{equation}
where
\begin{equation}
\frac{\delta P(x)}{\delta u_i} = \sum_\al (-\d_x)^\al 
\frac{\d P(x)}{\d u_i^{(\al)}} 
\end{equation}
is the variational derivative.

More generally, one usually considers a class of Poisson brackets of the form:
\begin{equation}
\{u_i(x), u_j(y)\} = B_{ij}(y, \d_y^\al) \delta(x-y),
\end{equation}
where $B_{ij}$ are differential operators in $\d_y^\al$ whose coefficients 
are polynomials in $u_k^{(\gamma)}(y)$. 
Then the $r\times r$ matrix $B = (B_{ij})$ is 
called a Hamiltonian operator, and, integrating by parts, formula \eqref{16.3} 
can be rewritten in its most familiar form:
\begin{equation}
\{I_P, I_Q\} = \int \Bigl( B \frac{\delta P(x)}{\delta u} \Bigr)
\frac{\delta Q(x)}{\delta u}\di x.
\end{equation}

Given a Hamiltonian $h = \int P(x)\di x$, 
we have the corresponding Hamiltonian system of 
evolutionary partial differential equations:
\begin{equation}\lbb{16.4}
\dot u_i = \{h,u_i\} \equiv \sum_j B_{ij} \frac{\delta P}{\delta u_j},
\end{equation}
so that if another Hamiltonian $h_1$ is in involution with $h$, i.e.,
$\{h, h_1\} = 0$, then $h_1$ is an integral of motion of \eqref{16.4}, i.e., 
$\dot h_1 = 0$.

It is shown in \cite{DN1} and \cite{M} that for $r\geq 2$, any
Poisson bracket of hydrodynamic type (i.e.,
linear in the derivatives) under certain nondegeneracy conditions can be
transformed into a linear Poisson bracket of hydrodynamic type by a
change of the field variables. The latter Poisson brackets have been studied
rather extensively (see \cite{Do}, \cite{DN2} and references
there, \cite{GD}, \cite{M}, \cite{Z}).

Let $H = \Cset[\frac{\d}{\d x_1}, \dots, \frac{\d}{\d x_N}]$ be the universal
 enveloping algebra of the $N$-dimensional abelian Lie algebra $\dd$. 
Let $F = \bigoplus_{i=1}^r H u_i$ be the free 
$H$-module of rank $r$ on generators  $u_i$.
Consider a linear Poisson bracket \eqref{16.1}. For any three subspaces
$A,B,C$ of $F$, we will use the notation $\{A,B\}\subset C$ if for any
$a\in A$, $b\in B$ all coefficients in front of $\d_y^\al\delta(x-y)$
in $\{a(x),b(y)\}$ belong to $C$. We call a {\em linear Poisson algebra\/}
any $H$-submodule $L$ of $F$ which is closed under the Poisson bracket,
i.e., such that $\{L,L\}\subset L$. By an isomorphism of two such algebras
we mean a $\Cset$--linear isomorphism preserving Poisson brackets.

If $L$ is a linear Poisson algebra, we
define the {\em $\lambda$-bracket\/} ($\lambda = (\lambda_1, \dots, 
\lambda_N)$) on $L$ as the Fourier transform of
the linear Poisson bracket \eqref{16.1}:
\begin{equation}\lbb{16.5}
[u_i\,_\lambda u_j] = \tsum_\al \, \lambda^\al B_{\al ij}.
\end{equation}
Then we get a Lie conformal algebra in $N$ (commuting) indeterminates 
(defined in the same way as for the $N=1$ case in the introduction). Thus, the 
classification of linear Poisson algebras follows from the classification 
of Lie
$\ue(\dd)$-conformal algebras, where $\dd$ is the $N$-dimensional abelian 
Lie algebra.

Recall that the structure of a Lie conformal algebra is equivalent to the
structure of a Lie \psalg\ (see \seref{shconf}). The relationship between the
linear Poisson bracket \eqref{16.1} and the pseudobracket can be described
explicitly as follows:
\begin{align}
\lbb{ui*uj}
[u_i * u_j] &= \tsum_k\, P_{ij}^k(\d\tt1,1\tt\d) \tt_H u_k,
\intertext{if}
\lbb{uixujy}
\{u_i(x) , u_j(y)\} &= \tsum_k\, P_{ij}^k(\d_x,\d_y) 
                       \bigl( u_k(y) \de(x-y) \bigr)
\intertext{for some polynomials $P_{ij}^k$.
Note that equation \eqref{16.1} can always be written in the form 
\eqref{uixujy}. Indeed, if}
\lbb{uixujy2}
\{u_i(x) , u_j(y)\} &= \tsum_k\, Q_{ij}^k(\d_y,\d_t) 
                       \bigl( u_k(t) \de(x-y) \bigr)\big|_{t=y}
\intertext{for some polynomials $Q_{ij}^k$, then we have \eqref{uixujy}
with $P_{ij}^k (z,w) = Q_{ij}^k (-z,z+w)$. 
In this case, the $\lambda$-bracket \eqref{16.5} is given by:}
\lbb{uilauj}
\{u_i(x) , u_j(y)\} &= \tsum_k\, Q_{ij}^k(\la,\d) u_k.
\end{align}

\begin{remark}\lbb{poissext}
The constant terms of $B_{\al ij}$ in \eqref{16.1} give a central 
extension of the linear Poisson algebra
corresponding to the $B_{\al  ij}$ with 
constant terms removed. In terms of the associated Lie \psalgs\
this corresponds to a central extension by $\Cset$ with a trivial action of 
$H$. By \thref{exts}, these central extensions are parameterized by
$\coh^2(L,\Cset)$.
\end{remark}
\begin{examples}[cf.\ \cite{M}, \cite{DN2}, \cite{K4}, \cite{BKV}]\lbb{16.7}
\hfill
\begin{enumerate}
\item {\em General Poisson algebra} $W_{r,N}$, where $1\leq r\leq N$ 
$(1 \leq i,j \leq r)$:
\begin{displaymath}
\{u_i(x), u_j(y)\} = \frac{\d u_j(y)}{\d y_i} \delta(x-y) + u_j(y) 
\frac{\d}{\d y_i} \delta(x-y) + u_i(y) \frac{\d}{\d y_j} \delta(x-y).
\end{displaymath}

\item {\em Special Poisson algebra} $S_{r,N,\chi},$ where $2\leq r\leq N$ and 
$\chi = (\chi_1, \dots, \chi_r) \in \Cset^r$, is the following
subalgebra of $W_{r,N}$:
\begin{displaymath}
\biggl\{ \sum_{i=1}^r P_i(\d_x) u_i(x) \; \bigg| \:
 \sum_{i=1}^r \Bigl(\frac{\d}{\d x_i} + \chi_i\Bigr) P_i(\d_x) = 0 \biggr\}.
\end{displaymath}
It is generated over $H = \Cset[\frac{\d}{\d x_1}, \dots, \frac{\d}{\d x_N}]$
by elements
\begin{displaymath}
u_{ij}(x) = \Bigl(\frac{\d}{\d x_i} + \chi_i\Bigr) u_j(x)
          - \Bigl(\frac{\d}{\d x_j} + \chi_j\Bigr) u_i(x).
\end{displaymath}

\item {\em Hamiltonian Poisson algebra} $H_{2s, N}$, $2\leq r=2s\leq N$:
\begin{displaymath}
\{u(x), u(y)\} = \sum_{i=1}^s 
\Bigl(\frac{\d u(y)}{\d y_i} \frac{\d}{\d y_{i+s}} 
\delta (x-y) - \frac{\d u(y)}{\d y_{i+s}} \frac{\d}{\d y_i} \delta (x-y)\Bigr).
\end{displaymath}
We have an inclusion $H_{2s,N} \subset W_{2s,N}$ by letting 
\begin{displaymath}
u(x) = \sum_{i=1}^s 
\Bigl( \frac{\d u_{i+s}(x)}{\d x_i} - \frac{\d u_i(x)}{\d x_{i+s}} \Bigr).
%
\end{displaymath}

\item {\em Current Poisson algebra} $\Cur_N\g$ associated to a simple 
$r$-dimensional
Lie algebra $\g$ with structure constants $c_{ij}^k$ $(1 \leq i,j,k \leq r)$:
\begin{displaymath}
\{v_i(x), v_j(y)\} = \sum_{k=1}^r c_{ij}^k v_k(y) \, \delta(x-y).
\end{displaymath}

\item Semidirect sum of $W_{r,N}$ or one of its subalgebras $S_{r,N,\chi}, 
H_{2s, N}$ with $\Cur_N\g$ defined by ($1 \leq i \leq r, \; v(x) \in 
\Cur_N\g$):
\begin{displaymath}
\{ u_i(x), v(y)\} = \frac{\d v(y)}{\d y_i} \delta(x-y) + v(y) 
\frac{\d}{\d y_i} \delta(x-y).
\end{displaymath}
\end{enumerate}
\end{examples}

A subspace $I$ of a Poisson algebra $L$ is called an {\em ideal\/}
if it is invariant under taking Poisson brackets with elements of $L$,
i.e., if $\{L,I\} \subset I$.
A Poisson algebra $L$ is called {\em simple\/} (respectively {\em semisimple})
if the Poisson bracket is not identically zero and $L$ contains no
nonzero $H$-invariant ideals $I$
such that $I\ne L$ (respectively $\{I,I\} \ne0$).
Note that the Poisson algebras that we consider here are {\em finite},
i.e., finitely generated as $H$-modules.

Then Theorems \ref{classify} and \ref{tsemisim} 
and \coref{csubwd} imply:

\begin{theorem}
{\rm(i)} 
Any simple linear Poisson algebra is isomorphic to one of the 
Poisson algebras $W_{r,N}, S_{r,N,\chi}, H_{2s, N}, \Cur_N\g$.

{\rm(ii)} 
Any semisimple linear Poisson algebra is a direct sum of simple ones and 
of the semidirect sums described in \exref{16.7}{\rm(5)}.
\end{theorem}

\begin{remark}
It follows from \reref{poissext} and the results of \seref{centrexts}
that all nontrivial
central extensions of simple Poisson algebras are described by the following
$2$-cocycles. 

For $\al \in \Cset^r$ let
\begin{displaymath}
\psi_\al(x,y) = \sum_{i=1}^r \al_i \frac{\d}{\d y_i} \delta(x-y).
\end{displaymath}

Then all nontrivial $2$-cocycles for $H_{r,N}$ are:
\begin{displaymath}
\ga_\al(u(x),u(y)) = \psi_\al(x,y).
\end{displaymath}

All nontrivial $2$-cocycles for $\Cur_N \g$ are:
\begin{displaymath}
\ga_\al(v_i(x), v_j(y)) = b_{ij} \psi_\al(x,y),
\end{displaymath}
where $b_{ij} = (v_i|v_j)$ is the invariant scalar product.

The Poisson algebra $W_{r,N}$ has a nontrivial central extension iff $r=1$,
and in the latter case it is given by the well-known Virasoro cocycle:
\begin{displaymath}
\ga(u(x),u(y)) = \left (\frac{\d}{\d y}\right )^3 \delta(x-y).
\end{displaymath}

The Poisson algebra $S_{r,N,\chi}$ has no nontrivial central extensions
if $r>2$ or $\chi \neq 0$, and $S_{2,N,0} \simeq H_{2,N}$ has nontrivial 
central extensions described above.
\end{remark}

\bibliographystyle{amsalpha}


\end{document}